\documentclass[a4paper, 12pt]{article}

\usepackage{amsmath} 
\usepackage{theorem}
\usepackage{amssymb}
\usepackage{amsfonts}
\usepackage{latexsym}
\usepackage{amscd}
\usepackage{diagramb}
\usepackage{yfonts}
\usepackage{hyperref}  
\usepackage{graphicx}
\usepackage{enumerate}
\usepackage{cancel}
\usepackage{extarrows}
\usepackage{upgreek}
\usepackage{gensymb}

\usepackage{pxfonts}

\usepackage{stmaryrd}
\usepackage{float}

\input xy
\xyoption{all}

   \newcommand{\axiom}[1]{\textbf{\hypertarget{#1}{(#1)}}}
    \newcommand{\axiomref}[1]{\hyperlink{#1}{(#1)}}

\DeclareMathAlphabet{\mathpzc}{OT1}{pzc}{m}{it}

\usepackage{xspace,colortbl}
\usepackage{color}

\definecolor{verde}{rgb}{0.,0.7,0.}
\definecolor{indigo}{rgb}{.18, .34, .78}
\definecolor{indigo1}{rgb}{.18, .24, .78}
\definecolor{indigo2}{rgb}{.18, .14, .78}
\definecolor{indigo3}{rgb}{.18, 0., .78}
\definecolor{rojo}{rgb}{1,0,0}
\definecolor{negro}{rgb}{0,0,0}
\definecolor{lila}{rgb}{.46, .16, .78}
\definecolor{lila1}{rgb}{.46, .16, .86}
\definecolor{lila2}{rgb}{.56, .16, .86}
	\definecolor{lila3}{rgb}{.63, .16, .78}
\definecolor{lila4}{rgb}{.7, .16, .78}
\definecolor{lila5}{rgb}{.78, .26, .78}
\definecolor{lila6}{rgb}{.6, 0., .78}
\definecolor{grey}{rgb}{0.7,0.7,0.7}
\definecolor{darkgrey}{rgb}{0.6,0.6,0.6}

\theoremstyle{plain}
\theoremheaderfont{\normalfont\bfseries}

\newtheorem{thm}{Theorem}[section]

\newtheorem{lma}[thm]{Lemma}
\newtheorem{cor}[thm]{Corollary}
\newtheorem{defn}[thm]{Definition}

\theorembodyfont{\rmfamily}

\newtheorem{rem}[thm]{Remark}
\newtheorem{prop}[thm]{Proposition}

\newcommand{\qed}{\hfill\quad\fbox{\rule[0mm]{0,0cm}{0,0mm}}  \par\bigskip}

\newcommand{\x}{\mbox{-}}

\newcommand{\w}{\hspace{-0,06cm}}
\newcommand{\s}{\hspace{0,06cm}}
\newcommand{\Cat}{\operatorname {Cat}}
\newcommand{\Set}{\operatorname {Set}}

\newcommand{\Dbl}{\operatorname{\mathbb D bl}}
\newcommand{\Ps}{{\rm Ps}}

\newcommand{\Fun}{{\rm Fun}}

\newcommand{\Int}{\operatorname {Int}}

\newcommand{\Aa}{{\mathbb A}}
\newcommand{\Bb}{{\mathbb B}}
\newcommand{\Cc}{{\mathbb C}}

\newcommand{\OO}{{\mathbb O}}
\newcommand{\I}{{\mathbb I}}

\newcommand{\Del}{\boxtimes}
\newcommand{\comp}{\circ}
\newcommand{\iso}{\cong}
\newcommand{\ot}{\otimes}

\newcommand{\C}{{\mathcal C}}

\newcommand{\M}{{\mathcal M}}
\newcommand{\D}{{\mathcal D}}
\newcommand{\F}{{\mathcal F}}
\newcommand{\G}{{\mathcal G}}
\newcommand{\Oo}{{\mathcal O}}

\newcommand{\A}{{\mathcal A}}
\newcommand{\B}{{\mathcal B}}
\newcommand{\E}{{\mathcal E}}
\newcommand{\V}{{\mathcal V}}

\def\ul{\underline}

\newcommand{\crta}{\overline}

\newcommand{\Id}{\operatorname {Id}}

\def\V{{\mathcal V}}  

\newcommand{\Lax}{\operatorname{\mathbb L ax}}
\newcommand{\Psd}{\operatorname{\mathbb P sd}}

\newcommand{\tlabel}[1]{\label{t:#1}}
\newcommand{\tref}[1]{T.~\ref{t:#1}}

\newcommand{\cref}[1]{C.~\ref{c:#1}}

\newcommand{\lelabel}[1]{\label{le:#1}}
\newcommand{\leref}[1]{Lemma~\ref{le:#1}}
\newcommand{\eqlabel}[1]{\label{eq:#1}}
\newcommand{\equref}[1]{(\ref{eq:#1})}
\newcommand{\thlabel}[1]{\label{th:#1}}
\newcommand{\thref}[1]{Theorem~\ref{th:#1}}
\newcommand{\delabel}[1]{\label{de:#1}}
\newcommand{\deref}[1]{Definition~\ref{de:#1}}
\newcommand{\prlabel}[1]{\label{pr:#1}}
\newcommand{\prref}[1]{Proposition~\ref{pr:#1}}
\newcommand{\colabel}[1]{\label{co:#1}}
\newcommand{\coref}[1]{Corollary~\ref{co:#1}}

\newcommand{\rmlabel}[1]{\label{rm:#1}}
\newcommand{\rmref}[1]{Remark~\ref{rm:#1}}
\newcommand{\selabel}[1]{\label{se:#1}}
\newcommand{\seref}[1]{Section~\ref{se:#1}}
\newcommand{\sslabel}[1]{\label{ss:#1}}
\newcommand{\ssref}[1]{Subsection~\ref{ss:#1}}
\newcommand{\ssslabel}[1]{\label{ss:#1}}
\newcommand{\sssref}[1]{Subsection~\ref{ss:#1}}

\newcommand*{\threefrac}[3]{%
  \begin{array}{@{\,}c@{\,}}%
    #1\\
    \hline
    #2\\
    \hline
    #3%
  \end{array}%
}

\newcommand*{\fourfrac}[4]{%
 \begin{array}{@{\,}c@{\,}c@{\,}}%
    #1\\
    \hline
    #2\\
    \hline
    #3\\
    \hline
		#4%
  \end{array}%
	}

\newdir{:=}{{}}
\newcommand{\fit}[3]{\ar@{:=}@/{#3}/[#1] |{\Downarrow #2} }

\begin{document}

\title{Gray (skew) multicategories: \\ double and Gray-categorical cases}

\author{Bojana Femi\'c \vspace{6pt} \\
{\small Mathematical Institite of  \vspace{-2pt}}\\
{\small Serbian Academy of Sciences and Arts } \vspace{-2pt}\\
{\small Kneza Mihaila 36,} \vspace{-2pt}\\
{\small 11 000 Belgrade, Serbia} \vspace{-2pt}\\
{\small femicenelsur@gmail.com}}

\date{}
\maketitle

\begin{abstract}
We construct in a unifying way skew-multicategories and multicategories of double and Gray-categories that we call Gray (skew) multicategories. We study their different versions depending on the types of functors and higher transforms. We construct Gray type products by generators and relations 
and prove that  Gray skew-multicategories are closed and representable on one side, and that the Gray multicaticategories taken with the strict type of functors are representable. We conclude that the categories of double and Gray-categories with strict functors underlying Gray (skew) multicategories are skew monoidal, respectively monoidal, depending on the type of the inner-hom and product considered. The described Gray (skew) multicategories we see as prototypes of general Gray (skew) multicategories, which correspond to (higher) categories of 
higher dimensional internal and enriched categories. 
\end{abstract}

{\small {\em Keywords:} double categories, Gray monoidal product, Gray-categories, multicategories}

{\em 2020 Mathematics Subject Classification: 18N20, 18N10}. \\

\section{Introduction}

The classical Bifunctor Theorem states that two functors $\A\to\C$ and $\B\to\C$ among categories give rise to a functor from the Cartesian product $\A\times\B\to\C$, provided that a certain commutativity equation is fulfilled, and the obvious converse holds, \cite{McL}. 
A 2-categorical version of the result was proved in \cite{FMS}. When working with 2-categories, one has the liberty to choose pseudo or 
(op)lax functors, the authors treated the lax case. The mentioned commutativity equation among two functors is substituted by a 2-cell complying with some axioms in the 2-categorical setting. Such a data they naturally called {\em distributive law}. Though, they were treated by Gray in \cite{Gray} under the name {\em quasi-functors}. Gray used quasi-functors to construct a tensor product among 2-categories.

In \cite{Fem} we proved the Bifunctor theorem for double categories and lax double functors and we explored the limitations of 
obtaining a Gray type tensor product for $Dbl_{lx}$, the category of double categories and lax double functors. 
Namely, 
we introduced a product $\Aa\ot\Bb$ among two double categories $\Aa$ and $\Bb$ based on the candidate 
$\Lax_{hop}(\Aa,\Bb)$ for inner-hom in $Dbl_{lx}$. The double category $\Lax_{hop}(\Aa,\Bb)$ has lax double functors $\Aa\to\Bb$ as objects, 
horizontal oplax transformations and vertical lax transformations as horizontal, respectively vertical 1-cells, and corresponding modifications for 2-cells. Following Gray's ideas, to construct $\ot$ we characterized a lax double functor $\Aa\to\Lax_{hop}(\Bb,\Cc)$, 
we described the corresponding {\em lax double quasi-functors} $\Aa\times\Bb\to\Cc$ (different from usual functors on Cartesian product!) 
and then determined $\Aa\ot\Bb$ by generators and relations reading off the axioms satisfied by the quasi-functors. However, 
$-\ot\Bb$ and $\Lax_{hop}(\Bb,-)$ can not be adjoint functors on $Dbl_{lx}$, because there is no way to determine any kind of a double functor $\Aa\ot\Bb\to\Cc$ out of a lax double quasi-functor $\Aa\times\Bb\to\Cc$ other than a strict one. Though, we do have that this product strictifies lax double quasi-functors. It is moreover so that $\Lax_{hop}(\Bb,-)$ is even not a functor on $Dbl_{lx}$. 
Open for a future research remained the question of what else can be said about the product that we introduced, and which changes 
are to be made in order to obtain a monoidal, or at least skew-monoidal structure on some version of the category of double categories. 



Recently, in some ways similar construction was carried out in \cite{BL}, where skew-monoidal category structures are obtained on the category of Gray-categories with strict Gray functors. To obtain this the authors used the result from \cite{BLack} that left representable skew-multicategories induce skew-monoidal categories. The multimaps that they used to construct a (left representable) skew-multicategory of Gray-categories are nothing but the quasi-functors from the construction of Gray in the context of Gray-categories. Following this multicategory approach, including Hermida's equivalence between representable multicategories and monoidal categories, we embarked on 
a research that would tackle by same strategies and in a unifying way the questions of monoidality and skew-monoidality of the categories of double categories (including 2-categories) and Gray-categories, with a hope that analogous proofs and findings can be carried out in higher dimensional enriched and internal categories. 

\medskip

We denote generically by $\Oo$ a collection of double categories, or Gray-categories with the motivation that 
$\Oo$ can stand also for their higher dimensional analogues. 
We start by introducing, for such a generic $\Oo$, a quasi-functor of two variables, or a binary quasi-functor. By definition it corresponds uniquely to a functor $\F:\A\to\Oo^*_\bullet(\B,\C)$, where $\Oo^*_\bullet(\B,\C)$ is a hom-set, but also a hom-object in $\Oo$. We recall our characterization of binary double categorical quasi-functors from \cite{Fem1} in 
\prref{char df}. 
Depending on the variations of hom-sets, where $*$ stands for the nature of functors (strict, lax, oplax, pseudo) and $\bullet$ for the nature of transformations (we neglect notations for higher transforms, but will usually preserve the same type on modifications as taken for transformations), we then describe how one generates what we call a {\em meta-product} $\ot^*_\bullet$ out of the obtained data for quasi-functors. Analyzing two quasi-functors having opposite $\bullet$-types, we apply Crans's type of argument to rule out the 
``strict functors-weak transformations'' possibility for the type of the inner-hom, {\em i.e.} of the quasi-functors and meta-products for Gray-categories. 
For double categories and Gray-categories we introduce transformations, modifications and perturbations on quasi-functors, obtaining double {\em i.e.} Gray-categories $q\x\Oo^*_\bullet(\B,\C)$ (see \ssref{tr for qf dim2} and \ssref{Gray-cat of qf}). After 
having introduced some notions and results for double categories in Section 2, Sections 3 and 4 are dedicated merely to double categories and Sections 5 and 6 basically to Gray-categories. In \seref{n-double cat} 
we define lax ternary quasi-functors and their oplax transformations for double categories and study the necessary ingredients for obtaining a multicategory of double categories. We also study duality among double categories by which horizontal and vertical 1-cells flip their directions. Though, the obtained results will not play the role in our multicategory. In 
\seref{meta-lax double} we study representability of the multicategory of double categories with respect to the meta-products. Here is where we show that only the strict functor case yields representability, and that weakening in the sense of studying mixed functor structures, that will appear in the construction of skew-multicategories, one gets one-sided representable skew-multicategories. This solution was suggested in \cite{Bour:SS} and employed for left representability in \cite{BL}. 

In Section 5 we turn to Gray-categories. We briefly review binary quasi-functors from \cite{BL} and characterize them in 
\prref{quasi-Gray} (this is a variation of the ``compact definition of binary maps'' from \cite[Section 4.5]{BL}). We introduce a concise, although somewhat loose, fraction notation for transversal composition of 3-cells. By it we omit writing the domain and codomain 2-cells 
(although, if one adds them, one recovers the unique possible way to write the full expression), and we also neglect identity and interchange 
3-cells and thus substantially simplify the annotation of the axioms expressed in terms of 3-cells. As announced before, we then introduce higher cells on binary quasi-functors for Gray-categories obtaining a Gray-category $q\x G\x\Cat^{lx}_{oplx}(\A,\B)$. In \ssref{incub} we recall 
ternary quasi-functors from \cite{BL}, we introduce their lax and oplax transformations, and we also introduce $n$-ary quasi-functors,  determined by a collection of 4-ary quasi-functors. Then we proceed to construct tight multimaps for the left skew-multicategory: they consist 
of quasi-functors whose left-most component functor is strict and the rest ones are weak, of type $*$. We construct evaluation quasi-functors, 
and closedness isomorphisms connecting binary and ternary quasi-functors, and ternary and 4-ary quasi-functors (these isomorphisms we call for short 
binary-to-ternary and ternary-to-4-ary). The latter we do by first observing that for $n=3,4$ one has isomorphisms 
$$(A_n)\quad q_{n-1}\x\Oo(\A_1\times...\times\A_{n-1},\Oo(\A_n,\C))\iso \Oo_{st}(\A_1,q_{n-1}(\A_2\times...\times\A_n,\C)),$$
and then by proving 1-1 correspondences 
$$\hspace{0,38cm}\text{transformations of binary quasi-functors}\quad \leftrightarrow \quad \text{ternary quasi-functors}\vspace{-0,2cm}$$
$$\text{transformations of ternary quasi-functors}\quad \leftrightarrow \quad \text{4-ary quasi-functors},$$
in \leref{gain incubator} and \leref{gain mecon}. Then by taking a quasi-functor $\F$ on the right in $(A_n)$ and a 1-cell $f$ in $\A_1$, by our characterization of quasi-functors we obtain a transformation of quasi-functors $\F(f)$ of $n\x 1$-variables on the left in $(A_n)$, which gives us an $n$-ary quasi-functor, by  the above correspondence. Thus we end up obtaining  
$$\hspace{-0,54cm} (B_n)\quad \qquad\qquad q_n\x\Oo(\A_1\times...\times\A_n,\C)\iso q_{n-1}\x\Oo(\A_1\times...\times\A_{n-1},\Oo(\A_n,\C)). $$
This trend seems to repeat for higher $n$, which may be useful for considering higher dimensions. This approach to the proof of the closedness of the multicategory is one of the differences of our study of the Gray-categorical case with respect to that of \cite{BL},  
which is enabled by our explicit definitions of transformations of quasi-functors. We finish Section 5 by proving the substitution for the (skew) multicategory. Here appears the second difference of our proofs with respect to \cite{BL}: while the authors use the isomorphisms $(B_n)$ and ``duality bijections'' 
$$d_n: q_n\x G\x\Cat^{ps}_{lx}(\Pi_{i=1}^n\A_i,\B)\to q_n\x G\x\Cat^{ps}_{lx}(\Pi_n^{i=1}\A_i^{co},\B^{co}),$$ 
where $(-)^{co}:G\x\Cat_{ps}\to G\x\Cat_{ps}$ is an isomorphism of categories that flips the 2-cells, instead of the duality bijections we use isomorphisms $(A_n)$. Here $G\x\Cat_{ps}$ is the category of Gray-categories with pseudo-Gray-functors. The use of pseudo-Gray-functors 
is crucial in the constructions of \cite{BL}, as they indeed need invertibility of the pseudofunctor structures in the duality bijections (so to prove substitution and obtain a skew-multicategory), but also they use them - via the pseudofunctor classifier for Gray-categories 
from \cite{Go} - to prove that $G\x\Cat^{ps}_{lx}(-,-):G\x\Cat_{ps}^{op}\times G\x\Cat_{ps}\to G\x\Cat_{ps}$ is a (bi)functor (and similarly $G\x\Cat^{ps}_{lx}(-,-)$). The third difference in our approach is that we prove that $G\x\Cat^*_\bullet(-,-)$, for any $*$ and $\bullet$, is a functor by direct inspection of axioms. The obtained (skew) multicategories for double and Gray-categories we call {\em Gray (skew) multicategories}, as the basis for all the necessary ingredients for their construction was indeed settled down by the Gray's construction for 2-categories in \cite{Gray}. 

In the final section we deduce properties of the Gray skew-multicategories $\OO^{(st,*)}_\bullet$: one-sided closedness and representability, and also representability 
of the multicategory $\OO^{st}_\bullet$. We also conclude on (skew) closedness and (skew) monoidality of the induced underlying categories 
of the (unary tight) multicategories, and we present in Table \ref{table:1} the classification of these findings for all the studied versions of categories of 2-categories, double and Gray-categories. In \ssref{conjecture} we conjecture the existence of a lax double functor classifier and use its Gray-categorical companion from \cite[Section 14.2]{Gur} to give other descriptions of the meta-products $\ot^*_\bullet$ for double and Gray-categories, for $*$ being $lax$ and $pseudo$. This way we add higher dimensional analogues to the result for 
2-categories from \cite{Nik} (see isomorphism (146) in {\em loc.cit.}). As announced above, we prove in \ssref{lax and strict} that 
$G\x\Cat^{st}_\bullet(-,-)$ and 
 $G\x\Cat^{ps}_\bullet(-,-)$ are functors, and using our meta-products, constructed by generators and relations, we obtain analogously as above, the (representable) multicategory $G\x\Cat^{st}_\bullet$. Here we see two more differences with respect to the construction in \cite{BL}. We conclude the representability of the (skew) multicategory by constructing explicitly the meta-product (in case of Gray-categories by giving direct indications for that). Given that the skew-multicategory is closed, the authors prove representability by proving the existence of a nullary map classifier and that each endo-hom functor has a left adjoint (see their Subsections 2.5 and 5.2). The existence of left adjoint they prove by showing that the category $G\x\Cat_{st}$ is locally finitely presentable. The other mentioned difference is that rather than proving directly that  $G\x\Cat^{st}_\bullet$ is a multicategory, it is concluded in \cite{BL} that 
$\Lax^\sharp=G\x\Cat^{st}_{lx}$ and $\Psd^\sharp=G\x\Cat^{st}_{ps}$ are multicategories by considering them as {\em sub-multicategories} of 
$\Lax=G\x\Cat^{ps}_{lx}$ and $\Psd=G\x\Cat^{ps}_{ps}$, respectively, according to their Proposition 7.2. So, in this approach one needs to pass via the larger multicategories for whose construction the duality bijections (for the proof of substitution) were used. In the last subsection we draw concluding remarks, we highlight the above exposed differences between the two studies of skew-multicategories for Gray-categories, and name a possible direction for future research.

\section{Quasi-functors of two variables}

In this work we are going to present a generalized version of Gray's procedure in \cite{Gray} by which he proved that the category 
$2\x\Cat$ of 2-categories is monoidal closed. To start, in this Section we will introduce ``general quasi-functors''. They will play the 
role of multimaps for the construction of a multicategory that we will call {\em Gray multicategory}. 
We start by an overview of some existing constructions of Gray products.

\subsection{Some existing Gray products and multicategories}

The internal hom that Gray used in his proof of monoidal closedness of $2\x\Cat$ is the 2-category of 2-functors, lax transformations 
and modifications. (It is the laxness of transformations 
that doesn't allow closedness with the Cartesian product, as explained nicely in \cite[Section 5.1]{GP}.) 

In \cite{Gabi} the author proved 
that $Dbl$, the category of double categories and double functors, is monoidal closed with Gray type of monoidal product. 
The internal hom she used has double functors for objects, horizontal and vertical (strong) pseudotransformations for horizontal 
and vertical 1-cells, respectively, and modifications. 

In the proof of monoidality, Gray introduced {\em quasi-functors} (of $n$-variables, $n\geq 2$), he defined the tensor product via generators and relations, and proved that it represents quasi-functors. 
From the representability by Yoneda one obtains the existence of monoidality constraints. For the same reason, any axiom that should 
hold among these constraints corresponds to axioms among arrows between domains of quasi-functors, which are Cartesian products, so the axioms 
trivially hold. In essence, Gray has constructed a multicategory of 2-categories, a binary map classifier (we recall these notions further below) and proved that $2\x\Cat$ is monoidal.  

B\"ohm took a different approach: she proved that the inner-hom is representable, which guaranteed the existence of a tensor product. She obtained monoidality constraints and the isomorphisms determining monoidal closedness of the category using Yoneda lemma. 
The pentagon and triangle axiom she proved by using Yoneda lemma again and with diagram chasing arguments. 

B\"ohm's result that $Dbl$ is closed and monoidal generalizes Gray's one for  $2\x\Cat$ in that double categories have an extra direction 
for 1-cells. However, if one forgets the vertical direction, her transformations being pseudo gives a stricter version of a result, compared to lax transformations in Gray's case. Studying her proof one easily sees that it can not be copied to generalize from pseudo 
to (op)lax transformations of double functors. Namely, while the extra natural transformation $l$ passes smoothly to the lax-oplax case, in the definition of $r$ one really needs invertibility of the structure 2-cells of (both horizontal and vertical) transformations. 

In \cite{Fem} we constructed candidates for a product $\ot$ and inner-hom for the category $Dbl_{lx}$ of double categories and lax 
double functors. 
We did this following \cite{Gray}: we introduced {\em double quasi-functors} and then $\ot$ by generators and relations. 
Our candidate for inner-hom was the double category $\Lax_{hop}(\Aa,\Bb)$, which has lax double functors 
$\Aa\to\Bb$ as objects, horizontal oplax transformations and vertical lax transformations as horizontal, respectively, vertical 1-cells, and corresponding modifications for 2-cells. However, we found: 
\begin{enumerate} [a)]
\item the laxness of double functors hinders $\Lax_{hop}(-,-)$ to be a functor: in the second entry one really needs pseudo double functors (so, one lacks the closedness for $Dbl_{lx}$); 
\item starting from a lax double functor $\Aa\to\Lax_{hop}(\Bb,\Cc)$ (one gets, as an intermediate step, a lax double 
quasi-functor $\Aa\times\Bb\to\Cc$), there is no way to determine any kind of a double functor 
$\Aa\ot\Bb\to\Cc$ other than a strict double one. (For this reason we say that this product strictifies lax double quasi-functors.)
\end{enumerate}
In conclusion, our $\ot$ (whose construction has $\Lax_{hop}(\Aa,\Bb)$ as a starting point) can not yield a ``monoidal and closed'' structure 
for $Dbl_{lx}$: to get a ``monoidal and closed'' category of double categories, one should change $Dbl_{lx}$ by $Dbl$. 
Still, one can ask the following two questions: 
\begin{enumerate} 
\item can $\ot$ equip $Dbl$ with a monoidal or skew-monoidal structure?, and 
\item changing 
from $\Lax_{hop}(\Aa,\Bb)$ to $\Dbl_{hop}(\Aa,\Bb)$ - where the functors are strict double, and transformations and modifications 
remain the same - can one get a closed monoidal structure on $Dbl$ (which will then be different than the one from \cite{Gabi}), and what 
relation is there between these two monoidal structures? 
\end{enumerate} 
We will address these questions by studying an even more general construction. 

\medskip

Lately, Bourke and Lobbia proved in \cite{BL} that the category $G\x\Cat$ of Gray-categories and strict Gray functors is left closed and 
skew-monoidal with four different versions of tensor products and inner-homs. In the first two inner-homs the functors are pseudo, but they differ in that the transformations and modifications are lax in one and pseudo in the other, and finally both have the corresponding perturbations. In the last two tensor products, studied in the last section, the inner-homs have strict Gray functors for objects, and the rest of the cells are as in the previous two inner-homs. 
The authors used the approach of multicategories: due to \cite[Section 6.2]{BLack}, left representable skew-multicategories give rise to skew-monoidal categories. Before we pursue, let us recall what closed (skew) multicategories are. For more on multicategories we refer the reader to \cite{Lein}.

\begin{defn}
A multicategory consists of:
\begin{itemize} 
\item a collection of objects, \vspace{-0,2cm}
\item for each list of objects $a_1,...,a_n$ for $n\geq 0$ and an object $b$, a set $\M_n(a_1,...,a_n;b)$, \vspace{-0,2cm}
\item for each object $a$ an element $1_a\in\M_1(a;a)$, \vspace{-0,2cm}
\item for all lists $\crta{a_i}, i=1,...,n, b_1,...,b_n$ and objects $c$ a function, called {\em substitution}:
$$\M_n(b_1,...,b_n,c)\times\Pi_{j=1}^{n}\M_{k_i}(\crta{a_i};b_i)\to\M_{\sum_{i=1}^n k_i}(\crta{a_1},...,\crta{a_n};c)$$
$$(g, f_1,...,f_n)\mapsto g\circ(f_1,...,f_n)$$
satisfying a natural associativity and two identity axioms. Here $\crta{a_i}$ is a short annotation for a list of objects, 
and the elements of $\M_n(a_1,...,a_n;b)$ are called {\em multimaps}. 
\end{itemize}
\end{defn}

\begin{rem} \rmlabel{multcat-cat}
An equivalent definition of a multicategory can be found in \cite[Proposition 1.1]{L}, by which a multicategory $\M$ consists of a 
category $\C$, functors $\M_n(-;-):(\C^n)^{op}\times\C\to\Set$ so that $\M_1(-;-)=\C(-,-):\C^{op}\times\C\to\Set$, and the analogous substitution functions that are natural in each variable $a_1,...,a_m,b_1,...,b_n,c$, so that $\C$ is the underlying category of the multicategory $\M$. 
\end{rem}

A skew-multicategory differs from a multicategory in that in it one differentiates {\em tight} and {\em loose} components 
$\M_n^t(a_1,...,a_n;b)\subseteq\M_n^l(a_1,...,a_n;b)$, respectively, for $n\geq 1$, so that the nullary components are loose, the unary identities in $\M_1^l(a;a)$ are tight, and the multimaps $g\circ(f_1,...,f_n)$ obtained by substitution are tight just when $g$ and $f_1$ 
are tight.

\begin{defn}
A multicategory $\M$ is said to be {\em left closed} if for all objects $b,c$ there is an object $[b,c]$ and a binary map 
$e_{b,c}:[b,c], b\to c$ so that the induced functions $\M_n(\crta a;[b,c])\to\M_{n+1}(\crta a,b;c)$ are bijections for all $n\geq 0$. 
If there exists an object $\{b,c\}$ and a binary map $\crta{e_{b,c}}:b,\{b,c\}\to c$ for which the induced functions 
$\M_n(\crta a;\{b,c\})\to\M_{n+1}(b,\crta a;c)$ are bijections for all $n\geq 0$, then $\M$ is said to be {\em right closed}. 
$\M$ is called {\em biclosed} if it is both left and right closed. 
\end{defn}

\begin{defn}
A skew-multicategory $\M$ is said to be {\em left/right closed} if its multicategory $\M^l$ of loose maps is left/right closed,  
the map $e_{b,c}:[b,c], b\to c$ (resp. $\crta{e_{b,c}}:b,\{b,c\}\to c$) is tight and the induced functions 
$\M_n^t(\crta a;[b,c])\to\M_{n+1}^t(\crta a,b;c)$ (resp. $\M_n^t(\crta a;\{b,c\})\to\M_{n+1}^t(b,\crta a;c)$) are bijections for $n\geq 1$. 
\end{defn}

\begin{defn} \delabel{lrep}
A multicategory $\M$ is said to be {\em left representable} if for any list of objects $\crta a=a_1,..,a_n$ there exists an object %
$m(a_1,..,a_n)$ and a multimap $j_{\crta a}:a_1,..,a_n\to m(a_1,..,a_n)$ that for every object $c$ and list $\crta b$ (of length $l$) induces bijections 
$$\M_1(m(a_1,..,a_n);c)\to\M_n(a_1,..,a_n;c)$$
 and 
$$\M_{1+l}(m(a_1,..,a_n),\crta b;c)\to\M_{n+l}(a_1,..,a_n,\crta b;c).$$
Right representability is similarly defined. 
\end{defn}

\begin{defn}
A multicategory $\M$ is said to be {\em representable} if for any list of objects $\crta a=a_1,..,a_n$ there exists an object %
$m(a_1,..,a_n)$ and a multimap $j_{\crta a}:a_1,..,a_n\to m(a_1,..,a_n)$ that for every object $c$ and lists $\crta x, \crta b$ (of lengths 
$k$ and $l$, respectively) induces bijections 
$$\M_1(m(a_1,..,a_n);c)\to\M_n(a_1,..,a_n;c)$$
 and 
$$\M_{k+1+l}(\crta x,m(a_1,..,a_n),\crta b;c)\to\M_{k+n+l}(\crta x,a_1,..,a_n,\crta b;c).$$
\end{defn}

The multimap $j_{\crta a}$ satisfying the first isomorphism in the above two definitions is called an {\em $n$-ary map classifier}. The second 
isomorphism in \deref{lrep} makes it {\em left universal}, while the second isomorphism in the latter definition makes it {\em universal}. 

\smallskip

The objective of our research is to: answer the above listed questions arisen regarding double categories, and to study the case of Gray-categories in the same style as we address the double categorical case, using multicategories, so to give a unifying approach to the study of (skew) monoidal (one-sided) closedness of the categories of internal or enriched categories in law dimensions. Let's begin.

\subsection{General quasi-functors}

Let $\V_{m+1}$ and $\I_{m+1}$ denote the categories of weak $m$-categories and internal $m$-categories, respectively. The latter are 
$(1\times m)$-categories in the sense of \cite[Section 1]{Shul}. The dimension of the categories in  $\V_{m+1}$ and $\I_{m+1}$ is $m$. 
For example, double categories are internal 1-categories and we consider that their dimension is 2. We choose the labeling $\V_{m+1}$ for 
weak $m$-categories in accordance with iterated enrichment, which generates strict $m$-categories, where $\V_{m+1}=\V_m\x\Cat$ starting with $\V_0\x\Cat=\Set$.

\medskip

To define general quasi-functors we will denote by $\Oo$ any of $\V_{m+1}$ or $\I_{m+1}$ for some fixed $m$. 
As motivating examples we think of the categories $2\x\Cat$ of 2-categories, $Dbl_*$ of double categories and any sort of double functors 
(strict, lax, oplax or pseudo ones) and similarly of the category $G\x\Cat_*$ of Gray-categories. 
We will use the following notational conventions. When we write $\Oo$ for the category, we are not specifying which kind of morphisms, {\em i.e.} functors we have in mind. When we do want to specify morphisms, we will denote by $\Oo_*$ the category with 
functors of type $*$. There are various candidates for inner-hom for $\Oo$. Namely, concerning {\em e.g.} Gray-categories, one chooses among strict, lax, oplax or pseudo functors, strict, lax, oplax or pseudo transformations, similarly for modifications, and then one takes the corresponding perturbations. For $\A,\B\in\Oo$ let us denote by $\Oo^*_\bullet(\A,\B)$ a generic inner-hom candidate for $\Oo$, where 
$*$ is a label for the type of functors and $\bullet$ for the type of transformations. Strictly speaking, we should also add labels for modifications and perturbations, but we will omit to do so for simplicity reasons. In the enriched case we will consider that for now it is irrelevant 
if lax or oplax version is meant, while in the internal case we fix that internal functors are lax, horizontal transformations are oplax, and vertical transformations are lax (as we used in \cite{Fem}). 

The starting point of our considerations is that $\Oo^*_\bullet(\A,\B)$ must be an object of $\Oo$. For the questions of closedness of 
$\Oo_*$ we will want $\Oo^\star_\bullet: (\Oo_*)^{op}\times\Oo_*\to\Oo_*$ to be a bifunctor. It will turn out that the ambient categories 
$\Oo_*$ for inner-homs are only $\Oo_{st}$ and $\Oo_{ps}$, but that the type $\star$ can differ from $*$.  

Take $\A,\B,\C\in\Oo$. For a fixed choice of parameters $*,\bullet$, in order to obtain the notion of a {\em quasi-functor}, we start by describing any 0-cell in $\Oo_*(\A,\Oo^*_\bullet(\B,\C))$. To make the reasoning simpler, 
let us fix the weakest version, namely that we consider lax functors among $\A$ and $\Oo^*_\bullet(\B,\C)$ and oplax(-lax) transformations. 
Being from $\Oo$, both $\A$ and $\B$ have $0-, 1-,.., m$-cells, and possibly additionally vertical $1-,.., m-1$-cells. Thus a lax functor 
\begin{equation} \eqlabel{F to be charact}
F:\A\to\Oo^{lx}_\bullet(\B,\C)
\end{equation} 
yields images $F(a)=(-,a)$ in $\Oo(\B,\C)$ where $a$ is any $k$-cell in $\A$, for $k=0,...,m$ (and possibly vertical cells). 
For a 0-cell $A\in\A$ one easily sees that $(-,A):\B\to\C$ is a lax functor. Similarly, for every 0-cell $B\in\B$ one gets that 
$(B,-):\A\to\C$ is a lax functor. Namely, the axioms of $F$ for being  a lax functor are given by equalities of transformations and modifications. Evaluating these at $B$ one gets axioms for $(B,-)$ to be a lax functor. 
The rest of the data that one has from $F$ being a lax functor are: 
\begin{itemize}
\item $(-,a_1)$ is a (horizontal) oplax transformation, 
\item ($(-,a_1')$ is a vertical lax transformation), 
\item compositor and unitor for $F$, 
\item $(-,a_2)$ is a modification of (horizontal) oplax transformations (and vertical lax transformations), 
\item ($(-,a_2')$ is a modification of horizontal lax transformations and vertical oplax transformations), 
\item $(-,a_3)$ is a perturbation,
\end{itemize}
for a $k$-cell $a_k$ in $\A$, $k=1,2,3$, and so on for higher dimensions. We will see examples for $m=2,3$ after the following definition.

\begin{defn} \delabel{gen-quasi}
Let $\Oo$ be any of the categories $\V_{m+1}, \I_{m+1}$ and let $\A,\B,\C\in\Oo$. 
A {\em general quasi-functor (of two variables)} $H: \A\times\B\to\C$ of type $(*,\bullet)$ in $\Oo$ consists of: 
\begin{enumerate}
\item two families 
$$(-,A):\B\to\C, \qquad (B,-):\A\to\C$$ 
of type $*$ of functors of dimension $m$ indexed by objects $A\in \A$ and $B\in\B$, 
\item a single 2-cell in the enriched case, {\em i.e.} four 2-cells in the internal case, 
we will refer to these 2-cells as ``structure cells of type $\bullet$'';
\item $k$-cells for $k=3,..,m$,
\end{enumerate}
which fulfill certain axioms so that the above data and axioms form an $m$-dimensional functor $\A\to\Oo^*_\bullet(\B,\C)$ of type $*$. 
\end{defn}

This generic definition is not precise, though we give now examples for double and Gray-categories. 
(The specific form of the axioms depends also on the type $\bullet$ of transformations and also of the higher transforms considered in the inner-hom.) 
For $\Oo=Dbl$ see \prref{char df} of Appendix B.1: apart from the two families of lax double functors,  
four families of 2-cells and 20 axioms characterize a lax double functor $F:\Aa\to \Lax_{hop}(\B,\C)$. For $\Oo=G\x\Cat$ and types 
$*=ps, \bullet=lx$ see Table 1 of \cite[Section 4]{BL} and Appendix D.1 thereof: apart from the two families of pseudo-Gray-functors,  
 one gets three families of 2-cells and four families of 3-cells satisfying 14 axioms and 12 degeneracy equations.
The four 3-cells in the pseudo functor case correspond to (appear as a weakening of) the four axioms \axiomref{$k'k,K$}, \axiomref{$k,K'K$}, \axiomref{$(k,K)-l-nat$} and \axiomref{$(k,K)-r-nat$} of \cite[Proposition 3.3]{Fem} (see Appendix B.1). The authors consider {\em normal, i.e. unital} functors, otherwise there would appear additional two 3-cells, which in the (op)lax case would come from the two axioms \axiomref{$1_B,K$} and \axiomref{$k,1_A$} from Appendix B.1. 
If one would consider internal pseudo unital case for $m=3$, additional six 3-cells would appear coming from the six axioms \axiomref{$u,K'K$}, \axiomref{$k'k,U$}, 
\axiomref{$u,\frac{U}{U'}$}, \axiomref{$\frac{u}{u'},U$}, \axiomref{$(u,U)-l-nat$} and \axiomref{$(u,U)-r-nat$}, while in the non-unital  case the additional two 3-cells would appear stemming from the two axioms \axiomref{$u,1_A$} and \axiomref{$1_B,U$}, all from Appendix B.1. (Here $k,k'$ are horizontal 1-cells in $\B$, $K,K'$ are horizontal 1-cells in $\A$, $u,u'$ are vertical 1-cells in $\B$ and $U,U'$ are vertical 1-cells in $\A$.)

\bigskip

By \deref{gen-quasi}, observe that if $\Oo^*_\bullet(\A,\B)$ is an object of $\Oo$, then 
there is a bijection 
\begin{equation} \eqlabel{quasi-set*}
q\x\Oo^*_\bullet(\A\times\B,\C)\iso\Oo_*(\A, \Oo^*_\bullet(\B,\C)),
\end{equation}
where on the left-hand side is the set of general quasi-functors of type $(*,\bullet)$. 
One can actually also introduce higher cells on quasi-functors so that the transformations are {\em possibly of another type $\#$} at $\A$ but of the same type $\bullet$ in $\B$, 
so that $q\x\Oo^*_\bullet\vert_{(\#,\bullet)}(\A\times\B,\C)$ is even an object in $\Oo$. Then one even has natural isomorphisms
\begin{equation} \eqlabel{quasi}
q\x\Oo^*_\bullet\vert_{(\#,\bullet)}(\A\times\B,\C)\iso\Oo^*_\#(\A, \Oo^*_\bullet(\B,\C))
\end{equation}
in $\Oo$ (to see this, the double categorical proof in Subsection 4.2 is insightful). 
We name now some examples studied in literature of the isomorphism \equref{quasi} and the description of the corresponding quasi-functors. 
As mentioned, for double categories \equref{quasi} is given in \cite[Proposition 3.3, Definition 3.4 and Section 4]{Fem}. 
For Gray-categories this isomorphism corresponds to \cite[Proposition 4.2]{BL}, where quasi-maps are taken for multimaps of a (skew) multicategory. For bicategories it can be found in \cite[Section 1.3]{Ver}.



\begin{rem} \rmlabel{tr in qu.f.}
Let $\F:\A\to\Oo^*_{\bullet}(\B,\C)$ denote a generic $*$-type functor. 
In particular, among the data and axioms in \deref{gen-quasi} we find those that come from $\F(x)=(-,x)$ being a $k$-cell in 
$\Oo^*_{\bullet}(\B,\C)$ for every $k$-cell $x$ in $\A$ and $k=1,2,..,n$. Moreover, as we recorded in Table 4 in the Appendix  
of \cite{Fem1} (see the last column of the present Appendix B.2) some of the obtained axioms in \deref{gen-quasi} have an additional meaning so that overall one gets that $(y,-)$ is a $k$-cell in $\Oo^*_{\#}(\A,\C)$ for any $k$-cell $y$ in $\B$. (For instance, 
compositor and unitor laws for $F$ together with a modification law for $(-,a_2)$ equivalently mean that $(b_1,-)$ is a 
(horizontal) oplax transformation. Also, one axiom of $(-,a_1)$ being a (horizontal) oplax transformation equivalently expresses an axiom 
for $(b_2,-)$ to be a modification.) One finds similar occurrences in higher dimensions. 
\end{rem}

For future use in this paper we recall \cite[Proposition 4.1]{Fem1}, it illustrates what we just said in the remark. 
It characterizes a lax double quasi-functor $H\colon \A\times\B\to\C$ that arises from a lax double functor $\A\to\Lax_{hop}(\B, \C)$. 
Here $\Lax_{hop}(\B, \C)$ is the double category consisting of lax double functors $\B\to\C$, horizontal oplax 
transformations \axiom{$h.o.t.$} as for 1h-cells, vertical lax transformations \axiom{$v.l.t.$} as for 1v-cells, and the corresponding modifications \axiom{$m.ho\x vl$} as for 2-cells. In the terminology of \deref{gen-quasi} we say that the above quasi-functor $H$ is of type $(lx,o-l)$.

\begin{prop} \prlabel{quasi-fun}
Let $\A,\B,\C$ be double categories. The following are equivalent:
\begin{enumerate} 
\item $H\colon \A\times\B\to\C$ is a double quasi-functor of type $(lx,o-l)$, meaning that there are two families of lax double functors 
$(-,A)\colon\B\to\C\quad\text{ and}\quad (B,-)\colon\A\to\C$ for objects $A\in\A, B\in\B$, 
such that $H(A,-)=(-, A), H(-, B)=(B,-)$ and $(-,A)\vert_B=(B,-)\vert_A=(B,A)$, and 
there are four families of 2-cells 
\begin{equation} \eqlabel{four cells}
\scalebox{0.78}{
\bfig
 \putmorphism(-150,50)(1,0)[(B,A)`(B', A)`(k, A)]{600}1a
 \putmorphism(450,50)(1,0)[\phantom{A\ot B}`(B', A') `(B', K)]{680}1a
\putmorphism(-180,50)(0,-1)[\phantom{Y_2}``=]{450}1r
\putmorphism(1100,50)(0,-1)[\phantom{Y_2}``=]{450}1r
\put(330,-190){\fbox{$(k,K)$}}
 \putmorphism(-150,-400)(1,0)[(B,A)`(B,A')`(B,K)]{600}1a
 \putmorphism(450,-400)(1,0)[\phantom{A\ot B}`(B', A') `(k, A')]{680}1a
\efig}
\end{equation}

$$
\scalebox{0.78}{
\bfig
\putmorphism(-150,50)(1,0)[(B,A)`(B,A')`(B,K)]{600}1a
\putmorphism(-150,-400)(1,0)[(\tilde B, A)`(\tilde B,A') `(\tilde B,K)]{640}1a
\putmorphism(-180,50)(0,-1)[\phantom{Y_2}``(u,A)]{450}1l
\putmorphism(450,50)(0,-1)[\phantom{Y_2}``(u,A')]{450}1r
\put(-20,-180){\fbox{$(u, K)$}}
\efig}
\quad
\scalebox{0.78}{
\bfig
\putmorphism(-150,50)(1,0)[(B,A)`(B',A)`(k,A)]{600}1a
\putmorphism(-150,-400)(1,0)[(B, \tilde A)`(B', \tilde A) `(k,\tilde A)]{640}1a
\putmorphism(-180,50)(0,-1)[\phantom{Y_2}``(B,U)]{450}1l
\putmorphism(450,50)(0,-1)[\phantom{Y_2}``(B',U)]{450}1r
\put(0,-180){\fbox{$(k,U)$}}
\efig}
$$

$$
\scalebox{0.78}{
\bfig
 \putmorphism(-150,500)(1,0)[(B,A)`(B,A) `=]{600}1a
\putmorphism(-180,500)(0,-1)[\phantom{Y_2}`(B, \tilde A) `(B,U)]{450}1l
\put(-20,50){\fbox{$(u,U)$}}
\putmorphism(-150,-400)(1,0)[(\tilde B, \tilde A)`(\tilde B, \tilde A) `=]{640}1a
\putmorphism(-180,50)(0,-1)[\phantom{Y_2}``(u,\tilde A)]{450}1l
\putmorphism(450,50)(0,-1)[\phantom{Y_2}``(\tilde B, U)]{450}1r
\putmorphism(450,500)(0,-1)[\phantom{Y_2}`(\tilde B, A) `(u,A)]{450}1r
\efig}
$$ 
in $\C$ determined by all 1h-cells $K\colon A\to A'$ and 1v-cells $U\colon A\to\tilde A$ in $\A$, and 1h-cells $k\colon B\to B'$ 
and 1v-cells $u\colon B\to\tilde B$ in $\B$, which satisfy 20 axioms from \cite[Proposition 3.3]{Fem}, and
\item there are two families of lax double functors 
$(-,A)\colon\B\to\C\quad\text{ and}\quad (B,-)\colon\A\to\C$ for objects $A\in\A, B\in\B$, 
such that $(-,A)\vert_B=(B,-)\vert_A=(B,A)$, and the following hold: 
\begin{enumerate}[(i)]
\item $(-,K)\colon (-,A)\to(-,A')$ is a horizontal oplax transformation for each 1h-cell $K\colon A\to A'$, 
$(-,U)\colon (-,A)\to(-,\tilde A)$ is a vertical lax transformation for each 1v-cell $U\colon A\to\tilde A$ in $\A$, and 
$(-,\zeta)$ is a modification with respect to horizontally oplax and vertically lax transformations for each 2-cell $\zeta$ in $\A$; 
\item $(k,-)\colon (B,-)\to (B', -)$ is a horizontal lax transformation for each 1h-cell $k\colon B\to B'$, 
$(u,-)\colon (B,-)\to(\tilde B,-)$ is a vertical oplax transformation for each 1v-cell $u\colon B\to\tilde B$ in $\B$, and 
$(\omega,-)$ is a modification with respect to horizontally lax and vertically oplax transformations for each 2-cell $\omega$ in $\B$;
\item for 1h-cells $K,k$ and 1v-cells $U,u$ the following 2-cell components of oplax resp. lax transformations coincide:
$(-,K)\vert_k=(k,-)\vert_K$, $(-,U)\vert_u=(u,-)\vert_U$, 
$(-,K)\vert_u=(u,-)\vert_K$ and $(-,U)\vert_k=(k,-)\vert_U$.
\end{enumerate}
\end{enumerate}
\end{prop}

\begin{rem} \rmlabel{dual quasi}
We make a few remarks here. 
\begin{enumerate} [a)]
\item Changing the type of $H$ in the above proposition from $(lx,o-l)$ into $(lx,l-o)$, one gets a lax double functor 
$\A\to\Lax_{hlx}(\B, \C)$, where the codomain double category differs from $\Lax_{hop}(\B, \C)$ in that its 1h-cells are horizontal 
{\em lax} transformations (h.l.t.), its 1v-cells are vertical {\em oplax} transformations (v.o.t.), and 2-cells are the corresponding modifications (m.hl-vo). In the generalized notation we write $\Lax_{hop}(\B, \C)$ as $\Oo^{lx}_{o-l}(\B, \C)$, and $\Lax_{hlx}(\B, \C)$ 
as $\Oo^{lx}_{l-o}(\B, \C)$ with $\Oo=Dbl$.  
\item Let $\#$-transformations be opposite to $\bullet$-transformations. Changing $\bullet$ in \equref{quasi-set*} to $\#$,
for general $\Oo$ we have bijections 
\begin{equation} \eqlabel{quasi-set-sharp}
q\x\Oo^*_{\#}(\A\times\B,\C)\iso\Oo_*(\A, \Oo^*_{\#}(\B,\C)),
\end{equation}
where $q\x\Oo^*_{\#}(\A\times\B,\C)$ denotes the set of quasi-functors of type $(*; \#)$. 
\item Inspecting \prref{quasi-fun} we see that the data of quasi-functors $\A\times\B\to\C$ of type $(lx,l-o)$ and $(lx,o-l)$ 
differ in that in the items (i) and (ii) of point 2. in the above proposition, transformations (h.o.t.) and (h.l.t.) change their roles, 
and likewise (v.l.t.) and (v.o.t.). 
\item Consider the axioms determining a lax double functor $\A\to\Lax_{hop}(\B, \C)$ in \prref{quasi-fun} that we included in 
the Appendix B.1, \prref{char df}. If we make the following two changes in these axioms: 1) change the roles of $\A$ and $\B$, 
and 2) change horizontal oplax transformations into horizontal lax transformations, and the same with vertical transformations, 
then the axioms remain precisely the same. Thus a quasi-functor $\A\times\B\to\C$ stemming from a lax double functor 
$\A\to\Lax_{hop}(\B, \C)$ is the same as a quasi-functor that stems from a lax double functor $\B\to\Lax_{hlx}(\A, \C)$. 
We thus may write $\Oo_{lx}(\A,\Oo^{lx}_{o-l}(\B, \C))\iso\Oo_{lx}(\B,\Oo^{lx}_{l-o}(\A, \C))$ for $\Oo=Dbl$. 
\end{enumerate}
\end{rem}

We find the same occurrence as in the part d) of the Remark in Gray-categories.  
To see that the same argument applied there is applicable also to 3-cells occurring in the structure of a 
3-dimensional quasi-functor, we will explain it in the example of the structure 3-cells  
$(k,-)^2_{K'K}$ and $(-,K)^2_{k'k}$ of a Gray-categorical quasi-functor, which were denoted by $F^{A_2,A_1}_B$ and $F^A_{B_2,B_1}$ 
in \cite[Section 4.3]{BL}, respectively. The 3-cell $(k,-)^2_{K'K}$ can be imagined as a 3-cell going from left to right in our double categorical axiom \axiomref{($k,K'K$)}, while $(-,K)^2_{k'k}$ as a 3-cell going from left to right in our axiom \axiomref{($k'k,K$)}. 
We already explained in the Remark that changing the roles of $\A$ and $\B$, and of horizontal oplax and horizontal lax transformations, 
the structure 2-cells $(k,K)$ appearing in the axioms remain unchanged (mind the difference that our transformations $(-,K)$ are oplax 
and $(k,-)$ lax, whereas to them corresponding transformations $F(A)$ of \cite{BL} are lax, and $F(B)$ are oplax, respectively). 
Observe that by doing the two announced changes: 
1) the domain 2-cell of $(k,-)^2_{K'K}$ becomes the codomain of the structure 3-cell $(-,K)^2_{k'k}$ (and similarly 
the codomain of the former becomes the domain of the latter), 2) the direction of $(k,-)^2_{K'K}$ was determined by $(k,-)$ being lax, 
and of $(-,K)^2_{k'k}$ by $(-,K)$ being oplax, and these roles, hence directions, are now flipped around.  
Thus by doing these two changes $(k,-)^2_{K'K}$ and $(-,K)^2_{k'k}$ change the roles. Similar flips of roles happen with the rest of 
the 3-cells, as the data in a quasi-functor structure are symmetric. 

\smallskip

This brings us to bijections 
\begin{equation} \eqlabel{quasi-right bij}
q\x\Oo^*_{\bullet}(\A\times\B,\C)\iso q\x\Oo^*_{\#}(\B\times\A,\C)\stackrel{\equref{quasi-set-sharp}}{\iso}\Oo_*(\B, \Oo^*_{\#}(\A,\C)),
\end{equation}
where now $\#$-transformations are opposite to $\bullet$-transformations, that lift to natural isomorphisms 
\begin{equation} \eqlabel{quasi-right}
q\x\Oo^*_{\bullet}\vert_{(\#,\bullet)}(\A\times\B,\C)\iso q\x\Oo^*_{\#}\vert_{(\bullet,\#)}(\B\times\A,\C)\iso\Oo^*_{\bullet}(\B, \Oo^*_{\#}(\A,\C)) 
\end{equation}
for $\Oo$ of dimension at least $m\leq 3$ (in \ssref{Gray-cat of qf} we show that in dimension 3 the left-hand side above is a Gray-category).  
Compare our \rmref{dual quasi} and the above two isomorphisms with 2-categorical results from \cite[Proposition I,4.13]{Gray}. 
%
We will denote by $ev:\Oo^*_{\bullet}(\B,\C)\times\B\to\C$ and $\crta{ev}:\B\times\Oo^*_{\#}(\B,\C)\to\C$ the quasi-functors so that $ev$ satisfies the usual universal property with respect to \equref{quasi}, and $\crta{ev}$ satisfies a similar universal property with respect to \equref{quasi-right}.

\subsection{Dimension matters for the possible types}

Crans has argued in the introduction of \cite{Cr} that there can not be a biclosed monoidal structure on the category of Gray-categories and strict Gray-functors so that the inner-homs have strict functors for objects and weak transformations for 1-cells. This argument was fine grained in \cite[Corollary 4.5]{BG}, where it was shown that under the assumption of biclosedness for $G\x\Cat$ (with either inner-hom having then necessarily strict functors for objects), there exist no proper natural transformations among strict Gray functors in the first place.
Although we will be concluding on biclosedness and monoidality of $G\x\Cat$ only in \ssref{I and II}, already at this point we can draw the following conclusion applying Crans's argument.

\begin{prop} \prlabel{Crans before}
For $m\geq 3$ if there exists a quasi-functor $\A\times\B\to\C$ of type $(st,\bullet)$, then either $\bullet=st$, or $\C$ is a proper 3-category (if $m=3$), {\em i.e.} has strict interchange among 2-cells (for $m>3$). 
\end{prop}

\begin{proof}
Joining \equref{quasi-right bij} and \equref{quasi-set*} we get 
$$\Oo_*(\B, \Oo^*_{\#}(\A,\C))\iso \Oo_*(\A, \Oo^*_{\bullet}(\B,\C)).$$
Take $H:\A\times\B\to\C$ from the left in \equref{quasi-right bij}, set $*=st$, and consider a transformation $\alpha:F\Rightarrow G$ in $\Oo^{st}_{\#}(\A,\C)$. Set $\B=2$, the interval 3-category (it has solely $F$ and $G$ for objects, and a single non-trivial 1-cell $\alpha$). As a 1-cell in $\Oo^{st}_{\#}(\A,\C)$, we can see $\alpha$ as a strict Gray-functor 
$2\to\Oo^{st}_{\#}(\A,\C)$, and by the above isomorphism further as a strict Gray-functor $\crta\alpha: \A\to 
\Oo^{st}_{\bullet}(2,\C)$. By the property of a strict functor, for 1-cells $f,g$ in $\A$ we have $\crta\alpha(gf)=
\crta\alpha(g)\crta\alpha(f)$, 
which is the following composition in $\C$ 
$$\scalebox{0.64}{
\bfig
\putmorphism(-530,210)(0,-1)[\phantom{Y_2}``=]{450}1r
\putmorphism(1030,210)(0,-1)[\phantom{Y_2}``=]{450}1r
 \putmorphism(450,210)(1,0)[F(C)`G(C) `\alpha(C)]{580}1a
  \putmorphism(-550,210)(1,0)[F(A)`\phantom{F(B)}`F(gf)]{1000}1a
 \putmorphism(-550,-265)(1,0)[F(A)`G(A)`\alpha(A)]{620}1a
 \putmorphism(80,-265)(1,0)[\phantom{F(B)}`G(C) `G(gf)]{970}1a
\put(100,-10){\fbox{$\alpha_{gf}$}}
\efig}= 
\scalebox{0.64}{
\bfig
 \putmorphism(450,360)(1,0)[F(B)`F(C) `F(g)]{680}1a
 \putmorphism(1120,360)(1,0)[\phantom{F(B)}`G(C) `\alpha(C)]{600}1a
\put(1000,130){\fbox{$\alpha_g$}}

  \putmorphism(-150,-90)(1,0)[F(A)` F(B) `F(f)]{600}1a
\putmorphism(450,-90)(1,0)[\phantom{F(A)}` G(B) `\alpha(B)]{680}1a
 \putmorphism(1120,-90)(1,0)[\phantom{F(A)}`G(C) ` G(g)]{620}1a

\putmorphism(450,360)(0,-1)[\phantom{Y_2}``=]{450}1l
\putmorphism(1710,360)(0,-1)[\phantom{Y_2}``=]{450}1r

 \putmorphism(-150,-540)(1,0)[F(A)`G(A)`\alpha(A)]{600}1a
 \putmorphism(450,-540)(1,0)[\phantom{F(B)}`G(B) `G(f)]{680}1a

\putmorphism(-180,-90)(0,-1)[\phantom{Y_2}``=]{450}1r
\putmorphism(1040,-90)(0,-1)[\phantom{Y_2}``=]{450}1r
\put(350,-320){\fbox{$\alpha_f$}}
\efig}=:\alpha_f*\alpha_g
$$ 
of the 2-cell components $\alpha_f$ of the transformation $\alpha$. 
Similarly, let us allow a third strict functor $H:\A\to\C$, and another transformation $\beta:G\Rightarrow H$ in 
$\Oo^{st}_{\#}(\A,\C)$, and let us 
consider $\alpha, \beta, \frac{\alpha}{\beta}$ as strict Gray-functors $3\to\Oo^{st}_{\#}(\A,\C)$, and the induced 
strict Gray-functors  $\crta\alpha, \crta\beta, \frac{\crta\alpha}{\crta\beta}:\A\to\Oo^{st}_{\bullet}(3,\C)$. As strict functors, both $\crta\beta$ and $\frac{\crta\alpha}{\crta\beta}$ should satisfy the above rule 
$\alpha_{gf}=\alpha_f*\alpha_g$, and the only natural way for the latter to be defined at $gf$ is $\frac{\crta\alpha}{\crta\beta}\vert_{gf}=\frac{\alpha_{gf}}{\beta_{gf}}$, which is the vertical composition of two instances of 2-cells as in the above diagram. So, it should at the same time be $\frac{\crta\alpha}{\crta\beta}\vert_{gf}=
\frac{\crta\alpha}{\crta\beta}\vert_f * \frac{\crta\alpha}{\crta\beta}\vert_g$, 
meaning that the two following expressions should be equal:
$$
\scalebox{0.64}{
\bfig
 \putmorphism(450,860)(1,0)[F(B)`F(C) `F(g)]{680}1a
 \putmorphism(1120,860)(1,0)[\phantom{F(B)}`G(C) `\alpha(C)]{600}1a
\put(1000,630){\fbox{$\alpha_g$}}

  \putmorphism(-150,430)(1,0)[F(A)` F(B) `F(f)]{600}1a
\putmorphism(450,430)(1,0)[\phantom{F(A)}` G(B) `\alpha(B)]{680}1a
 \putmorphism(1120,430)(1,0)[\phantom{F(A)}`G(C) ` G(g)]{620}1a

\putmorphism(450,860)(0,-1)[\phantom{Y_2}``=]{410}1l
\putmorphism(1710,860)(0,-1)[\phantom{Y_2}``=]{410}1r

 \putmorphism(-150,0)(1,0)[F(A)`G(A)`\alpha(A)]{600}1a
 \putmorphism(450,0)(1,0)[\phantom{F(B)}`G(B) `G(f)]{680}1a

\putmorphism(-180,430)(0,-1)[\phantom{Y_2}``=]{450}1r
\putmorphism(1040,430)(0,-1)[\phantom{Y_2}``=]{450}1r
\put(350,210){\fbox{$\alpha_f$}}

 \putmorphism(1120,0)(1,0)[\phantom{F(B)}`G(C) `G(g)]{600}1a
 \putmorphism(1720,0)(1,0)[\phantom{F(B)}`H(C) `\beta(C)]{600}1a
\put(1600,-250){\fbox{$\beta_g$}}

  \putmorphism(450,-470)(1,0)[G(A)` G(B) `G(f)]{600}1a
\putmorphism(1050,-470)(1,0)[\phantom{F(A)}` H(B) `\beta(B)]{680}1a
 \putmorphism(1720,-470)(1,0)[\phantom{F(A)}`H(C) ` H(g)]{620}1a

\putmorphism(1050,0)(0,-1)[\phantom{Y_2}``=]{450}1l
\putmorphism(2310,0)(0,-1)[\phantom{Y_2}``=]{450}1r

 \putmorphism(450,-900)(1,0)[G(A)`H(A)`\beta(A)]{600}1a
 \putmorphism(1050,-900)(1,0)[\phantom{F(B)}`H(B) `H(f)]{680}1a

\putmorphism(420,-470)(0,-1)[\phantom{Y_2}``=]{450}1r
\putmorphism(1640,-470)(0,-1)[\phantom{Y_2}``=]{450}1r
\put(950,-690){\fbox{$\beta_f$}}
\efig}
\quad=\quad  
\scalebox{0.64}{
\bfig
 \putmorphism(450,860)(1,0)[F(B)`F(C) `F(g)]{680}1a
 \putmorphism(1120,860)(1,0)[\phantom{F(B)}`G(C) `\alpha(C)]{600}1a
\put(1000,650){\fbox{$\alpha_g$}}

\putmorphism(450,430)(1,0)[F(B)` G(B) `\alpha(B)]{680}1a
 \putmorphism(1120,430)(1,0)[\phantom{F(A)}`G(C) ` G(g)]{620}1a

\putmorphism(450,860)(0,-1)[\phantom{Y_2}``=]{410}1l
\putmorphism(1710,860)(0,-1)[\phantom{Y_2}``=]{410}1r

  \putmorphism(-150,-20)(1,0)[F(A)` F(B) `F(f)]{600}1a
\putmorphism(-180,-20)(0,-1)[\phantom{Y_2}``=]{450}1r
\putmorphism(1040,-20)(0,-1)[\phantom{Y_2}``=]{450}1r
\put(350,-240){\fbox{$\alpha_f$}}
 \putmorphism(-150,-450)(1,0)[F(A)`G(A)`\alpha(A)]{600}1a
\putmorphism(450,-20)(1,0)[\phantom{F(B)}` G(B) `\alpha(B)]{680}1a
\putmorphism(1100,-20)(1,0)[\phantom{F(A)}` H(B) `\beta(B)]{620}1a

 \putmorphism(1720,450)(1,0)[\phantom{F(B)}`H(C) `\beta(C)]{600}1a

\put(1600,200){\fbox{$\beta_g$}}
 \putmorphism(1720,-20)(1,0)[\phantom{F(A)}`H(C) ` H(g)]{620}1a

  \putmorphism(450,-470)(1,0)[\phantom{F(A)}` G(B) `G(f)]{600}1a
\putmorphism(1050,-470)(1,0)[\phantom{F(A)}` H(B) `\beta(B)]{680}1a

\putmorphism(1050,450)(0,-1)[\phantom{Y_2}``=]{450}1l
\putmorphism(2310,450)(0,-1)[\phantom{Y_2}``=]{450}1r

 \putmorphism(450,-900)(1,0)[G(A)`H(A)`\beta(A)]{600}1a
 \putmorphism(1050,-900)(1,0)[\phantom{F(B)}`H(B). `H(f)]{680}1a

\putmorphism(420,-470)(0,-1)[\phantom{Y_2}``=]{450}1r
\putmorphism(1640,-470)(0,-1)[\phantom{Y_2}``=]{450}1r
\put(950,-690){\fbox{$\beta_f$}}
\efig}
$$ 
Such an equality is only possible in 3-categories, given that in Gray-categories a non-strict interchange should act between the two different compositions of the 2-cells $\alpha_f$ and $\beta_g$ in the middle of the above two diagrams. Hence $\C$ must be a proper 3-category, or all the 2-cells $\alpha_f, \beta_f$ are trivial, meaning that $\#=\bullet=st$. 
\qed\end{proof}

\subsection{Candidate for general Gray tensor product - meta-product}

For $\A,\B,\C\in\Oo$, any choice of a quasi-functor $H$ of type $(*,\bullet)$ in $\Oo$ carries a choice of: 
1) two families $(-,A):\B\to\C, \qquad (B,-):\A\to\C$ for $A\in \A$ and $B\in\B$ of functors of type $*$ in $\Oo$, 
2) 2- and possibly higher cells in $\C$ (depending on the dimension $m$ in $\Oo$), and 3) axioms in $\C$. It thus also fixes 
a meaning for cells $H(a,b)=(b,a)$ in $\C$, where $b$ from $\B$ and $a$ from $\A$ are cells of suitable degrees, which appear by the 
definition of quasi-functor (as a characterization of a functor $F:\A\to\Oo^*_\bullet(\B,\C)$ of type $*$). In other words, 
$F(a)(b)=H(a,b)=(b,a)$ in $\C$ for certain cells $a$ in $\A$ and $b$ in $\B$. 
Observe that the obtained axioms are relations holding between such cells $(b,a)$ in $\C$. 

If $\Oo^*_\bullet(-, -)$ were a functor $\Oo^{op}_*\times\Oo_*\to\Oo_*$ and if there were an adjunction $(-\times\B, \Oo^*_\bullet(\B, -))$ 
for $\B\in\Oo$, 
from all the possible choices for $F:\A\to\Oo^*_\bullet(\B,\A\times\B)$ if $F$ were the unit $\eta$ of this adjunction, then we would be 
able to define a candidate $\ot$ for a monoidal product that would equip $\Oo$ with a structure of a closed monoidal category. 
This goes by 
defining $\A\ot\B$ as a quotient of the Cartesian product $\A\times\B$, with generators and relations that we read off from the quasi-functor 
$H_\eta$ determined by $\eta$. Namely, by taking for generators $a\ot b:=H_\eta(a,b)$ for suitable cells $a$ in $\A$ and $b$ in $\B$, 
and by rewriting relations holding among cells $H_\eta(a,b)$ in $\A\times\B$ obtained from the axioms of a quasi-functor in terms of relations among $a\ot b$. Then one should finally also check if there exist associativity and unity constraints for $\ot$ obeying the 
pentagonal and triangular equation. 

However, as done in \cite{Gray} and \cite{Fem}, we can also define $\ot$ as a ``candidate for Gray tensor product'' 
without knowing if the above functoriality and adjunction conditions are fulfilled, and then we can also ask the question whether $\ot$ gives 
a (non-Cartesian) monoidal product on $\Oo$. In order to make the name more compact, for a ``candidate for a monoidal product''  
we will say {\em meta-product}. To formulate a general form of a meta-product, for every cell $H(a,b)\in\C=\A\times\B$ for suitable cells 
$a\in\A, b\in\B$, and every axiom among them that appear in \deref{gen-quasi}, one writes $a\ot b$ and the resulting axioms among them. Then the meta-product $\A\ot\B$ in $\Oo$ is generated by objects $A\ot B$ for objects $A\in\A, B\in\B$ and $k$-cells $A\ot b, a\ot B$ 
for all $k$-cells $a\in\A, b\in\B$ and $k=1,...,m$, and all the other cells originating from points 2. and 3. in \deref{gen-quasi}, 
and the above described relations. 

We will give the definition for $m=2$ and $\I_1=\Cat(\Cat)=Dbl$ and the type $(lx, o-l)$ in the next definition.  
A meta-product for 
higher dimensional cases is constructed similarly, but the complexity rises. 
%
For horizontal and vertical 1-cells in a double category we will say shorty 1h- and 1v-cells. 
Horizontal composition we will denote by concatenation, and the vertical one we will write as fractions. 


\begin{defn} \delabel{tensor}
Let $\A,\B\in Dbl$ and let $\A\ot\B$ be generated as an object in $Dbl$ by the following data: 
\ul{objects}: $A\ot B$ for objects $A\in\A, B\in\B$; \\ 
\ul{1h-cells}: $A\ot k, K\ot B$, where $k$ is a 1h-cell in $\B$ and $K$ a 1h-cell in $\A$; \\ 
%
\ul{1v-cells}: $A\ot u, U\ot B$ and vertical compositions of such obeying the following rules: 
$$\frac{A\ot u}{A\ot u'}=A\ot \frac{u}{u'}, \quad \frac{U\ot B}{U'\ot B}=\frac{U}{U'}\ot B, \quad A\ot 1^B=1^{A\ot B}=1^A\ot B$$
where $u,u'$ are 1v-cells of $\B$ and $U,U'$ 1v-cells of $\A$; \\
\ul{2h-cells}: $A\ot\beta, \alpha\ot B$ for 2h-cells $\beta$ in $\B$ and $\alpha$ in $\A$; \\
$$
\scalebox{0.86}{
\bfig
\putmorphism(-150,50)(1,0)[A\ot B`A\ot B'`A\ot k]{600}1a
\putmorphism(-150,-400)(1,0)[A\ot\tilde B`A'\ot\tilde B' `A\ot l]{640}1a
\putmorphism(-180,50)(0,-1)[\phantom{Y_2}``A\ot u]{450}1l
\putmorphism(450,50)(0,-1)[\phantom{Y_2}``A\ot v]{450}1r
\put(0,-180){\fbox{$A\ot \omega$}}
\efig}
\qquad
\scalebox{0.86}{
\bfig
\putmorphism(-150,50)(1,0)[A\ot B`A'\ot B`K\ot B]{600}1a
\putmorphism(-150,-400)(1,0)[\tilde A\ot B`\tilde A'\ot B `L\ot B]{640}1a
\putmorphism(-180,50)(0,-1)[\phantom{Y_2}``U\ot B]{450}1l
\putmorphism(450,50)(0,-1)[\phantom{Y_2}``V\ot B]{450}1r
\put(0,-180){\fbox{$\zeta\ot B$}}
\efig}
$$
where $\omega$ and $\zeta$ are as in \equref{omega-zeta}, \\
four 2h-cells from the laxity of functors $(-,A)$ and $(B,-)$: ``structure cells of type $*$''
\begin{equation} \eqlabel{4 laxity 2-cells}
(A\ot k')(A\ot k)\stackrel{(A\ot -)_{k'k}}{\Rightarrow} A\ot (k' k), \quad 
(K'\ot B)(K\ot B)\stackrel{(-\ot B)_{K'K}}{\Rightarrow}(K' K)\ot B
\end{equation}
$$1_{A\ot B}\stackrel{(A\ot -)_B}{\Rightarrow} A\ot 1_B, 
\quad 1_{A\ot B}\stackrel{(-\ot B)_A}{\Rightarrow} 1_A\ot B$$
which satisfy associativity and unitality laws, and 
where $k,k'$ are 1h-cells of $\B$ and $K, K'$ 1h-cells of $\A$, \\
and four types of 2-cells coming from the 2-cells of point 2 in \prref{char df}: ``structure cells of type $\bullet$''
a vertically globular 2-cell $K\ot k\colon (A'\ot k)(K\ot B) \Rightarrow(K\ot B')(A\ot k)$, 
a horizontally globular 2-cell $U\ot u\colon \frac{U\ot B}{\tilde A\ot u}\Rightarrow\frac{A\ot u}{U\ot\tilde B}$ 
(so that $1^A\ot 1^B=1_{A\ot B}$), and 
2-cells $K\ot u$ and $U\ot k$, 
subject to the rules induced by the rules of point 2 in \cite[Proposition 3.3]{Fem} \\
and the following ones: ``axioms for horizontal composition of type $*$ of 2-cells'', 
``axioms for vertical composition of 2-cells'' and ``unity axioms''. 
$$
\scalebox{0.86}{
\bfig
\putmorphism(-150,500)(1,0)[A\ot B`A\ot B' `A\ot k]{600}1a
 \putmorphism(450,500)(1,0)[\phantom{F(A)}`A\ot B'' `A\ot k' ]{620}1a

 \putmorphism(-150,50)(1,0)[A\ot \tilde B`A\ot \tilde B' ` A\ot l]{600}1a
 \putmorphism(470,50)(1,0)[\phantom{F(A)}`A\ot \tilde B'' `A\ot l']{620}1a

\putmorphism(-180,500)(0,-1)[\phantom{Y_2}``A\ot u]{450}1l
\putmorphism(450,500)(0,-1)[\phantom{Y_2}``]{450}1r
\putmorphism(300,500)(0,-1)[\phantom{Y_2}``A\ot v]{450}0r
\putmorphism(1080,500)(0,-1)[\phantom{Y_2}``A\ot w]{450}1r
\put(-40,280){\fbox{$A\ot\omega$}}
\put(620,280){\fbox{$A\ot\omega'$}}

\putmorphism(-150,-400)(1,0)[A\ot \tilde B`A\ot \tilde B'' `A\ot l'l]{1200}1a

\putmorphism(-180,50)(0,-1)[\phantom{Y_2}``=]{450}1l
\putmorphism(1080,50)(0,-1)[\phantom{Y_3}``=]{450}1r
\put(240,-160){\fbox{$(A\ot -)_{l'l}$}}

\efig}
=
\scalebox{0.86}{
\bfig
\putmorphism(-150,500)(1,0)[A\ot B`A\ot B' `A\ot k]{600}1a
 \putmorphism(450,500)(1,0)[\phantom{F(A)}`A\ot B'' `A\ot k' ]{620}1a
 \putmorphism(-150,50)(1,0)[A\ot B`A\ot B''`A\ot k'k]{1220}1a

\putmorphism(-180,500)(0,-1)[\phantom{Y_2}``=]{450}1r
\putmorphism(1080,500)(0,-1)[\phantom{Y_2}``=]{450}1r
\put(240,290){\fbox{$(A\ot -)_{k'k}$}}

\putmorphism(-150,-400)(1,0)[A\ot \tilde B`A\ot \tilde B'' `A\ot l'l]{1200}1a

\putmorphism(-180,50)(0,-1)[\phantom{Y_2}``A\ot u]{450}1l
\putmorphism(1080,50)(0,-1)[\phantom{Y_3}``A\ot w]{450}1r
\put(300,-180){\fbox{$A\ot\omega'\omega$}} 
\efig}
$$

$$
\scalebox{0.86}{
\bfig
\putmorphism(-150,500)(1,0)[A\ot B`A'\ot B`K\ot B]{600}1a
 \putmorphism(470,500)(1,0)[\phantom{F(A)}`A''\ot B `K'\ot B]{600}1a
 \putmorphism(-150,50)(1,0)[A\ot B`A''\ot B`K'K\ot B]{1220}1a

\putmorphism(-180,500)(0,-1)[\phantom{Y_2}``=]{450}1r
\putmorphism(1080,500)(0,-1)[\phantom{Y_2}``=]{450}1r
\put(240,290){\fbox{$(-\ot B)_{K'K}$}}

\putmorphism(-150,-400)(1,0)[\tilde A\ot B`\tilde A''\ot B `L'L\ot B]{1200}1a

\putmorphism(-180,50)(0,-1)[\phantom{Y_2}``U\ot B]{450}1l
\putmorphism(1080,50)(0,-1)[\phantom{Y_3}``U''\ot B]{450}1r
\put(300,-180){\fbox{$\zeta'\zeta\ot B$}} 
\efig}
=
\scalebox{0.86}{
\bfig
\putmorphism(-150,500)(1,0)[A\ot B`A'\ot B`K\ot B]{600}1a
 \putmorphism(470,500)(1,0)[\phantom{F(A)}`A''\ot B `K'\ot B]{600}1a

 \putmorphism(-150,50)(1,0)[\tilde A\ot B`\tilde A'\ot B`L\ot B]{600}1a
 \putmorphism(470,50)(1,0)[\phantom{F(A)}`\tilde A''\ot B `L' \ot B]{620}1a

\putmorphism(-180,500)(0,-1)[\phantom{Y_2}``U\ot B]{450}1l
\putmorphism(450,500)(0,-1)[\phantom{Y_2}``]{450}1r
\putmorphism(300,500)(0,-1)[\phantom{Y_2}``U'\ot B]{450}0r
\putmorphism(1080,500)(0,-1)[\phantom{Y_2}``U''\ot B]{450}1r
\put(-40,280){\fbox{$\zeta\ot B$}}
\put(620,280){\fbox{$\zeta' \ot B$}}

\putmorphism(-150,-400)(1,0)[\tilde A\ot B`\tilde A''\ot B `L'L\ot B]{1200}1a

\putmorphism(-180,50)(0,-1)[\phantom{Y_2}``=]{450}1l
\putmorphism(1080,50)(0,-1)[\phantom{Y_3}``=]{450}1r
\put(240,-160){\fbox{$(-\ot B)_{L'L}$}}

\efig}
$$

$$A\ot\frac{\omega}{\omega'}=\frac{A\ot\omega}{A\ot\omega'}, \quad \frac{\zeta}{\zeta'}\ot B=\frac{\zeta\ot B}{\zeta'\ot B},$$
$$A\ot\Id_k=\Id_{A\ot k}, \quad \Id_K\ot B=\Id_{K\ot B}, \quad A\ot\Id^u=\Id^{A\ot u}, \quad \Id^U\ot B=\Id^{U\ot B} .$$

The source and target functors $s,t$ on $\A\ot\B$ are defined as in $\A\times\B$, the composition functor $c$ is defined by 
horizontal juxtaposition of 
the corresponding 2-cells, and the unit functor $i$ is defined on generators as follows: 
$$i(A\ot B)=1_{A\ot B}, \, i(A\ot v)=1^A\ot v (=Id^{A\ot v}) \quad \text{and} \quad i(U\ot B)=U\ot 1^B (=Id^{U\ot B}).$$ 
\end{defn}

\bigskip

Since $\A\ot^*_\bullet\B$ is defined by generators and relations on $\A\times\B$, 
it is clear that there is a double quasi-functor 
\begin{equation} \eqlabel{J}
J^*_\bullet\colon\A\times\B\to\A\ot^*_\bullet\B
\end{equation} 
of type $*$ given by $J^*_\bullet(-, B)(a)=a\ot B, J^*_\bullet(A,-)(b)=A\ot b$ for cells $a$ in $\A$ and $b$ in $\B$ and with unique 
structure cells of types $*$ and $\bullet$ satisfying axioms for compositions and unities of  higher cells. 
An analogous quasi-functor $J^*_\bullet$ one finds for any $\Oo$.

\subsection{Higher cells for quasi-functors of two variables in dimension 2} \sslabel{tr for qf dim2}

In this subsection we first recall definitions of transformations and modifications for quasi-functors for double categories from 
\cite[Section 4.1]{Fem}. The definitions of transformations we write in the form for lax double quasi-functors. The version for 
oplax double quasi-functors is obtained by reformulating the axioms in the obvious way so that the globular structure 2-cells of 
type $\bullet$ ($(k,K)$ and $(u,U)$) are written in the opposite direction. 
The corresponding definitions for 2-categories are obtained by making all the 1v-cells identities. 

\begin{defn} \delabel{oplax tr cubical}
A (horizontal) oplax transformation $\theta\colon (-,-)_1\Rightarrow (-,-)_2$ between lax double quasi-functors $(-,-)_1,(-,-)_2\colon 
\A\times\B\to\C$ is given by: for each $A\in\A$ a (horizontal) oplax transformation $\theta^A\colon (-,A)_1\Rightarrow(-,A)_2$ and 
for each $B\in\B$ a (horizontal) oplax transformation $\theta^B\colon (B,-)_1\Rightarrow(B,-)_2$, both of lax double functors, such that 
$\theta^A_B=\theta^B_A$ and such that 

\medskip

\noindent \axiom{$HOT^q_1$} \vspace{-0,6cm}
$$\scalebox{0.8}{
\bfig
 \putmorphism(450,700)(1,0)[(B', A)_1 `(B', A')_1 `(B', K)_1]{680}1a
 \putmorphism(1140,700)(1,0)[\phantom{A\ot B}`(B', A')_2 ` \theta^{B'}_{A'}]{680}1a

 \putmorphism(-150,250)(1,0)[(B, A)_1 `(B', A)_1`(k, A)_1]{600}1a
 \putmorphism(450,250)(1,0)[\phantom{A\ot B}`(B', A)_2 `\theta^{B'}_A]{680}1a
 \putmorphism(1130,250)(1,0)[\phantom{A\ot B}`(B', A')_2 ` (B', K)_2]{680}1a

\putmorphism(450,700)(0,-1)[\phantom{Y_2}``=]{450}1r
\putmorphism(1750,700)(0,-1)[\phantom{Y_2}``=]{450}1r
\put(1020,470){\fbox{$ \theta^{B'}_K$}}

 \putmorphism(-150,-200)(1,0)[(B, A)_1`(B, A)_2 `\theta^A_B]{640}1a
 \putmorphism(460,-200)(1,0)[\phantom{A'\ot B'}`(B', A)_2 `(k, A)_2]{680}1a

\putmorphism(-180,250)(0,-1)[\phantom{Y_2}``=]{450}1l
\putmorphism(1120,250)(0,-1)[\phantom{Y_3}``=]{450}1r
\put(310,0){\fbox{$ \theta^A_k$}}

 \putmorphism(1170,-200)(1,0)[\phantom{A\ot B}`(B', A')_2 ` (B', K)_2]{650}1a
\putmorphism(450,-200)(0,-1)[\phantom{Y_2}``=]{450}1r
\putmorphism(1750,-200)(0,-1)[\phantom{Y_2}``=]{450}1r

 \putmorphism(460,-650)(1,0)[(B, A)_2 `(B, A')_2 `(B, K)_2]{650}1a
 \putmorphism(1260,-650)(1,0)[ `(B', A')_2 `(k, A')_2]{600}1a
\put(920,-440){\fbox{$ (k,K)_2$}}
\efig}
=
\scalebox{0.8}{
\bfig

 \putmorphism(-150,700)(1,0)[(B,A)_1`(B',A)_1`(k,A)_1]{600}1a
 \putmorphism(450,700)(1,0)[\phantom{A\ot B}`(B', A')_1 `(B',K)_1]{680}1a

 \putmorphism(-150,250)(1,0)[(B,A)_1`(B,A')_1`(B,K)_1]{600}1a
 \putmorphism(450,250)(1,0)[\phantom{A\ot B}`(B', A')_1 `(k,A')_1]{680}1a
 \putmorphism(1120,250)(1,0)[\phantom{A'\ot B'}`(B', A')_2 `\theta^{A'}_{B'}]{650}1a

\putmorphism(-180,700)(0,-1)[\phantom{Y_2}``=]{450}1r
\putmorphism(1100,700)(0,-1)[\phantom{Y_2}``=]{450}1r
\put(350,450){\fbox{$(k,K)_1$}}
\put(1000,0){\fbox{$\theta^{A'}_k$}}

 \putmorphism(-150,-200)(1,0)[(B,A)_1`(B,A')_1`(B,K)_1]{600}1a
 \putmorphism(450,-200)(1,0)[\phantom{A''\ot B'}` (B, A')_2 `\theta^{A'}_B]{680}1a
 \putmorphism(1100,-200)(1,0)[\phantom{A''\ot B'}`(B', A')_2 ` (k, A')_2]{660}1a

\putmorphism(450,250)(0,-1)[\phantom{Y_2}``=]{450}1l
\putmorphism(1750,250)(0,-1)[\phantom{Y_2}``=]{450}1r
\putmorphism(-180,-200)(0,-1)[\phantom{Y_2}``=]{450}1r
\putmorphism(1100,-200)(0,-1)[\phantom{Y_2}``=]{450}1r
 \putmorphism(-150,-650)(1,0)[(B,A)_1`(B,A)_2`\theta^B_A]{600}1a
 \putmorphism(580,-650)(1,0)[`(B,A')_2`(B,K)_2]{540}1a
\put(320,-420){\fbox{$\theta^{B}_K$}}
\efig} \vspace{-0,2cm}
$$
for every 1h-cells $K\colon A\to A'$ and $k\colon B\to B'$, 

\medskip

\noindent \axiom{$HOT^q_2$} \vspace{-0,3cm}
$$\scalebox{0.86}{
\bfig
\putmorphism(-150,500)(1,0)[(B,A)_1`(B,A')_1`(B,K)_1]{600}1a
 \putmorphism(450,500)(1,0)[\phantom{F(A)}`(B,A')_2 `\theta^{A'}_B]{640}1a

 \putmorphism(-150,50)(1,0)[(\tilde B,A)_1`(\tilde B,A')_1`(\tilde B,K)_1]{600}1a
 \putmorphism(450,50)(1,0)[\phantom{F(A)}`(\tilde B,A')_2 `\theta^{A'}_{\tilde B}]{640}1a

\putmorphism(-180,500)(0,-1)[\phantom{Y_2}``(u,A)_1]{450}1l
\putmorphism(450,500)(0,-1)[\phantom{Y_2}``]{450}1r
\putmorphism(300,500)(0,-1)[\phantom{Y_2}``(u,A')_1]{450}0r
\putmorphism(1100,500)(0,-1)[\phantom{Y_2}``(u,A')_2]{450}1r
\put(-40,280){\fbox{$(u,K)_1$}}
\put(700,280){\fbox{$\theta^{A'}_u$}}

\putmorphism(-150,-400)(1,0)[(\tilde B,A)_1`(\tilde B,A)_2 `\theta^{A}_{\tilde B}]{640}1a
 \putmorphism(490,-400)(1,0)[\phantom{F(A')}` (\tilde B,A')_2 `(\tilde B,K)_2]{640}1a

\putmorphism(-180,50)(0,-1)[\phantom{Y_2}``=]{450}1l
\putmorphism(1120,50)(0,-1)[\phantom{Y_3}``=]{450}1r
\put(320,-200){\fbox{$\theta^{\tilde B}_K$}}
\efig}
\quad=\quad
\scalebox{0.86}{
\bfig
\putmorphism(-150,500)(1,0)[(B,A)_1`(B,A')_1`(B,K)_1]{600}1a
 \putmorphism(450,500)(1,0)[\phantom{F(A)}`(B,A')_2 `\theta^{A'}_B]{640}1a
 \putmorphism(-150,50)(1,0)[(B,A)_1`(B,A)_2`\theta^{A}_B]{600}1a
 \putmorphism(470,50)(1,0)[\phantom{F(A)}`(B,A')_2 `(B,K)_2]{660}1a

\putmorphism(-180,500)(0,-1)[\phantom{Y_2}``=]{450}1r
\putmorphism(1100,500)(0,-1)[\phantom{Y_2}``=]{450}1r
\put(320,280){\fbox{$\theta^B_K$}}

\putmorphism(-150,-400)(1,0)[(\tilde B,A)_1`(\tilde B,A)_2 `\theta^{A}_{\tilde B}]{640}1a
 \putmorphism(490,-400)(1,0)[\phantom{F(A')}` (\tilde B,A')_2 `(\tilde B,K)_2]{640}1a

\putmorphism(-180,50)(0,-1)[\phantom{Y_2}``(u,A)_1]{450}1l
\putmorphism(450,50)(0,-1)[\phantom{Y_2}``]{450}1l
\putmorphism(610,50)(0,-1)[\phantom{Y_2}``(u,A)_2]{450}0l 
\putmorphism(1120,50)(0,-1)[\phantom{Y_3}``(u,A')_2]{450}1r
\put(-40,-180){\fbox{$\theta^A_u$}} 
\put(620,-180){\fbox{$(u,K)_2$}}

\efig}
$$
for every 1h-cell $K\colon A\to A'$ and 1v-cell $u\colon B\to \tilde B$, 

\medskip

\noindent \axiom{$HOT^q_3$} \vspace{-0,3cm}
$$\scalebox{0.86}{
\bfig
\putmorphism(-150,500)(1,0)[(B,A)_1`(B',A)_1`(k,A)_1]{600}1a
 \putmorphism(480,500)(1,0)[\phantom{F(A)}`(B',A)_2 `\theta^A_{B'}]{640}1a

 \putmorphism(-150,50)(1,0)[(B,\tilde A)_1`(B',\tilde A)_1`(k,\tilde A)_1]{600}1a
 \putmorphism(480,50)(1,0)[\phantom{F(A)}`(B',\tilde A)_2 `\theta^{\tilde A}_{B'}]{640}1a

\putmorphism(-180,500)(0,-1)[\phantom{Y_2}``(B,U)_1]{450}1l
\putmorphism(450,500)(0,-1)[\phantom{Y_2}``]{450}1r
\putmorphism(300,500)(0,-1)[\phantom{Y_2}``(B',U)_1]{450}0r
\putmorphism(1100,500)(0,-1)[\phantom{Y_2}``(B',U)_2]{450}1r
\put(-40,280){\fbox{$(k,U)_1$}}
\put(700,280){\fbox{$\theta_U^{B'}$}}

\putmorphism(-150,-400)(1,0)[(B,\tilde A)_1`(B,\tilde A)_2 `\theta^{\tilde A}_B]{640}1a
 \putmorphism(480,-400)(1,0)[\phantom{A'\ot B'}` (B',\tilde A)_2 `(k,\tilde A)_2]{680}1a

\putmorphism(-180,50)(0,-1)[\phantom{Y_2}``=]{450}1l
\putmorphism(1120,50)(0,-1)[\phantom{Y_3}``=]{450}1r
\put(320,-200){\fbox{$\theta^{\tilde A}_k$}}
\efig}
\quad=\quad
\scalebox{0.86}{
\bfig
\putmorphism(-150,500)(1,0)[(B,A)_1`(B',A)_1`(k,A)_1]{600}1a
 \putmorphism(480,500)(1,0)[\phantom{F(A)}`(B',A)_2 `\theta^{A}_{B'}]{640}1a
 \putmorphism(-150,50)(1,0)[(B,A)_1`(B,A)_2`\theta^{A}_B]{600}1a
 \putmorphism(480,50)(1,0)[\phantom{F(A)}`(B',A)_2 `(k,A)_2]{680}1a

\putmorphism(-180,500)(0,-1)[\phantom{Y_2}``=]{450}1r
\putmorphism(1100,500)(0,-1)[\phantom{Y_2}``=]{450}1r
\put(320,280){\fbox{$\theta^A_k$}}

\putmorphism(-150,-400)(1,0)[(B,\tilde A)_1`(B,\tilde A)_2 `\theta^{\tilde A}_{B}]{640}1a
 \putmorphism(490,-400)(1,0)[\phantom{F(A')}` (B',\tilde A)_2 `(k,\tilde A)_2]{640}1a

\putmorphism(-180,50)(0,-1)[\phantom{Y_2}``(B,U)_1]{450}1l
\putmorphism(480,50)(0,-1)[\phantom{Y_2}``]{450}1l
\putmorphism(570,50)(0,-1)[\phantom{Y_2}``(B,U)_2]{450}0l 
\putmorphism(1120,50)(0,-1)[\phantom{Y_3}``(B',U)_2]{450}1r
\put(-40,-180){\fbox{$\theta_U^B$}} 
\put(620,-180){\fbox{$(k,U)_2$}}
\efig}
$$
for every 1v-cell $U\colon A\to \tilde A$ and 1h-cell $k\colon B\to B'$, and 

\medskip
\noindent \axiom{$HOT^q_4$} \vspace{-0,3cm}
$$\scalebox{0.86}{
\bfig
 \putmorphism(-170,500)(1,0)[(B,A)_1`(B,A)_1 `=]{600}1a
 \putmorphism(510,500)(1,0)[\phantom{Y_2}` `\theta^{A}_B]{420}1a
\putmorphism(-200,500)(0,-1)[\phantom{Y_2}` `(B,U)_1]{450}1l
\putmorphism(-240,500)(0,-1)[\phantom{Y_2}`(B, \tilde A)_1 `]{450}0l
\put(-70,50){\fbox{$(u,U)_1$}}
\putmorphism(-190,-400)(1,0)[(\tilde B, \tilde A)_1` `=]{500}1a
\putmorphism(-200,50)(0,-1)[\phantom{Y_2}``(u,\tilde A)_1]{450}1l
\putmorphism(450,50)(0,-1)[\phantom{Y_2}`(\tilde B, \tilde A)_1`(\tilde B, U)_1]{450}1l
\putmorphism(450,500)(0,-1)[\phantom{Y_2}`(\tilde B, A)_1 `(u,A)_1]{450}1l
\put(600,260){\fbox{$\theta^A_u$}}
\putmorphism(430,50)(1,0)[\phantom{(B, \tilde A)}``\theta^{A}_{\tilde B}]{480}1a
\putmorphism(1070,50)(0,-1)[\phantom{(B, A')}`(\tilde B, \tilde A)_2`(\tilde B,U)_2]{450}1r
\putmorphism(1070,500)(0,-1)[(B, A)_2`(\tilde B, A)_2`(u,A)_2]{450}1r
\putmorphism(450,-400)(1,0)[\phantom{(B, \tilde A)}``\theta^{\tilde A}_{\tilde B}]{480}1a
\put(600,-170){\fbox{$\theta_U^{\tilde B}$}}
\efig}
\quad=\quad
\scalebox{0.86}{
\bfig
 \putmorphism(-150,500)(1,0)[(B,A)_1`(B,A)_2 `\theta^{A}_B]{600}1a
 \putmorphism(450,500)(1,0)[\phantom{(B,A)}` `=]{540}1a
\putmorphism(-180,500)(0,-1)[\phantom{Y_2}`(B, \tilde A)_1 `(B,U)_1]{450}1l
\put(0,280){\fbox{$\theta_U^B$}}
\putmorphism(-180,-400)(1,0)[(\tilde B, \tilde A)_1` `\theta^{\tilde A}_{\tilde B}]{500}1a
\putmorphism(-180,50)(0,-1)[\phantom{Y_2}``(u,\tilde A)_1]{450}1l
\putmorphism(450,50)(0,-1)[\phantom{Y_2}`(\tilde B, \tilde A)_2`(u,\tilde A)_2]{450}1r
\putmorphism(450,500)(0,-1)[\phantom{Y_2}`(B, \tilde A)_2 `(B,U)_2]{450}1r
\put(620,50){\fbox{$(u,U)_2$}}
\putmorphism(-180,50)(1,0)[\phantom{(B, \tilde A)}``\theta^{\tilde A}_B]{500}1a
\putmorphism(1130,50)(0,-1)[\phantom{(B, A')}`(\tilde B, \tilde A)_2`(\tilde B,U)_2]{450}1r
\putmorphism(1130,500)(0,-1)[(B, A)``(u,A)_2]{450}1r
\putmorphism(1170,500)(0,-1)[`(\tilde B, A)_2`]{450}0r
\putmorphism(450,-400)(1,0)[\phantom{(B, \tilde A)}``=]{520}1b
\put(0,-170){\fbox{$\theta^{\tilde A}_u$}}
\efig}
$$
for every 1v-cells $U\colon A\to \tilde A$ and $u\colon B\to\tilde B$. 
\end{defn}

For the definition of vertical transformations of quasi-functors for double categories we refer the reader to the Appendix B.3.

\begin{defn} \delabel{modif btw horiz}
Let horizontal oplax transformations $\theta, \theta'$ and vertical lax transformations $\theta_0, \theta'_0$ acting between lax double quasi-functors 
$H_1, H_2, H_3, H_4\colon \A\times\B\to\C$ be given as in the left diagram below. Denote by $(-,A)_i\colon\B\to\C, (B,-)_i\colon\A\to\C, i=1,2,3,4$ the pairs of 
lax double functors corresponding to $H_1, H_2, H_3, H_4$, respectively. A modification $\Theta$ (on the left below) is given by a pair of modifications 
$\tau^A, \tau^B$ acting between transformations among lax double functors:
\begin{equation} \eqlabel{q-modif}
\scalebox{0.86}{
\bfig
\putmorphism(-150,50)(1,0)[H_1` H_2` \theta]{400}1a
\putmorphism(-150,-270)(1,0)[H_3 ` H_4 ` \theta' ]{400}1b
\putmorphism(-170,50)(0,-1)[\phantom{Y_2}``\theta_0]{320}1l
\putmorphism(250,50)(0,-1)[\phantom{Y_2}``\theta_0']{320}1r
\put(-30,-140){\fbox{$\tau$}}
\efig}
\qquad\qquad
\scalebox{0.86}{
\bfig
\putmorphism(-180,50)(1,0)[(-,A)_1` (-,A)_2`\theta^A]{550}1a
\putmorphism(-180,-270)(1,0)[(-,A)_3`(-,A)_4 `\theta'^A]{550}1b
\putmorphism(-170,50)(0,-1)[\phantom{Y_2}``\theta_0^A]{320}1l
\putmorphism(350,50)(0,-1)[\phantom{Y_2}``\theta_0'^A]{320}1r
\put(0,-140){\fbox{$\tau^A$}}
\efig}
\qquad\qquad
\scalebox{0.86}{
\bfig
\putmorphism(-180,50)(1,0)[(B,-)_1` (B,-)_2`\theta^B]{550}1a
\putmorphism(-180,-270)(1,0)[(B,-)_3`(B,-)_4 `\theta'^B]{550}1b
\putmorphism(-170,50)(0,-1)[\phantom{Y_2}``\theta_0^B]{320}1l
\putmorphism(350,50)(0,-1)[\phantom{Y_2}``\theta_0'^B]{320}1r
\put(0,-140){\fbox{$\tau^B$}}
\efig}
\end{equation}
such that $\tau^A_B=\tau^B_A$ for every $A\in\A, B\in\B$. 
\end{defn}

\section{Quasi-functors of $n$-variables - lax double category case} \selabel{n-double cat}


Quasi-functors of $n$-variables for 2-categories were introduced in \cite{Gray}. 
In \cite[Definition 4.5]{Fem1} we introduced pseudodouble quasi-functors of 3-variables, which we called ``cubical functors'' there 
in accordance to \cite[Section 4]{GPS} for 2-categories. At the beginning of this section we will generalize this definition into an $n$-variable and lax double case. In the remainder of the section we will show how quasi-functors for lax double categories make a multicategory.

\subsection{Quasi-functors of $n$-variables and transformations} \sslabel{q.f.t}

These two is what we need to build a multicategory. Apart from giving their definitions, in this subsection we will also prove the substitution for the multicategory.

\subsubsection{Quasi-functors of $n$-variables}

\begin{defn} \delabel{double qf 3}
A {\em lax double quasi-functor of $n$-variables} $H:\A_1\times...\times\A_n\to\C$ for $n\geq 2$ consists of lax double quasi-functors 
of two variables 
$$H(A_1,...,A_{i-1},\, -,\, A_{i+1},...,A_{j-1}, \, -,\, A_{j+1},..., A_n): \A_i\times\A_j\to\C$$
for all $i<j$ and all choices of objects $A_l\in\A_l, l=1,...,n$, which agree on objects and give unambiguous lax double functors 
$$H(A_1,...,A_{i-1},\, -,\,A_{i+1},..., A_n): \A_i\to\C,$$ 
such that the structure 2-cells of the former (from point 2. of \deref{gen-quasi}) 
relate in the following way (we simplify the notation by only writing the three involved variables):
\begin{enumerate} [i)]
\item for all 1h-cells $(f,g,h):(A,B,C)\to(A',B', C')$ in $\A\times\B\times\C$ it is  
$$\scalebox{0.8}{
\bfig
 \putmorphism(450,700)(1,0)[` `(A,g,C')_2]{680}1a
 \putmorphism(1140,700)(1,0)[` ` (f,B',C')_1]{650}1a

 \putmorphism(-150,250)(1,0)[``(A,B,h)_3]{600}1a
 \putmorphism(450,250)(1,0)[` `(f,B,C')_1]{680}1a
 \putmorphism(1130,250)(1,0)[` ` (A',g, C')_2]{650}1a

\putmorphism(450,700)(0,-1)[\phantom{Y_2}``=]{450}1r
\putmorphism(1750,700)(0,-1)[\phantom{Y_2}``=]{450}1r
\put(940,480){\fbox{$ (f,g,C')_{12}$}}

 \putmorphism(-150,-200)(1,0)[``(f,B,C)_1]{640}1a
 \putmorphism(460,-200)(1,0)[` `(A',B,h)_3]{680}1a

\putmorphism(-160,250)(0,-1)[\phantom{Y_2}``=]{450}1l
\putmorphism(1120,250)(0,-1)[\phantom{Y_3}``=]{450}1r
\put(310,20){\fbox{$(f,B,h)_{13}$}}

 \putmorphism(1140,-200)(1,0)[` ` (A',g,C')_2]{630}1a
\putmorphism(450,-200)(0,-1)[\phantom{Y_2}``=]{450}1r
\putmorphism(1750,-200)(0,-1)[\phantom{Y_2}``=]{450}1r

 \putmorphism(460,-650)(1,0)[` `(A',g,C)_2]{700}1a
 \putmorphism(1140,-650)(1,0)[ ` `(A',B',h)_3]{630}1a
\put(860,-430){\fbox{$ (A',g,h)_{23}$}}
\efig}
=
\scalebox{0.8}{
\bfig

 \putmorphism(-150,700)(1,0)[``(A,B,h)_3]{600}1a
 \putmorphism(450,700)(1,0)[` `(A,g,C')_2]{680}1a

 \putmorphism(-150,250)(1,0)[``(A,g,C)_2]{600}1a
 \putmorphism(450,250)(1,0)[` `(A,B',h)_3]{680}1a
 \putmorphism(1120,250)(1,0)[` `(f,B',C')_1]{650}1a

\putmorphism(-180,700)(0,-1)[\phantom{Y_2}``=]{450}1r
\putmorphism(1100,700)(0,-1)[\phantom{Y_2}``=]{450}1r
\put(280,470){\fbox{$(A,g,h)_{23}$}}
\put(910,20){\fbox{$(f,B',h)_{13}$}}

 \putmorphism(-150,-200)(1,0)[``(A,g,C)_2]{600}1a
 \putmorphism(450,-200)(1,0)[` `(f,B',C)_1]{680}1a
 \putmorphism(1100,-200)(1,0)[` ` (A',B',h)_3]{660}1a

\putmorphism(450,250)(0,-1)[\phantom{Y_2}``=]{450}1l
\putmorphism(1750,250)(0,-1)[\phantom{Y_2}``=]{450}1r
\putmorphism(-180,-200)(0,-1)[\phantom{Y_2}``=]{450}1r
\putmorphism(1100,-200)(0,-1)[\phantom{Y_2}``=]{450}1r
 \putmorphism(-150,-650)(1,0)[``(f,B,C)_1]{600}1a
 \putmorphism(580,-650)(1,0)[`` (A',g,C)_2]{540}1a
\put(280,-420){\fbox{$(f,g,C)_{12}$}}
\efig} \vspace{-0,2cm}
$$
\item for all 1v-cells $(u,v,z):(A,B,C)\to(\tilde A, \tilde B, \tilde C)$ in $\A\times\B\times\C$ it is  
$$
\scalebox{0.8}{
\bfig
 \putmorphism(0,480)(1,0)[` `=]{470}1a
\putmorphism(0,500)(0,-1)[` `(A,B,z)_3]{450}1l
\put(40,200){\fbox{$(v,z)_{23}$}}
\putmorphism(0,-320)(1,0)[ ` `=]{400}1a 
\putmorphism(0,90)(0,-1)[ ``(A,v,\tilde C)_2]{450}1l
\putmorphism(470,90)(0,-1)[``]{450}1l 
\putmorphism(490,90)(0,-1)[``(A,\tilde B,z)_3]{450}0l 
\putmorphism(470,500)(0,-1)[``]{450}1l 
\putmorphism(450,90)(1,0)[``=]{520}1a 
\putmorphism(970,90)(0,-1)[ ``]{450}1r
\putmorphism(950,90)(0,-1)[ ``(u,\tilde B,C)_1]{450}0r
\putmorphism(970,-300)(0,-1)[``(\tilde A,\tilde B,z)_3]{450}1r
\putmorphism(470,-730)(1,0)[ ``=]{500}1a 
\putmorphism(470,-300)(0,-1)[` `(u,\tilde B,\tilde C)_1]{450}1l 
\put(530,-430){\fbox{$(u,z)_{13}$}}

 \putmorphism(980,480)(1,0)[` `=]{470}1a %
\putmorphism(970,500)(0,-1)[` `(A,v,C)_2]{450}1l
\putmorphism(1450,500)(0,-1)[` `(u,B,C)_1]{450}1r
\putmorphism(1450,90)(0,-1)[``(\tilde A,v,C)_2]{450}1r
\putmorphism(1010,-320)(1,0)[` `=]{400}1a 
\put(1040,200){\fbox{$(u,v)_{12}$}}
\efig}
=
\scalebox{0.8}{
\bfig
 \putmorphism(-170,480)(1,0)[` `=]{480}1a
\putmorphism(-180,500)(0,-1)[\phantom{Y_2}` `(A,B,z)_3]{450}1l
\put(-120,250){\fbox{$(u,z)_{13}$}}
\putmorphism(-180,-320)(1,0)[` `=]{480}1a
\putmorphism(-180,90)(0,-1)[\phantom{Y_2}``]{450}1l
\putmorphism(-150,50)(0,-1)[\phantom{Y_2}``(u,B,\tilde C)_1]{450}0l
\putmorphism(300,90)(0,-1)[\phantom{Y_2}``]{450}1r
\putmorphism(280,90)(0,-1)[\phantom{Y_2}``(\tilde A,B,z)_3]{450}0r
\putmorphism(300,500)(0,-1)[\phantom{Y_2}` `(u,B, C)_1]{450}1r
\putmorphism(-650,90)(1,0)[ ``=]{400}1a
\putmorphism(-670,90)(0,-1)[ ``(A,v,\tilde C)_2]{450}1l
\putmorphism(-670,-300)(0,-1)[``(u,\tilde B,\tilde C)_1]{450}1l
\putmorphism(-650,-730)(1,0)[ ``=]{440}1a 
\putmorphism(-180,-300)(0,-1)[`  `(\tilde A,v,\tilde C)_2]{450}1r
\put(-600,-500){\fbox{$(u,v)_{12}$}}

\putmorphism(330,90)(1,0)[` `=]{450}1a
\putmorphism(770,90)(0,-1)[``(\tilde A,v,C)_2]{450}1r
\putmorphism(300,-300)(0,-1)[`  `]{450}1r 
\putmorphism(770,-300)(0,-1)[`  `(\tilde A,\tilde B,z)_3]{450}1r
\putmorphism(350,-730)(1,0)[` `=]{440}1a 
\put(360,-500){\fbox{$(v,z)_{23}$}}
\efig}
$$
(where we simplified the notation by writing $(v,z)_{23}$ for the 2-cell $(A,v,z)_{23}$ and so on...), 
\item for $(f,v,h):(A,B,C)\to(A',\tilde B,C')$ \vspace{-0,3cm}
$$\scalebox{0.86}{
\bfig
\putmorphism(-150,500)(1,0)[``(A,B,h)_3]{600}1a
 \putmorphism(480,500)(1,0)[` `(f,B,C')_1]{640}1a
 \putmorphism(-150,50)(1,0)[``(f,B,C)_1]{600}1a
 \putmorphism(470,50)(1,0)[` `(A',B,h)_3]{660}1a

\putmorphism(-180,500)(0,-1)[\phantom{Y_2}``=]{450}1l
\putmorphism(1100,500)(0,-1)[\phantom{Y_2}``=]{450}1r
\put(230,280){\fbox{$(f,B,h)_{13}$}}

\putmorphism(-170,-400)(1,0)[` `(f,\tilde B,C)_1]{640}1b
 \putmorphism(470,-400)(1,0)[` `(A',\tilde B,h)_3]{640}1b

\putmorphism(-180,50)(0,-1)[\phantom{Y_2}``(A,v,C)_2]{450}1l %
\putmorphism(450,50)(0,-1)[\phantom{Y_2}``]{450}1l
\putmorphism(1100,50)(0,-1)[\phantom{Y_3}``(A',v,C')_2]{450}1r
\put(-90,-180){\fbox{$(f,v,C)_{12}$}} 
\put(540,-180){\fbox{$(A',v,h)_{23}$}}
\efig}
\quad
=
\quad
\scalebox{0.86}{
\bfig
\putmorphism(-150,500)(1,0)[``(A,B,h)_3]{600}1a
 \putmorphism(480,500)(1,0)[` `(f,B,C')_1]{640}1a

 \putmorphism(-150,50)(1,0)[``(A,\tilde B,h)_3)]{600}1a
 \putmorphism(450,50)(1,0)[``(f,\tilde B,C')_1]{640}1a

\putmorphism(-180,500)(0,-1)[\phantom{Y_2}``(A,v,C)_2]{450}1l
\putmorphism(450,500)(0,-1)[\phantom{Y_2}``]{450}1r
\putmorphism(1100,500)(0,-1)[\phantom{Y_2}``(A',v,C')_2]{450}1r
\put(-90,280){\fbox{$(A,v,h)_{23}$}}
\put(540,280){\fbox{$(f,v,C')_{12}$}}

\putmorphism(-170,-400)(1,0)[` `(f,\tilde B,C)_1]{640}1b
 \putmorphism(470,-400)(1,0)[` `(A',\tilde B,h)_3]{640}1b

\putmorphism(-180,50)(0,-1)[\phantom{Y_2}``=]{450}1l
\putmorphism(1100,50)(0,-1)[\phantom{Y_3}``=]{450}1r
\put(270,-200){\fbox{$(f,\tilde B,h)_{13}$}}
\efig}
$$

and 2 similar conditions, one for $(f,g,z):(A,B,C)\to(A',B',\tilde C)$ and the other for $(u,g,h):(A,B,C)\to(\tilde A,B',C')$, 
\item for $(u,v,h):(A,B,C)\to(\tilde A,\tilde B,C')$  \vspace{-0,3cm}
$$\scalebox{0.86}{
\bfig
 \putmorphism(-120,500)(1,0)[` `=]{550}1a
 \putmorphism(450,500)(1,0)[` `(A,B,h)_3]{550}1a
\putmorphism(-140,520)(0,-1)[` `(A,v,C)_2]{480}1l
\put(-80,50){\fbox{$(u,v,C)_{12}$}}
\putmorphism(-150,-380)(1,0)[` `=]{540}1a
\putmorphism(-140,80)(0,-1)[``(u,\tilde B,C)_1]{450}1l
\putmorphism(430,50)(0,-1)[` `(\tilde A,v,C)_2]{450}1l
\putmorphism(430,520)(0,-1)[` `(u,B,C)_1]{480}1l
\put(470,290){\fbox{$(u,B,h)_{13}$}}
\putmorphism(430,50)(1,0)[``(\tilde A,B,h)_3]{540}1a
\putmorphism(1000,80)(0,-1)[``(\tilde A,v,C')_2]{450}1r
\putmorphism(1000,520)(0,-1)[``(u,B,C')_1]{480}1r
\putmorphism(450,-380)(1,0)[``(\tilde A,\tilde B,h)_3]{540}1b
\put(500,-190){\fbox{$(\tilde A, v,h)_{23}$}}
\efig}
\quad=\quad
\scalebox{0.86}{
\bfig
 \putmorphism(-150,500)(1,0)[` `(A,B,h)_3]{600}1a
 \putmorphism(450,500)(1,0)[` `=]{540}1a
\putmorphism(-180,520)(0,-1)[` `(A,v,C)_2]{450}1l
\put(-120,280){\fbox{$(A,v,h)_{23}$}}
\putmorphism(-150,-380)(1,0)[` `(\tilde A,\tilde B,h)_3]{500}1b
\putmorphism(-180,80)(0,-1)[``(u,\tilde B,C)_1]{450}1l
\putmorphism(400,80)(0,-1)[``(u,\tilde B,C')]{450}1r
\putmorphism(400,520)(0,-1)[` `(A,v,C')_2]{450}1r
\putmorphism(-150,50)(1,0)[``(A,\tilde B,h)_3]{500}1a
\putmorphism(1000,80)(0,-1)[``(\tilde A,v,C')_2]{450}1r
\putmorphism(1000,520)(0,-1)[``(u,B,C')_1]{450}1r
\putmorphism(450,-380)(1,0)[``=]{520}1b
\put(-120,-170){\fbox{$(u, \tilde B,h)_{13}$}}
\put(460,50){\fbox{$(u,v,C')_{12}$}}
\efig}
$$

and 2 similar conditions, one for $(u,g,z):(A,B,C)\to(\tilde A,B',\tilde C)$ and the other for $(f,v,z):(A,B,C)\to(A',\tilde B,C')$, 
where $f,g,h$ are 1h-cells and $u,v,z$ are 1v-cells, as usual. 
\end{enumerate}
\end{defn}

Thus a lax double quasi-functor of $n$-variables amounts to lax double quasi-functors of three variables 
$\A_i\times\A_j\times\A_k\to\C$ for all $i,j,k\in\{1,...,n\}$.

\subsubsection{Substitution for quasi-functors - lax double category case}

The following is a lax and double categorical version of \cite[Theorem I,4.7]{Gray}.

\begin{prop} \prlabel{subst}
Given lax double quasi-functors $F_i:\A_{i1}\times...\times\A_{im_i}\to\B_i$ of $m_i$-variables with $i=1,...,n$ and 
a lax double quasi-functor $G:\B_1\times...\times\B_n\to\C$ of $n$-variables, then the composition 
$$\Pi_{j=1}^{m_1}\A_{1j}\times...\times\Pi_{j=1}^{m_n}\A_{nj}\stackrel{F_1\times...\times F_n}{\longrightarrow}\B_1\times...\times\B_n
\stackrel{G}{\to}\C$$
is a lax double quasi-functor of $m_1+..+m_n$-variables (all quasi-functors are meant to be of type $(lax, o-l)$). 
\end{prop}

\begin{proof}
We prove this in \prref{subst-gen} both for double categories and Gray-categories. 
\qed\end{proof}

\subsubsection{Transformations of quasi-functors of three variables}

One can form a double category of lax double quasi-functor of $n$-variables. We content ourselves for the time being 
with a corresponding category. 
In \cite[Definition 4.6]{Fem1} we introduced vertical strict transformations for lax double quasi-functors of three variables. 
Since there we relied heavily on the assumption that 1v-cell components of vertical transformations had companions and conjoints, 
and with this hypothesis vertical strict transformations lift to horizontal pseudotransformations, vertical strict transformations 
were more suitable and sufficient for us to work with. Here we may directly consider horizontal pseudotransformations.

\begin{defn} \delabel{tr of qf 3}
A {\em horizontal pseudonatural transformation} $\theta: H_1\Rightarrow H_2$ between lax double quasi-functors $H_1,H_2:\A\times\B\times\C\to\E$ consists of horizontal pseudonatural transformations 
$$\theta^A: H_1(A,-,-) \Rightarrow H_2(A,-,-),$$ 
$$\theta^B: H_1(-,B,-) \Rightarrow H_2(-,B,-),$$ 
$$\theta^C: H_1(-,-,C) \Rightarrow H_2(-,-,C)$$ 
of lax double quasi-functors, which give unambiguous horizontal pseudonatural transformations 
$$\theta^{A;B}: (A,B,-)_3^1\to (A,B,-)_3^2$$
$$\theta^{B;C}: (-,B,C)_1^1\to (-,B,C)_1^2$$ 
$$\theta^{A;C}: (A,-,C)_2^1\to (A,-,C)_2^2$$ 
of lax double functors for each $(A,B,C)\in\A\times\B\times\C$, so that twelve equalities between structure 2-cells of 
$\theta^A, \theta^B, \theta^C$, on one hand, 
and the twelve structure 2-cells of both $H_1$ and $H_2$ (as quasi-functors from \deref{double qf 3}), on the other hand, hold. We present these twelve equalities schematically as lists consisting of those structure 2-cells which are related in the only possible way by one equation: 
$$(\theta^{A;C})_g, (\theta^{B;C})_f, (f,g,C)^i_{12}, \quad (\theta^{A;B})_h, (\theta^{B;C})_f, (f,B,h)^i_{13}, \quad 
(\theta^{A;B})_h, (\theta^{A;C})_g, (A,g,h)^i_{23}$$
$$(\theta^{A;C})^v, (\theta^{B;C})_f, (f,v,C)^i_{12}, \quad (\theta^{A;B})^z, (\theta^{B;C})_f, (f,B,z)^i_{13}, \quad 
(\theta^{A;B})^z, (\theta^{A;C})_g, (A,g,z)^i_{23}$$
$$(\theta^{A;C})_g, (\theta^{B;C})^u, (u,g,C)^i_{12}, \quad (\theta^{A;B})_h, (\theta^{B;C})^u, (u,B,h)^i_{13}, \quad 
(\theta^{A;B})_h, (\theta^{A;C})^v, (A,v,h)^i_{23}$$
$$(\theta^{A;C})^v, (\theta^{B;C})^u, (u,v,C)^i_{12}, \quad (\theta^{A;B})^z, (\theta^{B;C})^u, (u,B,z)^i_{13}, \quad 
(\theta^{A;B})^z, (\theta^{A;C})^v, (A,v,z)^i_{23}$$
where $(f,g,C)_{12}^i$ for $i=1,2$ presents a structure 2-cell $(f,g,C)_{12}$ from \deref{double qf 3} for $H_1$ and $H_2$, respectively, and similarly for the remaining eleven 2-cells of that type. 
\end{defn}

Observe that it holds $\theta^{A;B}(C)=\theta^{B;C}(A)=\theta^{A;C}(B)$ for all $(A,B,C)\in\A\times\B\times\C$ in the above definition. 
We also note that 
although the above twelve axioms have the same labels as those of \cite[Definition 4.6]{Fem1} defining vertical strict transformations, 
the shape of the involved structure 2-cells of the transformations are different.

\subsection{Multicategory of double categories}


We define two multicategories for double categories. For one we set $\M_n(\A_1,...,\A_n;\C)=q_n\x\Fun^*_\bullet(\A_1\times...\times\A_n,\C)$ (where $*$ denotes some weak double functors) 
 and for the other $\M_n(\A_1,...,\A_n;\C)=q_n\x Dbl_\bullet(\A_1\times...\times\A_n,\C)$. For the unary identities in $\M_n(\A;\A)$ we set the identity double functors. Substitution holds by \prref{subst} (and its strict version). 
We denote these two multicategories by $\Fun^*_\bullet(-,-)$ and $Dbl_\bullet(-,-)$.

\subsubsection{Left and right evaluation} \ssslabel{left}

We now establish evaluation quasi-functors. The two sides of evaluation will lead to quasi-functors of two different types. 
For this reason we will stress the type of functors in question.

\begin{prop} \prlabel{ev}
There is a lax double quasi-functor $ev:Lax_{hop}(\B, \C)\times\B\to\C$ of type $(lx,o-l)$ such that given any lax double quasi-functor 
$H:\A\times\B\to\C$ of type $(lx,o-l)$ it is 
$$ev(H^t\times \Id_\B)=H$$
where $H^t:\A\to Lax_{hop}(\B, \C)$ is the lax double functor corresponding to $H$ by \equref{quasi}.
\end{prop}

\begin{proof}
For a 0-cell $F\in Lax_{hop}(\B, \C)$ and any cell $x$ in $\B$ set $ev(F,x)=F(x)$; for a 0-cell $B\in\B$ and 1h- or 1v-cell $\alpha$ in 
$Lax_{hop}(\B, \C)$ set $ev(\alpha,B)=\alpha(B)$, the 1h- respectively 1v-cell component of the transformation $\alpha$; for a 2-cell 
$b$ in $Lax_{hop}(\B, \C)$ set $ev(b,B)=b(B)$, the 2-cell component of the modification $b$; and finally, for a 1h-cell $g$ and 1v-cell 
$v$ in $\B$ set $ev(\alpha,g)=\alpha_g$ and $ev(\alpha,u)=\alpha^u$, the 2-cell components of the (horizontal or vertical) transformation 
$\alpha$ (see Definitions 2.2, 2.4 and 2.7 of \cite{Fem} for the definitions and notations for horizontal and vertical transformations 
and (double) modifications).  

First condition for $ev$ to be a lax quasi-functor is that $ev(F,-):\B\to\C$ and $ev(-,B):Lax_{hop}(\B, \C)\to\C$ be lax double functors. 
It is clear that $ev(F,-)$ is such a functor, as so is $F$. Since both compositions of 1-cells in $Lax_{hop}(\B, \C)$ are vertical compositions of (horizontal resp. vertical) transformations, we obtain that  
$ev(-,B)$ is a strict functor (see Lemmata 2.3 and 2.6 of \cite{Fem} for the vertical compositions of transformations).

Condition $iii)$ of \prref{quasi-fun} is trivially fulfilled. 
For a (horizontal oplax or vertical lax) transformation $\alpha$ and a modification 
$b$ it is clear that $ev(\alpha,-)$ is a (horizontal oplax, resp. vertical lax) transformation, and that $ev(b,-)$ is a corresponding 
modification. In Table \ref{table:1} we argument why $ev(-,g)$ and $ev(-,v)$ are a horizontal lax and a 
vertical oplax transformation, respectively. 

\begin{table}[h!]
\begin{center}
\begin{tabular}{ c c } 
transformation axiom & its meaning on cells in $Lax_{hop}(\B, \C)$ on which evaluated \\ [0.5ex]
\hline
$(h.l.t.\x 1)$ for  $ev(-,g)$ & \cite[Lemma 2.3 part 3]{Fem} (vertical composition of h.o.t. $\alpha,\beta$) 
\\ [1ex]   
$(v.o.t.\x 1)$ for  $ev(-,v)$ & \cite[Lemma 2.3 part 2]{Fem} (vertical composition of h.o.t. $\alpha,\beta$) 
\\ [1ex]   
$(h.l.t.\x 3)$ for  $ev(-,g)$ & \cite[Lemma 2.6 part i)]{Fem} (vertical composition of v.l.t. $\alpha_0,\beta_0$) 
\\ [1ex]   
$(v.o.t.\x 3)$ for  $ev(-,v)$ & \cite[Lemma 2.6 part ii)]{Fem} (vertical composition of v.l.t. $\alpha_0,\beta_0$) 
\\ [1ex]   
$(h.l.t.\x 2), (h.l.t.\x 4)$ for  $ev(-,g)$ & hold clearly \\ [1ex]    
$(v.o.t.\x 2), (v.o.t.\x 4)$ for  $ev(-,v)$ & hold clearly \\ [1ex]   
$(h.l.t.\x 5)$ for  $ev(-,g)$ & $(m.ho\x vl.\x 1)$ for $ev(b,-)$ \\ [1ex]   
$(v.o.t.\x 5)$ for  $ev(-,v)$ & $(m.ho\x vl.\x 2)$ for $ev(b,-)$ \\ [1ex]   
\end{tabular}
\caption{Why $ev(-,g)$ and $ev(-,v)$ are transformations}
\label{table:1}
\end{center}
\end{table}
Finally, the modification axioms $(m.hl\x vo.\x 1)$ at a 1h-cell $\alpha$ in $Lax_{hop}(\B, \C)$ and $(m.hl\x vo.\x 2)$ 
at a 1v-cell $\alpha_0$ in $Lax_{hop}(\B, \C)$ for $ev(-,\beta)$ and a 2-cell $\beta$ in $\B$ coincide with the transformation axioms 
$(h.o.t.\x 5)$ for $\alpha$ and $(v.l.t.\x 5)$ for $\alpha_0$, respectively. This finishes the proof that $ev$ is a lax double 
quasi-functor by \prref{quasi-fun}, and the last statement is clear. 
\qed\end{proof}

\begin{thm} \thlabel{left closed multicat}
\begin{enumerate}
\item
Given a lax double quasi-functor of $n+1$-variables $H:\Pi_{i=1}^n \A_i\times\B\to\C$ of type $(lx,o-l)$, 
there is a unique lax double quasi-functor of $n$-variables $H^t:\Pi_{i=1}^n \A_i\to Lax_{hop}(\B,\C)$ of type $(lx,o-l)$ 
such that $ev(H^t\times \Id_\B)=H$. 
\item 
The above correspondence extends to a natural isomorphism of double categories: 
$$q_{n+1}\x Lax_{hop}(\Pi_{i=1}^n \A_i\times\B,\C)\iso q_n\x Lax_{hop}(\Pi_{i=1}^n \A_i, Lax_{hop}(\B,\C)).$$
\end{enumerate} 
\end{thm}

\begin{proof} 
%
For the first part the reader is referred to the proof of \thref{mixed-quasi n}: the only difference is that there mixed 
``strict-weak'' 
quasi-functors were treated instead of entirely lax ones. This does not affect the proof: one can simply substitute the only strict 
component there, which is leftmost, and take ``lax'' for the weakness, $*=lx$. We omit the proof of the second part, it is straightforward. 
\qed\end{proof}


\bigskip

From the proof of \prref{ev} it is clear that the analogous result holds in the strict functor case, so that we have: 

\begin{thm} \thlabel{strict-left closed} \tlabel{strict-left closed}
\begin{enumerate}
\item There is a double quasi-functor $ev:Dbl_{hop}(\B, \C)\times\B\to\C$ of type $(lx,o-l)$ such that given any double quasi-functor 
$H:\A\times\B\to\C$ of type $(lx,o-l)$ it is 
$$ev(H^t\times \Id_\B)=H$$
where $H^t:\A\to Dbl_{hop}(\B, \C)$ is the double functor corresponding to $H$ by \equref{quasi}.
\item Given a double quasi-functor of $n+1$-variables $H:\Pi_{i=1}^n \A_i\times\B\to\C$ of type $(lx,o-l)$, there is a unique 
double quasi-functor of $n$-variables $H^t:\Pi_{i=1}^n \A_i\to Dbl_{hop}(\B,\C)$ of type $(lx,o-l)$ such that $ev(H^t\times \Id_\B)=H$. 
\item There is a double category isomorphism for $n\geq 2$ 
$$q_{n+1}\x Dbl_{hop}(\Pi_{i=1}^n \A_i\times\B,\C)\iso q_n\x Dbl_{hop}(\Pi_{i=1}^n \A_i, Dbl_{hop}(\B,\C)).$$
\end{enumerate}
\end{thm}

We have that the multicategories $Lax_{hop}(-,-)$ and $Dbl_{hop}(-,-)$ are left closed.

\begin{rem}
Observe that the isomorphisms in \thref{left closed multicat} correspond to the isomorphisms $\lambda_{n+1}^{-1}$ from \cite{BL},  
and that in our case it is also true that ``$\lambda_n$ restrict to tight maps, {\em i.e.} they preserve tightness''. 
\end{rem}

\medskip

We now turn to right evaluation. Apart from the evaluation from \prref{ev} there are three more evaluation maps: 
 $ev_2:Lax_{hlx}(\B, \C)\times\B\to\C$, \,\, $ev_3:\B\times Lax_{hop}(\B, \C)\to\C$ and $ev_4:\B\times Lax_{hlx}(\B, \C)\to\C$. 
Analogously as in the proof of \prref{ev} one readily sees that, if $\beta$ denotes a horizontal transformation {\em i.e.} 
a 1h-cell in the corresponding inner-hom, then $ev_2(\beta,-)$ and $ev_4(-,\beta)$ will be transformations of type $\#=lx$ and 
$ev_3(-,\beta)$ of type $\bullet=oplx$, because that is the type of their respective inner-homs. Then according to \prref{quasi-fun}, $ev_2$ and $ev_3$ will be quasi-functors of type $(lx,l-o)$, while $ev_4$ will be a quasi-functor of type $(lx,o-l)$. In what follows we will be interested in the evaluation $ev_4$ and this is the one to which we will refer to as the {\em right evaluation}. Observe that $ev_2$ and $ev_4$ are related via the flip map, so that the rightmost triangle in the diagram below commutes. The remaining two triangles commute due to \equref{quasi-right bij}, where $H_\bullet$ denotes a quasi-funtor of type $(lx,o-l)$, and $H_\#$ a one of type $(lx,l-o)$. 
$$\scalebox{0.8}{
\bfig
 \putmorphism(400,700)(1,0)[\A\times\B `Lax_{hlx}(\B, \C)\times\B ` H^t\times 1]{960}1a
 \putmorphism(-100,380)(1,1)[``flip]{240}1l
 \putmorphism(-150,250)(1,0)[\B\times\A `\C`H_\bullet]{960}1b
\putmorphism(400,700)(1,-1)[\phantom{Y_2}``H_\#]{460}1l
\putmorphism(1300,700)(-1,-1)[\phantom{Y_2}``ev_2]{450}1r
 \putmorphism(780,250)(1,0)[\phantom{A\ot B}`\B\times Lax_{hlx}(\B, \C)` \crta{ev}=ev_4]{1180}{-1}b
\putmorphism(1350,700)(1,-1)[\phantom{Y_2}``flip]{450}1r
\put(750,490){\fbox{\equref{quasi-right bij}}}
\efig}
$$
The above arguments are sufficient to prove the part 1. in the next theorem.

\begin{thm} \thlabel{lax-right closed}
\begin{enumerate}
\item There is a lax double quasi-functor $\crta{ev}:\B\times Lax_{hlx}(\B, \C)\to\C$ of type $(lx,o-l)$ such that given any lax double quasi-functor $H:\B\times\A\to\C$ of type $(lx,o-l)$ it is 
$$\crta{ev}(\Id_\B\times H^t)=H$$
where $H^t:\A\to Lax_{hlx}(\B, \C)$ is the lax double functor corresponding to $H$ by 
\equref{quasi-right bij}. 
\item Given a lax double quasi-functor of $n+1$-variables $H:\B\times\Pi_{i=1}^n \A_i\to\C$ of type $(lx,o-l)$, there is a unique 
lax double quasi-functor of $n$-variables $H^t:\Pi_{i=1}^n \A_i\to Lax_{hlx}(\B,\C)$ of type $(lx,o-l)$ such that 
$\crta{ev}(\Id_\B\times H^t)=H$. 
\item There is a double category isomorphism 
$$q_{n+1}\x Lax_{hlx}(\B\times\Pi_{i=1}^n \A_i,\C)\iso q_n\x Lax_{hlx}(\Pi_{i=1}^n \A_i, Lax_{hlx}(\B,\C)).$$
\item There is a double category isomorphism 
$$q_{n+1}\x Dbl_{hlx}(\B\times\Pi_{i=1}^n \A_i,\C)\iso q_n\x Dbl_{hlx}(\Pi_{i=1}^n \A_i, Dbl_{hlx}(\B,\C)).$$
\end{enumerate}
\end{thm}

\begin{proof}
We only comment the part 2. For it, one substitutes $\A$ in the diagrams above by 
$\Pi_{i=1}^n\A_i$, and $H^t:\Pi_{i=1}^n \A_i\to Lax_{hlx}(\B,\C)$ should now be a quasi-functor of type $(lx,o-l)$. Observe 
that since the types of quasi-functors $H_\#$ and $ev_3$ in the middle triangle \equref{quasi-right bij} are $(lx,l-o)$, the 
quasi-functor $H^t:\Pi_{i=1}^n \A_i\to Lax_{hlx}(\B,\C)$ on its top arrow is also of that type, but it is the composite $flip\circ
H^t\circ flip$ that is of type $(lx,o-l)$, as we wish, and it is this composite that obeys the identity in point 2.
\qed\end{proof}

We thus have that the multicategories $\OO^{lx}_{hlx}$ and $\OO^{st}_{hlx}$ with $\Oo=Dbl$ are right closed.


\subsection{Dualities} \sslabel{dual}

In this subsection we establish some dualities. For a double category $\A$ let $\A^t$ denote the {\em transpose double category} which is the dual double category to $\A$ in which both 1h-cells and 1v-cells are taken in the opposite direction, \cite[Section 3.2.2]{MG}. To visualize this we 
show a transpose $\omega^t$ of a 2-cell $\omega$ below 
$$
\scalebox{0.86}{
\bfig
\putmorphism(-150,190)(1,0)[A`B`f]{500}1a
\putmorphism(-150,-180)(1,0)[\tilde A`\tilde B`g]{500}1a
\putmorphism(-150,190)(0,-1)[\phantom{Y_2}``u]{370}1l
\putmorphism(360,190)(0,-1)[\phantom{Y_2}``v]{370}1r
\put(30,30){\fbox{$\omega$}}
\efig}
\qquad\qquad
\scalebox{0.86}{
\bfig
\putmorphism(-150,190)(1,0)[\tilde B`\tilde A`g^t]{500}1a
\putmorphism(-150,-180)(1,0)[B`A.`f^t]{500}1a
\putmorphism(-150,190)(0,-1)[\phantom{Y_2}``v^t]{370}1l
\putmorphism(360,190)(0,-1)[\phantom{Y_2}``u^t]{370}1r
\put(30,30){\fbox{$\omega^t$}}
\efig}
$$
Observe that reading the data of a lax functor $F:\A\to\B$ from the double categories $\A^t$ and $\B^t$ one obtains a lax 
double functor $F^t:\A^t\to\B^t$ defined by $F^t(a^t)=F(a)$ for any cell $a$ in $\A$. Namely, for a 1h-cell composition $gf$ in $\A$ we have the same structure 2-cell $F^t(f^t)\comp^t F^t(g^t)=F(g)F(f)\Rightarrow F(gf)=F^t(f^t\comp^t g^t)$ as for $F$. 

Now, a horizontal oplax transformation $\alpha:F\Rightarrow G:\A\to\B$ is given via structure 2-cells 
$$
\scalebox{0.86}{
\bfig
 \putmorphism(-170,190)(1,0)[F(A)`F(B)`F(f)]{540}1a
 \putmorphism(360,190)(1,0)[\phantom{F(f)}`G(B) `\alpha(B)]{560}1a
 \putmorphism(-170,-190)(1,0)[F(A)`G(A)`\alpha(A)]{540}1a
 \putmorphism(360,-190)(1,0)[\phantom{G(B)}`G(A) `G(f)]{560}1a
\putmorphism(-180,190)(0,-1)[\phantom{Y_2}``=]{380}1r
\putmorphism(940,190)(0,-1)[\phantom{Y_2}``=]{380}1r
\put(280,0){\fbox{$\alpha_f$}}
\efig}
\qquad\text{and}\qquad
\scalebox{0.86}{
\bfig
\putmorphism(-150,190)(1,0)[F(A)`G(A)`\alpha(A)]{560}1a
\putmorphism(-150,-180)(1,0)[F(\tilde A)`G(\tilde A)`\alpha(\tilde A)]{600}1a
\putmorphism(-180,190)(0,-1)[\phantom{Y_2}``F(u)]{370}1l
\putmorphism(410,190)(0,-1)[\phantom{Y_2}``G(u)]{370}1r
\put(30,30){\fbox{$\alpha^u$}}
\efig}
$$
in $\B$ for every 1h-cell $f\colon A\to B$ and 1v-cell $u: A\to\tilde A$ in $\A$. When read from $\A^t$ and $\B^t$, together with the axioms for a \axiomref{$h.o.t.$} for $\alpha$, one obtains a \axiom{$h.l.t.$} $\alpha^t:G^t\Rightarrow F^t:\A^t\to\B^t$. Similarly, a 
\axiomref{$v.l.t.$} $\alpha_0: F\Rightarrow \tilde F:\A\to \B$ corresponds to a \axiom{$v.o.t.$} $\alpha_0^t:\tilde F^t\Rightarrow F^t:\A^t\to\B^t$. Finally, a 2-cell in 
$Lax_{hop}(\A,\B)$, {\em i.e.} a modification \axiomref{$m.ho\x vl$} among transformations $F,G,\tilde F,\tilde G:\A\to\B$ 
corresponds to a modification \axiom{$m.hl\x vo$} among transformations $\tilde G^t,\tilde F^t,G^t,F^t:\A^t\to\B^t$. We may conclude 
that there is a (strict) double category isomorphism 
\begin{equation} \eqlabel{Verity-t}
Lax_{hop}(\A,\B)\iso(Lax_{hlx}(\A^t,\B^t))^t.
\end{equation}
From here we get 
$$Lax_{hop}(\A,Lax_{hop}(\B,\C))\iso\large(Lax_{hlx}(\A^t,Lax_{hlx}(\B^t, \C^t))\large)^t.$$
Applying \equref{quasi} on the left with both transformation types being $\bullet$, and on the right with both being $\#$, we further get 
\begin{equation} \eqlabel{quasi-t}
q\x Lax_{hop}\vert_{(hop,hop)}(\A\times\B,\C)\iso( q\x Lax_{hlx}\vert_{(hlx,hlx)}(\A^t\times\B^t,\C^t))^t.
\end{equation}
On the other hand, by further iterating the above relation between inner-homs and setting $[\A,\B]=Lax_{hop}(\A,\B)$ 
and $\{\A,\B\}=Lax_{hlx}(\A,\B)$ for short, we obtain
\begin{equation} \eqlabel{inner-iter}
Lax_{hop}(\A_1, [\A_2,...[\B,\C]...])\iso(Lax_{hlx}(\A_1^t, \{\A_2^t,...\{\B^t,\C^t\}...\}))^t.
\end{equation}
Now, by iterative applications of \thref{left closed multicat}, on one hand, and of part 3. of \thref{lax-right closed}, on the other hand, we get 
$$q_n\x Lax_{hop}(\Pi_{i=1}^n \A_i, Lax_{hop}(\B,\C))\iso q\x Lax_{hop}(\A_1\times\A_2, [\A_3,...,[\A_n,[\B,\C]]...])$$
and 
\begin{align*}
q_n\x Lax_{hlx}(\Pi_{i=n}^1 \A_i^t, Lax_{hlx}(\B^t,\C^t)) 
& \iso q\x Lax_{hlx}(\A_2^t\times\A_1^t, \{\A_3^t,...,\{\A_n^t, \{\B^t,\C^t\}\}...\}) \\
&\stackrel{\equref{quasi-t}}{\iso} \large(q\x Lax_{hop}(\A_2\times\A_1, \{\A_3^t,...,\{\A_n^t, \{\B^t,\C^t\}\}...\}^t) \large)^t \\
&\stackrel{\equref{inner-iter}}{\iso} \large(q\x Lax_{hop}(\A_2\times\A_1, [\A_3,...,[\A_n,[\B,\C]]...]) \large)^t, 
\end{align*}
respectively. Then among quasi-functors of $n$-variables of opposite types we find the relation
$$q_n\x Lax_{hop}(\Pi_{i=1}^n \A_i, Lax_{hop}(\B,\C))\iso \large(q_n\x Lax_{hlx}(\A_n^t\times\Pi_{i=n-1}^3 \times \A_1^t\times\A_2^t, 
Lax_{hlx}(\B^t,\C^t))\large)^t$$
and thus also
\begin{equation} \eqlabel{d-n}
q_{n+1}\x Lax_{hop}(\Pi_{i=1}^n \A_i\times\B,\C)\iso( q_{n+1}\x Lax_{hlx}(\B^t\times\A_n^t\times...\times\A_3^t\times \A_1^t\times\A_2^t,\C^t)^t.
\end{equation}


\begin{rem}
In analogy to the isomorphisms in \thref{left closed multicat} we could denote the isomorphisms in the part 3. of 
\thref{lax-right closed} by $\crta{\lambda}_{n+1}^{-1}$. The authors in \cite{BL} do not use these morphisms, as they do not study right closedness. Instead and opposed to our \equref{d-n}, they use dualities $d^{n+1}$ among multimaps of the same type 
(see Propositions 4.4, 4.7, 4,11 thereof). Mind the difference between  the opposite Gray-categories $\A^{co}$ of \cite{BL} and our ``transpose duality'' for double categories: in the former 2-cells are reversed and the 3-dimensional functors are pseudo, which is necessary to be able to define the isomorphism 
functor $(-)^{co}:G\x\Cat_p\to G\x\Cat_p$ in \cite[Section 3.5.]{BL}. In our duality 
\equref{d-n} the (double) functors are lax. Observe that it is derived from \equref{Verity-t}, which generalizes to double categories the duality of Verity $Bicat^{op}(\B,\C)\iso (Bicat(\B^{op},\C^{op}))^{op}$ for bicategories $\B,\C$, where $op$ stands for reversing only the 1-cells, $Bicat(\B,\C)$ is the bicategory of lax functors, lax transformations and modifications, whereas $Bicat^{op}(\B,\C)$ is the bicategory of lax functors, oplax transformations and modifications (see page 21 of \cite{Ver}). 
\end{rem}

\section{Relations to the meta-product - lax double category case} \selabel{meta-lax double}

In \cite[Proposition 6.3 and (15)]{Fem} we obtained the following isomorphisms of double categories: 
$$q\x Lax_{hop}(\A\times\B,\C)\iso Dbl_{hop}(\A\ot^{lx}_{o-l}\B,\C)$$ 
and 
$$Dbl_{hop}(\A\ot^{lx}_{o-l}\B,\C)\iso Lax_{hop}(\A,Lax_{hop}(\B,\C)),$$
where the latter one is obtained as a consequence of \equref{quasi} with $\Oo=Dbl$ and type $(lx,o-l)$. 
Here $\ot^{lx}_{o-l}$ indicates that the meta-product is obtained from an inner-hom whose cells are 
lax double functors, 1h-cells horizontal oplax transformations and 1v-cells vertical lax transformations. 
From the proof of \cite[Proposition 6.3]{Fem}, the first isomorphism above, one sees immediately that 
\begin{equation} \eqlabel{any-strict}
q\x \Fun^*_{hop}(\A\times\B,\C)\iso Dbl_{hop}(\A\ot^*_{o-l}\B,\C)
\end{equation} 
holds for any type $*$ of double functors. Observe that the type $*$ of functors used in the construction of meta-product is reflected in the type of functors that one obtains on the left-hand side. However, the only type of functors from the meta-product on the right hand-side 
that one may get out of any type of functors on the left-hand side is the strict one. 
One has that for every double quasi-functor $H\colon\A\times\B\to\C$ of type $*$ there is a unique  
strict double functor $\crta H\colon \A\ot\B\to\C$ such that $H=\crta{H}J^*_{o-l}$ with $J^*_{o-l}$ from \equref{J}. 
The isomorphism \equref{any-strict} is crucial for checking if $\ot^*_{o-l}$ is a monoidal product. 
Namely, if we can prove the marked isomorphisms in the following sequence of double category isomorphisms:
\begin{align}
Dbl_{hop}((\A\ot^*_{o-l}\B)\ot^*_{o-l}\C,\D) & \iso q\x \Fun^*_{hop}((\A\ot^*_{o-l}\B)\times\C,\D) \label{eq1}\\
& \stackrel{!}{\iso} q_3\x \Fun^*_{hop}(\A\times \B\times\C,\D) \label{eq2}\\
& \stackrel{!}{\iso}  q\x \Fun^*_{hop}(\A\times (\B\ot^*_{o-l}\C),\D) \label{eq3}\\
& \iso Dbl_{hop}(\A\ot^*_{o-l}(\B\ot^*_{o-l}\C),\D), \nonumber
\end{align}
then by The Yoneda lemma we would obtain that there is a double category isomorphism 
$(\A\ot^*_{o-l}\B)\ot^*_{o-l}\C\to \A\ot^*_{o-l}(\B\ot^*_{o-l}\C)$. However, a $*$-type quasi-functor in the right-hand side of 
(\ref{eq1}) consists of $*$-type functors $\A\ot^*_{o-l}\B\to\D$ and $\C\to\D$, and we should relate them bijectively with 
three quasi-functors of two variables appearing in (\ref{eq2}). For the same reason as explained above, 
the functors from $\A\ot^*_{o-l}\B$ in (\ref{eq1}) can only be strict, {\em i.e.} in (\ref{eq1}) it should be 
$q\x \Fun^*_{hop}(-,-)=q\x Dbl_{hop}(-,-)$. This determines $*=st$ in the above reasoning 
for obtaining an invertible associativity constraint for $(Dbl, \ot^*_{o-l})$. 

On the other hand, if one does not have at least one of the two isomorphisms marked above, {\em i.e.} either left or right representability, one 
is not able to obtain a possibly non-invertible constraint $\alpha$ going in either direction. Namely, setting $\D=\A\ot^*_{o-l}(\B\ot^*_{o-l}\C)$, 
a desired double functor $(\A\ot^*_{o-l}\B)\ot^*_{o-l}\C\to \A\ot^*_{o-l}(\B\ot^*_{o-l}\C)$ would correspond to 
$J_{\A\ot^*_{o-l}\B,\C}(J_{\A,\B}\times\Id_\C)$ in (\ref{eq2}), with $J=J^*_{o-l}$. 
Similarly, we would obtain a double functor 
$\A\ot^*_{o-l}(\B\ot^*_{o-l}\C)\to (\A\ot^*_{o-l}\B)\ot^*_{o-l}\C$ as the correspondent to $J_{\A,\B\ot^*_{o-l}\C}(\Id_\A\times J_{\B,\C})$. 
Thus, having representability on one side yields an associativity constraint in one direction. 

Under the condition $*=st$, necessary in the above reasoning for obtaining an associativity constraint for $(Dbl, \ot^*_{o-l})$, 
we can prove the results of the following subsection.

\subsection{The case of strict (double) functors} \sslabel{strict case}

By left closedness of $Dbl_{hop}(-,-)$ and right closedness of $Dbl_{hlp}(-,-)$ one can prove 
separately left representability of $\ot^{st}_{o-l}$ and right representability of $\ot^{st}_{l-o}$:
%
\begin{multline*} 
q_n\x Dbl_{hop}((\A\ot^{st}_{o-l}\B)\times(\Pi_{i=1}^{n-1}\C_i),\D)  
\iso q_{n+1}\x Dbl_{hop}(\A\times \B\times(\Pi_{i=1}^{n-1}\C_i),\D) \\
 \text{left representability of} \,\, \ot^{st}_{o-l} 
\end{multline*}  
\begin{multline*} 
q_n\x Dbl_{hlx}((\Pi_{i=1}^{n-1}\A_i)\times (\B\ot^{st}_{l-o}\C),\D)
\iso q_{n+1}\x Dbl_{hlx}((\Pi_{i=1}^{n-1}\A_i)\times \B\times\C,\D)  \\
  \text{right representability of} \,\, \ot^{st}_{l-o}.
	\end{multline*}  
Namely, by multicategory-closedness \thref{strict-left closed} we can pass the $\Pi_{i=1}^{n-1}\C_i$'s to the codomain to get an equivalent question. Thus for left representability it is enough to prove: 
$Dbl_{hop}((\A\ot^{st}_{o-l}\B),\C)\iso q\x Dbl_{hop}(\A\times \B, \C)$ for any double category $\C$, but this is true by \equref{any-strict}. 
For the other isomorphism use \thref{lax-right closed} (to drop out the $\Pi_{i=1}^{n-1}\A_i$'s) and the ``hlx'' version of \equref{any-strict}.

\bigskip

Let 
\begin{equation} \eqlabel{acumul left}
\crta\ot^n_{o-l}\A_i=(..(\A_1\ot\A_2)\ot...\ot\A_{n-1})\ot\A_n
\end{equation} 
whereby $\ot=\ot^{st}_{o-l}$.

\begin{cor} \colabel{strict-rep} [{\em representability}]
There are bijections:
\begin{enumerate}
\item 
$q_n\x Dbl_{hop}((\Pi_{i=1}^{r}\A_i)\times(\B\ot^{st}_{o-l}\C)\times(\Pi_{1}^{s}\D_i),\E)  \iso$ \\
$$q_{n+1}\x Dbl_{hop}((\Pi_{i=1}^{r}\A_i)\times(\B\times\C)\times(\Pi_{i=1}^{s}\D_i),\E) $$
natural in $\A_i$'s, $\D_i$'s and $\E$, with $n=r+1+s$;
%
\item 
$q_n\x Dbl_{hop}((\Pi_{i=1}^{r}\A_i)\times(\crta\ot^k_{o-l}\B_i)\times(\Pi_{1}^{s}\D_i),\E)  \iso$ \\
$$q_p\x Dbl_{hop}((\Pi_{i=1}^{r}\A_i)\times(\Pi_{i=1}^k\B_i)\times(\Pi_{i=1}^{s}\D_i),\E) $$
natural in $\A_i$'s, $\D_i$'s and $\E$, with $n=r+1+s$ and $p=r+k+s$;
\item $Dbl_{hop}(\crta\ot^n_{o-l}\A_i,\B)\iso q_n\x Dbl_{hop}(\Pi_{i=1}^n\A_i,\B)$ \hspace{0,2cm} natural in $\B$. 
\end{enumerate}
\end{cor}

\begin{proof}
For the first part, take an $n$-variable quasi-functor $\F$ from the left-hand side. To prove that it yields a quasi-functor 
$\crta\F$ of $n+1$-variables on the right hand-side, whereby $\crta\F$ is induced in the obvious way using the universal property 
\equref{any-strict}, 
we should check the quasi-functor property of $\crta\F$ for any three variables. We do this by analyzing the known quasi-functor 
property of $\F$. If three variables of $\F$ fall into $\B\ot^{st}_{o-l}\C\times\Pi_{i=1}^{s}\D_i$, we are done by the above-proved 
left representability of $\ot^{st}_{o-l}$. Next we study the case when three variables live in 
$\A_i\times\A_j\times (\B\ot^{st}_{o-l}\C)$. We have: 
\begin{align}
 q_3\x Dbl_{hop}(\A_i\times\A_j\times (\B\ot^{st}_{o-l}\C),\E) & \stackrel{\tref{strict-left closed}}{\iso}
q\x Dbl_{hop}(\A_i\times\A_j, Dbl_{hop}(\B\ot^{st}_{o-l}\C,\E)) \nonumber\\ 
& \stackrel{\equref{any-strict},\equref{quasi}}{\iso} q\x Dbl_{hop}(\A_i\times\A_j, Dbl_{hop}(\B, Dbl_{hop}(\C,\E)))  \nonumber\\
& \stackrel{\tref{strict-left closed}}{\iso} q_3\x Dbl_{hop}(\A_i\times\A_j \times\B, Dbl_{hop}(\C,\E))  \nonumber\\
& \stackrel{\tref{strict-left closed}}{\iso}q_4\x Dbl_{hop}(\A_i\times\A_j \times\B\times\C,\E) \nonumber
\end{align}
thus the desired quasi-functor property for $\crta\F$ in this case is fulfilled. The last case that we need to check is when three variables 
live in $\A_i\times(\B\ot^{st}_{o-l}\C)\times\D_j$. Then the eight axioms of \deref{double qf 3} are fulfilled, whereby in the first one 
the 1h-cells $(f,g,h):(A,B, C)\to(A', B', C')$ are now substituted by $(f,g\ot C,h):(A,B\ot C, D)\to(A', B'\ot C', D')$ in 
$\A_i\times(\B\ot^{st}_{o-l}\C)\times\D_j$. Similar adjustments apply to the remaining seven axioms. The axioms express compatibilities 
of certain 2-cells living in $\E$. By the universal property the restriction of $\F$ to $\B\ot^{st}_{o-l}\C$ determines a quasi-functor 
$\B\times\C\to\E$ of type $(st, o-l)$, 
which consists of two strict functors $\B\to\E, \C\to\E$. By definition these two strict functors together with the remaining strict 
functors from the rest of variables of $\F$, on one hand, and all the quasi-functors of two variables making $\F$, on the other hand, 
determine $\crta\F$. The latter is indeed a quasi-functor: the desired eight axioms in $\E$ are precisely the ones named earlier. 

The second part follows from the first part by iteration, and the third one 
follows from the second part with $r=s=0$.  
(The third part is also a consequence of first applying alternating \equref{any-strict} and \equref{quasi} iteratively until getting 
to $q\x Db_{hop}(\A_1\times\A_2, I)$ where $I$ is an accumulatively obtained inner-hom, and then applying ``backwards'' 
\thref{strict-left closed}.)
\qed\end{proof}

\medskip

The point 2. above means that the multicategory $\OO^{st}_{hop}$ with $\Oo=Dbl$ is representable.

\medskip



Completely analogously to \coref{strict-rep} one has 
\begin{multline*} 
q_n\x Dbl_{hlx}((\Pi_{i=1}^{r}\A_i)\times(\B\ot^{st}_{l-o}\C)\times(\Pi_{1}^{s}\D_i),\E)  \iso \\
 q_{n+1}\x Dbl_{hlx}((\Pi_{i=1}^{r}\A_i)\times(\B\times\C)\times(\Pi_{i=1}^{s}\D_i),\E) 
\end{multline*} 
and 
$$Dbl_{hlx}(\crta\ot^n_{l-o}\A_i,\B)\iso q_n\x Dbl_{hlx}(\Pi_{i=1}^n\A_i,\B),$$
where $\crta\ot^n_{l-o}$ is the l-o and right-hand side analogue of \equref{acumul left}. 

\medskip

From the discussion below (\ref{eq2}) and part 1. of \coref{strict-rep} (once with $r=0$ and then with $s=0$) we have that there is a double category isomorphism 
$$\alpha^{st}:(\A\ot^{st}_{o-l}\B)\ot^{st}_{o-l}\C\to \A\ot^{st}_{o-l}(\B\ot^{st}_{o-l}\C).$$
Joining \equref{quasi} and \equref{any-strict} one obtains 
\begin{equation} \eqlabel{closedness}
Dbl_{hop}(\A\ot^*_{o-l}\B,\C)\iso\Fun^*_{hop}(\A,\Fun^*_{hop}(\B,\C)).
\end{equation} 
From here we see that if $*\not= st$, then $(Dbl, \ot^*_{o-l})$ can be neither monoidal (we got the associativity constraint 
$\alpha^{st}$ via Yoneda lemma only in the case $*=st$), nor closed.

\begin{thm} \thlabel{Dbl left closed mon}
The category $(Dbl, \ot^{st}_{o-l})$ is biclosed monoidal. 
\end{thm}


Instead of proving this theorem here, it will be covered in \ssref{I and II} using \coref{strict-rep}. 
It generalizes Gray's result from 2-categories to double categories, but also B\"ohm's one. Namely, recall that Gray's 
inner-hom consisted of 2-functors, lax transformations and modifications, while B\"ohm's inner-hom had {\em pseudo} transformations both for 1h- and 1v-cells.

\medskip

\begin{thm}
The category $(Dbl, \ot^{st}_{o-l}, \ot^{st}_{l-o}, 1)$ is strong duoidal. 
\end{thm}

\begin{proof}
For the interchange take the strict double functor induced by $1\times flip\times 1$ and the two projections $J^{st}_\bullet$. The two monoidal categories are then braided and isomorphic as braided monoidal categories by \cite[Prop.6.11]{ASw}. 
\qed\end{proof}

\subsection{The case of lax (double) functors}

For meta-products constructed from non-strict double functors one could still study skew-monoidality of the corresponding category of double categories. 
For the meta-product $\ot^{lx}_{hop}$ of type $(lx, hop)$ we still have \equref{any-strict} and \equref{closedness}
with $*=lx$: 
$$Dbl_{hop}(\A\ot^{lx}_{o-l}\B,\C)\iso q\x Lax_{hop}(\A\times\B,\C)$$
and 
\begin{equation}  \eqlabel{Nikol}
Dbl_{hop}(\A\ot^{lx}_{o-l}\B,\C)\iso Lax_{hop}(\A, Lax_{hop}(\B,\C)).
\end{equation}
However, the first of these isomorphisms shows that $\ot^{lx}_{hop}$ fails to be a ``tight binary map classifier'', since the right-hand side does not consist of strict double quasi-functors. Thus the multicategory $Lax^{lx}_{hop}(-,-)$ fails to be left representable and we can not obtain an associativity constraint $\alpha$ in either direction 
as we fail to have any of the isomorphisms marked in (\ref{eq2}) and (\ref{eq3}).

This problem can be waved by passing to skew-multicategories. This will be clarified in \ssref{sk-m rep}.

\section{Gray skew-multicategories}

For the construction of Gray skew-multicategories for general $\Oo$ the first step to begin with is 
that for two objects $\A,\B\in\Oo$ the set $\Oo^*_\bullet(\A,\B)$ must be an object in $\Oo$. \\
The sort of a skew-multicategory that suits to our construction originates from \cite[Section 3.2]{Bour:SS} and \cite[Section 4]{BL}. 
Namely, in the type of multimaps observed there one finds the least necessary condition to obtain the (one-sided) representability of the 
skew-multicategories that we are going to study. 
Before turning to our tight multimaps, we announce that for the category of loose multimaps we will set 
$q_n\x\Oo^*_\bullet(\Pi_{i=1}^n\A_i,\C)$, the category of multimaps for 
general Gray multicategories where $*\not=st$. In this section all definitions and results will be for $\Oo=G\x\Cat$, 
if otherwise not stated. Nevertheless, we will still use the notation $\Oo$ with the idea that the results that we will expose 
could well work also for higher dimensions and also internal categories. 

We start by defining the category of multimaps for Gray-categories in the first two subsections.

\subsection{Quasi-functors of two variables for Gray-categories} \sslabel{qf 2 Gray}

Quasi-functors for Gray-categories of type $(ps, lx, lx)$ are the ``loose binary maps'' of \cite[Section 4.3]{BL}. 
Namely, these quasi-maps are obtained from pseudo-Gray-functors $\A\to\Ps_{lx}(\B,\C)$ for Gray-categories $\A,\B,\C$, 
where the inner-hom $\Ps_{lx}(\B,\C)$ consists of pseudo-maps, {\em lax} transformations, lax modifications and the corresponding perturbations. 
Pseudo-Gray-functors were defined in \cite{Go} and were called ``pseudo-maps'' in \cite{BL}. 
The inner-hom $\Ps_{lx}(\B,\C)$ was proved to be a Gray-category in \cite{Go} and was denoted by $\bf{Lax}(\B,\C)$ in \cite{BL}. 

\smallskip

We list here the structural data and the labels of axiom they satisfy for binary maps, {\em i.e.} quasi-functors of two variables $H: \A\times\B\to\C$ of type $(*=ps, lx, lx)$ for Gray-categories from \cite[Section 4.3]{BL}. 
They consist of:
\begin{itemize}
\item pseudo Gray functors $H(-,B):\A\to\C$ and $H(A,-):\B\to\C$ s.t. $H(A,-)\vert_B=H(-,B)\vert_A$, 
\item 2-cells $H(f,g)$ that we write $(f,g)$ for short, 
for 1-cells $f:A\to A'$ in $\A$ and $g:B\to B'$ in $\B$, 
\item invertible 3-cells \axiom{$f'f,g$}=$(-,g)^2_{f',f}$, \axiom{$f,g'g$}=$(f,-)^2_{g',g}$ (they correspond to axioms \axiomref{($k,K'K$)} and \axiomref{($k'k,K$)} in double categories, and mean the transformation components at compositions of 1-cells, we call them {\em cocycles}) satisfying the following type of axioms:
\begin{itemize}
\item \axiom{$A_1,A_2,A_3:B$} and \axiom{$A:B_1,B_2,B_3$}: associators for $(f'f,g)$ and $(f,g'g)$, 
\item \axiom{$A:B_1,\beta$}, \axiom{$A:\beta, B_2$}, \axiom{$A_1,\alpha:B$}, \axiom{$\alpha,A_2:B$}: left and right whiskerings ({\em i.e.} naturality with respect to 2-cells) for $(f'f,g)$ and $(f,g'g)$, and 
\item \axiom{$A_1, A_2:B_1,B_2$}: compatibility relating the 3-cells $(f'f,g),(f'f,g'), (f,g'g),(f',g'g)$;
\end{itemize}
\item 3-cells \axiom{$a,g$}=$(-,g)_a$ and \axiom{$f,b$}=$(f,-)_b$ for 2-cells $a$ in $\A$ and $b$ in $\B$ (they correspond to axioms \axiomref{$(k,K)$-r-nat} and \axiomref{$(k,K)$-l-nat} in double categories, and mean the transformation components at 2-cells), satisfying the following type of axioms:
\begin{itemize}
\item $(A_1, A_2:\beta)=(f'f,b)$ and $(\alpha:B_1, B_2)=(a,g'g)$: 
compatibility with the composition of 1-cells, 
\item $(A:\beta_1,\beta_2)=(f,\frac{b}{b'})$ and $(\alpha_1,\alpha_2:B)=(\frac{a}{a'},g)$: 
compatibility with the vertical composition of 2-cells, 
\item $(\Lambda:B)=(\Lambda,g)$ and $(A:\theta)=(f,\Theta)$: 
compatibility with 3-cells $\Lambda$ in $\A$ and $\Theta$ in $\B$,
\item $(\alpha:\beta)=(a,b)$: compatibility between the 3-cells $(a,g)$ and $(f,b)$;
\end{itemize}
\end{itemize}
and satisfy in total 16 degeneracy equations. One should keep in mind that \cite{BL} use pseudo Gray functors which are {\em unitary} or {\em normal}, that is to say that the unitor 2-cells $id_{\F(A)}\Rightarrow\F(id_A)$ are identities. 

\medskip

Since the annotation in diagrams expressing equations among 3-cells is space and time consuming, when annotating 
the axioms 
from the definitions from \cite{BL} and in our definitions of higher cells for binary quasi-functors of Gray-categories, which we will introduce in the next subsection, to annotate equations holding among 3-cells we will be using a concise notation of fractions that we now introduce. (In case of too long fractions, we write them as a horizontal concatenation separated by dots.)  
All these axioms express equalities of transversal compositions of 3-cells. 
Each 3-cell in this composition is indeed a vertical composition of horizontal compositions of ``atomic'' 3-cells, whereby there is a single ``non-trivial'' atomic 3-cell: the rest are identity 3-cells or interchange laws. Let us refer to this single non-trivial 3-cell as 
{\em essential}. At a single level of every fraction we will write a single essential 3-cell, so that the whole fraction replaces the mentioned transversal composition of 3-cells. This way we omit writing identity 3-cells and interchange laws and we simplify substantially the expression. It is important to remark that once one gives a sequence of essential 3-cells (appearing in the axioms and fractions), there is a unique way to recover the full form of the expression, {\em i.e.} to recover the pasting diagrams.

\medskip

Now we may state a Gray-categorical version of \prref{quasi-fun}. 
Observe that in quasi-functors for Gray-categories the following types of structure cells appear, where $i,j\in\{0,1,2,3\}$ denote 
the degree of cells: 
$$(0,0), \quad (1,0), (0,1), \quad (2,0), (1,1), (0,2), \quad (2,1), (1,2), \quad (3,0), (0,3)$$
so that a cell of type $(i,j)$ is a $i+j$-cell. 

\begin{prop} \prlabel{quasi-Gray}
Let $\A,\B,\C$ be Gray-categories. The following are equivalent:
\begin{enumerate} 
\item $H\colon \A\times\B\to\C$ is a quasi-functor of type $(*,\bullet,\bullet)$; 
 \item there are two families of Gray functors of type $*$: 
$$(A,-)\colon\B\to\C\quad\text{ and}\quad (-,B)\colon\A\to\C$$ 
for objects $A\in\A, B\in\B$, such that $(A,-)\vert_B=(-,B)\vert_A=(A,B)$, and so that the following hold: 
\begin{enumerate}[(i)]
\item $(f,-)\colon (A,-)\to(A',-)$ is a $\bullet$-type transformation (of $*$-type Gray functors) for each 1-cell $f\colon A\to A'$,  
$(a,-)$ is a $\bullet$-type modification (of $\bullet$-type transformations) for each 2-cell $a$ in $\A$, while 
$(-,-)^2_{f'f}:\frac{(f,-)}{(f',-)}\Rightarrow(f'f,-)$ is an invertible modification, and 
$(\Lambda,-):(a,-)\Rrightarrow(a',-)$ is a $\bullet$-type perturbation for every 3-cell $\Lambda$ in $\A$; 
\item $(-,g)\colon (-,B)\to (-, B')$ is a $\#$-type transformation (of $*$-type Gray functors) for each 1-cell $g\colon B\to B'$, 
$(-,b)$ is a $\#$-type modification (of $\#$-type transformations) for each 2-cell $b$ in $\B$, while 
$(-,-)^2_{g'g}:\frac{(-,g)}{(-,g')}\Rightarrow(g'g,-)$ is an invertible modification, and 
$(-,\Sigma):(-,b)\Rrightarrow(-,b')$ is a a $\#$-type perturbation for every 3-cell $\Sigma$ in $\B$;
\item the functor values at 2-cells coincide with the modification 2-cell components at objects, and 
the functor values at 3-cells coincide with the perturbation 3-cell components at objects:  
$$(A, -)\vert_b=(-,b)_A, \qquad (-,B)_a=(a,-)_B,$$ 
$$(A, -)\vert_\Sigma=(-,\Sigma)_A, \qquad (-,B)_\Lambda=(\Lambda,-)_B;$$ 
\item the 2-cell components of the respective transformations coincide: 
$$(f,-)\vert_g=(-,g)\vert_f;$$
\item the following 3-cell components of transformations and modifications coincide: 
$$(f,-)\vert_b=(-,b)\vert_f, \qquad (-,g)\vert_a=(a,-)\vert_g$$ 
$${(-,-)^2_{f',f}}\vert_g=(-,g)^2_{f',f} \qquad {(-,-)^2_{g',g}}\vert_f=(f,-)^2_{g',g}.$$
\end{enumerate}
\end{enumerate}
\end{prop}

\begin{proof}
A version of this claim can be found in \cite[Section 4.5]{BL}, with the difference that we added the point (v) above and 
the statement that for all 3-cells 
$\Lambda$ in $\A$ and $\Sigma$ in $\B$ there are perturbations $(\Lambda,-):(a,-)\Rrightarrow(a',-)$ and $(-,\Sigma):(-,b)\Rrightarrow(-,b')$.  
That $(\Lambda,-)$ is a perturbation is the consequence of it being the image of a 3-cell by a pseudo-Gray-functor 
$\F:\A\to G\x\Cat^*_\bullet(\B,\C)$ determining the quasi-functor. In the convention of \cite[Section 4.3]{BL}, where the inner-hom 
in the codomain of $\F$ is $\Ps_{lx}(\B,\C)$, $(\Lambda,-)$ is a {\em lax} perturbation and $(-,\Sigma):(-,b)\Rrightarrow(-,b')$ is an 
{\em oplax} perturbation. The axiom that the latter satisfies corresponds to the axiom \axiomref{$A:\Theta$} from \cite[Appendix D.1]{BL}, which is the part (ii) of $(f,-)$ being a lax transformation (see Appendix C.1). In our above announced fraction notation for transversal composition of 3-cells it reads 
$$\frac{(f,-)_b}{(A',-)_\Lambda}=\frac{(A,-)_\Lambda}{(f,-)_{b'}} \quad \Leftrightarrow \quad 
\frac{(-,b)_f}{(-,\Lambda)_{A'}}=\frac{(-,\Lambda)_A}{(-,b')_f}.$$
Observe that the latter is the axiom \equref{op-pert} 
of an oplax perturbation of oplax modifications in Appendix C.1.
\qed\end{proof}

\bigskip

In \cite{BL} quasi-functors for Gray-categories were considered of the type $(*,\bullet,\bullet)=(ps, lax, lax)$, 
whereas we considered double categorical quasi-functors in \cite{Fem} and in the previous sections 
of the type $(lx, oplax-lax)$, meaning that horizontal transformations were taken oplax and the vertical ones lax. 



In view of \prref{quasi-Gray}, the part d) of our \rmref{dual quasi} and the discussion following it become evident.

\subsection{Gray-category of binary quasi-functors for Gray-categories} \sslabel{Gray-cat of qf}

We now define higher cells for quasi-functors and build a Gray-category $q\x G\x\Cat^{lx}_{oplx}(\A, \B)$. Its objects are lax quasi-functors, and higher cells are oplax transformations, oplax modifications and perturbations thereof. The other versions of this 
Gray-category are constructed in a similar way. The 1-, 2- and 3-cells of $q\x G\x\Cat^{lx}_{oplx}(\A, \B)$ will be induced by pairs of the respective cells in the Gray-category $G\x\Cat^{lx}_{oplx}(\A, \B)$. The proof that the latter is indeed a Gray-category we 
leave for \ssref{lax and strict} (it is for organizational reasons, to make the reading smoother, but the curious reader may jump directly to that proof).

Throughout this subsection the fractions will stand for 
our convention described in the previous subsection for annotation of transversal composition of 3-cells.  
To define oplax transformations of quasi-functors for Gray-categories, the axiom \axiomref{$HOT^q_1$} of \deref{oplax tr cubical} is substituted by a 3-cell fulfilling certain axioms. (The axioms \axiomref{$HOT^q_2$}-\axiomref{$HOT^q_4$} become irrelevant, as they only appear in the internal category case.) We give an oplax quasi-functor version of the definition. This definition is in accordance with the unitarity/normality assumption on Gray functors. Observe that to give the two structure 3-cells \axiomref{$a,g$}=$(-,g)_a$ and \axiomref{$f,b$}=$(f,-)_b$ of a quasi-functor is equivalent to giving a single 3-cell 
$$
\scalebox{0.86}{
\bfig
 \putmorphism(-150,500)(1,0)[(A,B)`(A,B')`(A,g)]{600}1a
 \putmorphism(450,500)(1,0)[\phantom{A\ot B}`(A',B') `(f,B')]{680}1a

 \putmorphism(-150,50)(1,0)[(A, B)`( A,B')`(A,g')]{600}1a
 \putmorphism(450,50)(1,0)[\phantom{A\ot B}`( A',B') `( f',B')]{680}1a
\putmorphism(-180,500)(0,-1)[\phantom{Y_2}``]{450}1l
\putmorphism(-160,500)(0,-1)[\phantom{Y_2}``=]{450}0l
\putmorphism(450,500)(0,-1)[\phantom{Y_2}``]{450}1l
\putmorphism(610,500)(0,-1)[\phantom{Y_2}``=]{450}0l 
\putmorphism(1120,500)(0,-1)[\phantom{Y_3}``=]{450}1r
\put(-40,270){\fbox{$(A,b)$}} 
\put(650,270){\fbox{$(a,B')$}}
\putmorphism(-150,-400)(1,0)[(A, B)`(A',B) `(f', B)]{640}1a
 \putmorphism(450,-400)(1,0)[\phantom{A'\ot B'}` ( A',B') `(A',g')]{680}1a

\putmorphism(-180,50)(0,-1)[\phantom{Y_2}``=]{450}1l
\putmorphism(1120,50)(0,-1)[\phantom{Y_3}``=]{450}1r
\put(300,-200){\fbox{$(f',g')$}}

\efig}
\quad\stackrel{\axiom{$(a,b)\x lr\x nat$}}{\Rrightarrow}
\quad
\scalebox{0.86}{
\bfig
 \putmorphism(-150,500)(1,0)[(A,B)`(A,B')`(A,g)]{600}1a
 \putmorphism(450,500)(1,0)[\phantom{A\ot B}`(A', B') `(f,B')]{680}1a
 \putmorphism(-150,50)(1,0)[(A,B)`(A',B)`(f,B)]{600}1a
 \putmorphism(450,50)(1,0)[\phantom{A\ot B}`(A', B') `(A',g)]{680}1a

\putmorphism(-180,500)(0,-1)[\phantom{Y_2}``=]{450}1r
\putmorphism(1100,500)(0,-1)[\phantom{Y_2}``=]{450}1r
\put(350,260){\fbox{$(f,g)$}}

\putmorphism(-150,-400)(1,0)[(A, B)`(A',B) `(f',B)]{640}1a
 \putmorphism(450,-400)(1,0)[\phantom{A'\ot B'}` (A',B'). `(A',g')]{680}1a

\putmorphism(-180,50)(0,-1)[\phantom{Y_2}``=]{450}1l
\putmorphism(450,50)(0,-1)[\phantom{Y_2}``]{450}1r
\putmorphism(300,50)(0,-1)[\phantom{Y_2}``=]{450}0r
\putmorphism(1100,50)(0,-1)[\phantom{Y_2}``]{450}1r
\putmorphism(1080,50)(0,-1)[\phantom{Y_2}``=]{450}0r
\put(-20,-180){\fbox{$(a,B)$}}
\put(660,-180){\fbox{$(A',b)$}}
\efig}
$$

\begin{defn} \delabel{tr Gray qu-f 2}
Let $\A,\B,\C$ be Gray-categories. 
An oplax transformation $\theta\colon (-,-)_1\Rightarrow (-,-)_2$ between lax quasi-functors $(-,-)_1,(-,-)_2\colon 
\A\times\B\to\C$ of oplax type is given by: for each $A\in\A$ an oplax transformation $\theta^A\colon (A,-)_1\Rightarrow(A,-)_2$ and 
for each $B\in\B$ an oplax transformation $\theta^B\colon (-,B)_1\Rightarrow(-,B)_2$, both of lax functors, such that 
$\theta^A_B=\theta^B_A$ and a 3-cell 
$$\scalebox{0.8}{
\bfig
 \putmorphism(450,700)(1,0)[(A,B')_1 `(A', B')_1 `(f, B')_1]{680}1a
 \putmorphism(1140,700)(1,0)[\phantom{A\ot B}`(A', B')_2 ` \theta^{B'}_{A'}]{680}1a

 \putmorphism(-150,250)(1,0)[(A,B)_1 `(A,B')_1`(A,g)_1]{600}1a
 \putmorphism(450,250)(1,0)[\phantom{A\ot B}`(A,B')_2 `\theta^{B'}_A]{680}1a
 \putmorphism(1130,250)(1,0)[\phantom{A\ot B}`(A', B')_2 ` (f,B')_2]{680}1a

\putmorphism(450,700)(0,-1)[\phantom{Y_2}``=]{450}1r
\putmorphism(1750,700)(0,-1)[\phantom{Y_2}``=]{450}1r
\put(1020,470){\fbox{$ \theta^{B'}_f$}}

 \putmorphism(-150,-200)(1,0)[(A, B)_1`(A, B)_2 `\theta^A_B]{640}1a
 \putmorphism(460,-200)(1,0)[\phantom{A'\ot B'}`(A, B')_2 `(A,g)_2]{680}1a

\putmorphism(-180,250)(0,-1)[\phantom{Y_2}``=]{450}1l
\putmorphism(1120,250)(0,-1)[\phantom{Y_3}``=]{450}1r
\put(310,0){\fbox{$ \theta^A_g$}}

 \putmorphism(1170,-200)(1,0)[\phantom{A\ot B}`(A', B')_2 ` (f, B')_2]{650}1a
\putmorphism(450,-200)(0,-1)[\phantom{Y_2}``=]{450}1r
\putmorphism(1750,-200)(0,-1)[\phantom{Y_2}``=]{450}1r

 \putmorphism(460,-650)(1,0)[(A, B)_2 `(A',B)_2 `(f, B)_2]{650}1a
 \putmorphism(1260,-650)(1,0)[ `(A', B')_2 `(A',g)_2]{600}1a
\put(920,-440){\fbox{$ (f,g)_2$}}
\efig}
\stackrel{\axiom{$\theta, f,g$}}{\Rrightarrow}
\scalebox{0.8}{
\bfig

 \putmorphism(-150,700)(1,0)[(A,B)_1`(A, B')_1`(A,g)_1]{600}1a
 \putmorphism(450,700)(1,0)[\phantom{A\ot B}`(A', B')_1 `(f, B')_1]{680}1a

 \putmorphism(-150,250)(1,0)[(A,B)_1`(A',B)_1`(f, B)_1]{600}1a
 \putmorphism(450,250)(1,0)[\phantom{A\ot B}`(A', B')_1 `(A',g)_1]{680}1a
 \putmorphism(1120,250)(1,0)[\phantom{A'\ot B'}`(A', B')_2 `\theta^{A'}_{B'}]{650}1a

\putmorphism(-180,700)(0,-1)[\phantom{Y_2}``=]{450}1r
\putmorphism(1100,700)(0,-1)[\phantom{Y_2}``=]{450}1r
\put(350,450){\fbox{$(f,g)_1$}}
\put(1000,0){\fbox{$\theta^{A'}_g$}} 

 \putmorphism(-150,-200)(1,0)[(A,B)_1`(A',B)_1`(f, B)_1]{600}1a
 \putmorphism(450,-200)(1,0)[\phantom{A''\ot B'}` (A',B)_2 `\theta^{A'}_B]{680}1a
 \putmorphism(1100,-200)(1,0)[\phantom{A''\ot B'}`(A', B')_2 ` (A',g)_2]{660}1a

\putmorphism(450,250)(0,-1)[\phantom{Y_2}``=]{450}1l
\putmorphism(1750,250)(0,-1)[\phantom{Y_2}``=]{450}1r
\putmorphism(-180,-200)(0,-1)[\phantom{Y_2}``=]{450}1r
\putmorphism(1100,-200)(0,-1)[\phantom{Y_2}``=]{450}1r
 \putmorphism(-150,-650)(1,0)[(A,B)_1`(A,B)_2`\theta^B_A]{600}1a
 \putmorphism(580,-650)(1,0)[`(A',B)_2`(f, B)_2]{540}1a
\put(320,-420){\fbox{$\theta^{B}_f$}}
\efig} \vspace{-0,2cm}
$$
for all 1-cells $f\colon A\to A'$ in $\A$ and $g\colon B\to B'$ in $\B$ satisfying the following axioms 
$$\threefrac{(\theta,f',g')}{(a,b)_1\x lr\x nat}{(\frac{\theta^{A'}_b}{\theta^B_a})}=
\threefrac{(\frac{\theta^{B'}_a}{\theta^A_b})}{(a,b)_2\x lr\x nat}{(\theta,f,g)}$$
for 2-cells $a:f\Rightarrow f'$ and $b:g\Rightarrow g'$, and 
\begin{multline*}
(\theta,f'f,g'g) \cdot \big(\frac{(\theta^{B''})^2_{f'f}}{(\theta^A)^2_{g'g}}\big) \cdot {(f'f,-)^2_2}_{g'g} \cdot 
{(-,g')^2_2}_{f'f} \cdot{(-,g)^2_2}_{f'f}=\\
{(f'f,-)^2_1}_{g'g} \cdot {(-,g')^2_1}_{f'f} \cdot{(-,g)^2_1}_{f'f} \cdot \big(\frac{(\theta^{A''})^2_{g'g}}{(\theta^B)^2_{f'f}}\big) 
\cdot (\theta,f',g) \cdot \big(\frac{(\theta,f',g')}{(\theta,f,g)}\big) \cdot (\theta,f,g') 
\end{multline*}
for composable 1-cells $A\stackrel{f}{\to} A'\stackrel{f'}{\to} A''$ and $B\stackrel{g}{\to} B'\stackrel{g'}{\to} B''$, where the 
$\cdot$ denotes transversal composition of 3-cells and is read from right to left (our fraction notation would occupy too much space for this axiom), and two degeneracy equations: $(\theta, id_A,g)=\Id_{\theta^A_g}$ and $(\theta, f, id_B)=
\Id_{\theta^B_f}$.
\end{defn}

The above two non-degeneracy axioms are equivalent to the following four axioms
\begin{equation} \eqlabel{teta-f-g-left}
\threefrac{(\theta,f',g)}{(a,g)_1\x l\x nat}{\theta^B_a}=\threefrac{\theta^{B'}_a}{(a,g)_2\x l\x nat}{(\theta,f,g)}
\qquad\qquad \threefrac{{(-,g)^2_2}_{f'f}}{(\theta^{B'})^2_{f'f}}{(\theta,f'f,g)}=
\fourfrac{(\theta,f,g)}{(\theta,f',g)}{(\theta^B)^2_{f'f}}{{(-,g)^2_1}_{f'f}}
\end{equation}

\begin{equation} \eqlabel{teta-f-g-right}
\threefrac{(\theta,f,g')}{(f,b)_1\x r\x nat}{\theta^{A'}_b}=\threefrac{\theta^A_b}{(f,b)_2\x r\x nat}{(\theta,f,g)}
\qquad\qquad \threefrac{{(f,-)^2_2}_{g'g}}{(\theta^{A})^2_{g'g}} {(\theta,f,g'g)}=
\fourfrac{(\theta,f,g')}{(\theta,f,g)}{(\theta^{A'})^2_{g'g}}{{(f,-)^2_1}_{g'g}}
\end{equation}
(with the same meanings of $a,b, f, f',g, g'$ as above) 
where we wrote $(a,g)_1\x l\x nat, (a,g)_2\x l\x nat$ for the 3-cell \axiomref{$a,g$} and 
$(f,b)_1\x r\x nat,(f,b)_2\x r\x nat$ for the 3-cell \axiomref{$f,b$} of the two quasi-functors, respectively. 

\smallskip

We record that the cocycles \axiomref{$f'f,g$}=$(-,g)^2_{f',f}$, \axiomref{$f,g'g$}=$(f,-)^2_{g',g}$ in the case of oplax transformations and that are used above, have the following form: 
\begin{equation} \eqlabel{oplax coc}
\scalebox{0.78}{
\bfig
 \putmorphism(450,650)(1,0)[(A,B') `(A,B'') `(A,g')]{680}1a
 \putmorphism(1140,650)(1,0)[\phantom{A\ot B}`(A',B'') ` (f,B'')]{680}1a

 \putmorphism(-150,200)(1,0)[(A,B) `(A,B')`(A,g)]{600}1a
 \putmorphism(450,200)(1,0)[\phantom{A\ot B}`(A', B') `(f,B')]{680}1a
 \putmorphism(1130,200)(1,0)[\phantom{A\ot B}`(A',B'') ` (A',g')]{680}1a

\putmorphism(450,650)(0,-1)[\phantom{Y_2}``=]{450}1r
\putmorphism(1750,650)(0,-1)[\phantom{Y_2}``=]{450}1r
\put(1000,420){\fbox{$ (f,g')$}}

 \putmorphism(-150,-250)(1,0)[(A,B)`(A',B) `(f,B)]{640}1a
 \putmorphism(480,-250)(1,0)[\phantom{A'\ot B'}`(A', B') `(A',g)]{680}1a

\putmorphism(-180,200)(0,-1)[\phantom{Y_2}``=]{450}1l
\putmorphism(1120,200)(0,-1)[\phantom{Y_3}``=]{450}1r
\put(310,-50){\fbox{$ (f,g)$}}

 \putmorphism(1170,-250)(1,0)[\phantom{A\ot B}`(A',B'') ` (A',g')]{650}1a
\putmorphism(450,-250)(0,-1)[\phantom{Y_2}``=]{450}1r
\putmorphism(1750,-250)(0,-1)[\phantom{Y_2}``=]{450}1r

 \putmorphism(480,-700)(1,0)[(A',B) `(A'', B'') `(A',g'g)]{1320}1a
\put(920,-470){\fbox{$ (A',-)_{g'g}$}}
\efig}\quad
\stackrel{(f,-)^2_{g'g}}{\Rrightarrow}
\quad
\scalebox{0.78}{
\bfig
 \putmorphism(-150,500)(1,0)[(A, B) `(A,B')`(A,g)]{580}1a
 \putmorphism(450,500)(1,0)[\phantom{(B, A)} `(A, B'') `(A,g')]{660}1a
\putmorphism(-180,500)(0,-1)[\phantom{Y_2}``=]{450}1r
\put(240,270){\fbox{$ (A,-)_{g'g}$}}

 \putmorphism(-150,50)(1,0)[(A,B)` `(A,g'g)]{1080}1a
 \putmorphism(1080,50)(1,0)[(A,B'')`(A',B'') ` (f,B'')]{680}1a

\putmorphism(1050,500)(0,-1)[\phantom{Y_2}``=]{450}1r

 \putmorphism(-150,-400)(1,0)[(A,B)`(A',B) `(f,B)]{640}1a
 \putmorphism(530,-400)(1,0)[\phantom{Y_2X}`(A',B'') `(A',g'g)]{1220}1a
\put(570,-200){\fbox{$ (f,g'g)$}}

\putmorphism(-180,50)(0,-1)[\phantom{Y_2}``=]{450}1l
\putmorphism(1700,50)(0,-1)[\phantom{Y_3}``=]{450}1r
\efig}
\end{equation}
and
$$
\scalebox{0.78}{
\bfig
 \putmorphism(-150,450)(1,0)[(A,B)`(A,B')`(A,g)]{600}1a
 \putmorphism(450,450)(1,0)[\phantom{A\ot B}`(A', B') `(f,B')]{680}1a

 \putmorphism(-150,0)(1,0)[(A,B)`(A',B)`(f,B)]{600}1a
 \putmorphism(450,0)(1,0)[\phantom{A\ot B}`(A', B') `(A',g)]{680}1a
 \putmorphism(1120,0)(1,0)[\phantom{A'\ot B'}`(A'',B') `(f',B')]{660}1a

\putmorphism(-180,450)(0,-1)[\phantom{Y_2}``=]{450}1r
\putmorphism(1100,450)(0,-1)[\phantom{Y_2}``=]{450}1r
\put(350,210){\fbox{$(f,g)$}}
\put(1000,-250){\fbox{$(f',g)$}}

 \putmorphism(450,-450)(1,0)[\phantom{A''\ot B'}` (A'',B) `(f', B )]{680}1a
 \putmorphism(1100,-450)(1,0)[\phantom{A''\ot B'}`(A'',B') ` (A'',g)]{660}1a

\putmorphism(450,0)(0,-1)[\phantom{Y_2}``=]{450}1l
\putmorphism(1750,0)(0,-1)[\phantom{Y_2}``=]{450}1r
 \putmorphism(-150,-450)(1,0)[(A,B)`(A',B)`(f,B)]{600}1a
\putmorphism(-180,-450)(0,-1)[\phantom{Y_2}``=]{450}1r
\putmorphism(1100,-450)(0,-1)[\phantom{Y_2}``=]{450}1r
 \putmorphism(-150,-900)(1,0)[(A,B)`(A'',B)`(f'f, B)]{1280}1a
\put(260,-670){\fbox{$(-,B)_{f'f}$}}
\efig}\quad
\stackrel{(-,g)^2_{f'f}}{\Rrightarrow}\quad
\scalebox{0.78}{
\bfig
 \putmorphism(450,500)(1,0)[(A,B') `(A', B') `(f,B')]{680}1a
 \putmorphism(1140,500)(1,0)[\phantom{A\ot B}`(A'',B') ` (f',B')]{680}1a

 \putmorphism(-150,50)(1,0)[(A, B) `(A,B')`(A,g)]{600}1a
 \putmorphism(450,50)(1,0)[\phantom{A\ot B}`(A'',B') `(f'f,B')]{1350}1a

\putmorphism(450,500)(0,-1)[\phantom{Y_2}``=]{450}1r
\putmorphism(1750,500)(0,-1)[\phantom{Y_2}``=]{450}1r
\put(880,290){\fbox{$ (-,B')_{f'f}$}}

 \putmorphism(-150,-400)(1,0)[(A, B)`(A'',B) `(f'f, B)]{980}1a
 \putmorphism(780,-400)(1,0)[\phantom{A'\ot B'}`(A'',B') `(A'',g)]{980}1a

\putmorphism(-180,50)(0,-1)[\phantom{Y_2}``=]{450}1l
\putmorphism(1750,50)(0,-1)[\phantom{Y_3}``=]{450}1r
\put(560,-200){\fbox{$ (f'f,g)$}}
\efig}
$$

\medskip

Given three quasi-functors and composable oplax transformations 
$(-,-)_1\stackrel{\theta}{\Rightarrow}(-,-)_2\stackrel{\theta'}{\Rightarrow}(-,-)_3$ their vertical composition 
$\frac{\theta}{\theta'}$ is an oplax transformation $(-,-)_1\Rightarrow(-,-)_3$ with $(\frac{\theta}{\theta'})_A^B=
(\frac{\theta}{\theta'})_B^A=(A,B)_1\stackrel{\theta^B_A}{\to}(A,B)_2\stackrel{\theta'^B_A}{\to}(A,B)_3$ and 
\begin{equation} \eqlabel{comp teta}
(\frac{\theta}{\theta'}, f,g)=\frac{(\theta', f,g)}{(\theta, f,g)}.
\end{equation}
This composition is clearly strictly associative and unital.

\begin{rem} \rmlabel{oplax q.f.-lax tr.}
For later use we record that a {\em lax transformation} of oplax quasi-functors of two variables (with the same structure 2-cells 
$(k,K)_i$ as in the above definition, which now we will denote as $(f,g)^i$, because of later usage) 
but the 2-cells $\theta$ go in the opposite direction) are defined similarly via structure 3-cells 
$$\threefrac{(f,g)^2}{\theta_f}{\theta_g}\stackrel{(\theta,f,g)}{\Rrightarrow}\threefrac{\theta_g}{\theta_f}{(f,g)^1}.$$
Here the fractions present a loose way of denoting a composition of 2-cells as in the above definition. 
\end{rem}

\begin{defn} \delabel{oplax modif}
An oplax modification $\Lambda$ between oplax transformations $\theta\Rrightarrow\theta'\colon (-,-)_1\Rightarrow (-,-)_2$ between lax quasi-functors 
$(-,-)_1,(-,-)_2\colon\A\times\B\to\C$ of Gray-categories is given by a pair of oplax modifications 
$\Lambda^A:\theta^A\Rrightarrow\theta'^A$ and $\Lambda^B:\theta^B\Rrightarrow\theta'^B$ 
acting between oplax transformations 
of functors of Gray-categories such that $\Lambda^A_B=\Lambda^B_A$ for every $A\in\A, B\in\B$ and the following axiom for all 1-cells 
$f:A\to A'$ in $\A$ and $g:B\to B'$ in $\B$ holds 
\begin{equation} \eqlabel{modif qf 2} 
\threefrac{\Lambda^{A}_g}{(\theta,f,g)}{\Lambda^B_f}=\threefrac{\Lambda^{B'}_f}{(\theta',f,g)}{\Lambda^{A'}_g}.
\end{equation}
\end{defn}

Modifications for Gray-categories were defined in \cite{Go}, see also \cite[Section 3.3]{BL}. Observe that there are two types of modifications 
depending on the direction of the defining 3-cell components. We will call {\em lax modifications} among lax transformations the ones with 3-cell components directed as in 
\cite[Section 3.3]{BL}. Namely, the defining 3-cell components of a lax modification $\Lambda:\alpha\Rrightarrow\beta:F\Rightarrow G:\A\to\B$ of 
Gray-categories have the form 
$$
\scalebox{0.82}{
\bfig
 \putmorphism(-150,500)(1,0)[F(A)`G(A)`\alpha(A)]{600}1a

 \putmorphism(-150,50)(1,0)[F(A)`G(A)`\beta(A)]{600}1a
 \putmorphism(450,50)(1,0)[\phantom{A\ot B}`G( B) `G(f)]{680}1a
\putmorphism(-180,500)(0,-1)[\phantom{Y_2}``]{450}1l
\putmorphism(-160,500)(0,-1)[\phantom{Y_2}``=]{450}0l
\putmorphism(450,500)(0,-1)[\phantom{Y_2}``]{450}1l
\putmorphism(610,500)(0,-1)[\phantom{Y_2}``=]{450}0l 
\put(0,270){\fbox{$\Lambda(A)$}} 
\putmorphism(-150,-400)(1,0)[F(A)`F(B) `F(f)]{640}1a
 \putmorphism(450,-400)(1,0)[\phantom{A'\ot B'}` G(B) `\beta(B)]{680}1a

\putmorphism(-180,50)(0,-1)[\phantom{Y_2}``=]{450}1l
\putmorphism(1120,50)(0,-1)[\phantom{Y_3}``=]{450}1r
\put(350,-180){\fbox{$\beta_f$}}
\efig}
\quad\stackrel{\Lambda_f}{\Rrightarrow}\quad
\scalebox{0.82}{
\bfig
 \putmorphism(-150,500)(1,0)[F(A)`G(A)`\alpha(A)]{600}1a
 \putmorphism(450,500)(1,0)[\phantom{A\ot B}`G(B) `G(f)]{680}1a
 \putmorphism(-150,50)(1,0)[F(A)`F(B)`F(f)]{600}1a
 \putmorphism(450,50)(1,0)[\phantom{A\ot B}`G(B) `\alpha(B)]{680}1a

\putmorphism(-180,500)(0,-1)[\phantom{Y_2}``=]{450}1r
\putmorphism(1100,500)(0,-1)[\phantom{Y_2}``=]{450}1r
\put(350,260){\fbox{$\alpha_f$}}

 \putmorphism(450,-400)(1,0)[F(B)` G(B) `\beta(B)]{680}1a

\putmorphism(450,50)(0,-1)[\phantom{Y_2}``]{450}1r
\putmorphism(300,50)(0,-1)[\phantom{Y_2}``=]{450}0r
\putmorphism(1100,50)(0,-1)[\phantom{Y_2}``]{450}1r
\putmorphism(1080,50)(0,-1)[\phantom{Y_2}``=]{450}0r
\put(660,-160){\fbox{$\Lambda(B)$}}
\efig}
$$
for every 1-cell $f:A\to B$ in $\A$. In contrast, oplax modifications $\Lambda:\alpha\Rrightarrow\beta:F\Rightarrow G:\A\to\B$ among oplax transformations of Gray-categories (appearing in \deref{oplax modif}) are given via their defining 3-cell components that have the form 
$$
\scalebox{0.82}{
\bfig
 \putmorphism(-150,500)(1,0)[F(A)`F(B)`F(f)]{600}1a
 \putmorphism(450,500)(1,0)[\phantom{A\ot B}`G(B) `\alpha(B)]{680}1a
 \putmorphism(-150,50)(1,0)[F(A)`G(A)`\alpha(A)]{600}1a
 \putmorphism(450,50)(1,0)[\phantom{A\ot B}`G(B) `G(f)]{680}1a

\putmorphism(-180,500)(0,-1)[\phantom{Y_2}``=]{450}1r
\putmorphism(1100,500)(0,-1)[\phantom{Y_2}``=]{450}1r
\put(350,260){\fbox{$\alpha_f$}}

 \putmorphism(-150,-400)(1,0)[F(B)` G(B) `\beta(A)]{630}1a

\putmorphism(450,50)(0,-1)[\phantom{Y_2}``]{450}1r
\putmorphism(300,50)(0,-1)[\phantom{Y_2}``=]{450}0r
\putmorphism(-180,50)(0,-1)[\phantom{Y_2}``=]{450}1r
\putmorphism(-160,50)(0,-1)[\phantom{Y_2}``]{450}0r
\put(-30,-160){\fbox{$\Lambda(A)$}}
\efig}\quad\stackrel{\Lambda_f}{\Rrightarrow}\quad
\scalebox{0.82}{
\bfig
 \putmorphism(450,500)(1,0)[F(B)`G(B)`\alpha(B)]{600}1a

 \putmorphism(-150,50)(1,0)[F(A)`F(B)`F(f)]{600}1a
 \putmorphism(420,50)(1,0)[\phantom{A\ot B}`G( B) `\beta(B)]{700}1a
\putmorphism(450,500)(0,-1)[\phantom{Y_2}``=]{450}1l
\putmorphism(1100,500)(0,-1)[\phantom{Y_2}``=]{450}1l
\put(660,270){\fbox{$\Lambda(B)$}} 
\putmorphism(-150,-400)(1,0)[F(A)`G(A) `\beta(A)]{640}1a
 \putmorphism(450,-400)(1,0)[\phantom{A'\ot B'}` G(B). `G(f)]{680}1a

\putmorphism(-180,50)(0,-1)[\phantom{Y_2}``=]{450}1l
\putmorphism(1120,50)(0,-1)[\phantom{Y_3}``=]{450}1r
\put(370,-180){\fbox{$\beta_f$}}
\efig}
$$
For the sake of completeness of the exposition we refer the reader to Appendix C.1 where we 
give in our fraction notation for 3-cells the defining axioms of the notions of a lax transformation, (lax) modification and (lax) perturbation from \cite[Section 3.3]{BL}, but also our definitions of an oplax transformation, oplax modification, and oplax perturbation. 

\medskip

The horizontal composition of 2-cells in $q\x G\x\Cat^{lx}_{oplx}(\A, \B)$ is induced on 2-cell components by the horizontal composition of the corresponding 2-cells: 
$$[\Lambda\vert\Lambda'](A)=
\scalebox{0.86}{
\bfig
\putmorphism(-100,250)(1,0)[F(A)`G(A)`\theta_1(A)]{550}1a
 \putmorphism(430,250)(1,0)[\phantom{F(A)}`H(A) `\theta_2(A)]{580}1a

 \putmorphism(-100,-200)(1,0)[F(A)`G(A)`\theta_1'(A)]{550}1a
 \putmorphism(450,-200)(1,0)[\phantom{F(A)}`H(A) `\theta_2'(A)]{580}1a

\putmorphism(-100,250)(0,-1)[\phantom{Y_2}``=]{450}1l
\putmorphism(420,250)(0,-1)[\phantom{Y_2}``]{450}1r
 \putmorphism(400,250)(0,-1)[\phantom{Y_2}``=]{450}0r
\putmorphism(1020,250)(0,-1)[\phantom{Y_2}``=]{450}1r
\put(50,40){\fbox{$\Lambda_A$}}
\put(680,40){\fbox{$\Lambda_A' $}}
\efig}
$$
and on component 3-cells by the transversal composition of the corresponding 3-cells:
$[\Lambda\vert\Lambda']_f=\frac{\Lambda_f}{\Lambda'_f}$.

Similarly, the vertical composition of modifications is induced on component 2-cells by the vertical composition of the corresponding 2-cells:
$$(\frac{\Lambda}{\Sigma })(A)=
\scalebox{0.86}{
\bfig
 \putmorphism(-150,470)(1,0)[F(A)`G(A)  `\theta(A)]{440}1a
\putmorphism(-160,470)(0,-1)[\phantom{Y_2}`F(A) `=]{400}1l
\putmorphism(-160,80)(0,-1)[\phantom{Y_2}`F(A)`=]{400}1l
\putmorphism(300,80)(0,-1)[\phantom{Y_2}`G(A)`=]{400}1r
\putmorphism(300,470)(0,-1)[\phantom{Y_2}`G(A) `=]{400}1r
\put(-20,280){\fbox{$\Lambda_A$}}

\putmorphism(-180,70)(1,0)[\phantom{F(A)}``\theta'(A)]{380}1a
\putmorphism(-100,-300)(1,0)[\phantom{Y_2}` `\theta''(A)]{300}1b
\put(-20,-130){\fbox{$\Sigma_A$}}
\efig}
$$ 
and on component 3-cells by the transversal composition of the corresponding 3-cells:
$(\frac{\Lambda}{\Sigma})_f=\frac{\Lambda_f}{\Sigma_f}$.

It is clear that the associativity and unitality of 2-cells in both horizontal and vertical direction hold strictly.

\medskip

We also define perturbations for quasi-functors on Gray-categories.

\begin{defn}
A perturbation $\Theta:\Lambda\to\Lambda':\theta\Rrightarrow\theta'\colon (-,-)_1\Rightarrow (-,-)_2$ for quasi-functors 
$(-,-)_1,(-,-)_2\colon\A\times\B\to\C$ of type $(*,\bullet)$ of Gray-categories, where 
$\Lambda=(\Lambda^A, \Lambda^B)_{A\in\A, \\ B\in\B}, \theta=(\theta^A, \theta^B)_{A\in\A, \\ B\in\B}$ and similarly 
for $\Lambda'$ and $\theta'$, 
is given by a pair of perturbations $\Theta^A:\Lambda^A\Rrightarrow\Lambda'^A$ and $\Theta^B:\Lambda^B\Rrightarrow\Lambda'^B$ 
acting between modifications 
for Gray-categories, such that $\Theta^A_B=\Theta^B_A$ for every $A\in\A, B\in\B$. 
\end{defn}

Three types of compositions of perturbations, as 3-cells in $q\x G\x\Cat^{lx}_{oplx}(\A, \B)$, are induced by the respective 
compositions of the component perturbations in the Gray-category $G\x\Cat^{lx}_{oplx}(\A, \B)$.

An interchange in $q\x G\x\Cat^{lx}_{oplx}(\A, \B)$ among modifications $\Lambda=(\Lambda^A, \Lambda^B)_{A\in\A, \\ B\in\B}$ and 
$\Sigma=(\Sigma^A, \Sigma^B)_{A\in\A, \\ B\in\B}$ is given by the perturbation $\Int_{\Lambda,\Sigma}:=(\Int_{\Lambda_A,\Sigma_A}, 
\Int_{\Lambda_B,\Sigma_B})_{A\in\A, \\ B\in\B}$, where $\Int_{\Lambda_A,\Sigma_A}$ and $\Int_{\Lambda_B,\Sigma_B}$ are perturbation  
3-cell components from \rmref{interch}. This finishes the proof that $q\x G\x\Cat^{lx}_{oplx}(\A, \B)$is a Gray-category.

\subsection{Quasi-functors of $n$-variables for Gray-categories} \sslabel{incub}

Because of the rise in dimension, before defining quasi-functors of $n$-variables for Gray-categories, unlike for double categories, 
one has to know what quasi-functors of three variables $H: \A\times\B\times\C\to\D$ for Gray-categories are. They were defined in the type 
$\bullet=lx$ in \cite[Section 4.10]{BL} as ternary multimaps. We give a brief description of them in the type $\bullet=oplx$. For this reason 
we will accordingly differentiate {\em lax incubators} and {\em oplax incubators}, the latter we show below. 
Quasi-functors of three variables of type $(*,\bullet)$ consist of quasi-functors of two variables ({\em i.e. binary maps}) of 
type $(*,\bullet)$ 
$$H(A,-,-): \B\times\C\to\D, \quad H(-,B,-): \A\times\C\to\D, \quad H(-,-,C): \A\times\B\to\D \vspace{-0,3cm}$$ 
s.t. \vspace{-0,3cm}
\begin{enumerate} [(i)]
\item 
\begin{equation} \eqlabel{obj}
H(A,-,-)\vert_B=H(-,B,-)\vert_A, H(A,-,-)\vert_C=H(-,-,C)\vert_A, H(-,B,-)\vert_C=H(-,-,C)\vert_B, 
\end{equation}
\item 
the structure 2-cells of type $\bullet$ (we draw them for the type $\bullet=oplax$) of the binary maps relate via the {\em oplax incubator} 3-cells 
$$\scalebox{0.8}{
\bfig
 \putmorphism(450,700)(1,0)[` `(A,g,C')_2]{680}1a
 \putmorphism(1140,700)(1,0)[` ` (f,B',C')_1]{650}1a

 \putmorphism(-150,250)(1,0)[``(A,B,h)_3]{600}1a
 \putmorphism(450,250)(1,0)[` `(f,B,C')_1]{680}1a
 \putmorphism(1130,250)(1,0)[` ` (A',g, C')_2]{650}1a

\putmorphism(450,700)(0,-1)[\phantom{Y_2}``=]{450}1r
\putmorphism(1750,700)(0,-1)[\phantom{Y_2}``=]{450}1r
\put(940,480){\fbox{$ (f,g,C')_{12}$}}

 \putmorphism(-150,-200)(1,0)[``(f,B,C)_1]{640}1a
 \putmorphism(460,-200)(1,0)[` `(A',B,h)_3]{680}1a

\putmorphism(-160,250)(0,-1)[\phantom{Y_2}``=]{450}1l
\putmorphism(1120,250)(0,-1)[\phantom{Y_3}``=]{450}1r
\put(310,20){\fbox{$(f,B,h)_{13}$}}

 \putmorphism(1140,-200)(1,0)[` ` (A',g,C')_2]{630}1a
\putmorphism(450,-200)(0,-1)[\phantom{Y_2}``=]{450}1r
\putmorphism(1750,-200)(0,-1)[\phantom{Y_2}``=]{450}1r

 \putmorphism(460,-650)(1,0)[` `(A',g,C)_2]{700}1a
 \putmorphism(1140,-650)(1,0)[ ` `(A',B',h)_3]{630}1a
\put(860,-430){\fbox{$ (A',g,h)_{23}$}}
\efig}
\stackrel{\axiom{f,g,h}}{\Rrightarrow}
\scalebox{0.8}{
\bfig

 \putmorphism(-150,700)(1,0)[``(A,B,h)_3]{600}1a
 \putmorphism(450,700)(1,0)[` `(A,g,C')_2]{680}1a

 \putmorphism(-150,250)(1,0)[``(A,g,C)_2]{600}1a
 \putmorphism(450,250)(1,0)[` `(A,B',h)_3]{680}1a
 \putmorphism(1120,250)(1,0)[` `(f,B',C')_1]{650}1a

\putmorphism(-180,700)(0,-1)[\phantom{Y_2}``=]{450}1r
\putmorphism(1100,700)(0,-1)[\phantom{Y_2}``=]{450}1r
\put(280,470){\fbox{$(A,g,h)_{23}$}}
\put(910,20){\fbox{$(f,B',h)_{13}$}}

 \putmorphism(-150,-200)(1,0)[``(A,g,C)_2]{600}1a
 \putmorphism(450,-200)(1,0)[` `(f,B',C)_1]{680}1a
 \putmorphism(1100,-200)(1,0)[` ` (A',B',h)_3]{660}1a

\putmorphism(450,250)(0,-1)[\phantom{Y_2}``=]{450}1l
\putmorphism(1750,250)(0,-1)[\phantom{Y_2}``=]{450}1r
\putmorphism(-180,-200)(0,-1)[\phantom{Y_2}``=]{450}1r
\putmorphism(1100,-200)(0,-1)[\phantom{Y_2}``=]{450}1r
 \putmorphism(-150,-650)(1,0)[``(f,B,C)_1]{600}1a
 \putmorphism(580,-650)(1,0)[`` (A',g,C)_2]{540}1a
\put(280,-420){\fbox{$(f,g,C)_{12}$}}
\efig} 
$$
in $\C$,  
for all 1h-cells $(f,g,h):(A,B,C)\to(A',B', C')$ in $\A\times\B\times\C$, which satisfy six axioms that can symbolically be written as follows
$$(f'f,g,h), \,\ (f, g'g,h), \,\, (f,g, h'h),$$
$$(\alpha,g,h),  \,\ (f,\beta,h), \,\, (f,g, \gamma)$$
and three degeneracy equations 
$$(1_A,g,h)=1_{(A,g,h)}, \,\ (f, 1_B,h)=1_{(f,B,h)}, \,\, (f,g, 1_C)=1_{(f,g,C)}.$$
The first six axioms express compatibility of the incubators with compositions of 1-cells in each variable and with 2-cells 
$(\alpha,\beta,\gamma)$ in $\A\times\B\times\C$. 
\end{enumerate}
(The 1-, 2- and 3-cells $(a,b,c)_l, (a,b,c)_{pq}, l,p,q\in\{1,2,3\}$ and $(f,g,h)$ in $\C$ appearing in our diagram above, for 1-, 2- and 3-cells $(a,b,c)\in\A\times\B\times\C$ respectively, are short annotations for $H(a,b,c)$.) 
The above nine axioms are analogous to those in \cite[Appendix D.2]{BL}, where they are drawn for the type $(lx,lx,lx)$). There they 
were denoted as follows 
$$\axiom{$A_1, A_2:B:C$}, \,\ \axiom{$A:B_1, B_2:C$}, \,\ \axiom{$A:B:C_1,C_2$}$$ 
$$\axiom{$\alpha:B:C$}, \,\ \axiom{$A:\beta: C$}, \,\ \axiom{$A:B:\gamma$}$$ 
$$\axiom{$1_a:B:C$}, \,\ \axiom{$A:1_b:C$}, \,\ \axiom{$A:B:1_c$}.$$ 
The first six axioms of this lax version of incubators we list in the Appendix C.2 using our fractions for transversal composition of 3-cells. Here we write them for oplax incubators for later use: 

\axiom{$f'f,g,h$} \qquad $\threefrac{(-,B,h)^2_{f',f}}{(-,g,C')^2_{f',f}}{(f'f,g,h)}=
\fourfrac{(f',g,h)}{(f,g,h)}{(-,g,C)^2_{f',f}}{(-,B',h)^2_{f',f}}$ \qquad\qquad 
\axiom{$a,g,h$} \qquad $\threefrac{(-,g,C')_a}{(-,B,h)_a}{(f,g,h)}=\threefrac{(f',g,h)}{(-,B',h)_a}{(-,g,C)_a}$

\axiom{$f, g'g,h$} \qquad $\threefrac{(A',-,h)^2_{g',g}}{(f,-,C')^2_{g',g}}{(f,g'g,h)}=
\fourfrac{(f,g,h)}{(f,g',h)}{(f,-,C)^2_{g',g}}{(A,-,h)^2_{g',g}}$ \qquad\qquad 
\axiom{$f,b,h$} \qquad $\threefrac{(f,-,C')_b}{(A',-,h)_b}{(f,g,h)}=\threefrac{(f,g',h)}{(A,-,h)_b}{(f,-,C)_b}$

\axiom{$f,g, h'h$}  \qquad $\threefrac{(A',g,-)^2_{h',h}}{(f,B,-)^2_{h',h}}{(f,g,h'h)}= 
\fourfrac{(f,g,h')}{(f,g,h)}{(f,B',-)^2_{h',h}}{(A,g,-)^2_{h',h}}$ \qquad\qquad 
\axiom{$f,g, c$} \qquad $\threefrac{(f,g,h')}{(A,g,-)_c}{(f,B',-)_c}=\threefrac{(f,B,-)_c}{(A',g,-)_c}{(f,g,h)}$. 

\smallskip

\noindent Observe that for oplax incubators we use oplax cocycles \equref{oplax coc}.

\medskip

We say that $H:\A\times\B\times\C\to\D$ is of type $((st,*),\bullet, \star)$ if so are the binary quasi-functors 
$H(-,B,-): \A\times\C\to\D$ and $H(-,-,C): \A\times\B\to\D$.

\smallskip

Observe that according to \prref{quasi-Gray}, 
in each constituting 2-cell $(f,g,C), (f,B,h), (A,g,h)$ in the (co)domain of the incubator 3-cell $(f,g, h)$ in the above diagram one has that $(f,-,C), (f,B,-), (A,g,-)$ are oplax transformations, while $(-,g,C), (-,B,h), (A,-,h)$ are lax transformations.

\smallskip

We will also need transformations among quasi-functors of 3-variables for Gray-categories.

\begin{defn} \delabel{tr of qf 3}
A transformation $\theta: H_1\Rightarrow H_2$ of type $\bullet$ between quasi-functors $H_1,H_2:\A\times\B\times\C\to\E$ of type 
$(*,oplax)$ of Gray-categories consists of transformations of type $\bullet$
$$\theta^A: H_1(A,-,-) \Rightarrow H_2(A,-,-),$$ 
$$\theta^B: H_1(-,B,-) \Rightarrow H_2(-,B,-),$$ 
$$\theta^C: H_1(-,-,C) \Rightarrow H_2(-,-,C)$$ 
of quasi-functors of two variables of type $(*,oplax)$, which satisfy the axiom on the left below if the transformations are 
$\bullet=oplx$ (as in \deref{tr Gray qu-f 2}), and the axiom on the right below if the transformations are $\bullet=lx$ (this is obtained by reading the diagrams for (co)domain 2-cells of our oplax transformations in \deref{tr Gray qu-f 2} from bottom to top, and then the structure 3-cell goes from right to the left):
\begin{equation} \eqlabel{tr of qf 3}
\fourfrac{(\theta^C,f,g)}{(\theta^B,f,h)}{(\theta^A,g,h)}{(f,g,h)^1_{oplx-inc}}=
\fourfrac{(f,g,h)^2_{oplx-inc}}{(\theta^A,g,h)}{(\theta^B,f,h)}{(\theta^C,f,g)} 
\qquad \qquad
\fourfrac{(f,g,h)^1_{oplx-inc}}{(\theta^C,f,g)}{(\theta^B,f,h)}{(\theta^A,g,h)}=
\fourfrac{(\theta^A,g,h)}{(\theta^B,f,h)}{(\theta^C,f,g)}{(f,g,h)^2_{oplx-inc}}.
\end{equation}
%
\end{defn}

(The above axioms are written in our loose manner, omitting identity 3-cells and interchange laws.) 

\begin{rem} \rmlabel{obs frac}
Mind that since here we take quasi-functors of $\bullet=oplax$ type, then the incubator 3-cells $(f,g,h)^1_{oplx-inc}$ of $H_1$ and 
$(f,g,h)^2_{oplx-inc}$ of $H_2$ are oplax as in \axiomref{f,g,h}. However, {\em lax} transformations between quasi-functors of type 
$(*,lax)$ will satisfy {\em an axiom of the exactly the same form} as the right one in \equref{tr of qf 3}. Namely, although in that case the incubator cells will be of lax type, and one might expect a permutation of the fraction in \equref{tr of qf 3}, what happens is that in the 
(co)domain 2-cells among the 3-cells in the fraction 
the 2-cells $\theta$  
will compose in a different order, but the order of the {\em actions} dictated by the 3-cells in the fractions will remain the same.  
\end{rem}

\smallskip

Observe that it is $\theta^{A;B}(C)=\theta^{B;C}(A)=\theta^{A;C}(B)$ for all $(A,B,C)\in\A\times\B\times\C$ in the above definition.

\medskip

Other versions of transformations for quasi-functors include: 
\begin{enumerate} 
\item {\em oplaxs} transformations between quasi-functors of type $(*,lx,\star)$ - when 2-cells $(f,g,C)^i$, and so on, go in the opposite direction; 
\item {\em laxs} transformations between quasi-functors of type $(*,lx,\star)$ - when both 2-cells $(f,g,C)^i$ and $\theta^{A;C}_g$, and so on, go in the opposite direction; 
\item pseudo and strict versions of transformations and quasi-functors - in the obvious way, 
\item but also mixed version where univalent components of the transformation are of different types: one speaks about $(x,y,z)$-type where 
$x,y,z$ are any of the adjectives. 
\end{enumerate}

\bigskip

We finally introduce quasi-functors of $n$-variables for Gray-categories.

\begin{defn} 
Let $\A_1, ...,\A_n,\C$ be Gray-categories. 
A {\em quasi-functor of $n$-variables} $H:\A_1\times...\times\A_n\to\C$ for $n\geq 4$ of type $(*,\bullet)$ consists of 
quasi-functors of three variables of type $(*,\bullet, \star)$
$$H(A_1,...,A_{i-1},\, -,\, A_{i+1},...,A_{j-1}, \, -,\, A_{j+1},...,A_{k-1}, \, -,\, A_{k+1},..., A_n): \A_i\times\A_j\times\A_k\to\C$$
for all $i<j<k$ and all choices of objects $A_l\in\A_l, l=1,...,n$, 
so that 
\begin{itemize}
\item omitting the irrelevant variables and reducing to four variables $A_i=A,A_j=B,A_k=C,A_l=D$, the induced binary maps coincide:
\begin{equation} \eqlabel{obj 4-ary}
H(A,-,-,-)\vert_B=H(-,B,-,-)\vert_A, H(A,-,-,-)\vert_C=H(-,-,C,-)\vert_A,  \vspace{-0,32cm}
\end{equation}
$$H(-,B,-,-)\vert_C=H(-,-,C,-)\vert_B, H(A,-,-,-)\vert_D=H(-,-,-,D)\vert_A, \vspace{-0,16cm}$$
$$H(-,B,-,-)\vert_D=H(-,-,-,D)\vert_B, H(-,-,-,D)\vert_C=H(-,-,-,D)\vert_C. $$
and 
\item the {\em mecon} axiom among the incubator 3-cells of type $\bullet$ holds, which states that the only two possible ways to compose incubator 3-cells determined by 1-cells $(f,g,h,k):(A,B,C,D)\to(A',B',C',D')\in\A_i\times\A_j\times\A_k\times\A_l$ coincide. 
\end{itemize}
\end{defn}

In a loose manner we may write the mecon axiom as follows
\begin{equation} \eqlabel{mecon}
\fourfrac{(f,g,h,D)}{(f,g,C',k)}{(f,B,h,k)}{(A',g,h,k)}=\fourfrac{(A,g,h,k)}{(f,B',h,k)}{(f,g,C,k)}{(f,g,h,D')}.
\end{equation}
Recall, as in \ssref{qf 2 Gray}, that fractions here present the transversal composition of 3-cells, whereas some identity 3-cells and interchange laws are omitted to simplify the expression. Note that 
on the left-hand side the fixed objects along the composition are placed in the counter-diagonal, while on the right-hand side they form a diagonal. In Appendix D.3 of \cite{BL} one can find the diagrams describing the mecon axiom for incubators of type $\bullet=lx$. Interestingly, 
the loose form of the mecon axiom above is exactly the same both for lax and oplax incubators, a similar occurrence we saw in \rmref{obs frac}.

\medskip

Thus a quasi-functor of $n$-variables for Gray-categories amounts to quasi-functors of four variables 
$\A_i\times\A_j\times\A_k\times\A_l\to\C$ for all $i,j,k,l\in\{1,...,n\}$. 

\smallskip

Given a quasi-functor of $n$-variables $\F$ its component quasi-functors of $k$-variables with $k<n$ we will sometimes call 
{\em restrictions of $\F$}.

\subsection{Tight multimaps: binary case} \sslabel{tight}

We will differentiate left and right skew-multicategories, depending on whether the tightness of multimaps is required on the left 
or on the right. In the case of left skew-multicategories, 
for tight multimaps we will take those quasi-functors that are {\em strict only in the leftmost entry and loose (of type $*$) 
in all the remaining entries}. We will denote the category of such tight multimaps as below on the left, while the isomorphism on 
the right will follow by iteration from the result of \thref{mixed-quasi n}.
$$
q_n\x\Oo^{(st,*)}_\bullet(\Pi_{i=1}^n\A_i,\C)=q_n\x\Oo^{(st,*)}_\bullet(\A_1{}^{st}\w\times^*\Pi_{i=2}^n\A_i,\C)
\iso\Oo^{st}_\bullet(\A_1,\Oo^*_\bullet(\A_2,...,\Oo^*_\bullet(\A_n,\C)...).
$$
Here $*$ stands for some loose type of functors. The meaning of the supra-index $(st,*)$ alludes to the combination of a strict 
functor (leftmost entry) and functors of loose type $*$ in the rest of the entries making a quasi-functor of $n$ variables of a mixed type. 
For $n=2$ we can already say that such a quasi-functor of two variables of mixed type is 
obtained from a strict functor $\A\to\Oo^*_\bullet(\B,\C)$ into the inner-hom taken with loose functors of type $*$, transformations and modifications of types $\bullet$ and $\star$, respectively. This induces bijections 
\begin{equation} \eqlabel{mixed-quasi-bij}
q\x\Oo^{(st,*)}_\bullet(\A{}^{st}\w\times^*\B,\C) \iso\Oo_{st}(\A,\Oo^*_\bullet(\B,\C))\iso\Oo_*(\B,\Oo^{st}_\#(\A,\C))
\end{equation} 
that lifts to natural isomorphisms in $\Oo$
\begin{equation} \eqlabel{mixed-quasi}
q\x\Oo^{(st,*)}_\bullet\vert_{(\#,\bullet)}(\A{}^{st}\w\times^*\B,\C) \iso\Oo^{st}_\#(\A,\Oo^*_\bullet(\B,\C))
\iso\Oo^*_\bullet(\B,\Oo^{st}_\#(\A,\C)).
\end{equation} 

In other words, we have:

\begin{defn} \delabel{skew-guasi}
Let $\Oo$ be any of the categories $\V_m, \I_{m-1}$ with $m\geq 2$ or $G\x\Cat$, and let $\A,\B,\C\in\Oo$. 
A {\em general quasi-functor (of two variables)} $H: \A\times\B\to\C$ of $(st,*)$-mixed type, {\em i.e.} of type $((st,*),\bullet)$, 
in $\Oo$ consists of: 
\begin{enumerate}
\item two families 
$$(-,A):\B\to\C, \qquad (B,-):\A\to\C$$ 
of functors of dimension $m$, where $(-,A)$ is of type $*$ and $(B,-)$ is strict, indexed by objects $A\in \A$ and $B\in\B$, 
\item a single 2-cell in the enriched and Gray-category case, {\em i.e.} four 2-cells in the internal case, 
we will refer to these 2-cells as ``structure 2-cells of type $\bullet$'';
\item if $m=3$, 
four 3-cells in the case of pseudo functors, {\em i.e.} six 3-cells in the case of (op)lax functors, both in the enriched case, 
and in the internal case additional six 3-cells in the pseudo case {\em i.e.} additional eight 3-cells in the case of (op)lax functors 
(the additional 3-cells in the internal case involve vertical 1-cells); \\ 
for $m>3$ corresponding $k$-cells for $k=4,..,m$, 
\end{enumerate}
so that the above data fulfill the corresponding axioms making an $m$-dimensional strict functor $\A\to\Oo^*_\bullet(\B,\C)$. 
(The axioms depend on the type $\bullet$ of transformations and higher morphisms considered in the inner-hom.) 
\end{defn}

\medskip

The following should hold for general $\Oo$, we prove it for $\Oo=G\x\Cat$ (for $\Oo=Dbl$ it follows easily from \prref{ev} 
and  \thref{strict-left closed}).

\begin{prop} \prlabel{ev-mixed}
Let $\Oo=G\x\Cat$. 
There is a quasi-functor $ev:\Oo^*_{\bullet}(\B, \C){}^{st}\w\times^*\B\to\C$ of type $((st,*),\bullet)$ such that given 
any quasi-functor $H:\A{}^{st}\w\times^*\B\to\C$ of type $((st,*),\bullet)$ it is 
$$ev((H^t){}^{st}\w\w\times^*\Id_\B)=H$$
where $H^t:\A\to \Oo^*_{\bullet}(\B, \C)$ is the functor corresponding to $H$ by \equref{mixed-quasi-bij}.
\end{prop}

\begin{proof}
The definition of $ev$ is analogous as in \prref{ev}. We recall the notation of cells and complete it for $m=3$ in Table 
\ref{table:1,5}
\begin{table}[H]
\begin{center}
\begin{tabular}{ c c c } 
degree of cells in: &  $\Oo^*_{\bullet}(\B, \C)$ & $\B$ \\ [0.5ex]
\hline
0 & $F$ & $B$ \\ [0.2ex]   
1 & $\alpha$ & $g$ \\ [0.2ex]   
2 & $b$ & $\beta$ \\ [0.2ex]   
3 & $\Sigma$ & $\Lambda$ \\ [0.1ex]   
\end{tabular}
\caption{Notations for cells}
\label{table:1,5}
\end{center}
\end{table}
\noindent and we add that $ev(\alpha,\beta)=\alpha_\beta, ev(b,g)=b_g, ev(\Lambda, B)=\Lambda(B), 
ev(F,\Sigma)=F(\Sigma)$ are the corresponding structure 3-cells of $\C$. The equation relating $H$ and $H^t$ is clear: 
$H^t(a)(b)=H(a,b)$ for cells $a\in\A, b\in\B$ of degrees $deg(a)+deg(b)\leq 3$. 

We already observed in the proof of \prref{ev} that the left evaluation is strict in the left variable. 
This holds also in the Gray-categorical case: in \cite{BL} strict Gray functors $ev_B=ev(-,B): G\x\Cat^{ps}_{lx}(\B, \C)\to\C$ 
were used. For $F\in\Oo^*_{\bullet}(\B, \C)$ we find that $ev(F,-)$ is clearly a Gray functor of type $*$, as so is $F$. 
Point (i) of \prref{quasi-Gray} clearly holds. 

For (ii) of \prref{quasi-Gray}, we present in Table \ref{table:2} why $ev(-,g)$ is a lax transformation, 
in Table \ref{table:3} why $ev(-,\beta)$ is a lax modification, and in Table \ref{table:4} why $ev(-,\Lambda)$ is a perturbation, 
respectively. 

\begin{table}[H]
\begin{center}
\begin{tabular}{ c c } 
transf. axiom for $ev(-,g)$ & its meaning on cells in $\Oo^*_{\bullet}(\B, \C)$ on which evaluated \\ [0.5ex]
\hline
$(i)$ deg. & clearly holds - identity transf. $ev(1_B,-)$ \\ [1ex]   
$(ii)$  & naturality of 3-cells $b_g$ in $b$ {\em i.e.} pert. ax. for $ev(\Sigma,-)$ \\ [1ex]    
$(iii)$ & composition of modification 3-cells $(b'\cdot b)_g$ at $g$ \\ [1ex]   
$(iv)$ deg. & clearly holds - identity modif. 3-cell $ev(\Id_\alpha,g)=\Id_\alpha\w\vert_g$ at $g$ \\ [1ex]   
$(v)$ & associativity of composition of transf. $\alpha,\beta,\gamma$ at $g$ \\ [1ex]   
$(vi)$ deg. & clearly holds - unitality of comp. of transformations at $g$ \\ [1ex]   
$(vii)$ & naturality in both coord. of comp. of transformations at $g$ \\ [1ex]   
\end{tabular}
\caption{Why $ev(-,g)$ is a transformation}
\label{table:2}
\end{center}
\end{table}

\begin{table}[H]
\begin{center}
\begin{tabular}{ c c } 
modif. axiom for $ev(-,\beta):$ & \\ [0.5ex]
$ev(-,g)\Rrightarrow ev(-,g')$ & its meaning on cells in $\Oo^*_{\bullet}(\B, \C)$ on which evaluated \\ [0.5ex]
\hline
$(i)$ deg. & deg. for transf. $\alpha$ at $1_B$ \\ [1ex]   
$(ii)$  & how transf. 3-cells $\alpha_{\frac{\beta}{\beta'}}$ are given at 1-cells $\beta,\beta'$ \\ [1ex]    
 & {\em i.e.} (iii) of transf. $\alpha$ \\ [1ex] 
$(iii)$ & naturality of modif. 3-cell $b_g$ in $g$ {\em i.e.} (iii) for modif. $ev(b,-)$ \\ [1ex]   
\end{tabular}
\caption{Why $ev(-,\beta)$ is a modification}
\label{table:3}
\end{center}
\end{table}

\begin{table}[h!]
\begin{center}
\begin{tabular}{ c c } 
perturb. axiom for $ev(-,\Lambda):$ & \\ [1ex]
$ev(-,\beta)\Rrightarrow ev(-,\beta')$ & its meaning on cells in $\Oo^*_{\bullet}(\B, \C)$ on which evaluated \\ [0.5ex]
\hline
single axiom & axiom (ii) for transf. $\alpha$ w.r.t. $\Lambda:\beta\Rrightarrow\beta'$ \\ [1ex]   
\end{tabular}
\caption{Why $ev(-,\Lambda)$ is a perturbation}
\label{table:4}
\end{center}
\end{table}

\medskip

The coincidences of the 3-cells in (iii) and (iv) of \prref{quasi-Gray} obviously hold by the definition of $ev$, it only remains to comment 
the axiom in (iii). Observe that the structure 3-cell $ev(-,g)^2_\alpha$ for the lax transformation $ev(-,g)$ is given by the identity: 
the functors $ev(-,B), ev(-,B')$ are strict and the 3-cell $ev(-,g)^2_\alpha$ is then just the identity expressing how the composition component $\frac{\alpha}{\alpha'}\vert_g$ is given in terms of components $\alpha_g$ and $\alpha'_g$. The axiom that we need to prove 
in (iii) then becomes much simpler and it melts down to the identity saying how the transformation structure 3-cell 
$(\frac{\alpha}{\alpha'})^2_{g',g}$ is given in terms of  transformation structure 3-cells $\alpha_{g',g}^2$ and $\alpha^{'2}_{g',g}$, 
which clearly holds. 
\qed\end{proof}

In particular, left evaluation, being strict in the left variable, is easily seen to extend to a quasi-functor $ev: \Oo^*_\bullet(\B,\C)\times\B\to C$ of type $(*,\bullet)$, so that given any quasi-functor $H:\A\times\B\to\C$ of type $(*,\bullet)$ it is 
$
ev(H^t\times \Id_\B)=H
$ 
where $H^t:\A\to \Oo^*_{\bullet}(\B, \C)$ is the functor corresponding to $H$ by \equref{quasi}. We will get back to this quasi-functor in \thref{mixed-quasi n}.

\subsection{Nullary and unary multimaps for the multicategory} \sslabel{nullary/unary}

For the multicategory and skew-multicategory $\OO$ that we are going to construct corresponding to $\Oo$ the class of 2-categories, double categories or Gray-categories, for the class of nullary maps $\OO_0^l(-;\A)$ we set the set of objects of $\A$, 
and for tight and loose unary maps we set $\OO_1^t(\A;\B)=\Oo_{st}(\A,\B)$ and $\OO_1^l(\A;\B)=\Oo_*(\A;\B)$. For a quasi-functor 
$F_n:\A_1\times...\times\A_n\to\B$ and an object $A_k\in\A_k$ we set $F_n\circ_k A_k = F_n(...,A_k,...)$, which is the $k$-th $n\x 1$-ary multimap constituting $F_n$, where $A_k$ denotes the nullary multimap $(-)\to\A_k$ picking up the object $A_k$. Similarly, substituting $A_k$ by any unary multimap $F_1: \A\to\A_k$, we define $F_n\circ_k F_1 = F_n(...,F_1(-),...)$. At an object $A\in\A$ 
it is defined as above, whereas at higher cells its definition is induced according to \prref{char df} and \prref{quasi-Gray}. 

When only nullary maps, or a nullary and a unary multimap, are composed into an $n$-ary multimap at any two variables, the associativity of the substitution holds by \equref{obj} and \equref{obj 4-ary}. For a further unary map $G_1:\B\to\C$ we define
$$(G_1\circ F_n)\circ_k A:=G_1\circ(F_n\circ_k A) \qquad\text{and}\qquad (G_1\circ F_n)\circ_k F_1:=G_1\circ(F_n\circ_k F_1)$$ 
assuring associativity in these cases. They mean 
$(G_1\circ F_n)(...,A_k,...)=G_1\circ F_n(...,A_k,...)$ and $(G_1\circ F_n)(...,a_k,...)=G_1\circ F_n(...,a_k,...)$ for a higher cell $a_k$ in $\A_k$. In particular, for $F_n:\A_1\times...\times\A_n\to\Oo(\B,\C)$ and $G_1=ev_B$ for an object $B\in\B$, 
we have 
\begin{equation}\eqlabel{first-last}
(ev_B\circ F_n)(...,A_k,...)=ev_B\circ F_n(...,A_k,...).
\end{equation}

\subsection{Properties of tight multimaps: from binary to ternary} \sslabel{2-3}

Iterating \equref{mixed-quasi-bij}, {\em i.e.} for $n=3$, we obtain
\begin{equation} \eqlabel{towards ternary} 
q\x\Oo^{(st,*)}_\bullet(\A_1{}^{st}\w\times^*\A_2,\Oo^*_\bullet(\B,\C)) 
\stackrel{\equref{mixed-quasi-bij}}{\iso}\Oo_{st}(\A_1, \Oo^*_\bullet(\A_2,\Oo^*_\bullet(\B,\C)))\stackrel{\equref{quasi}}{\iso}
\Oo_{st}(\A_1, q\x\Oo^*_\bullet\vert_{(\bullet,\bullet)}(\A_2\times\B,\C)).
\end{equation} 
A $(st,*)$-mixed quasi-functor from the left $\F:\A_1{}^{st}\w\times^*\A_2\to\Oo^*_\bullet(\B,\C)$ 
consists of a strict functor $\A_1\to\Oo^*_\bullet(\B,\C)$ and a $*$-type functor 
$\A_2\to\Oo^*_\bullet(\B,\C)$. In return, the former consists of a strict functor $\A_1{}^{st}\to\C$ and a $*$-type functor 
$\B\to\C$, whereas the latter consists of $*$-type functors $\A_2\to\C$ and the same $\B\to\C$, so that in total we have  
one strict and two $*$-type functors. The analogous three functors we find in any $(st,*)$-mixed quasi-functor 
$\G:\A_1{}^{st}\w\times^*\A_2\times^*\B\to\C$. 

To prove that there is an isomorphism \equref{ternary} below we should verify \equref{obj} (both for double and Gray-categories), check that $\F$, seen as a functor of three variables $(-,-,-):=:\crta\F(-,-,-)=\F(-,-)(-)$, and any $\G$ living in the left-hand side of \equref{ternary} satisfy the same axioms, and identify an inverse assignment in the opposite direction. Let $A_1\in\A_1, A_2\in\A_2, B\in\B$. 
Since $\F$ is a quasi-functor (on the left in \equref{towards ternary}), we have $\F(A_1,-)\vert_{A_2}=\F(-,A_2)\vert_{A_1}$ and thus 
$(A_1,-,-)\vert_{A_2}=(-,A_2,-)\vert_{A_1}$. On the other hand, since $\F(A_1)$ is a quasi-functor (on the right in \equref{towards ternary}, by abuse of notation), we have $\F(-)(A_2,-)\vert_B=\F(-)(-,B)_{A_2}$, and thus $(-,A_2,-)\vert_{B}=(-,-,B)\vert_{A_2}$ (by identifying $(-,-,-)=\F(-)(-,-)$).The identity $(A_1,-,-)\vert_{B}=(-,-,B)\vert_{A_1}$ follows by \equref{first-last} with $k=1$, as $\F(A_1)(-,-)\vert_B=\F(-)(-,B)\vert_{A_1}$. 
To check the axioms, let $f:A\to A', g:B\to B', h:C\to C'$ be 1h-cells in $\A_1, \A_2$ and $\B$, respectively. By the composite isomorphism in \equref{towards ternary}, we get that $\F(f)(-,-)=(f,-,-)$ is a transformation of 
type $(\bullet,\bullet)$ among quasi-functors of two variables $(A,-,-)\Rightarrow(A',-,-):\A_2\times\B\to\C$. 
In dimension 2, the axiom \axiomref{$HOT^q_1$} for $(f,-,-)=\theta$ is then satisfied for all $g=K$ and $h=K$ as indicated above. This is precisely the axiom $i)$ of \deref{double qf 3}. In the internal case the proof is similar for the 
1v-cells and the rest of the axioms (observe that for a 1v-cell $U:A\to\tilde A$ we have that $\F(U)(-,-)=(U,-,-)$ is a vertical lax transformation and as such it satisfies axioms \axiomref{$VLT^q_1$}-\axiomref{$VLT^q_4$}). Thus the eight axioms of 
\deref{double qf 3} in the internal case are recovered from the axioms \axiomref{$HOT^q_1$}-\axiomref{$HOT^q_4$} and \axiomref{$VLT^q_1$}-\axiomref{$VLT^q_4$} of horizontal oplax and vertical lax transformations between lax double quasi-functors of two variables. 

In dimension 3, for Gray-categories, $\F(f)(-,-)=(f,-,-)$ being a transformation of type $\bullet=oplx$, it means that there is a 3-cell 
$((f,-,-), g,h)$ for all 1-cells $(g,h)\in(\A_2,\B)$ as in \deref{tr Gray qu-f 2} (with $g=K$ and $h=k$), satisfying compatibility axioms 
\equref{teta-f-g-left} and \equref{teta-f-g-right} and two degeneracy equations: $(\theta, id_A,g)=\Id_{\theta^A_g}$ and 
$(\theta, f, id_B)=\Id_{\theta^B_f}$. Relabeling in those axioms first $f,g\mapsto g,h$ (including their (co)domains) and then applying that the transformation $\theta:(-,-)_1\Rightarrow(-,-)_2$ is now $(f,-,-):(A,-,-)\to(A',-,-)$, 
the cocycles ${(f,-)^2_1}_{g'g}, \, {(f,-)^2_2}_{g'g}, \, {(-,g)^2_1}_{f'f}$ and ${(-,g)^2_2}_{f'f}$ become 
$(A,g,-)^2_{h'h}, (A',g,-)^2_{h'h}, \, (A,-,h)^2_{g'g}$ and $(A',-,h)^2_{g'g}$, respectively. Furthermore, the 3-cells 
$(\theta,f,g)$ become $(f,g,h)$, while the cocycles $(\theta^{A})^2_{g'g}$ and $(\theta^{B})^2_{f'f}$ become $(f,B,-)^2_{h'h}$ 
and $(f,-,C)^2_{g'g}$, respectively. With these changes, the right-hand-sides of \equref{teta-f-g-left} and \equref{teta-f-g-right} 
become exactly the incubator axioms \axiomref{$f, g'g,h$} and \axiomref{$f,g, h'h$}. Similarly, for the left-hand sides in 
\equref{teta-f-g-left} and \equref{teta-f-g-right} the 3-cells $(a,g)_1, (a,g)_2, \, (f,b)_1, (f,b)_2$ and $\theta^A_b, \theta^B_a$ 
become $(A,-,h)_b, (A',-,h)_b, (A,g,-)_c, (A',g,-)_c$ and $(f,B,-)_c, (f,-,C)_b$, respectively, so that the two mentioned equations 
become precisely the incubator axioms \axiomref{$f,b,h$} and \axiomref{$f,g, c$}. The two degeneracy equations correspond to the 
incubator axioms \axiomref{$A:1_b:C$} and \axiomref{$A:B:1_c$}.
For the remaining three axioms from 
\cite[Appendix D.2]{BL} that would make out of $\F$ a quasi-functor of three variables we have the following. 
The relation between the oplax transformation $\F(f'f)(-,-)$ and the composition of transformations $\F(f)(-,-)$ and $\F(f')(-,-)$ is expressed 
by the axiom \equref{comp teta}, whereby passing from the generic transformation $\theta$ to $\F(f)(-,-)$ (and the structure 3-cell $(\theta, K,k)$ to the 3-cell $(f,g,h)$), in order for the expression of the composition to make sense, one needs to add into account the cocycles 
$(-,g,C)^2_{f'f}=\F(f'f,g)(C)$ (of the quasi-functor $\F(-,-)(C):\A_1{}^{st}\w\times^*\A_2\to\C$) and $(-,B,h)^2_{f'f}=\F(f'f,B)(h)$ (of the quasi-functor $\F(-,B)(-):\A_1\times\B\to\C$), where the object $C$ and 1-cell $h$ are from $\B$ and $B\in\A_2$. Then this recovers the axiom \axiomref{$A_1, A_2:B:C$}. 
The modification axiom \equref{modif qf 2} of $\F(a)(-,-)$, for any 2-cell $a:f\Rightarrow f'$ in $\A_1$, corresponds to the axiom \axiomref{$\alpha:B:C$}: the 3-cells $\Lambda^A_g$ and $\Lambda^B_f$ becomes 3-cells $(a,B,h)=(-,B,h)_a$ and $(a,g,C)=(-,g,C)_a$, respectively.  
Finally, $\F(id_A)(-,-)$ is evidently the identity transformation, which corresponds to the degeneracy axiom \axiomref{$1_a:B:C$}. 

There is another way to prove the same. Namely, 
observe that for a fixed object $B\in\B$ we have $\F(-,-)(B)=(-,-,B)\in q\x\Oo^{(st,*)}_\bullet(\A_1{}^{st}\w\times^*\A_2,\C)$, a quasi-functor of two variables. Thus for each 1-cell $h:B\to B'$ there is a transformation $(-,-,h)$ of type $(\#,\#)$ of such quasi-functors (the type $\#$ we see by considering $\F$ as an object of the right-hand side of \equref{towards ternary}). This is lax type of transformation. Then there is a 3-cell 
$((-,-,h), f,g)$ as in \rmref{oplax q.f.-lax tr.} 
and the same arguments as above show that $((-,-,h), f,g)$ determines an incubator 3-cell. As a matter of fact, they both are equal 
to $\F(f,g)(h)=\crta\F(f,g,h)$, which is indeed an oplax type of an incubator 3-cell. Observe that this 3-cell is also the modification 
3-cell component $\F(f,g)(h)$, obtained by evaluating the modification $\F(f,g)$ at the 1h-cell $h$. 

\smallskip

Applying exactly the same arguments backwards and with identifications as above: $\G(-,*,\bullet)=(-,*,\bullet)=\F(-,*)(\bullet)=\F(-)(*,\bullet)$, we 
obtain the assignment in the other direction, proving finally the bijection  
\begin{equation} \eqlabel{ternary}
q_3\x\Oo^{(st,*)}_\bullet(\A_1{}^{st}\w\times^*\A_2\times^*\B,\C)\iso 
q\x\Oo^{(st,*)}_\bullet(\A_1{}^{st}\w\times^*\A_2,\Oo^*_\bullet(\B,\C))
\end{equation}
both for double categories and Gray-categories. 

Observe that for every fixed $B\in\B$ mixed quasi-functors $\F$ and $\G$ restrict to mixed quasi-functors 
$\A_1{}^{st}\w\times^*\A_2\to\C$, and similarly, for every fixed object in $\A_2$ they restrict to  mixed quasi-functors 
$\A_1{}^{st}\w\times^*\B\to\C$, and this does not depend on the dimension of $\Oo$. For further use we highlight:

\begin{cor} \colabel{deducing q.f.}
A quasi-functor $\A_1{}^{st}\w\times^*\A_2\to\Oo^*_\bullet(\B,\C)$ yields quasi-functors $\A_1{}^{st}\w\times^*\A_2\to\C$ 
and $\A_1{}^{st}\w\times^*\B\to\C$ for general $\Oo$. 
\end{cor}

\medskip

We also may highlight a property that we used at the beginning of proving \equref{ternary}, which also holds for general $\Oo$ 
(compare to \rmref{tr in qu.f.}).

\begin{cor} \colabel{tr of qf of 2v}
Given a quasi-functor of three variables $\G:\A{}^{st}\w\times^*\B\times^*\C\to\E$ of type $((st,*),\bullet)$ in general $\Oo$. 
Let $f:A\to A', g:B\to B', h:C\to C'$ be 1h-cells and $\alpha, \beta, \gamma$ 2-cells in $\A, \B$ and $\C$, respectively. 
Then $\G(f,-,-)$ is a transformation of type $(\bullet, \bullet)$, $\G(-,g,-)$ is of type $(\#, \bullet)$, and $\G(-,-,h)$ is 
a transformation of type $(\#, \#)$ of quasi-functors of two variables, $\G(\alpha,-,-), \G(-,\beta,-), \G(-,-,\gamma)$ are the corresponding modifications, and similarly higher $k$-cells, $k=3,...,m$, induce $k$-cells in $\Oo^*_\bullet(\C,\E)$.  
\end{cor}

A further very useful tool that we proved while proving \equref{ternary} and that we will use in two further instances is the following:

\begin{lma} \lelabel{gain incubator}
Let $\F$ be from 
$q\x\Oo^{(st,*)}_\bullet(\A_1{}^{st}\w\times^*\A_2,\Oo^*_\bullet(\B,\C))\iso\Oo_{st}(\A_1, q\x\Oo^*_\bullet\vert_{(\bullet,\bullet)}(\A_2\times\B,\C))$ and $\G\in q_3\x\Oo^{(st,*)}_\bullet(\A_1{}^{st}\w\times^*\A_2\times^*\B,\C)$. 
Oplax incubator 3-cells $(f,g,h)=\G(f,g,h)$ for the quasi-functor $\G$ are in 1-1 correspondence with either of: 
\begin{itemize}
\item the structure 3-cell $((f,-,-), g,h)$ of the oplax transformation $\F(f)(-,-)=(f,-,-)$ corresponding to the structure 2-cells $(g,h)_i$ 
for $i=1,2$ of the corresponding quasi-functors of type $\bullet=oplax$, 
\item the structure 3-cell $((-,-,h), f,g)$ of the lax transformation $\F(-,-)(h)=(-,-,h)$ corresponding to the structure 2-cells and $(f,g)_i$ 
for $i=1,2$ of the corresponding quasi-functors of type $\bullet=oplax$. 
\end{itemize}
\end{lma}

\medskip

\begin{rem}
We comment in more details the types of transformations appearing in \equref{ternary}. 
In \equref{towards ternary} we indeed have more generally 
\begin{equation} \eqlabel{more gen}
q\x\Oo^{(st,*)}_\bullet(\A_1{}^{st}\w\times^*\A_2,\Oo^*_\star(\B,\C)) 
\stackrel{\equref{mixed-quasi-bij}}{\iso}\Oo_{st}(\A_1, \Oo^*_\bullet(\A_2,\Oo^*_\star(\B,\C)))\stackrel{\equref{quasi}}{\iso}
\Oo_{st}(\A_1, q\x\Oo^*_\bullet\vert_{(\bullet, \star)}(\A_2\times\B,\C)))
\end{equation}
but it was to obtain \equref{ternary} that we chose $\star=\bullet$. 
Namely, for a quasi-functor $\G$ from the left of  \equref{ternary} and 1h-cells $(f,g,h)\in(\A_1, \A_2, \B)$ we have by 
the definition of a ternary map, \prref{quasi-fun} and \prref{quasi-Gray} that 
\begin{itemize}
\item $\G(f,A_2,-)$ and $\G(A_1,g,-)$ are transformations of type $\bullet$  - see them as $(f,-), (g,-)$ in the cited propositions, and 
\item $\G(-,A_2,h)$ and $\G(A_1,-,h)$ are transformations of type $\#$  - see them both as $(-,h)$.  
\end{itemize}
On the other hand, for a quasi-functor $\F$ on the right of \equref{ternary}, {\em i.e.} on the left of \equref{more gen}, 
we have
\begin{itemize}
\item $\F((f,A_2),-)$ and $\F((A_1,g),-)$ are transformations of type $\star$ (this is why we wrote $\star=\bullet$ in \equref{towards ternary}, and 
\item $\F((-,-),h)$ is a transformation of type $(\#,\#)$, that is, both $\F((-,A_2),h)$ and $\F((A_1,-),h)$ are of type $\#$.   
\end{itemize}
\end{rem}

\subsection{Properties of tight multimaps: from ternary to 4-ary and any $n$-multimap}

In order to complete the proof for Gray-categories of the relation between $n$-ary and $(n+1)$-ary multimaps, we need one step more, which we do in this subsection. For that purpose we do the analogous as in \ssref{2-3}. We have 
\begin{align}
q_3\x\Oo^{(st,*)}_\bullet(\A_1{}^{st}\w\times^*\A_2^*\times^*\A_3,\Oo^*_\bullet(\B,\C)) 
& \stackrel{\equref{ternary}}{\iso}q\x\Oo^{(st,*)}(\A_1\times\A_2, \Oo^*_\bullet(\A_3,\Oo^*_\bullet(\B,\C))) \nonumber \\
& \stackrel{\equref{towards ternary}}{\iso}
\Oo_{st}(\A_1, q\x\Oo^*_\bullet\vert_{(\bullet,\bullet)}(\A_2\times\A_3,\Oo^*_\bullet(\B,\C))) \nonumber \\
& \stackrel{\equref{ternary}}{\iso}
\Oo_{st}(\A_1, q_3\x\Oo^*_\bullet\vert_{(\bullet,\bullet,\bullet)}(\A_2\times\A_3\times\B,\C)). \nonumber \label{towards 4-ary}\\
\end{align}
Using this we now want to show the isomorphism 
\begin{equation} \eqlabel{4-ary}
q_4\x\Oo^{(st,*)}_\bullet(\A_1{}^{st}\w\times^*\A_2^*\times^*\A_3^*\times^*\B,\C) \iso 
q_3\x\Oo^{(st,*)}_\bullet(\A_1{}^{st}\w\times^*\A_2^*\times^*\A_3,\Oo^*_\bullet(\B,\C)). 
\end{equation}
Take $\F$ from the right-hand side above and let $\crta\F$ denote the induced map $\A_1{}^{st}\w\times^*\A_2^*\times^*\A_3^*\times^*\B\to\C$. 
For the identities \equref{obj 4-ary} for $\crta\F$: the cases involving $(A,B,C)\in(\A_1\times\A_2\times\A_3)$ hold, as they already hold for $\F$, while similarly the cases involving $(B,C,D)\in(\A_2\times\A_3\times\B)$ hold because the corresponding identities hold in the right-hand side in the third line in (\ref{towards 4-ary}). The only remaining case is the one involving $(A,D)\in(\A_1\times\B)$: it holds by \equref{first-last}, similarly as in the previous subsection. 
For a 1-cell $f:A\to A'$ in $\A_1$ by the right-hand side of (\ref{towards 4-ary}) we have that $\F(f)(-,-,-)$ is a transformation of type 
$(\bullet, \bullet,\bullet)$, which is oplax, of quasi-functors of three variables $\F(A)(-,-,-)\to\F(A')(-,-,-)$, which are oplax, too. Then 
$\F(f)(-,-,-)$ satisfies the left-hand side axiom in \equref{tr of qf 3}.  
To avoid confusion, let us shift the lettering in this axiom from $(f,g,h):(A,B,C)\to(A',B',C')$ into 
$(g,h,k):(B,C,D)\to(B',C',D')$, so that we have 
$$\fourfrac{(\theta^D,g,h)}{(\theta^C,g,k)}{(\theta^B,h,k)}{(g,h,k)^1_{oplx-inc}}=
\fourfrac{(g,h,k)^2_{oplx-inc}}{(\theta^B,h,k)}{(\theta^C,g,k)}{(\theta^D,g,h)}.$$ 
The transformations for two variables $\theta^B, \theta^C, \theta^D$ appearing here in the present context are 
$$(f,B,-,-):(A,B,-,-)\Rightarrow(A',B,-,-)$$ 
$$(f,-,C,-):(A,-,C,-)\Rightarrow(A',-,C,-)$$ 
$$(f,-,-,D):(A,-,-,D)\Rightarrow(A',-,-,D)$$ 
so that the 3-cell $(\theta^D,g,h)$ here presents the structure 3-cell of the oplax transformation 
$\theta=\F(f)(-,-,-)$ corresponding to the oplax quasi-functors of two variables $\theta^D=\F(f)(-,-,D)$ and the 2-cells 
$(g,h)\in\A_2\times\A_3$. According to \leref{gain incubator} it is the oplax incubator 3-cell $(f,g,h,D)=\crta\F(f,g,h,D)=\F(f,g,h)(D)$. 
By continuing this interpretation for the rest of the 3-cells in the fractions above containing $\theta$ and observing that 
$(g,h,k)^1_{oplx-inc}=(A,g,h,k)$ 
and $(g,h,k)^2_{oplx-inc}=(A',g,h,k)$ as the oplax incubator 3-cells of the oplax quasi-functors $\F(A)(-,-,-)$ and $\F(A')(-,-,-)$ 
from above, we end up obtaining precisely the mecon axiom \equref{mecon} for $\crta\F$. This proves that $\crta\F$  
is a 4-ary quasi-functor in the left-hand side of \equref{4-ary}. 

\medskip

As in the proof of \equref{ternary}, there is one more way to prove \equref{4-ary}. 
For objects $D, D'$ and a 1-cell $k:D\to D'$ from $\B$ we know that $\F(-,-,-)(D)$ and $\F(-,-,-)(D')$ 
are quasi-functors of three variables, and that $\F(-,-,-)(k)$ 
is a lax transformation between them. 
Let $(f,g,h)^1_{oplx-inc}=(f,g,h,D)_{oplx-inc}$ and $(f,g,h)^2_{oplx-inc}=(f,g,h,D')_{oplx-inc}$ denote the incubator 3-cells of $\F(-,-,-)(D)$ and $\F(-,-,-)(D')$, respectively, for 1-cells $(f,g,h)\in\A_i\times\A_j\times\A_k$. Given that $\F$ is a quasi-functor of type 
$\bullet=oplax$, then so are $\F(-,-,-)(D)$ and $\F(-,-,-)(D')$, hence their incubator 3-cells are of oplax type. Then the right-hand side of \equref{tr of qf 3} applies for the transformation $\F(-,-,-)(k)$. 
The 3-cells $(\theta^C,f,g)$ appearing there in the present context are structure 3-cells of the lax transformation $\theta=\F(-,-,-)(k)$ 
that correspond to the quasi-functors of two variables $\theta^C=\F(-,-,C)(k)$ and the 2-cells $(f,g)\in\A_i\times\A_j$. 
According to \leref{gain incubator}
the structure 3-cells $(\theta^C,f,g)$ are the incubator 3-cells $(f,g,C,k)=\crta\F(f,g,C,k)=\F(f,g,C)(k)$. 
By continuing this interpretation for the rest of the 
similar 3-cells in the right-hand side of \equref{tr of qf 3} we end up getting precisely the mecon axiom \equref{mecon} for $\crta\F$.

\bigskip

As in \leref{gain incubator} we may formulate:

 \begin{lma} \lelabel{gain mecon}
Let $\F$ be from 
$q_3\x\Oo^{(st,*)}_\bullet(\A_1{}^{st}\w\times^*\A_2\w\times^*\A_3,\Oo^*_\bullet(\B,\C))\iso\Oo_{st}(\A_1, q_3\x\Oo^*_\bullet\vert_{(\bullet,\bullet,\bullet)}(\A_2\times\A_3\times\B,\C))$ and $\G\in q_4\x\Oo^{(st,*)}_\bullet(\A_1{}^{st}\w\times^*\A_2\w\times^*\A_3\times^*\B,\C)$. 
The mecon axiom for 1-cells $(f,g,h,k)\in(\A_1\times\A_2\times\A_3\times\B)$ for the quasi-functor $\G$ is obtained as either of: 
\begin{itemize}
\item the oplax transformation axiom for $\F(f)(-,-,-)=(f,-,-,-)$ for quasi-functors of three variables corresponding to the 1-cells $(g,h,k)\in(\A_2\times\A_3\times\B)$, 
\item the lax transformation axiom for $\F(-,-,-)(k)=(-,-,-,k)$ for quasi-functors of three variables corresponding to the 1-cells $(f,g,h)\in(\A_1\times\A_2\times\A_3)$. 
\end{itemize}
\end{lma}

\subsection{Closedness isomorphisms} \sslabel{closedness}

We now may prove:

\begin{thm} \thlabel{mixed-quasi n}
Let $\Oo$ stand for the categories of 2-categories, double categories, or Gray-categories so that $\Oo^*_\bullet(\A,\B)$ is an object 
of $\Oo$ for any objects $\A,\B$ of $\Oo$, and let $n\geq 2$. 
\begin{enumerate}
\item There is a bijection
$$
q_n\x\Oo^{(st,*)}_\bullet(\Pi_{i=1}^n\A_i,\C) \iso q_{n-1}\x\Oo^{(st,*)}_\bullet(\Pi_{i=1}^{n-1}\A_i,\Oo^*_\bullet(\A_n,\C))
$$
so that given a quasi-functor of $n$-variables $H$ on the left there is a unique quasi-functor of $n-1$-variables $H^t$ 
on the right such that $ev((H^t){}^{st}\w\times^*\Id_{\A_n})=H$ with $ev$ from \prref{ev-mixed}. 
\item 
There is a quasi-functor $ev:\Oo^*_{\bullet}(\B, \C)\times\B\to\C$ of type $(*,\bullet)$ that induces bijections  
$$q_{n}\x\Oo^*_\bullet(\Pi_{i=1}^{n-1} \A_i\times\B,\C)\iso q_{n-1}\x \Oo^*_\bullet(\Pi_{i=1}^{n-1} \A_i, \Oo^*_\bullet(\B,\C)).$$
\end{enumerate}
\end{thm}

\begin{proof}
We only prove the first part. Set $\B=\A_n$. 
Let $\F:\Pi_{i=1}^{n-1}\A_i\to\Oo^*_\bullet(\B,\C)$ be from the right, let $\crta\F:\Pi_{i=1}^n\A_i\to\C$ be induced by $\F$ in the obvious way, that is $\crta\F(a_1,...,a_n)=\F(a_1,...,a_{n-1})(a_n)$ with $a_i\in\A_i, i=1,..,n$, where $\sum_{i=1}^n deg(a_i)\leq m$, so that at 
most three variables of $a_1,...,a_{n-1}$ have degrees greater than $0$. Let us consider the 
quasi-functor property of $\F$ at arbitrarily chosen three variables $i,j,k$ of $\Pi_{l=1}^{n-1}\A_l$, and let 
$A_i\in\A_i, A_j\in\A_j,A_k\in\A_k, B\in\B$. Then $\crta\F(-,-,-,B)=\F(-,-,-)(B):\A_i\times\A_j\times\A_k\to\C$ is a quasi-functor, which is a restriction of both $\crta\F$ and $\F$. On the other hand, 
we have in particular that $\F(A_i,-,-)\w:\A_j\times\A_k\to\Oo^*_\bullet(\B,\C)$ is a quasi-functor, then by \equref{ternary} 
it determines a quasi-functor $\crta\F(A_i,-,-,-)\w:\A_j\times\A_k\times\B\to\C$, which is a restriction of $\crta\F$. Instead of $i$ 
we could have done the same reasoning with $j$ and $k$. This way we cover all the possible triples of variables in $\Pi_{i=1}^n\A_i$ and 
thus we proved that $\crta\F$ is a quasi-functor $\Pi_{i=1}^n\A_i\to\C$ for double categories, while for Gray-categories this proves that $\crta\F$ is a collection of ternary quasi-functors. It remains to prove the identities \equref{obj 4-ary} 
and the mecon axiom for Gray-categories. Observe that it is sufficient to prove them for variables from $\A_i\times\A_j\times\A_k\times\B$, 
as for the remaining variables those properties pass along from holding for $\F$. But this now holds as in the proof of \equref{4-ary}, and we are done. 
The converse is clear. 
\qed\end{proof}

\begin{rem} \rmlabel{Gray-cat isos}
\begin{enumerate}
\item For every $n$ the sets $q_n\x\Oo^{(st,*)}_\bullet(\Pi_{i=1}^n\A_i,\C)$ are functorial in $\A_1,..,\A_n,\C$. Namely, we know that the right-hand side of \equref{quasi} is functorial in $\A,\B,\C$, thus the left-hand side, that is $q\x\Oo^*_\bullet(\A\times\B,\C)$,  becomes so, too. Then the right-hand side of \equref{ternary} is functorial in $\A_1,\A_2,\B,\C$, then so becomes 
$q_3\x\Oo^{(st,*)}_\bullet(\A_1{}^{st}\w\times^*\A_2\times^*\B,\C)$. Similarly, $q_4\x\Oo^{(st,*)}_\bullet(\A_1{}^{st}\w\times^*\A_2^*\times^*\A_3^*\times^*\B,\C)$ from \equref{4-ary} becomes functorial in $\A_1,\A_2,\A_3,\B,\C$, and for general $n$ the functoriality is obtained iteratively from the above Theorem. 
\item Analogously as in \ssref{Gray-cat of qf}, $q_n\x\Oo^{(st,*)}_\bullet(\Pi_{i=1}^n\A_i,\C)$ can be made into a Gray-category and the isomorphisms in the above Theorem are indeed in $G\x\Cat$. 
(For general $\Oo$ one should be able to form a category of dimension $m$, that is, an object in $\Oo$, of quasi-functors of $n$-variables 
(see Appendix B.3 for an insight on how vertical transformations and modifications in the internal case are defined).) 
\end{enumerate}
\end{rem}

\bigskip

The result of \thref{mixed-quasi n} unifies on one hand a version of our \thref{left closed multicat} and \thref{strict-left closed} 
for double categories, and on the other hand Propositions 4.2, 4.4 and 4.6 of \cite{BL} for Gray-categories and $n=2,3,4$. Apart from this, 
our proof of \equref{ternary} and  \equref{4-ary}, which uses the isomorphisms \equref{towards ternary} and (\ref{towards 4-ary}), 
respectively, provides an insightful hint for the proof of the analogues of \equref{ternary} and \equref{4-ary} in higher dimensions $m$, 
{\em i.e.} for any $\Oo$. Namely, our \leref{gain incubator} and \leref{gain mecon} prove 1-1 correspondences 
\begin{equation}\eqlabel{1-1}
\hspace{0,38cm}\text{transformations of binary quasi-functors}\quad \leftrightarrow \quad \text{ternary quasi-functors}\vspace{-0,2cm}
\end{equation}
$$\text{transformations of ternary quasi-functors}\quad \leftrightarrow \quad \text{4-ary quasi-functors},$$
which seem to mark a trend that repeats at higher dimensions.


There is more we can say. For $m>3$ let $n\geq m+1$ and set 
$$\axiom{$A_n$} \qquad q_{n-1}\x\Oo(\A_1\times...\times\A_{n-1},\Oo(\A_n,\C))\iso \Oo_{st}(\A_1,q_{n-1}(\A_2\times...\times\A_n,\C))$$
$$\hspace{-0,54cm}\axiom{$B_n$} \qquad\qquad q_n\x\Oo(\A_1\times...\times\A_n,\C)\iso q_{n-1}\x\Oo(\A_1\times...\times\A_{n-1},\Oo(\A_n,\C)).$$
In dimension 2, for double categories, n-ary quasi-functors are defined by ternary quasi-functors, 
while in dimension 3, for Gray-categories, n-ary case is determined by the 4-ary one. Then for dimension $m$ one should check successive pairs of isomorphisms $(A)_n$ and $(B)_n$ up to $n=m+1$ variables. Applying induction, knowing that $(A)_{n-1}$ and $(B)_{n-1}$ hold, one may prove $(A)_n$: 
\begin{align}
q_{n-1}\Oo(\A_1\times...\times\A_{n-1},\Oo(\A_n,\C)) & \stackrel{(B)_{n-1}}{\iso}
q_{n-2}\Oo(\A_1\times...\times\A_{n-2},\Oo(\A_{n-1},\Oo(\A_n,\C))) \nonumber \\
& \stackrel{(A)_{n-1}}{\iso}\Oo_{st}(\A_1,q_{n-2}(\A_2\times...\times\A_{n-1},\Oo(\A_n,\C)))\nonumber  \\
& \stackrel{(B)_{n-1}}{\iso} \Oo_{st}(\A_1,q_{n-1}(\A_2\times...\times\A_n,\C)).\nonumber 
\end{align}
To prove $(B)_n$, take $\F$ from its right-hand side, we want to prove that the induced map $\crta\F:\A_1\times...\times\A_n\to\C$ 
is an $n$-ary quasi-functor. 
For the identities ``on objects'' (as in \equref{obj 4-ary}) for $\crta\F$, the cases involving variables in $\A_1\times...\times\A_{n-1}$ hold, as they already hold for $\F$, while similarly the cases involving  variables in $\A_2\times...\times\A_n$ hold because the corresponding identities hold in the right-hand side in $(A)_n$. The only remaining case is the one involving the variables in $\A_1$ and 
$\A_n$, which holds by \equref{first-last} with $k=1$ (in the case of dimension $m$). 
For a 1-cell $f:A\to A'$ in $\A_1$ by the right-hand side of $(A)_n$ we have that $\F(f)(-,...,-)$ is a transformation of type 
$\bullet$, which is oplax, of quasi-functors of $n-1$ variables $\F(A)(-,...,-)\to\F(A')(-,...,-)$, which are oplax, too. By the above 
mentioned trend this should correspond to an $n$-ary quasi-functor, so that the isomorphism $(B)_n$ rest proved.

\bigskip

Quasi-functors of $(st,*)$-mixed type as above - more precisely, of the type $((st,*),\bullet)$ - induce meta-product 
$\ot^{(st,*)}_\bullet$ (they originate from a strict functor $\A\to\Oo^*_\bullet(\B,\A\times\B)$).

\subsection{Tight multimaps of right skew-multicategories}

For tight multimaps of a right skew-multicategory we will take quasi-functors that are {\em strict only in the rightmost entry} and loose 
(of the same type $*$) in all the remaining entries. Similarly as in \prref{ev-mixed}, there are 
evaluation functors $\crta{ev}:\A\s^*\w\times^{st}\Oo^*_\#(\A,\C)\to\C$ for quasi-functors of type $((*,st), \bullet)$ as in the bijections 
\begin{equation} \eqlabel{mixed-right-bij}
q\x\Oo^{(*,st)}_\bullet(\A\s^*\w\times^{st}\B,\C) \iso\Oo_{st}(\B,\Oo^*_\#(\A,\C))\iso\Oo_*(\A,\Oo^{st}_\bullet(\B,\C)).
\end{equation}
For the argument why $\crta{ev}:\A\s^*\w\times^{st}\Oo^*_\#(\A,\C)\to\C$ is a quasi-functor of type $((*,st), \bullet)$, recall 
\thref{lax-right closed} and the discussion preceding it. (Observe that for a quasi-functor $H$ and a 1-cell $f$ in $\A$ 
the transformation $H(f,-)$ is of type $\bullet$, then we know that consequently the transformation $H(-,g)$ is of type $\#$, for a 1-cell $g$ in $\B$; this explains the type $\#$ appearing in \equref{mixed-right-bij}). The above bijections lift to natural isomorphisms in $\Oo$  
\begin{equation} \eqlabel{mixed-quasi-right}
q\x\Oo^{(*,st)}_\bullet\vert_{(\#,\bullet)}(\A\s^*\w\times^{st}\B,\C) \iso\Oo^{st}_\bullet(\B,\Oo^*_\#(\A,\C))\iso
\Oo^*_\#(\A,\Oo^{st}_\bullet(\B,\C)).
\end{equation}
%
Observe that \equref{mixed-quasi-right} is a tight version of \equref{quasi-right}. 
These quasi-functors induce the meta-product $\ot^{(*,st)}_\bullet$. 

\smallskip

For the record we state that the bijections \equref{quasi-set-sharp} induce meta-products $\ot^{(st,*)}_\#$ and $\ot^{(*,st)}_\#$. 

\medskip

We comment the right-hand side version of \thref{mixed-quasi n}. For this purpose, observe that 
the right-hand side version of \equref{towards ternary} reads 
$$q\x\Oo^{(*,st)}_\bullet\vert_{(\bullet,\bullet)}(\A_2^*\times^{st}\A_1,\Oo^*_\#(\B,\C)) \stackrel{\equref{mixed-quasi-right}}{\iso}
\Oo^*_\bullet(\A_2, \Oo^{st}_\bullet(\A_1,\Oo^*_\#(\B,\C)))\stackrel{\equref{mixed-quasi-right}}{\iso}
\Oo^*_\bullet(\A_2, q\x\Oo^{(*,st)}_\bullet\vert_{(\bullet,\bullet)}(\B^*\times^{st}\A_1,\C))$$
and completely analogously one gets the right-hand side version of \equref{ternary} 
$$q_3\x\Oo^{(*,st)}_\bullet(\B^*\w\times^*\A_2^*\times^{st}\A_1,\C)\iso q\x\Oo^{(*,st)}_\bullet(\A_2^*\times^{st}\A_1, \Oo^*_\#(\B,\C)).$$
Consequently, analogously as in \thref{mixed-quasi n} we get 
\begin{multline} \label{right closed}
q_n\x\Oo^{(*,st)}_\bullet(\Pi_{i=1}^n\A_i,\C)=q_n\x\Oo^{(*,st)}_\bullet(\Pi_{i=1}^{n-1}\A_i^*\w\times^{st}\A_n,\C)
\iso\Oo^{st}_\bullet(\A_n,\Oo^*_\bullet(\A_{n-1},...,\Oo^*_\#(\A_1,\C)...) \\
  \iso q_{n-1}\x\Oo^{(*,st)}_\bullet(\Pi_{i=2}^n\A_i,\Oo^*_\#(\A_1,\C)),
\end{multline}
which extends to 
$$q_n\x\Oo^{*}_\bullet(\Pi_{i=1}^n\A_i,\C)  \iso q_{n-1}\x\Oo^{*}_\bullet(\Pi_{i=2}^n\A_i,\Oo^*_\#(\A_1,\C)).$$
%

\begin{rem}
Observe from \equref{mixed-quasi-bij} and \equref{mixed-right-bij} that in skew cases one does not encounter the collision 
of the types of transformations and the dimension as in \prref{Crans before}. For example in \equref{mixed-right-bij} above, 
a transformation 1-cell $\alpha$ in $\Oo^*_\#(\A,\C)$ can be seen as a strict Gray functor $2\to\Oo^*_\#(\A,\C)$, which is equivalently a loose $*$-type functor $\crta\alpha: \A\to\Oo^{st}_\bullet(2,\C)$, so there is no strictness of $\crta\alpha$ 
that would force the questionable equation that occurred in \prref{Crans before}. The situation in the case of \equref{mixed-quasi-bij} is symmetric. 
\end{rem}

\subsection{Substitution for quasi-functors (for double and Gray-categories)} \sslabel{substitution Gray}

We prove now the substitution for quasi-functors both for double categories (announced in \prref{subst}) and for Gray-categories. 
 

\begin{prop} \prlabel{subst-gen}
Given quasi-functors $F_i:\A_{i1}\times...\times\A_{im_i}\to\B_i$ of $m_i$-variables with $i=1,...,n, m_i\geq 2, n\geq 2$ and 
a quasi-functor $G:\B_1\times...\times\B_n\to\C$ of $n$-variables, all of type $(*,\bullet)$, for 
$\Oo\in\{2\x\Cat,Dbl,G\x\Cat\}$. 
Then 
\begin{enumerate}
\item 
the composition 
$$\Pi_{j=1}^{m_1}\A_{1j}\times...\times\Pi_{j=1}^{m_n}\A_{nj}\stackrel{F_1\times...\times F_n}{\longrightarrow}\B_1\times...\times\B_n
\stackrel{G}{\to}\C$$
is a quasi-functor of $m_1+..+m_n$-variables of type $(*,\bullet)$;
\item the above composition is of type $((st,*),\bullet)$ if so are $F_1$ and $G$, and of type $((*,st),\bullet)$ if so are $F_n$ and $G$.  
\end{enumerate}
\end{prop}

\begin{proof}
We should check the quasi-functor property for any three fixed variables. 
The non-trivial case is when they belong to the domains of two or three different quasi-functors $F_i$ 
(if they all belong to some $\A_i=\A_{i1}\times...\times\A_{im_i}$, then the quasi-functor property is fulfilled by $F_i$ and then preserved in the single variable 
in $\B_i$ by $G$). In the latter case $G$ has at least three variables, so we are in a situation  
$G(F_i(...),F_j(...), F_k(...))$, without loss of generality, and the quasi-functor property for 1h-cells $f,g,h$ in $\A_i,\A_j,\A_k$, respectively, corresponds to the quasi-functor property of $G$ for 1h-cells $g_i, g_j, h_k$ of $\B_i,\B_j,\B_k$, which we know is fulfilled, since $G$ is a quasi-functor. 

If $f,g,h$ belong to the domains of two quasi-functors $F_i, F_j$, we have to check, without loss of generality, whether 
$G(F_i(-,-),F_j(-,-))$ is a quasi-functor in any of the four choices of three variables. 
%
%
%
%
In the double categorical case, 
the first axiom of \deref{double qf 3} for $G(F_i(-,-),F_j(-,D))$, for an object $D$, to be a quasi-functor is equivalent to the condition that $G(F_i(f,g),-)=G(a,-)$, for 1h-cells $f,g$ and the corresponding 2-cell $a$, is a modification of the composition of transformations 
$G(F_i(f,B),-)=G(\crta f, -)$ and $G(F_i(A',g),-)=G(\crta g, -)$ (see \cite[Lemma 2.3]{Fem} for definition), which we know is true by 
\prref{quasi-fun}. 
The axioms of modifications of transformations of double functors can be found in \deref{modif-hv} in the Appendix A. 
(Namely, in the modification {\em i.e.} quasi-functor condition in one side of the relevant axiom there appear the following 2-cells: 
$G(F_i(f,g),F_j(C,D))$ as a modification component, and $G(F_i(f,B),F_j(h,D))=G(\crta f, \crta h)$ and $G(F_i(A',g),F_j(h,D))=
G(\crta g, \crta h)$ as transformation 2-cell components at the 1h-cell $F_j(h,D)=\crta h$.) The situation for $G(F_i(-,-),F_j(C,-))$ 
is similar. For $G(F_i(A,-),F_j(-,-))$ and $G(F_i(-,B),F_j(-,-))$ the axiom for a quasi-functor corresponds to the axiom for $G(-, F_j(g,h))$ 
to be a modification of the similar sort as above, which is again true due to \prref{quasi-fun}. 

The second axiom for a quasi-functor of three variables is a vertical analogue of the first one. The remaining four axioms are proved in a 
similar way, using proper modifications of double categories (the first two are globular modifications).

For Gray-categories, the proof for the first part is analogous, though the analogy for the ternary quasi-functor property works only partially. Namely, it
now involves 
a 3-cell $\axiomref{f,g,h}$ and nine axioms (from \cite[Appendix D.2]{BL}), while the modification condition for $G(a,-)$ (with respect to the composition of transformations $G(\crta f, -)$ and $G(\crta g, -)$), consists of a 3-cell and three axioms (see \cite[Appendix A.3]{BL}, 
and our Appendix C.2). Hereby, the former nine axioms involve composability conditions for 1-cells, compatibility with 2-cells and the identity 1-cells in all the three entries of $(f,g,h)$, while the latter three modification axioms express the same three types of compatibilities only in the single variable in $G(a,-)=G(F_i(f,g),-)$. The rest of the conditions should follow from $G(F_i(-,g)^2_{f',f},-)$ and $G(F_i(f,-)^2_{g',g},-)$ being 
perturbations, for the cocycle 3-cells $F_i(-,g)^2_{f',f}$ and $F_i(f,-)^2_{g',g}$, and from $G(-,c)$, for a 2-cell $c$ in $\C$, being a 
modification. However, some technical restrictions appear, and similar occurres in the checking of the 4-ary quasi-functor condition 
in the cases (ii) and (iii) below. For this reason, we recur to a more conceptual proof. 
 

As a matter of fact, the only cases that we have to work out in both the double categorical and Gray-category case, namely, to prove that 
$G(F_i(-,-),F_j(-,-))$ is a quasi-functor, we can also manage by seeing $G(F_i(-,-),F_j(C,-))$ and 
$G(F_i(-,-),F_j(-,D))$ as objects of the right-hand side of \equref{ternary}, and $G(F_i(A,-),F_j(-,-))$ and $G(F_i(-,B),F_j(-,-))$ 
as objects of the right-hand side of \equref{towards ternary} (abstracting the strictification in the leftmost argument, if necessary), 
which they are, since $F_i$, respectively $F_j$, are quasi-functors.  

This proves the ternary quasi-functor property for $G(F_1\times...\times F_n)$. For the 4-ary property, we distinguish the four cases 
of distributions of four variables and 1-cells $f,g,k,l$:
\begin{enumerate}[i)]
\item $G(F_1(-,-), F_2(-,-), F_3(-,-), F_4(-,-))$ with single distribution $(f,g,h,k)$ (up to choosing the left or right position in each $F_i$), 
\item $G(F_1(-,-), F_2(-,-), F_3(-,-))$ with  distributions $((f,g),h,k), (f,(g,h),k),(f,g,(h,k))$,
\item $G(F_1(-,-), F_2(-,-))$ with single distribution $((f,g),(h,k))$, 
\item $G(F_1(-,-,-), F_2(-,-,-))$ with distributions $((f,g,h),k)), (f,(g,h,k))$ (with three subcases depending on the positions of $k$ in $F_2$ and $f$ in $F_1$, respectively).
\end{enumerate}
The case (i) immediately holds by the 4-ary quasi-functor property of $G$ for the induced 1-cells $(g_1,g_2,g_3,g_4)$ in $\B_i$. 
The rest of the cases is proved using (\ref{towards 4-ary}) and \equref{4-ary}, analogously as we proved the ternary property in the paragraph above, where in (ii) and (iii) we use the isomorphisms \vspace{-0,3cm}
$$
q\x\Oo^{(st,*)}(\A_1\times\A_2, \Oo^*_\bullet(\A_3,\Oo^*_\bullet(\B,\C)))  \stackrel{\equref{quasi}}{\iso} 
q\x\Oo^{(st,*)}(\A_1\times\A_2, q\x\Oo^*_\bullet(\A_3\times\B,\C)) \vspace{-0,2cm}
$$
and \vspace{-0,2cm}
$$
\Oo_{st}(\A_1, q\x\Oo^*_\bullet\vert_{(\bullet,\bullet)}(\A_2\times\A_3,\Oo^*_\bullet(\B,\C)))  \stackrel{(\ref{towards 4-ary})}{\iso}
\Oo_{st}(\A_1, \Oo_{st}(\A_2, q\x\Oo^*_\bullet\vert_{(\bullet,\bullet)}(\A_3\times\B,\C))) \vspace{-0,2cm}
$$
respectively (observe that the left-hand sides appear in (\ref{towards 4-ary})). Namely, 
\begin{itemize}
\item in (ii):
\begin{itemize}
\item the case $((f,g),h,k)$ is to be understood as $G(F_1(f,g), F_2(h,X),F_3(k,Y))$ (with variations $F_2(X,h)$ and $F_3(Y,k)$),
see $G(F_1(-,-), F_2(-,X), F_3(-,Y))$ as living in 
$q\x\Oo^{(st,*)}(\A_1\times\A_2, \Oo^*_\bullet(\A_3,\Oo^*_\bullet(\B,\C)))$;
\item the case $(f,(g,h),k)$ is to be understood as $G(F_1(f,X), F_2(g,h),F_3(k,Y))$ (with variations $F_1(X,f)$ and $F_3(Y,k)$),
see $G(F_1(-,X), F_2(-,-), F_3(-,Y))$ as living in 
$\Oo_{st}(\A_1, q\x\Oo^*_\bullet\vert_{(\bullet,\bullet)}(\A_2\times\A_3,\Oo^*_\bullet(\B,\C)))$;
\item the case $(f,g,(h,k))$ is to be understood as $G(F_1(f,X), F_2(g,Y),F_3(h,k))$ (with variations $F_1(X,f)$ and $F_2(Y,g)$),
see $G(F_1(-,X), F_2(-,Y), F_3(-,-))$ as living in 
$\Oo_{st}(\A_1, \Oo_{st}(\A_2, q\x\Oo^*_\bullet\vert_{(\bullet,\bullet)}(\A_3\times\B,\C)))$;
\end{itemize}
\item in (iii) the single case $((f,g),(h,k))$ is to be understood as $G(F_1(f,g), F_2(h,k))$,
see $G(F_1(-,-), F_2(-,-))$ as living in $q\x\Oo^{(st,*)}(\A_1\times\A_2, q\x\Oo^*_\bullet(\A_3\times\B,\C))$;
\item in (iv): 
\begin{itemize}
\item for $((f,g,h),k)$, in the subcase $G(F_1(f,g,h), F_2(k,X,Y))$ (with variations $F_2(X,k,Y)$ and $F_2(X,Y,k)$) recognize 
$G(F_1(-, -,-), F_2(-,X,Y))$ as living in $q_3\x\Oo^{(st,*)}_\bullet(\A_1{}^{st}\w\times^*\A_2^*\times^*\A_3,\Oo^*_\bullet(\B,\C))$; 
\item for $(f,(g,h,k))$, in the subcase $G(F_1(f,X,Y), F_2(g,h,k))$ (with variations $F_1(X,f,Y)$ and $F_1(X,Y,f)$) recognize 
$G(F_1(-,X,Y), F_2(-,-,-))$ as living in 
$\Oo_{st}(\A_1, q_3\x\Oo^*_\bullet\vert_{(\bullet,\bullet,\bullet)}(\A_2\times\A_3\times\B,\C))$. 
\end{itemize}
\end{itemize}

\medskip

In order for $G(F_1\times...\times F_n)$ to be of type $((st,*),\bullet)$, it should be strict in $\A_{11}$. This entry is controlled firstly by 
$F_1:\Pi_{j=1}^{m_1}\A_{1j}\to\B_1$ - thus $F_1$ should be strict in $\A_{11}$ - and then we further need strictness in $\B_1$, which is 
controlled by $G:\B_1\times...\times\B_n\to\C$, thus both $F_1$ and $G$ should be of type $((st,*),\bullet)$. The case 
$((',st),\#)$ is now clear. 
\qed\end{proof}

For the substitution for the multicategory and skew-multicategory $\OO$ that we are constructing, additionally to the above proposition 
we also need to study the cases when some of the $F_1,.., F_n,G$ above instead of being quasi-functors, are nullary or unary maps. 
We find the following cases: 
\begin{itemize} 
\item all $F_i$'s are nullary: then the composition is nullary and there is nothing to prove, \vspace{-0,2cm}
\item some of $F_i$'s, say $k\geq 0$ many, are nullary and the rest are identities: this is the component quasi-functor $G_{n-k}$ of $G$ of $n-k$ variables, 
\item any unary map composed into any variable of a quasi-functor $G$ solely alters (functorially) the component quasi-functor $G_{n-1}$, 
so that the composition is a quasi-functor; 
\item some of $F_i$'s, say $k\geq 0$ many, are nullary and the rest are identities and/or unary maps: in the component quasi-functor $G_{n-k}$ of $G$ the variables where unary maps are composed change functorially - the composite is again a quasi-functor; 
\item some of $F_i$'s, say $k\geq 0$ many, are nullary and the rest are quasi-functors: 
in the above proposition replace $G$ by $G_{n-k}$, 
\vspace{-0,2cm}
\item if there is a single quasi-functor $F$: in the proof of the above proposition at the places where quasi-functor properties of non-trivial $F_i$ were used, replace the quasi-functor property of $F_i$ with the corresponding (restricted) quasi-functor property of $G$; \vspace{-0,2cm}
\item 
if both $F$ and $G$ are unary, the composition is clearly a map $GF$. \vspace{-0,2cm}
\end{itemize}
It is obvious that in order for the composition to be tight it is necessary that the leftmost non-nullary $F_i$ and $G$ be tight. 

The associativity of substitution involving nullary and unary maps we discussed in \ssref{nullary/unary}. The associativity of the rest of combinations of substitutions clearly holds, as all quasi-functors of $k$-variables consists of $k$ one variable functors.

\subsection{Gray skew-multicategories} \sslabel{skew mut gen}

We define left Gray skew-multicategory for $\Oo\in\{2\x\Cat, Dbl, G\x\Cat\}$ as follows:
\begin{itemize} 
\item take $\Oo$ for the collection of objects, \vspace{-0,2cm}
\item let $\M_0^l(-;\A)$ be the set of objects of $\A$, \vspace{-0,2cm}
\item set $\M_1^t(\A;\B)=\Oo_{st}(\A,\B)$ and $\M_1^l(\A;\B)=\Oo_*(\A;\B)$, \vspace{-0,2cm}
\item set $\M_n^t(\A_1,...,\A_n;\C)=q_n\x\Oo^{(st,*)}_\bullet(\A_1\times...\times\A_n,\C)$ for tight multimaps, \vspace{-0,2cm}
\item set $\M_n^l(\A_1,...,\A_n;\C)=q_n\x\Oo^*_\bullet(\A_1\times...\times\A_n,\C)$ for loose multimaps, \vspace{-0,2cm}
\item the unary identities in $\M_1^l(\A;\A)$ are obviously tight, \vspace{-0,2cm}
\item substitution is guaranteed by \prref{subst-gen}. \vspace{-0,2cm}
\end{itemize}
We denote left Gray skew-multicategory for $\Oo$ by $\OO^{(st,*)}_\bullet$.  

A right Gray skew-multicategory for $\Oo$ differs from the left one only in the tight $n$-ary part: now $\M_n^t(\A_1,...,\A_n;\C)=
q_n\x\Oo^{(*,st)}_\bullet(\A_1\times...\times\A_n,\C)$.  
The right Gray skew-multicategory for $\Oo$ we denote by $\OO^{(*,st)}_\bullet$.

\section{Properties of Gray skew-multicategories}

In this section we draw some conclusions about the above constructed skew-multicategories. 


\subsection{Skew-multicategories are closed on one side} \sslabel{skew m are closed}

Let  $\Oo\in\{2\x\Cat, Dbl, G\x\Cat\}$ as above. \\
Now that we have skew-multicategories $\OO^{(st,*)}_\bullet$ and 
$\OO^{(*,st)}_\bullet$ we may say, due to \thref{mixed-quasi n} that $\OO^{(st,*)}_{\bullet}$ is left closed. 

\medskip

We saw that quasi-functors $\crta{ev}:\B^*\w\times^{st}\Oo^*_\#(\B, \C)\to\C$ of type $((*,st),\bullet)$ induce 
natural isomorphisms \equref{mixed-quasi-right} and (\ref{right closed}). 
This makes the skew-multicategory $\OO^{(*,st)}_\bullet$ right closed with inner-hom $\Oo^*_\#(-,-)$. 

\medskip

We may conclude that there is no proper Gray skew-multicategory that is biclosed, but that the multicategory $\OO^{st}_{\bullet}$ is closed on both sides, {\em i.e.} biclosed. 
Since by definition the associated multicategory of loose maps of a left/right closed skew-multicategory is left/right closed, 
we have more generally that the multicategory $\OO^*_{\bullet}$ is biclosed. However, to these considerations we should include the restrictions on dimension that we saw in \prref{Crans before}. Then we may state:

\begin{prop} \prlabel{dim-biclosed}
The following multicategories are biclosed: 
\begin{itemize}
\item $\OO^{st}_{\bullet}$ for $m\leq 2$ and $\Oo=3\x\Cat$, 
\item $G\x\Cat^{st}_{st}$, and 
\item $G\x\Cat^*_{\bullet}$ for $*\not= st$ 
\end{itemize}
with left and right inner-homs $\Oo^*_\bullet(-,-)$ and $\Oo^*_\#(-,-)$ in each respective case, so that $\#$ is opposite to $\bullet$, 
and if $\bullet=\#=ps$, then the two inner-homs coincide. 
\end{prop}

Our findings are in accordance with Theorems 5.1 and 5.2 of \cite{BL}, where it was proved that the skew-multicategory $\Lax=\OO^{(*,st)}_\bullet$ is left-closed and that the multicategory $\Lax^l=\OO^{ps}_{lx}$ is biclosed, where $\Oo=G\x\Cat$. They are also in accordance with 
the statement from the last paragarph of \cite[Section 7]{BL}, saying that the multicategory $\Lax^\sharp=G\x\Cat^{st}_{ps}$ is not biclosed, 
but this is so already by the Crans's argument that we employed in \prref{Crans before}.

\subsection{Skew meta-product}

As announced, quasi-functors of types $((st,*),\bullet)$ and $((*,st),\bullet)$ induce meta-products $\ot^{(st,*)}_\bullet$ 
and $\ot^{(*,st)}_\bullet$, respectively. Theese meta-products are obtained by tightening the meta-product $\ot^*_\bullet$ 
in the left respectively right entry. In the case of double categories and type $((st,lx),o-l)$, by tightening \deref{tensor} the four 2-cells 
``structure 2-cells of type $*=lx$'' in \equref{4 laxity 2-cells} become two 2-cells and two axioms 
\begin{equation} \eqlabel{2 laxity 2-cells}
(A\ot k')(A\ot k)\stackrel{(A\ot -)_{k'k}}{\Rightarrow} A\ot (k' k), \qquad 1_{A\ot B}\stackrel{(A\ot -)_B}{\Rightarrow} A\ot 1_B
\end{equation}
\begin{equation} \eqlabel{2 strictness 2-cells} 
(K'\ot B)(K\ot B)=(K' K)\ot B,\quad 1_{A\ot B}=1_A\ot B
\end{equation}
and the second of the two diagrams for ``axioms for horizontal composition of type $*$ of 2-cells'' in \deref{tensor} simplifies into the axiom 
$$\zeta'\zeta\ot B=[\zeta\ot B \vert\zeta' \ot B].$$

\medskip

The list of the defining data and relations for the meta-product of Gray-categories is much longer and we omit it. It consists of the data for 
$A\ot-$ and $-\ot B$ to be pseudo-Gray-functors, for $A\in\A, B\in\B$ (see Appendix A.1 of \cite{BL}) and three families of 2-cells and for families of 3-cells coming from those from Section 4.3, and which satisfy axioms analogous to those of Appendix D.1. 

\medskip

The meta-product $\ot^{(st,*)}_\bullet$ satisfies the following universal property.

\begin{thm} \thlabel{univ}
Let  $\Oo\in\{2\x\Cat, Dbl, G\x\Cat\}$. There is an isomorphism in $\Oo$ 
$$q\x \Oo^{(st,*)}_{\bullet}(\A{}^{st}\w\times^*\B,\C)\iso \Oo^{st}_\bullet(\A\ot^{(st,*)}_\bullet\B,\C).$$
\end{thm}

\begin{proof}
The proof is analogous to that of \cite[Proposition 6.3]{Fem}. We start the proof for double categories. Namely,  
for $H\in q\x \Oo^{(st,*)}_{\bullet}(\A{}^{st}\w\times^*\B,\C)$ define $\crta{H}\colon\A\ot^{(st,*)}_\bullet\B\to\C$ by 
$\crta{H}(A\ot y)\colon\hspace{-0,2cm}=H(A,y)=(y,A)$ and $\crta{H}(x\ot B)\colon\hspace{-0,2cm}=H(x,B)=(B,x)$ 
for cells $x$ in $\A$ and $y$ in $\B$ of all degrees (for general $\Oo$ it would be $m+1$ of them in the enriched and $2m$ 
in the internal case). Extend $\crta H$ to a strict functor (in particular, $\crta H(1_{A\ot B})
=1_{H(A,B)}=1_{(B,A)}$). Define $\crta H((A\ot-)_B)\colon\hspace{-0,2cm}=(-,A)_B, \crta H((A\ot-)_{k'k})\colon\hspace{-0,2cm}=(-,A)_{k'k}$ 
and observe that these are $*$-loose 2-cells, whereas for the other entry we have that $(-\ot B)_A$ and $(-\ot B)_{K'K}$ are identity 2-cells, 
so by strictness of $\crta H$ their images by it are identities, too. At last, images by $\crta H$ of the 2-cells $K\ot k, K\ot u, U\ot k$ 
and $U\ot u$ define by $(k,K), (u,K), (k,U)$ and $(u,U)$, respectively. 

Conversely, given $G\in \Oo^{st}_\bullet(\A\ot^{(st,*)}_\bullet\B,\C)$, define 
$\crta{(-,A)}\colon\B\to\C$ by $\crta{(y,A)}\colon\hspace{-0,2cm}=G(A\ot y)$, and by the two globular 2-cells: 
$\crta{(-,A)}_{k'k}\colon G(A\ot k')G(A\ot k)=G((A\ot k')(A\ot k))\stackrel{G((A\ot-)_{k'k})}{\Rightarrow}G(A\ot k'k)$ and 
$\crta{(-,A)}_B\colon\hspace{-0,2cm}= G(A\ot -)_B$, and analogously for $\crta{(B,-)}\colon\Aa\to\Cc$. 
This equips $\crta{(-,A)}$ with a structure of a loose functor (here written in the form of lax), and $\crta{(B,-)}$ 
with a structure of a strict functor. Define the 2-cells 
$\crta{(k,K)}, \crta{(u,K)}, \crta{(k,U)}$ and $\crta{(u,U)}$ in the obvious way, then the laws from \deref{skew-guasi} for 
$\crta{(-,A)}$ and $\crta{(B,-)}$ to make a $(st,*)$-mixed quasi-functor 
pass {\em mutatis mutandi} from the defining relations of $\A\ot^{(st,*)}_\bullet\B$, since $G$ is a strict functor.  

\smallskip

The proof for Gray-categories is analogous. 
By the way how  $\A\ot \B$ was constructed (see \deref{gen-quasi} and the discussion above it) it is clear that $\crta H$ defined 
through $H$ analogously to the above described double-categorical case makes a strict Gray functor. 
That the converse works in an analogous way is readily seen. 

\smallskip

The 1-1 correspondence on higher cells follows through, accordingly, in the analogous way as it was proved in the double 
categorical case with lax double functors in \cite[Proposition 6.3]{Fem}. For reader's convenience we only list how the necessary 
structure cells for transformations and modifications are defined. Given a transformation (of whatever type) 
$\theta=(\theta^A, \theta^B)_{A\in\A, \\B\in\B}$ between quasi-functors $H\Rightarrow H'$ 
we define a transformation (of the same type) $\Sigma\colon\crta{H}\Rightarrow\crta{H'}$ by setting at objects 
$\Sigma(A\ot B)=\theta^A_B=\theta^B_A$, at 1-cells $\Sigma_{A\ot k}=\theta^A_k, \Sigma_{K\ot B}=\theta^B_K$, at 2-cells 
$\Sigma_{K\ot k}=\theta_{(k,K)}, \Sigma_{\zeta\ot B}=\theta^B_\zeta, \Sigma_{A\ot \omega}=\theta^A_\omega$, and at composable 
1-cells: $\Sigma_{A\ot k', A\ot k}=\theta^A_{k',k}, \Sigma_{K'\ot B, K\ot B}=\theta^B_{K',K}$. For the converse, given a transformation $\tilde\Sigma\colon G\Rightarrow G'$, define $\zeta=(\zeta^A, \zeta^B)_{A\in\A, \\B\in\B}$ in the obvious (converse) way. Finally, given a modification $\tau=(\tau^A, \tau^B)_{A\in\A, \\B\in\B}$ on the quasi-functor side, define 
$\Theta(A\ot B)=\tau^A_B=\tau^B_A$ at objects, and $\Theta_{A\ot k}=\tau^A_k, \Theta_{K\ot B}=\tau^B_K$ at 1-cells. 
For the converse, formulate the obvious (converse) definition. All the necessary axioms are straightforwardly seen to hold. 
\qed\end{proof}

\bigskip

Observe that in the above proof the mixed type $(st,*)$ of the quasi-functors and meta-product $\ot^{(st,*)}_\bullet$ can be changed into fully weak type $*$, so that we also have the following isomorphism in $\Oo$ 
\begin{equation} \eqlabel{any-strict-Gray}
q\x \Oo^*_\bullet(\A\times\B,\C)\iso \Oo^{st}_\bullet(\A\ot^*_\bullet\B,\C).
\end{equation}

\bigskip

Analogously, the universal property that the meta-product $\ot^{(*,st)}_\bullet$ satisfies is involved in the following isomorphisms in $\Oo$
\begin{equation} \eqlabel{right univ ot}
q\x \Oo^{(*,st)}_{\bullet}(\A^*\w\times^{st}\B,\C)\iso \Oo^{st}_\bullet(\A\ot^{(*,st)}_\bullet\B,\C).
\end{equation}

\subsection{Skew-multicategories are representable on one side} \sslabel{sk-m rep}

Let $m\leq 3$ of $\Oo$. \\
We saw before that because of \equref{any-strict} the isomorphism in (\ref{eq2}) needed for left representability of $\ot^*_{o-l}$ was possible only for $*=st$. For the mixed types of meta-products we obtain that $\ot^{(st,*)}_\bullet$ is left representable and 
that $\ot^{(*,st)}_\bullet$ is right representable.

\begin{prop} \prlabel{left rep}
There is an isomorphism in $\Oo$ 
$$q_n\x \Oo^{(st,*)}_{\bullet}((\A\ot^{(st,*)}_\bullet\B){}^{st}\w\times^*(\Pi_{i=1}^{n-1}\C_i),\D)  
\iso q_{n+1}\x \Oo^{(st,*)}_{\bullet}(\A{}^{st}\w\times^* \B\times^*(\Pi_{i=1}^{n-1}\C_i),\D).$$
\end{prop}

\begin{proof}
To both sides of the isomorphism that is to be proved we apply iteratively \thref{mixed-quasi n} 
and at the end on the left-hand side we apply \equref{mixed-quasi} 
to get to an equivalent statement $\Oo^{st}_\bullet(\A\ot^{(st,*)}_\bullet\B, \Oo^*_\bullet(\C_1, I))\iso
\Oo^{(st,*)}_\bullet(\A{}^{st}\w\times^*\B, \Oo^*_\bullet(\C_1, I))$, where $I$ is an inner-hom generated by accumulation. 
By \thref{univ} we know this is true, so the claim follows. 
\qed\end{proof}

Set 
\begin{equation} \eqlabel{acu left G}
\crta\ot^n_\bullet\A_i=(..(\A_1\ot\A_2)\ot...\ot\A_{n-1})\ot\A_n
\end{equation} 
whereby $\ot=\ot^{(st,*)}_\bullet$, and observe that by iterative applications of the above isomorphism we obtain 
$$q_{l+1}\x \Oo^{(st,*)}_{\bullet}((\crta\ot^k_\bullet\A_i){}^{st}\w\times^*(\Pi_{i=1}^{l}\B_i),\C)  
\iso q_{k+l}\x \Oo^{(st,*)}_{\bullet}(\A_1{}^{st}\w\times^* \Pi_{i=2}^k\A_i\times^*(\Pi_{i=1}^l\B_i),\C).$$

\smallskip

\begin{thm} \thlabel{left cl-rep}
The left skew-multicategory $\OO^{(st,*)}_{\bullet}$ is left closed and left representable. 
\end{thm}

\begin{proof}
We clearly have $q_{n+1}\x\Oo^{(st,*)}_{\bullet}(\I\times\Pi_{i=1}^n\A_i,\B) \iso q_n\x\Oo^*_{\bullet}(\Pi_{i=1}^n\A_i,\B)$ 
where $\I$ is the initial object in $\Oo$ (a trivial double or Gray-category). Left closedness follows 
by \thref{mixed-quasi n}, and left representability by \thref{univ} and \equref{acu left G}.  
\qed\end{proof}

\bigskip

Analogously as in \prref{left rep}, due to (\ref{right closed}), \equref{mixed-quasi-right} and \equref{right univ ot} we obtain 
$$q_n\x \Oo^{(*,st)}_{\bullet}((\Pi_{i=1}^{n-1}\C_i)^*\w\times^{st}(\A\ot^{(*,st)}_\bullet\B),\D)  
\iso q_{n+1}\x \Oo^{(*,st)}_{\bullet}((\Pi_{i=1}^{n-1}\C_i)^*\w\times\A^*\w\times^{st} \B,\D).$$
Accordingly, we get that the right skew-multicategory $\OO^{(*,st)}_{\bullet}$ is (right closed and) right representable.

\begin{rem}
Observe that we can not carry out the same computation as around (\ref{eq1}) to relate bijectively $\Oo^?((\A\ot\B)\ot\C,\D)$ with 
$\Oo^?(\A\ot(\B\ot\C),\D)$ using meta-product induced by mixed quasi-functors. This is because in (\ref{eq1}) we need to use 
left representability (from \prref{left rep}), which we have for $\ot^{(st,*)}_{\bullet}$, while in (\ref{eq3}) we need to use 
right representability (from the last isomorphism above), which we have for $\ot^{(*,st)}_{\bullet}$, so we do not have a suitable tensor product to complete this reasoning. Though, this is how far we may go: 
\begin{align}
\Oo^{(st,*)}_{\bullet}((\A\ot^{(st,*)}_{\bullet}\B)\ot^{(st,*)}_{\bullet}\C,\D) & \stackrel{\thref{univ}}{\iso}
q\x\Oo^{(st,*)}_{\bullet}((\A\ot^{(st,*)}_{\bullet}\B)^{st}\w\times^*\C,\D) \nonumber \\
& \stackrel{\prref{left rep}}{\iso} q_3\x\Oo^{(st,*)}_{\bullet}(\A^{st}\w\times^*\B\times^*\C,\D). \nonumber
\end{align}
Now setting $\D=\A\ot^{(st,*)}_{\bullet}(\B\ot^{(st,*)}_{\bullet}\C)$ consider $J_{\A,\B\ot^{(st,*)}_{\bullet}\C}\comp J_{\B,\C}$ in 
$q_3\x\Oo^{(st,*)}_{\bullet}(\A^{st}\times^*\B\times^*\C,\A\ot^{(st,*)}_{\bullet}(\B\ot^{(st,*)}_{\bullet}\C))$, where $J$ is analogous as in \equref{J}. The corresponding functor living in 
$\Oo^{(st,*)}_{\bullet}((\A\ot^{(st,*)}_{\bullet}\B)\ot^{(st,*)}_{\bullet}\C,\A\ot^{(st,*)}_{\bullet}(\B\ot^{(st,*)}_{\bullet}\C))$
is then the associator for the skew-monoidal structure, according to \cite[Section 6.2]{BLack}. 
\end{rem}

\subsection{The induced closed skew-monoidal categories $\Oo_{st}$} \sslabel{I and II}

In this subsection we draw conclusion about closedness and (skew-)monoidality of the categories obtained from the above studied (skew-)multicategories. Recall that as before the starting assumption is that for $\A,\B\in\Oo$ the classes $\Oo^*_\bullet(\A,\B)$ and 
$\Oo^*_\#(\A,\B)$ are objects in $\Oo$.

\subsubsection{Closedness of $\Oo_{st}$} \ssslabel{closedness}

As argued in \cite[Section 2.6]{BL} for any skew-multicategory $\M$, for those dimensions $m$ of $\Oo$ for which Gray skew-multicategories have a tight binary map classifier (but also a loos one) and are left (right) closed, one finds that the category $\Oo_{st}$ is left 
(right) closed with inner-hom $\Oo^*_\bullet(-,-)$, respectively $\Oo^*_\#(-,-)$, when they are bifunctors. Namely, we have 
$$\M_1^{t/l}(\A\ot\B;\C)\iso\M^{t/l}_2(\A,\B;\C)\iso\M^{t/l}_1(\A;[\B,\C])$$
for left skew-multicategories, and similarly $\M_1^{t/l}(\A\ot\B;\C)\iso\M^{t/l}_1(\B;\{\A,\C\})$ for right ones, with apriori a  different inner-hom. 
Clearly, if $\M$ is both left and right closed, $\Oo_{st}$ is too.

We have a tight binary map classifier and left (right) closed Gray skew-multicatego-\\ries: 
$$(\OO^{(st,*)}_\bullet, \ot^{(st,*)}_\bullet), \,\,\,\,  
(\OO^{(*,st)}_\bullet, \ot^{(*,st)}_\bullet), \,\, 
\quad\text{and hence also}\quad (\OO^{st}_\bullet, \ot^{st}_\bullet),$$ 
with $m\leq 3$, and by \equref{quasi-set*} and \equref{any-strict-Gray} also the loose case 
$$(\OO^*_\bullet, \ot^*_\bullet).$$
For double categories, from \cite[Section 3.1]{Fem} it may be deduced that $Dbl^{*}_{\bullet}(-,-): (Dbl_{st})^{op}\times Dbl_{st}\to Dbl_{st}$ is a bifunctor for any $*$. 
From \cite[Section 6]{Go} and \cite[Subsection 6.2]{BL} we have that 
$$G\x\Cat^{ps}_{lx}(-,-), \,\, G\x\Cat^{ps}_{ps}(-,-): (G\x\Cat_{st})^{op}\times G\x\Cat_{st}\to G\x\Cat_{st}$$
are bifunctors. As a matter of fact, as we will argue in \ssref{lax and strict}, the Gray-categorical bifunctors 
$G\x\Cat^{ps}_{lx}(-,-)$ and $G\x\Cat^{ps}_{ps}(-,-)$ can be substituted by more generally constructed $G\x\Cat^{*}_{lx}(-,-)$ and 
$G\x\Cat^{*}_{ps}(-,-)$, as in the double categorical case. For $\Oo$ of any dimension one has that $-\ot^{st,*}_\bullet\B:
\Oo_{st}\to\Oo_{st}$ is a functor (even $\Oo_{ps}\to\Oo_{ps}$), in virtue of \equref{2 strictness 2-cells}. Now, given that 
$\M_1^t(\A;\B)=\Oo_{st}(\A,\B)$, 
and substituting the tensor and inner-hom that we developed, the above isomorphisms read 
$$\Oo_{st}(\A\ot^{(st,*)}_\bullet\B, \C)\iso \Oo_{st}(\A, \Oo^{*}_\bullet(\B,\C)), \qquad 
\Oo_{st}(\A\ot^{(*,st)}_\bullet\B, \C)\iso \Oo_{st}(\B, \Oo^{*}_\#(\A,\C)).$$
We even proved in \equref{mixed-quasi} and \thref{univ} (left case) and \equref{mixed-quasi-right} and \equref{right univ ot} (right case) 
that we have natural isomorphisms in $\Oo$ 
$$\Oo^{st}_{\bullet}(\A\ot^{(st,*)}_\bullet\B,\C)\iso\Oo^{st}_{\bullet}(\A,\Oo^*_{\bullet}(\B,\C)), \qquad
\Oo^{st}_{\bullet}(\A\ot^{(*,st)}_\bullet\B,\C)\iso\Oo^{st}_{\bullet}(\B, \Oo^{*}_\#(\A,\C)).$$
%

\subsubsection{Adding (skew-)monoidality} \ssslabel{skew}

On the other hand, 
in \cite[Section 6.2]{BLack} it was shown that left representable skew-multicategories (and similarly right ones) yield skew-monoidal categories, 
and by Hermida's \cite[Theorem 9.8]{Her} there is a 1-1 correspondence between representable multicategories and monoidal categories. Then for the left-sided case we conclude (I) in:  
\begin{enumerate}[(I)]
 \item
  \begin{enumerate}[(a)]
  \item $(\Oo_{st}, \ot^{(st,*)}_\bullet, \Oo^{*}_\bullet(-,-))$ is left closed skew-monoidal, and 
  \item $(\Oo_{st}, \, \ot^{(*,st)}_\bullet, \,\, \Oo^*_\#(-,-))$ is right closed skew-monoidal, 
  \end{enumerate}
\item 
  \begin{enumerate}[(a)]
  \item $(\Oo_{st}, \ot^{st}_\bullet, \Oo^{st}_\bullet(-,-))$ is left closed monoidal, 
  \item $(\Oo_{st}, \, \ot^{st}_\bullet, \,\, \Oo^{st}_\#(-,-))$ right closed monoidal. 
  \end{enumerate}
\end{enumerate}
Here (II) comes from (I) with $*=st$ and for double categories use \coref{strict-rep}, whereas for Gray-categories see 
\thref{rep O} (the reader may jump to \ssref{lax and strict} to understand the protagonist of the theorem). It seems that in the literature the case when $*=lx$ has not been studied. From the two (II)'s we conclude that the product $\ot^{st}_{\bullet}$ gives a biclosed monoidal structure on $\Oo$. To (II) belong the following examples from the literature. 
The Gray's construction for $\Oo=2\x\Cat$ from \cite[Theorem I,4.14]{Gray}: his inner-homs are made out of strict 2-functors, 
oplax transformations (resp. lax) and modifications, thus $*=st$ and $\bullet=oplx$ (resp. $\bullet=lx$). 
And B\"ohm's construction for $\Oo=Dbl$ from \cite{Gabi}: her inner-hom is made out of strict double functors, (horizontal and vertical) pseudonatural transformations and modifications, thus $*=st$ and $\bullet=ps$. Both Gray's and B\"ohm's constructions are biclosed (observe that by having $\bullet=ps$, the latter construction has the same inner-hom in left and right sided closedness). 

To (I) belong the four constructions for $\Oo=G\x\Cat$ from \cite{BL}. Their first inner-hom is made out of Gray-pseudo-functors, lax transformations, lax modifications, and 
perturbations, thus $*=ps$ and $\bullet=lx,lx$ (see \cite[Theorem 5.8]{BL}). Their second inner-hom is made out of Gray-pseudo-functors, pseudo transformations, pseudo modifications, and perturbations, thus $*=ps$ and $\bullet=ps,ps$ (see \cite[Theorem 6.3]{BL}). The third and fourth inner-hom differ from the latter two in that they come from the full sub-Gray-categories $G\x\Cat^{st}_{lx}(\A,\B)$ and $G\x\Cat^{st}_{ps}(\A,\B)$ whose objects are strict Gray-functors, and the rest of the cells are the same (see Section 7 and Theorem 7.6 of \cite{BL}). Their proof that these induce skew-multicategories 
uses truly pseudo nature of the functors at some places (when constructing substitution, the duality functors are used). For this reason so constructed skew-multicategories are denoted with a superscript $\sharp$ rather than $st$.

\smallskip

In view of \prref{dim-biclosed}, that relies on Crans's argument from the introduction of \cite{Cr}, we see that for Gray-categories the possibility of $\ot^{st}_\bullet$ making a biclosed monoidal product with $\bullet\not=st$ is out ruled, 
and it is reserved 
%
for dimensions $m=1,2$ and 3-categories. They appear in the first column of Table (\ref{table:1}). 

\smallskip

We may now turn to the questions risen in \cite{Fem} and discuss the rest of the versions of categories of double categories. 
From (I) and (II) we may deduce that $(Dbl_{st}, \ot^{st}_{o-l})$ is left closed monoidal with inner-hom 
 $Dbl^{st}_{o-l}(-,-)$ (our \thref{Dbl left closed mon}), and that $(Dbl_{st}, \ot^{(st,ps)}_\bullet)$ is left closed 
skew-monoidal with inner-hom $\Fun^{ps}_\bullet(-,-)$. 
The Gray-categorical version of the former result we will prove in \ssref{lax and strict}. 
%
Recall that we can not obtain a monoidal category using the meta-product $\ot^{lx}_{o-l}$ constructed in \cite{Fem},
nor the meta-product $\ot^{ps}_{\bullet}$: 
by Hermida's \cite[Theorem 9.8]{Her} we need representability of the multicategory with respect to the meta-product $\ot^*_\bullet$, which is possible only when $*=st$, see \seref{meta-lax double}.    
This holds actually for $\Oo$ of any dimension $m$.

\begin{table}[h!]
\begin{center}
\begin{tabular}{ c c c c} 
m & biclosed mon. cat. &  left closed skew-mon. cat. & reference \\ [0.5ex]
\hline
2 & $(2\x\Cat, \, \ot^{st}_{oplx/lx}, \,\, 2\x\Cat^{st}_{oplx/lx}(-,-))$ & & \cite{Gray} \\ [1ex]   
2i & $(Dbl, \, \ot^{st}_{ps}, \,\, Dbl^{st}_{ps}(-,-))$ & & \cite{Gabi} \\ [1ex]   
2i & $(Dbl, \, \ot^{st}_{o-l/l-o}, \,\, Dbl^{st}_{o-l/l-o}(-,-))$ & & \thref{Dbl left closed mon} \\ [1ex]   
2i & & $(Dbl, \, \ot^{(st,*)}_{\bullet}, \,\, \Fun^*_{\bullet}(-,-))$ &  \sssref{skew} \\ [1ex]   
3 & & $(G\x\Cat, \, \ot^{(st,ps)}_{lx}, \,\,  G\x\Cat^{ps}_{lx}(-,-))$  & \cite[Section 5]{BL} \\ [1ex]   
3 & & $(G\x\Cat, \, \ot^{(st,ps)}_{ps}, \,\, G\x\Cat^{ps}_{ps}(-,-))$  & \cite[Section 6]{BL} \\ [1ex]   
3 & & $(G\x\Cat, \, \ot^{\sharp}_{lx}, \,\, G\x\Cat^{\sharp}_{lx}(-,-))$  & \cite[Section 7]{BL} \\ [1ex]   
3 & & $(G\x\Cat, \, \ot^{\sharp}_{ps}, \,\, G\x\Cat^{\sharp}_{ps}(-,-))$  & \cite[Section 7]{BL} \\ [1ex]   
3  & $(3\x\Cat, \, \ot^{st}_{\bullet}, \,\, 3\x\Cat^{st}_{\bullet}(-,-))$  & & \thref{rep O}  \\ [1ex]   
3  & $(G\x\Cat, \, \ot^{st}_{st}, \,\, G\x\Cat^{st}_{st}(-,-))$  & & \thref{rep O}  \\ [1ex]   
\end{tabular}
\caption{Closed skew monoidal structures  for $m\leq 3$} 
\label{table:1}
\end{center}
\end{table}
The symbol $\sharp$ in Table (\ref{table:1}) means that in  $\ot^{\sharp}_{\bullet}$ and $G\x\Cat^{\sharp}_{\bullet}(-,-)$ strict Gray functors are considered, but the construction of this meta-product and inner-hom passes through 'larger'  $\ot^{ps}_{\bullet}$ and $G\x\Cat^{ps}_{\bullet}(-,-)$. Our construction of the last two lines in Table (\ref{table:1}) will be treated in \ssref{lax and strict}.

\subsection{Lax functor classifier conjecture} \sslabel{conjecture}

From \equref{quasi} (with $\bullet=\#$) and \equref{any-strict} with $*=lx$ we obtain 
$$Dbl_{hop}(\A\ot^{lx}_{o-l}\B,\C)\iso Lax_{hop}(\A,Lax_{hop}(\B,\C))$$
and with $*=st$ we get $(C2)$ below. The study of a lax double functor classifier is out of the scope of this paper, for now we conjecture that there is a bijection $(C1)$ in  
\begin{flushleft} 
$(C1) \quad Lax_{hop}(\A,\B)\iso Dbl_{hop}(Q(\A),\B)$

$(C2) \quad Dbl_{hop}(\A\ot^{st}_{o-l}\B,\C)\iso Dbl_{hop}(\A, Dbl_{hop}(\B,\C))$.
\end{flushleft} 
Then from the following sequence of isomorphisms 
\begin{multline*}
Dbl_{hop}(\A\ot^{lx}_{o-l}\B,\C)\iso Lax_{hop}(\A,Lax_{hop}(\B,\C))\stackrel{(C1)}{\iso} \\ 
Dbl_{hop}(Q(\A),Dbl_{hop}(Q(\B),\C))\stackrel{(C2)}{\iso} Dbl_{hop}(Q(\A)\ot^{st}_{o-l}Q(\B),\C)
\end{multline*}
we would obtain 
\begin{equation} \eqlabel{lax d classif}
\A\ot^{lx}_{o-l}\B\iso Q(\A)\ot^{st}_{o-l}Q(\B).
\end{equation} 
The meta-product $\ot^{st}_{o-l}$ in the right-hand side is a biclosed monoidal product for the category $Dbl_{st}$ due to (II). 
It is the double categorical version of the Gray's biclosed monoidal product for 2-categories (see the first line in Table (\ref{table:1})). Yet, $\ot^{lx}_{o-l}$ is not a monoidal product for the category $Dbl_{lx}$, as we argued above. The isomorphism \equref{lax d classif} 
would be merely an alternative description of the product $\A\ot^{lx}_{o-l}\B$, it would be a double categorical version of the 2-categorical isomorphism $\C\Del\D\iso
\C^\dagger\ot_l\D^\dagger$ (146) of \cite{Nik}. Here 
$\C\Del\D$ is the 2-category that strictifies lax functors (in the sense that we proved in \equref{Nikol} 
for the double category $\A\ot^{lx}_{o-l}\B$), and $\ot_l$ is the Gray tensor product for 2-categories (for the inner-hom of strict functors, lax natural transformations and modifications).

\bigskip

For Gray-categories, by \equref{quasi} (with $\bullet=\#$) and \equref{any-strict-Gray} we have 
\begin{equation} \eqlabel{Gray*}
G\x\Cat^{st}_{\bullet}(\A\ot^{*}_{\bullet}\B,\C)\iso G\x\Cat^{*}_{\bullet}(\A,G\x\Cat^{*}_{\bullet}(\B,\C)),
\end{equation} 
so that for $*=ps$ we have the first isomorphism in 
\begin{multline*}
G\x\Cat^{st}_{\bullet}(\A\ot^{ps}_{\bullet}\B,\C)\iso G\x\Cat^{ps}_{\bullet}(\A,G\x\Cat^{ps}_{\bullet}(\B,\C))\\
\iso G\x\Cat^{st}_{\bullet}(Q(\A),G\x\Cat^{st}_{\bullet}(Q(\B),\C))\stackrel{\equref{Gray*}}{\iso}
G\x\Cat^{st}_{\bullet}(Q(\A)\ot^{st}_\bullet Q(\B),\C).
\end{multline*}
The pseudo Gray-functor classifier is constructed in \cite[Section 2]{Go}, see \cite[Remark 3.1]{BL} for the argument that the comonad $Q$ of Gohla plays the role of a pseudomap classifier, we applied this twice in the penultimate isomorphism above. 
In the last isomorphism it is $*=st$, and the bifunctor $G\x\Cat^{st}_{\bullet}(-,-)$ was implicitly studied in \cite[Section 7]{BL} (with $\bullet$ being $ps$ and $lx$). From the composition of isomorphisms we now deduce 
$$\A\ot^{ps}_{\bullet}\B\iso Q(\A)\ot^{st}_\bullet Q(\B).$$
In \cite[Section 14.2]{Gur} the existence of lax Gray functor classifier $Q^l$ was proved (we will comment the bifuntcor $G\x\Cat^{lx}_{\bullet}(-,-)$ in the next subsection). Then the same arguments for $*=lx$ imply: $\A\ot^{lx}_{\bullet}\B\iso Q^l(\A)\ot^{st}_\bullet Q^l(\B).$ 

We also record that in \cite{Mir} there are strictification results for tricategory homomorphisms. 



\subsection{Lax and strict Gray functor cases} \sslabel{lax and strict}

At last, we come to study the case of an inner-hom bifunctor for Gray-categories that was not studied in \cite{BL}. 
From \cite[Proposition 3.6 and Subsection 6.2]{BL} we know that 
$$G\x\Cat^{ps}_{lx}(-,-), \,\, G\x\Cat^{ps}_{ps}(-,-): (G\x\Cat_{st})^{ps}\times G\x\Cat_{ps}\to G\x\Cat_{ps}$$ 
are bifunctors and that they restrict to 
$(G\x\Cat_{st})^{op}\times G\x\Cat_{st}\to G\x\Cat_{st}$. The action of $G\x\Cat^{ps}_{lx}(-,-)$ is given in Subsection 3.4 thereof, 
and to avoid tedious verifications the authors carried out a conceptual proof of the bifunctor property by using the pseudomap classifier. 
The definition of the bifunctor action, in the same way as for double categories (see \cite[Section 2.3]{Gabi} and \cite[3.1]{Fem}), 
can equally well be given for $G\x\Cat^{*}_{lx}(-,-)$ for any $*$. Actually, for $*=lx$ the definition remains the same: Bourke and Lobbia 
expressed their definition in terms of the lax Gray functor structures. Let us see that $G\x\Cat^{lx}_{lx}(-,-)$ is a bifunctor: 1) on objects - $G\x\Cat^{lx}_{lx}(\A,\B)$ should be a Gray-category, and 2) on 1-cells - $G\x\Cat^{lx}_{lx}(F,G)$ should be a pseudo functor if $F,G$ are pseudo, {\em i.e.} a strict Gray functor if $F,G$ are strict. 
From the arguments that we will be using it should become clear that even if one takes any $*$ and $\bullet$ in places of $lx$, a bifunctor $G\x\Cat^{*}_{\bullet}(-,-)$ can be obtained. 

To 1): in Appendix B of \cite{BL} 
the Gray-category structure of $G\x\Cat^{ps}_{lx}(\A,\B)$ is described. It is readily seen that in it only the lax functor structure 
of its objects was used, and it is in the composition of its 1-cells, the lax transformations. Concretely, lax functor structures were used to define an invertible cocycle 3-cell $(\beta\cdot\alpha)^2_{f_2,f_2}$ at a composition of 1-cells $f_1,f_2$, for the composition of transformations $\alpha$ and $\beta$. It is given via the invertible cocycle 3-cells $\beta^2_{f_2,f_2}$ and $\alpha^2_{f_2,f_2}$. 
To prove that $\beta\cdot\alpha$ is well defined, the only axioms to be checked that contain structure cells of the functors are the axioms (v) and (vii). In them again only the lax functor and the lax transformations structures appear, so they hold also in the case of 
$G\x\Cat^{lx}_{lx}(\A,\B)$, and we can conclude that $G\x\Cat^{lx}_{lx}(\A,\B)$ is a Gray-category. 
(If colax instead of lax functors are used, the definitions accommodate accordingly, and the 
proof can equally be carried out even to prove that $G\x\Cat^{*}_{\bullet}(\A,\B)$, for any $*$ and $\bullet$, is a Gray-category.)

\begin{rem} \rmlabel{interch}
We highlight that the interchange in $G\x\Cat^*_{\bullet}(\A, \B)$ is given the same way as it is given in \cite[Appendix B.4]{BL} 
for the Gray-category $G\x\Cat^{ps}_{lx}(\A, \B)$. Namely, given modifications $\Gamma$ and $\Lambda$ an interchange among them is a 
perturbation $\Int_{\Gamma,\Lambda}$, which at an object $A\in\A$ is the interchange 3-cell in $\B$ among the 2-cell components 
$\Gamma(A)$ and $\Lambda(A)$.   
\end{rem}

To 2): consider $G\x\Cat^{*}_{lx}(-,-): (\Oo_{st/ps})^{op}\times\Oo_{st/ps}\to\Oo_{st/ps}$, where $\Oo$ stands for $G\x\Cat$ and we consider separately strict and pseudo functor case on the (co)domain category for $G\x\Cat^{*}_{lx}(-,-)$. On objects $\A$ and $\B$ we now know that 
 $G\x\Cat^{*}_{\bullet}(\A,\B)$ is again an object of $\Oo$. For morphisms $F:\A\to\A'$ and $G:\B\to\B'$ in $\Oo_{st/ps}$ we should get a morphism $G\x\Cat^{*}_{\bullet}(F,G): G\x\Cat^{*}_{\bullet}(\A',\B)\to G\x\Cat^{*}_{\bullet}(\A,\B')$ in $\Oo_{st/ps}$. 
On objects: take a $*$-type functor $H:\A'\to\B$ and we map it into $GHF:\A\to\B'$ as explained in Appendix C of \cite{BL} (via a left and right whiskering). It is immediately seen that both if $F$ and $G$ are strict, and if $F$ and $G$ are pseudo, the composition $GHF$ is a functor of type $*$. So, is $H$ is lax, so is $GHF$. By the definitions in Appendix C one can easily see that due to right whiskering via $G$, the pseudofunctor structure of $G$ is heavily used, and it is only at a single place where the pseudofunctor structure of $F$ is used. 
It is when defining $G\x\Cat^{*}_{\bullet}(F,G)$ at a 1-cell, lax transformation $\beta$, concretely, its invertible cocycle 3-cell 
$(G\beta F)^2_{f',f}$. The latter is given via $G(\beta^2_{F(f'),F(f)})$, which is clearly an invertible 3-cell, and $G(\beta_{F^2_{f',f}})$. 
The latter is invertible, by axiom (iii) of lax transformations, if the 2-cell $F^2_{f',f}$ is invertible, which is the case, since $F$ is pseudo. So we used that both $F$ and $G$ are pseudo (or strict), but this is indeed so, as we took them from $\Oo_{st/ps}$ from the beginning. Summing up, when checking $G\x\Cat^{*}_{lx}(-,-)$ on morphisms, a pair $(F,G)$, it was only when defining 
$G\x\Cat^{*}_{\bullet}(F,G)$ on objects (functors $H$) that $*$-type, or laxity of $H$ played a role.

With this we wish to point out that one can equally consider the bifunctor 
$G\x\Cat^{*}_{\bullet}(-,-): (\Oo_{st/ps})^{op}\times\Oo_{st/ps}\to\Oo_{st/ps}$, where for $*$ in particular one can consider $lx$ or $st$. 
As we commented in the previous subsection, the lax Gray functor classifier $Q^l$ was studied in \cite[Section 14.2]{Gur}, it strictifies the functors {\em i.e.} objects of $G\x\Cat^{lx}_{\bullet}(\A,\B)$ for Gray-categories $\A,\B$. 
More importantly, taking $*=st$, that is, imploying the bifunctor $G\x\Cat^{st}_{lx}(-,-)$, 
one is in position to have a simpler construction of the multicategory of {\em sharp} maps from Section 7 of \cite{BL}, without having to pass via the larger skew-multicategories $G\x\Cat^{ps}_{lx}$ and $G\x\Cat^{ps}_{ps}$, as in {\em loc. cit.}. Namely, the authors obtained 
multicategories $\Lax^\sharp=G\x\Cat^{st}_{lx}$ and $\Psd^\sharp=G\x\Cat^{st}_{ps}$ as sub-multicategories of $\Lax=G\x\Cat^{ps}_{lx}$ and 
$\Psd=G\x\Cat^{ps}_{ps}$, respectively, according to their Proposition 7.2. The difficulty that the authors encountered that obstructed 
the direct construction of these smaller multicategories is the way how they defined a part of the substitutions, for which they use duality bijections $d_n: q_n\x\Oo^{ps}_{lx}(\Pi_{i=2}^n\A_i,\B)\to q_n\x\Oo^{ps}_{lx}(\Pi_n^{i=2}\A_i^{co},\B^{co})$. Here 
$(-)^{co}:G\x\Cat_{ps}\to G\x\Cat_{ps}$ is an isomorphism of categories that flips the direction of 2-cells. It is these dualities that force the necessity to pass via the larger skew-multicategories $G\x\Cat^{ps}_{\bullet}$. Concretely, the troublesome substitutions are those that the authors denote by $\circ_2$, they correspond to what we denoted in a) below, which in \cite{BL} correspond respectively to b) below: 
\begin{enumerate} [a)]
\item $G(F_i(A,-),F_j(-,-))$ and $G(F_i(-,B),F_j(-,-))$ (when checking ternary condition) and $(f,(g,h,k))$ (when checking 4-ary condition, second subcase of case (iv)) in the proof of \prref{subst-gen}, 
\item third variable in Subsection 4.24 regarding ``binary into ternary'', second variable in Subsection 4.15 regarding ``binary into binary'', and second variable in Subsection 4.22 regarding ``ternary into binary''. 
\end{enumerate}
In our proof for the listed cases we used the isomorphisms \equref{towards ternary} and (\ref{towards 4-ary}), which correspond to 
\axiomref{$A_n$} from \ssref{closedness}, whereas \cite{BL} used the duality bijections. For the rest of the cases 
where ternary and 4-ary conditions needed to be checked, we used the isomorphisms \equref{ternary} and \equref{4-ary}, {\em i.e.} 
\axiomref{$B_n$}, analogously as was done in \cite{BL}.

Moreover, not only are the multicategories $G\x\Cat^{st}_{lx}$ and $G\x\Cat^{st}_{ps}$ left representable, as proven for 4-ary multicategories in Theorem 7.5, and that $(G\x\Cat, \ot^{st}_{lx}, G\x\Cat^{st}_{lx}(-,-))$ and $(G\x\Cat, \ot^{st}_{ps}, 
G\x\Cat^{st}_{ps}(-,-))$ are skew-monoidal categories, as concluded in Theorem 7.6 of \cite{BL}. Rather, we get indeed representable multicategories and monoidal induced categories, as we show next. 


\begin{thm} \thlabel{rep O}
The multicategory $G\x\Cat^{st}_{\bullet}$ is representable. 
\end{thm}

\begin{proof}
We prove a Gray-categorical version of \coref{strict-rep}, following the same lines of its proof. Set $\Oo=G\x\Cat$. Let $\F\in
q_n\x\Oo^{st}_\bullet((\Pi_{i=1}^{r}\A_i)\times(\B\ot^{st}_\bullet\C)\times(\Pi_{1}^{s}\D_i),\E)$ and we analyze the same three cases 
according to the places where three and four variables of $\F$ belong. The first case, $\B\ot^{st}_\bullet\C\times\Pi_{i=1}^{s}\D_i$, was proved by left representability, which we now also have, by \prref{left rep}. The second case $\A_i\times\A_j\times (\B\ot^{st}_\bullet\C)$ 
was proved by using \thref{strict-left closed} and \equref{any-strict}, \equref{quasi}, whose Gray-categorical versions are 
\thref{mixed-quasi n}, \thref{univ} and \equref{quasi}. 
The third case $\A_i\times(\B\ot^{st}_\bullet\C)\times\D_j$ needs a little bit of direct treatment. 
For 1-cells $(f,g,h):(A,B, D)\to(A', B', D')$ in $\A_i\times\B\times\D$ we have an incubator 3-cell $\F(f,g\ot C,h)$ in $\E$. 
Due to \thref{univ} (with $*=st$) we have $\F(f,g\ot C,h)=\crta\F(f,g,C,h)$, where 
$\crta\F: (\Pi_{i=1}^{r}\A_i)\times(\B\times\C)\times(\Pi_{i=1}^{s}\D_i)\to\E$ is the functor that induces $\F$ via the universal property. 
Thus, there is an incubator 3-cell for $\crta\F$ satisfying the required axioms, proving that $\crta\F$ is a quasi-functor. 
Because of the universal property the converse is similarly proved. The mecon axiom for four 1-cells is proved applying the analogous reasoning. Thus we have a bijection \vspace{-0,22cm}
$$q_n\x\Oo^{st}_\bullet((\Pi_{i=1}^{r}\A_i)\times(\B\ot^{st}_\bullet\C)\times(\Pi_{1}^{s}\D_i),\E)  \iso
q_{n+1}\x\Oo^{st}_\bullet((\Pi_{i=1}^{r}\A_i)\times(\B\times\C)\times(\Pi_{i=1}^{s}\D_i),\E). \vspace{-0,12cm}$$
As in \coref{strict-rep}, by iteration we get \vspace{-0,3cm}
$$q_n\x\Oo^{st}_\bullet((\Pi_{i=1}^{r}\A_i)\times(\crta\ot^k_\bullet\B_i)\times(\Pi_{1}^{s}\D_i),\E)\iso
q_p\x\Oo^{st}_\bullet((\Pi_{i=1}^{r}\A_i)\times(\Pi_{i=1}^k\B_i)\times(\Pi_{i=1}^{s}\D_i),\E)\vspace{-0,12cm}$$
and from here $\Oo^{st}_\bullet(\crta\ot^n_{o-l}\A_i,\B)\iso q_n\x\Oo^{st}_\bullet(\Pi_{i=1}^n\A_i,\B)$. 
\qed\end{proof}

\bigskip

\begin{cor} 
The category 
 $(G\x\Cat_{st}, \ot^{st}_{st}, G\x\Cat^{st}_{st}(-,-))$ is biclosed monoidal. 
\end{cor}

\begin{proof}
By (II) of \sssref{skew}, Hermida's \cite[Theorem 9.8]{Her} and the above theorem, the category 
$(G\x\Cat_{st}, \ot^{st}_{\bullet}, G\x\Cat^{st}_{\bullet}(-,-), G\x\Cat^{st}_{\#}(-,-))$ is a candidate for a monoidal category 
closed on any side, thus biclosed. Though, based on \prref{dim-biclosed} the only one that survives is when $\bullet=st$, and we have the claim. 
\qed\end{proof}

\subsection{Conclusions}

Summing up, we have studied monoidal and skew-monoidal closed structures on the categories of double categories (including 2-categories) 
and Gray-categories. We did this by constructing the corresponding skew-multicategories and multicategories applying the same reasonings in every step of the two constructions, the difference being only in the dimension. In this way we hope that analogous questions of 
constructing monoidal-like structures on categories (or higher categories) of higher dimensional enriched and internal categories can 
be tackled by following and continuing the arguments that we presented in this work. 

Various Gray-categorical results that we obtained  have been proved in \cite{BL}, but some of our approaches differ. We now want to highlight these differences. 

First of all, our meta-products present a constructive way of introducing monoidal or skew-monoidal products, whereas \cite{BL} 
obtain their skew-monoidal products by proving left representability, which in the presence of closedness is done by proving that the skew-multicategory admits a nullary map classifier and that each endo-hom functor has a left adjoint (as explained in their Subsection 2.5 and carried out in Subsection 5.2). The existence of left adjoint they prove in a conceptual way by showing that the category $G\x\Cat_{st}$ is locally finitely presentable.

The next important difference is that we do not use duality bijections 
$$d_n: q_n\x\Oo^{ps}_{lx}(\Pi_{i=1}^n\A_i,\B)\to q_n\x\Oo^{ps}_{lx}(\Pi_n^{i=1}\A_i^{co},\B^{co}),$$ 
where $(-)^{co}:G\x\Cat_{ps}\to G\x\Cat_{ps}$ is an isomorphism of categories 
that flips the 2-cells. Bourke and Lobbia use these dualities in the following occasions: 
\begin{itemize}
\item to prove 
some of the substitutions for the skew-multicategories, 
\item in order to prove right closedness, and thus biclosedness, of 
$\Lax^l=\OO^{ps}_{lx}$ in their Theorem 5.2, 
\item to consider the multicategories $\Lax^\sharp=\OO^{\sharp}_{lx}$ and $\Psd^\sharp=\OO^{\sharp}_{ps}$ in their Section 7 (and the induced skew-monoidal products $\ot^{\sharp}_{lx}$ and $\ot^{\sharp}_{ps}$): they construct them 
as sub-multicategories of $\Lax=\OO^{ps}_{lx}$ and $\Psd=\OO^{ps}_{ps}$ (because by their direct proof compositions of some tight maps result being really pseudo and not tight maps).
\end{itemize}
(Because of the use of duality, their multicategory $\Lax^\sharp=\OO^{\sharp}_{lx}$ is not right closed, and hence neither biclosed, 
as the authors argue at the end of Section 7. Though, it can not be biclosed already because of the Crans's argument.)

Instead of using duality, we cover the proof of all the substitutions in \prref{subst-gen} by using our \equref{towards ternary} and 
(\ref{towards 4-ary}) (these are isomorphisms \axiomref{$A_n$}, which do not seem to appear in \cite{BL}) and \equref{ternary} and 
\equref{4-ary} (these are isomorphisms \axiomref{$B_n$}, which they also use). This enables us to construct directly the multicategories 
$\OO^*_\bullet$ for any $*$ and $\bullet$. We do present some duality results, see \ssref{dual}, 
but of a different sort, by which lax functors remain lax, and these results do not form part of our constructions of (skew) multicategories. 
We use the symmetry present in quasi-functors (recall part d) of our \rmref{dual quasi}, \thref{lax-right closed} and the discussion preceding it) to prove right closedness of our multicategories and categories (without a need to pass to a larger construction in which $*=ps$).  

Bourke and Lobbia use {\em pseudo map classifier} to conceptually prove that $\bf{Lax}(-,-)=
G\x\Cat^{ps}_{lx}(-,-)$ is a bifunctor (their Proposition 3.6). We prove by direct checking that $G\x\Cat^*_\bullet(-,-)$ is a bifunctor. 



At last, the closedness isomorphisms proving the closedness of the (skew) multicategories, we prove using a different approach. 
While the ``binary-to-ternary'' and the ``ternary-to-4-ary'' isomorphisms in \cite{BL} are proved by directly comparing the axioms of the respective multimaps, we use explicit notions of transformations of quasi-functors of two and three variables that we introduced in 
\ssref{Gray-cat of qf} and \deref{tr of qf 3}, and we prove the 1-1 correspondences \equref{1-1}, which are based on \leref{gain incubator} and \leref{gain mecon}. 

In this research we have addressed with a same approach the study of monoidality and skew-monoidality of the categories of double categories 
and Gray-categories, so that the ideas in the proofs can be copied in larger dimensions. For 2-categories, Gray observed that $n$-ary multimaps (that is, quasi-functors of $n$-variables) are determined by ternary multimaps, which are collections of binary multimaps with compatibilities. For Gray-categories, because of the rise in dimension, $n$-ary multimaps are determined by 4-ary multimaps, which are collections of ternary multimaps satisfying the mecon axiom. Then in order to construct a multicategory corresponding to tetra-categories one would expect that one needs to describe 5-ary multimaps. However, in \cite[Section 2.6]{BL} it is proved that left representable 
{\em 4-ary} skew-multicategories give rise to skew-monoidal structures on the underlying category of tight unary maps (these ideas were 
extended in \cite{L}). This indicates that the notion of (skew) multicategories is possibly not enough to serve as a tool to describe a monoidal-like structure of the category of enriched or internal categories for dimensions higher than $m=3$. For double categories and Gray-categories we showed (in \cite[Subsection 4.1]{Fem} and \ssref{Gray-cat of qf}, see also \rmref{Gray-cat isos}, respectively) that the sets of multimaps 
$q_n\x\Oo^*_\bullet(\A,\B)$ are objects in $\Oo$. Thus, Gray multicategories $\M$ indeed seem to have this additional feature: 
that $\M_n(-,-): (\Oo^n)^{op} \times\Oo\to\Oo$ is a functor (compare to \rmref{multcat-cat}). 
This opens a space to explore higher dimensional analogues of multicategories as well as Hermida's equivalence result. 

\bigskip
\bigskip

{\bf Acknowledgments.} 
The author was supported by the Science Fund of the Republic of Serbia, Grant No. 7749891, Graphical Languages - GWORDS.

\pagebreak

{\bf\large{ Appendix A}}

\bigskip

{\bf Transformations and modifications of double categories}

\bigskip

\begin{defn} \delabel{hor nat tr} \cite[Definition 2.2]{Fem}
A {\em horizontal oplax transformation} $\alpha$ between lax double functors $F,G\colon \Aa\to\Bb$ consists of the following:
\begin{enumerate}
\item for every 0-cell $A$ in $\Aa$ a 1h-cell $\alpha(A)\colon F(A)\to G(A)$ in $\Bb$,
\item for every 1v-cell $u\colon A\to A'$ in $\Aa$ a 2-cell in $\Bb$:
$$
\scalebox{0.86}{
\bfig
\putmorphism(-150,50)(1,0)[F(A)`G(A)`\alpha(A)]{560}1a
\putmorphism(-150,-320)(1,0)[F(A')`G(A')`\alpha(A')]{600}1a
\putmorphism(-180,50)(0,-1)[\phantom{Y_2}``F(u)]{370}1l
\putmorphism(410,50)(0,-1)[\phantom{Y_2}``G(u)]{370}1r
\put(30,-110){\fbox{$\alpha^u$}}
\efig}
$$
\item for every 1h-cell $f\colon A\to B$  in $\Aa$ 
a 2-cell in $\Bb$:
$$
\scalebox{0.86}{
\bfig
 \putmorphism(-170,500)(1,0)[F(A)`F(B)`F(f)]{540}1a
 \putmorphism(360,500)(1,0)[\phantom{F(f)}`G(B) `\alpha(B)]{560}1a
 \putmorphism(-170,120)(1,0)[F(A)`G(A)`\alpha(A)]{540}1a
 \putmorphism(360,120)(1,0)[\phantom{G(B)}`G(A) `G(f)]{560}1a
\putmorphism(-180,500)(0,-1)[\phantom{Y_2}``=]{380}1r
\putmorphism(940,500)(0,-1)[\phantom{Y_2}``=]{380}1r
\put(280,310){\fbox{$\delta_{\alpha,f}$}}
\efig}
$$
\end{enumerate}
so that the following are satisfied: 

\begin{itemize}
\item (coherence with compositors for $\delta_{\alpha,-}$): for any composable 1h-cells $f$ and $g$ in $\Aa$ the 2-cell 
$\delta_{\alpha,gf}$ satisfies: \\
{\em \axiom{h.o.t.-1}} 
$$\scalebox{0.82}{
\bfig
  \putmorphism(-750,200)(1,0)[F(A)`\phantom{F(B)}`F(f)]{600}1a
\putmorphism(-130,200)(1,0)[F(A)`F(C)`F(g)]{580}1a
\putmorphism(-730,200)(0,-1)[\phantom{Y_2}``=]{400}1r
\putmorphism(420,200)(0,-1)[\phantom{Y_2}``=]{400}1r
\putmorphism(-730,-210)(0,-1)[\phantom{Y_2}``=]{400}1r
\putmorphism(1030,-210)(0,-1)[\phantom{Y_2}``=]{400}1r
 \putmorphism(450,-210)(1,0)[F(C)`G(C) `\alpha(C)]{580}1a
  \putmorphism(-750,-210)(1,0)[F(A)`\phantom{F(B)}`F(gf)]{1200}1a
\put(-270,20){\fbox{$F_{gf}$}}
 \putmorphism(-750,-615)(1,0)[F(A)`G(A)`\alpha(A)]{620}1a
 \putmorphism(-120,-615)(1,0)[\phantom{F(B)}`G(C) `G(gf)]{1170}1a
\put(-250,-420){\fbox{$\delta_{\alpha,gf}$}}
\efig}= 
\scalebox{0.82}{
\bfig
 \putmorphism(450,150)(1,0)[F(B)`F(C) `F(g)]{680}1a
 \putmorphism(1120,150)(1,0)[\phantom{F(B)}`G(C) `\alpha(C)]{600}1a
\put(1000,-100){\fbox{$\delta_{\alpha,g}$}}

  \putmorphism(-150,-300)(1,0)[F(A)` F(B) `F(f)]{600}1a
\putmorphism(450,-300)(1,0)[\phantom{F(A)}` G(B) `\alpha(B)]{680}1a
 \putmorphism(1120,-300)(1,0)[\phantom{F(A)}`G(C) ` G(g)]{620}1a

\putmorphism(450,150)(0,-1)[\phantom{Y_2}``=]{450}1l
\putmorphism(1710,150)(0,-1)[\phantom{Y_2}``=]{450}1r

 \putmorphism(-150,-750)(1,0)[F(A)`G(A)`\alpha(A)]{600}1a
 \putmorphism(450,-750)(1,0)[\phantom{F(B)}`G(B) `G(f)]{680}1a
 \putmorphism(1120,-750)(1,0)[\phantom{F(B)}`G(C) `G(g)]{620}1a

\putmorphism(-180,-300)(0,-1)[\phantom{Y_2}``=]{450}1r
\putmorphism(1040,-300)(0,-1)[\phantom{Y_2}``=]{450}1r
\put(350,-540){\fbox{$\delta_{\alpha,f}$}}
\put(1000,-960){\fbox{$G_{gf}$}}

 \putmorphism(450,-1200)(1,0)[G(A)` G(C) `G(gf)]{1300}1a

\putmorphism(450,-750)(0,-1)[\phantom{Y_2}``=]{450}1l
\putmorphism(1750,-750)(0,-1)[\phantom{Y_2}``=]{450}1r
\efig}
$$ 
(coherence with unitors for $\delta_{\alpha,-}$): for any object $A\in\Aa$: 
$$\text{{\em \axiom{h.o.t.-2}}}  \qquad\quad
\scalebox{0.86}{
\bfig
 \putmorphism(-150,420)(1,0)[F(A)`F(A)`=]{500}1a
\putmorphism(-180,420)(0,-1)[\phantom{Y_2}``=]{370}1l
\putmorphism(320,420)(0,-1)[\phantom{Y_2}``=]{370}1r
 \putmorphism(-150,50)(1,0)[F(A)`F(A)`F(1_A)]{500}1a
 \put(-80,250){\fbox{$F_A$}} 
\putmorphism(330,50)(1,0)[\phantom{F(A)}`G(A) `\alpha(A)]{560}1a
 \putmorphism(-170,-350)(1,0)[F(A)`G(A)`\alpha(A)]{520}1a
 \putmorphism(350,-350)(1,0)[\phantom{F(A)}`G(A) `G(1_A)]{560}1a

\putmorphism(-180,50)(0,-1)[\phantom{Y_2}``=]{400}1l
\putmorphism(910,50)(0,-1)[\phantom{Y_2}``=]{400}1r
\put(240,-150){\fbox{$\delta_{\alpha,1_A}$}}
\efig}
\quad=\quad
\scalebox{0.86}{
\bfig
 \putmorphism(-150,420)(1,0)[F(A)`G(A)`\alpha(A)]{500}1a
\putmorphism(-180,420)(0,-1)[\phantom{Y_2}``=]{370}1l
\putmorphism(320,420)(0,-1)[\phantom{Y_2}``=]{370}1r
  \put(-100,230){\fbox{$\Id_{\alpha(A)}$}} 
\putmorphism(-150,50)(1,0)[F(A)` \phantom{Y_2} `\alpha(A)]{450}1a

\putmorphism(350,50)(1,0)[G(A)` G(A) `=]{470}1a
\putmorphism(330,-300)(1,0)[G(A)` G(A) `G(1_A)]{480}1b
\putmorphism(330,50)(0,-1)[\phantom{Y_2}``=]{350}1l
\putmorphism(800,50)(0,-1)[\phantom{Y_2}``=]{350}1r
\put(470,-150){\fbox{$G_A$}}
\efig}
$$

\item (coherence with vertical composition and identity for $\alpha^\bullet$): for any composable 1v-cells $u$ and $v$ in $\Aa$:
$$\text{{\em \axiom{h.o.t.-3}}} \label{h.o.t.-3} \qquad\alpha^{\frac{u}{v}}=\frac{\alpha^u}{\alpha^v}\quad\qquad\text{ and}\quad\qquad
\text{{\em \axiom{h.o.t.-4}}} \label{h.o.t.-4} \qquad\alpha^{1^A}=\Id_{\alpha(A)};$$

\item (oplax naturality of 2-cells):
for every 2-cell in $\Aa$
$\scalebox{0.86}{
\bfig
\putmorphism(-150,50)(1,0)[A` B`f]{400}1a
\putmorphism(-150,-270)(1,0)[A'`B' `g]{400}1b
\putmorphism(-170,50)(0,-1)[\phantom{Y_2}``u]{320}1l
\putmorphism(250,50)(0,-1)[\phantom{Y_2}``v]{320}1r
\put(0,-140){\fbox{$a$}}
\efig}$ 
the following identity in $\Bb$ must hold:\\
$\text{{\em \axiom{h.o.t.-5}}}$ 
$$
\scalebox{0.86}{
\bfig
\putmorphism(-150,500)(1,0)[F(A)`F(B)`F(f)]{600}1a
 \putmorphism(450,500)(1,0)[\phantom{F(A)}`G(B) `\alpha(B)]{640}1a

 \putmorphism(-150,50)(1,0)[F(A')`F(B')`F(g)]{600}1a
 \putmorphism(450,50)(1,0)[\phantom{F(A)}`G(B') `\alpha(B')]{640}1a

\putmorphism(-180,500)(0,-1)[\phantom{Y_2}``F(u)]{450}1l
\putmorphism(450,500)(0,-1)[\phantom{Y_2}``]{450}1r
\putmorphism(300,500)(0,-1)[\phantom{Y_2}``F(v)]{450}0r
\putmorphism(1100,500)(0,-1)[\phantom{Y_2}``G(v)]{450}1r
\put(0,260){\fbox{$F(a)$}}
\put(700,270){\fbox{$\alpha^v$}}

\putmorphism(-150,-400)(1,0)[F(A')`G(A') `\alpha(A')]{640}1a
 \putmorphism(450,-400)(1,0)[\phantom{A'\ot B'}` G(B') `G(g)]{680}1a

\putmorphism(-180,50)(0,-1)[\phantom{Y_2}``=]{450}1l
\putmorphism(1120,50)(0,-1)[\phantom{Y_3}``=]{450}1r
\put(320,-200){\fbox{$\delta_{\alpha,g}$}}

\efig}
\quad=\quad
\scalebox{0.86}{
\bfig
\putmorphism(-150,500)(1,0)[F(A)`F(B)`F(f)]{600}1a
 \putmorphism(450,500)(1,0)[\phantom{F(A)}`G(B) `\alpha(B)]{680}1a
 \putmorphism(-150,50)(1,0)[F(A)`G(A)`\alpha(A)]{600}1a
 \putmorphism(450,50)(1,0)[\phantom{F(A)}`G(B) `G(f)]{680}1a

\putmorphism(-180,500)(0,-1)[\phantom{Y_2}``=]{450}1r
\putmorphism(1100,500)(0,-1)[\phantom{Y_2}``=]{450}1r
\put(350,260){\fbox{$\delta_{\alpha,f}$}}
\put(650,-180){\fbox{$G(a)$}}

\putmorphism(-150,-400)(1,0)[F(A')`G(A') `\alpha(A')]{640}1a
 \putmorphism(490,-400)(1,0)[\phantom{F(A')}` G(B'). `G(g)]{640}1a

\putmorphism(-180,50)(0,-1)[\phantom{Y_2}``F(u)]{450}1l
\putmorphism(450,50)(0,-1)[\phantom{Y_2}``]{450}1l
\putmorphism(610,50)(0,-1)[\phantom{Y_2}``G(u)]{450}0l 
\putmorphism(1120,50)(0,-1)[\phantom{Y_3}``G(v)]{450}1r
\put(40,-180){\fbox{$\alpha^u$}} 
\efig}
$$
\end{itemize}
A {\em horizontal strict transformation} is a  horizontal oplax transformation for which the 2-cells $\delta_{\alpha,f}$ in item 3 are identities. 
\end{defn}

\medskip

A {\em horizontal lax transformation} $\alpha$ between lax double functors $F,G\colon \Aa\to\Bb$ differs from its oplax counterpart 
in that the globular 2-cells $\delta_{\alpha,f}$ for any 1h-cell $f$ in $\Aa$ goes in the other direction, namely 
$$
\scalebox{0.86}{
\bfig
 \putmorphism(-170,500)(1,0)[F(A)`G(A)`\alpha(A)]{540}1a
 \putmorphism(360,500)(1,0)[\phantom{G(B)}`G(A) `G(f)]{560}1a
 \putmorphism(-170,120)(1,0)[F(A)`F(B)`F(f)]{540}1a
 \putmorphism(360,120)(1,0)[\phantom{F(f)}`G(B) `\alpha(B)]{560}1a
\putmorphism(-180,500)(0,-1)[\phantom{Y_2}``=]{380}1r
\putmorphism(940,500)(0,-1)[\phantom{Y_2}``=]{380}1r
\put(280,310){\fbox{$\sigma_{\alpha,f}$}}
\efig}
$$
and the axioms \axiomref{h.o.t.-1}-\axiomref{h.o.t.-5} are accordingly changed by the analogous axioms that we will refer to as to 
\axiom{h.l.t.-1}-\axiom{h.l.t.-5}. Indeed, note that only the three axioms \axiomref{h.o.t.-1}, \axiomref{h.o.t.-2} and \axiomref{h.o.t.-5} 
are changed into \axiom{h.l.t.-1}, \axiom{h.l.t.-2} and \axiom{h.l.t.-5}. 

\medskip

\begin{defn} \delabel{modif-hv}
A modification $\Theta$ between two horizontal oplax transformations $\alpha$ and $\beta$ and two vertical lax transformations $\alpha_0$ and $\beta_0$ 
depicted below on the left, where the lax double functors $F, G, F\s', G'$ act between $\Aa\to\Bb$, is given 
by a collection of 2-cells in $\Bb$ depicted below on the right:
\begin{equation} \eqlabel{modification cells}
\scalebox{0.86}{
\bfig
\putmorphism(-150,50)(1,0)[F` G`\alpha]{400}1a
\putmorphism(-150,-270)(1,0)[F'`G' `\beta]{400}1b
\putmorphism(-170,50)(0,-1)[\phantom{Y_2}``\alpha_0]{320}1l
\putmorphism(250,50)(0,-1)[\phantom{Y_2}``\beta_0]{320}1r
\put(-30,-140){\fbox{$\Theta$}}
\efig}
\qquad\qquad
\scalebox{0.86}{
\bfig
\putmorphism(-180,50)(1,0)[F(A)` G(A)`\alpha(A)]{550}1a
\putmorphism(-180,-270)(1,0)[F\s'(A)`G'(A) `\beta(A)]{550}1b
\putmorphism(-170,50)(0,-1)[\phantom{Y_2}``\alpha_0(A)]{320}1l
\putmorphism(350,50)(0,-1)[\phantom{Y_2}``\beta_0(A)]{320}1r
\put(0,-140){\fbox{$\Theta_A$}}
\efig}
\end{equation}
which satisfy the following rules: 

\medskip

\noindent {\em \axiom{m.ho-vl.-1}}  for every 1h-cell $f$, we have  
$$
\scalebox{0.86}{
\bfig
\putmorphism(-150,500)(1,0)[F(A)`F(B)`F(f)]{600}1a
 \putmorphism(450,500)(1,0)[\phantom{F(A)}`G(B) `\alpha(B)]{620}1a

 \putmorphism(-150,50)(1,0)[F\s'(A)`F\s'(B)`F\s'(f)]{600}1a
 \putmorphism(450,50)(1,0)[\phantom{F(A)}`G'(B) `\beta(B)]{620}1a

\putmorphism(-180,500)(0,-1)[\phantom{Y_2}``\alpha_0(A)]{450}1l
\putmorphism(450,500)(0,-1)[\phantom{Y_2}``]{450}1r
\putmorphism(300,500)(0,-1)[\phantom{Y_2}``\alpha_0(B)]{450}0r
\putmorphism(1080,500)(0,-1)[\phantom{Y_2}``\beta_0(B)]{450}1r
\put(0,280){\fbox{$(\alpha_0)_f$}}
\put(670,280){\fbox{$\Theta_B$}}

\putmorphism(-150,-400)(1,0)[F\s'(A)`G'(A) `\beta(A)]{600}1a
\putmorphism(510,-400)(1,0)[\phantom{Y_2}`G'(B) `G'(f)]{580}1a

\putmorphism(-180,50)(0,-1)[\phantom{Y_2}``=]{450}1l
\putmorphism(1080,50)(0,-1)[\phantom{Y_3}``=]{450}1r
\put(320,-180){\fbox{$\delta_{\beta,f}$}}

\efig}
\quad=\quad
\scalebox{0.86}{
\bfig
\putmorphism(-150,500)(1,0)[F(A)`F(B)`F(f)]{600}1a
 \putmorphism(450,500)(1,0)[\phantom{F(A)}`G(B) `\alpha(B)]{620}1a
\putmorphism(-150,50)(1,0)[F(A)`G(A) `\alpha(A)]{600}1a
\putmorphism(510,50)(1,0)[\phantom{Y_2}`G(B) `G(f)]{580}1a

\putmorphism(-180,500)(0,-1)[\phantom{Y_2}``=]{450}1r
\putmorphism(1080,500)(0,-1)[\phantom{Y_2}``=]{450}1r
\put(350,280){\fbox{$\delta_{\alpha,f}$}}

\putmorphism(-180,50)(0,-1)[\phantom{Y_2}``\alpha_0(A)]{450}1l
\putmorphism(1080,50)(0,-1)[\phantom{Y_3}``\beta_0(B)]{450}1r
\put(20,-180){\fbox{$\Theta_A$}}
\put(670,-180){\fbox{$(\beta_0)_f$}}

\putmorphism(450,50)(0,-1)[\phantom{Y_2}``]{450}1r
\putmorphism(300,50)(0,-1)[\phantom{Y_2}``\beta_0(A)]{450}0r

\putmorphism(-150,-400)(1,0)[F\s'(A)`G'(A) `\beta(A)]{600}1a
\putmorphism(510,-400)(1,0)[\phantom{Y_2}`G'(B) `G'(f)]{580}1a
\efig}
$$
and

\noindent {\em \axiom{m.ho-vl.-2}}   for every 1v-cell $u$, we have
$$
\scalebox{0.86}{
\bfig
 \putmorphism(-150,500)(1,0)[F(A)`F(A) `=]{600}1a
 \putmorphism(550,500)(1,0)[` `\alpha(A)]{400}1a
\putmorphism(-180,500)(0,-1)[\phantom{Y_2}`F\s'(A) `\alpha_0(A)]{450}1l
\put(30,50){\fbox{$\alpha_0^u$}}
\putmorphism(-150,-400)(1,0)[F\s'(\tilde A)` `=]{480}1a
\putmorphism(-180,50)(0,-1)[\phantom{Y_2}``F\s'(u)]{450}1l
\putmorphism(450,50)(0,-1)[\phantom{Y_2}`F\s'(\tilde A)` \alpha_0(\tilde A)]{450}1l
\putmorphism(450,500)(0,-1)[\phantom{Y_2}`F(\tilde A) `F(u)]{450}1l
\put(660,280){\fbox{$\alpha^u$}}
\putmorphism(450,50)(1,0)[\phantom{(B, \tilde A)}``\alpha(\tilde A)]{500}1a
\putmorphism(1070,50)(0,-1)[\phantom{(B, A')}`G'(\tilde A)`\beta_0(\tilde A)]{450}1r
\putmorphism(1070,500)(0,-1)[G(A)`G(\tilde A)`G(u)]{450}1r
\putmorphism(450,-400)(1,0)[\phantom{(B, \tilde A)}``\beta(\tilde A)]{500}1a
\put(640,-170){\fbox{$ \Theta_{\tilde A}$ } } 
\efig}\quad=\quad
\scalebox{0.86}{
\bfig
 \putmorphism(-150,500)(1,0)[F(A)`G(A) `\alpha(A)]{600}1a
 \putmorphism(450,500)(1,0)[\phantom{(B,A)}` `=]{460}1a
\putmorphism(-180,500)(0,-1)[\phantom{Y_2}`F\s'(A) `\alpha_0(A)]{450}1l
\put(650,50){\fbox{$\beta_0^u$}}
\putmorphism(-180,-400)(1,0)[F\s'(\tilde A)` `\beta(\tilde A)]{500}1a
\putmorphism(-180,50)(0,-1)[\phantom{Y_2}``F\s'(u)]{450}1l
\putmorphism(450,50)(0,-1)[\phantom{Y_2}`G'(\tilde A)`G'(u)]{450}1r
\putmorphism(450,500)(0,-1)[\phantom{Y_2}`G'(A) `\beta_0(A)]{450}1r
\put(0,260){\fbox{$\Theta_A$}}
\putmorphism(-180,50)(1,0)[\phantom{(B, \tilde A)}``\beta(A)]{500}1a
\putmorphism(1030,50)(0,-1)[\phantom{(B, A')}` G'(\tilde A). ` \beta_0(\tilde A)]{450}1r
\putmorphism(1030,500)(0,-1)[G(A)`G(\tilde A)` G(u)]{450}1r
\putmorphism(430,-400)(1,0)[\phantom{(B, \tilde A)}``=]{480}1b
\put(70,-170){\fbox{$\beta^u$}}
\efig}
$$
\end{defn}

\bigskip

{\bf\large{ Appendix B}}

\bigskip

{\bf Lax double quasi-functors and their vertical transformations}

\medskip

{\bf Appendix B.1}

\bigskip

\begin{prop} \prlabel{char df} \cite[Proposition 3.3]{Fem}
A lax double functor $\F\colon\Aa\to\Lax_{hop}(\Bb, \Cc)$ of double categories consists of the following: \\
1. lax double functors 
$$(-,A)\colon\Bb\to\Cc\quad\text{ and}\quad (B,-)\colon\Aa\to\Cc$$ 
such that $(-,A)\vert_B=(B,-)\vert_A=(B,A)$, 
for objects $A\in\Aa, B\in\Bb$, \\
2. 2-cells
$$
\scalebox{0.86}{
\bfig
 \putmorphism(-150,50)(1,0)[(B,A)`(B', A)`(k, A)]{600}1a
 \putmorphism(450,50)(1,0)[\phantom{A\ot B}`(B', A') `(B', K)]{680}1a
\putmorphism(-180,50)(0,-1)[\phantom{Y_2}``=]{450}1r
\putmorphism(1100,50)(0,-1)[\phantom{Y_2}``=]{450}1r
\put(350,-190){\fbox{$(k,K)$}}
 \putmorphism(-150,-400)(1,0)[(B,A)`(B,A')`(B,K)]{600}1a
 \putmorphism(450,-400)(1,0)[\phantom{A\ot B}`(B', A') `(k, A')]{680}1a
\efig}
$$

$$
\scalebox{0.86}{
\bfig
\putmorphism(-150,50)(1,0)[(B,A)`(B,A')`(B,K)]{600}1a
\putmorphism(-150,-400)(1,0)[(\tilde B, A)`(\tilde B,A') `(\tilde B,K)]{640}1a
\putmorphism(-180,50)(0,-1)[\phantom{Y_2}``(u,A)]{450}1l
\putmorphism(450,50)(0,-1)[\phantom{Y_2}``(u,A')]{450}1r
\put(0,-180){\fbox{$(u, K)$}}
\efig}
\quad
\scalebox{0.86}{
\bfig
\putmorphism(-150,50)(1,0)[(B,A)`(B',A)`(k,A)]{600}1a
\putmorphism(-150,-400)(1,0)[(B, \tilde A)`(B', \tilde A) `(k,\tilde A)]{640}1a
\putmorphism(-180,50)(0,-1)[\phantom{Y_2}``(B,U)]{450}1l
\putmorphism(450,50)(0,-1)[\phantom{Y_2}``(B',U)]{450}1r
\put(0,-180){\fbox{$(k,U)$}}
\efig}
$$

$$
\scalebox{0.86}{
\bfig
 \putmorphism(-150,500)(1,0)[(B,A)`(B,A) `=]{600}1a
\putmorphism(-180,500)(0,-1)[\phantom{Y_2}`(B, \tilde A) `(B,U)]{450}1l
\put(-20,50){\fbox{$(u,U)$}}
\putmorphism(-150,-400)(1,0)[(\tilde B, \tilde A)`(\tilde B, \tilde A) `=]{640}1a
\putmorphism(-180,50)(0,-1)[\phantom{Y_2}``(u,\tilde A)]{450}1l
\putmorphism(450,50)(0,-1)[\phantom{Y_2}``(\tilde B, U)]{450}1r
\putmorphism(450,500)(0,-1)[\phantom{Y_2}`(\tilde B, A) `(u,A)]{450}1r
\efig}
$$ 
in $\Cc$ for every 1h-cells $A\stackrel{K}{\to} A'$ and $B\stackrel{k}{\to} B'$ and 1v-cells $A\stackrel{U}{\to} \tilde A$ and 
$B\stackrel{u}{\to} \tilde B$ which satisfy: 

\noindent $\bullet$ \quad \axiom{($1_B,K$)}  
$$
\scalebox{0.86}{
\bfig
 \putmorphism(-210,420)(1,0)[(B,A)`(B,A)`=]{550}1a
\putmorphism(-210,50)(1,0)[(B,A)`(B,A) `(1_B,A)]{560}1a
\putmorphism(350,400)(1,0)[\phantom{F(A)}` (B,A') `(B,K)]{600}1a
\putmorphism(360,50)(1,0)[\phantom{F(A)}`(B,A') `(B,K)]{600}1a

\putmorphism(-170,420)(0,-1)[\phantom{Y_2}``=]{350}1l
\putmorphism(320,420)(0,-1)[\phantom{Y_2}``]{370}1r
\putmorphism(300,420)(0,-1)[\phantom{Y_2}``=]{370}0r
\put(-120,250){\fbox{$(-,A)_B$}}

\putmorphism(860,420)(0,-1)[\phantom{Y_2}``=]{350}1r
\put(450,240){\fbox{$\Id_{(B,K)}$}}
\putmorphism(-180,50)(0,-1)[\phantom{Y_2}``=]{350}1r
\putmorphism(860,50)(0,-1)[\phantom{Y_2}``=]{350}1r
 \putmorphism(-150,-300)(1,0)[(B,A)`(B,A')`(B,K)]{500}1a
 \putmorphism(350,-300)(1,0)[\phantom{A\ot B}`(B, A') `(1_B, A')]{600}1a
\put(170,-110){\fbox{$(1_B,K)$}}
\efig}
\quad=\quad
\scalebox{0.86}{
\bfig
 \putmorphism(-260,200)(1,0)[(B,A)`\phantom{F(A)} `(B,K)]{550}1a
\putmorphism(330,200)(1,0)[(B,A')`(B,A')`=]{600}1a
\putmorphism(-210,200)(0,-1)[\phantom{Y_2}``=]{370}1l
\putmorphism(320,200)(0,-1)[\phantom{Y_2}``]{370}1l
\putmorphism(340,200)(0,-1)[\phantom{Y_2}``=]{370}0l
\putmorphism(960,200)(0,-1)[\phantom{Y_2}``=]{370}1r
 \putmorphism(-260,-170)(1,0)[(B,A)`\phantom{Y_2}`(B,K)]{520}1a
 \put(420,30){\fbox{$(-,A')_{B}$}} 
\putmorphism(360,-170)(1,0)[(B,A')`(B,A') `(1_B,A')]{620}1a
\put(-160,30){\fbox{$\Id_{(B,K)}$}}
\efig}
$$


\noindent $\bullet$ \quad \axiom{($k,1_A$)} \vspace{-0,9cm}
$$ 
\scalebox{0.86}{
\bfig
 \putmorphism(-150,200)(1,0)[(B,A)`(B,A)`=]{500}1a
\putmorphism(360,200)(1,0)[\phantom{F(A)}`(B,A') `(k,A)]{500}1a
\putmorphism(830,200)(0,-1)[\phantom{Y_2}``=]{350}1r
\put(460,30){\fbox{$\Id_{(k,A)}$}}

\putmorphism(-180,200)(0,-1)[\phantom{Y_2}``=]{370}1l
\putmorphism(320,200)(0,-1)[\phantom{Y_2}``]{370}1r
\putmorphism(300,200)(0,-1)[\phantom{Y_2}``=]{370}0r
 \putmorphism(-150,-170)(1,0)[(B,A)`(B,A)`(B,1_A)]{500}1a
 \put(-140,30){\fbox{$(B,-)_A$}} 
\putmorphism(350,-170)(1,0)[\phantom{F(A)}`(B',A) `(k,A)]{560}1a
\efig}
\quad
=
\quad
\scalebox{0.86}{
\bfig
\putmorphism(-150,420)(1,0)[(B,A)` \phantom{Y_2}`(k,A)]{450}1a
\putmorphism(370,420)(1,0)[(B',A)` (B',A) `=]{500}1a
\putmorphism(-170,420)(0,-1)[\phantom{Y_2}``=]{350}1l
\putmorphism(350,420)(0,-1)[\phantom{Y_2}``=]{350}1l
\putmorphism(860,420)(0,-1)[\phantom{Y_2}``=]{350}1r
 \putmorphism(420,50)(1,0)[\phantom{Y_2}`(B',A)`(B',1_A)]{520}1a
\putmorphism(-150,50)(1,0)[(B,A)` (B',A) `(k,A)]{500}1a
 \put(410,250){\fbox{$(B',-)_A$}} 
\put(-120,250){\fbox{$\Id_{(k,A)}$}}

\putmorphism(-180,50)(0,-1)[\phantom{Y_2}``=]{350}1r
\putmorphism(860,50)(0,-1)[\phantom{Y_2}``=]{350}1r
\putmorphism(-150,-300)(1,0)[(B,A)`(B,A)`(B,1_A)]{500}1a
 \putmorphism(350,-300)(1,0)[\phantom{A\ot B}`(B', A) `(k, A)]{560}1a
\put(200,-120){\fbox{$(k,1_A)$}}
\efig}
$$

\noindent where the 2-cells $(-,A)_{B}$ and $(B,-)_A$ come from laxity of the lax double functors $(-,A)$ and $(B,-)$ 

\noindent $\bullet$ \quad  \axiom{($u,1_A$)} 
$$
\scalebox{0.86}{
\bfig
\putmorphism(-250,500)(1,0)[(B,A)`(B,A)` =]{550}1a
 \putmorphism(-250,50)(1,0)[(\tilde B,A)`(\tilde B,A)` =]{550}1a
 \putmorphism(-250,-400)(1,0)[(\tilde B,A)`(\tilde B,A)` (\tilde B,1_A)]{550}1a

\putmorphism(-280,500)(0,-1)[\phantom{Y_2}``(u,A)]{450}1l
 \putmorphism(-280,70)(0,-1)[\phantom{F(A)}` `=]{450}1l

\putmorphism(300,500)(0,-1)[\phantom{Y_2}``(u,A)]{450}1r
\putmorphism(300,70)(0,-1)[\phantom{Y_2}``=]{450}1r
\put(-150,270){\fbox{$Id^{(u,A)}$}}
\put(-190,-150){\fbox{$(\tilde B,-)_A$}}
\efig}
=
\scalebox{0.86}{
\bfig
\putmorphism(-250,500)(1,0)[(B,A)`(B,A)` =]{550}1a
 \putmorphism(-250,50)(1,0)[(B,A)`(B,A)` (B,1_A)]{550}1a
 \putmorphism(-250,-400)(1,0)[(\tilde B,A)`(\tilde B,A)` (\tilde B,1_A)]{550}1a

\putmorphism(-280,500)(0,-1)[\phantom{Y_2}``= ]{450}1l
 \putmorphism(-280,70)(0,-1)[\phantom{F(A)}` `(u,A)]{450}1l

\putmorphism(300,500)(0,-1)[\phantom{Y_2}``=]{450}1r
\putmorphism(300,70)(0,-1)[\phantom{Y_2}``(u,A)]{450}1r
\put(-160,290){\fbox{$(B,-)_A$}}
\put(-140,-150){\fbox{$(u,1_A)$}}
\efig}
$$

\noindent $\bullet$ \quad \axiom{($1_B,U$)} 
$$\scalebox{0.86}{
\bfig
\putmorphism(-250,500)(1,0)[(B,A)`(B,A)` =]{550}1a
 \putmorphism(-250,50)(1,0)[(B,A)`(B,A)` (1_B,A)]{550}1a
 \putmorphism(-250,-400)(1,0)[(B,\tilde A)`(B,\tilde A)` (1_B,\tilde A)]{550}1a

\putmorphism(-280,500)(0,-1)[\phantom{Y_2}``= ]{450}1l
 \putmorphism(-280,70)(0,-1)[\phantom{F(A)}` `(B,U)]{450}1l

\putmorphism(300,500)(0,-1)[\phantom{Y_2}``=]{450}1r
\putmorphism(300,70)(0,-1)[\phantom{Y_2}``(B,U)]{450}1r
\put(-170,290){\fbox{$(-,A)_B$}}
\put(-150,-160){\fbox{$(1_B,U)$}}
\efig}\quad
=
\scalebox{0.86}{
\bfig
\putmorphism(-250,500)(1,0)[(B,A)`(B,A)` =]{550}1a
 \putmorphism(-250,50)(1,0)[(B,\tilde A)`(B,\tilde A)` =]{550}1a
 \putmorphism(-250,-400)(1,0)[(B,\tilde A)`(B,\tilde A)` (\tilde B,1_{\tilde A})]{550}1a

\putmorphism(-280,500)(0,-1)[\phantom{Y_2}``(B,U)]{450}1l
 \putmorphism(-280,70)(0,-1)[\phantom{F(A)}` `=]{450}1l

\putmorphism(300,500)(0,-1)[\phantom{Y_2}``(B,U)]{450}1r
\putmorphism(300,70)(0,-1)[\phantom{Y_2}``=]{450}1r
\put(-170,270){\fbox{$Id^{(B,U)}$}}
\put(-180,-160){\fbox{$(-, \tilde A)_B$}}
\efig}
$$

\noindent $\bullet$ \qquad \axiom{($1^B,K$)}  \quad $(1^B,K)=Id_{(B,K)}$ \qquad\hspace{-0,16cm}\text{and}\qquad 
$\bullet$ \quad  \axiom{($k,1^A$)} \qquad $(k,1^A)=Id_{(k,A)}$ 

\noindent $\bullet$ \qquad  \axiom{($1^B,U$)} \quad $(1^B,U)=Id^{(B,U)}$ \qquad\text{and}\qquad \hspace{-0,22cm} 
$\bullet$ \quad  \axiom{($u,1^A$)} \qquad $(u,1^A)=Id^{(u,A)}$

\noindent $\bullet$ \quad  \axiom{($k'k,K$)} 
$$
\scalebox{0.78}{
\bfig
 \putmorphism(450,650)(1,0)[(B', A) `(B'', A) `(k', A)]{680}1a
 \putmorphism(1140,650)(1,0)[\phantom{A\ot B}`(B'', A') ` (B'', K)]{680}1a

 \putmorphism(-150,200)(1,0)[(B, A) `(B', A)`(k, A)]{600}1a
 \putmorphism(450,200)(1,0)[\phantom{A\ot B}`(B', A') `(B', K)]{680}1a
 \putmorphism(1130,200)(1,0)[\phantom{A\ot B}`(B'', A') ` (k', A')]{680}1a

\putmorphism(450,650)(0,-1)[\phantom{Y_2}``=]{450}1r
\putmorphism(1750,650)(0,-1)[\phantom{Y_2}``=]{450}1r
\put(1000,420){\fbox{$ (k',K)$}}

 \putmorphism(-150,-250)(1,0)[(B, A)`(B, A') `(B,K)]{640}1a
 \putmorphism(480,-250)(1,0)[\phantom{A'\ot B'}`(B', A') `(k, A')]{680}1a

\putmorphism(-180,200)(0,-1)[\phantom{Y_2}``=]{450}1l
\putmorphism(1120,200)(0,-1)[\phantom{Y_3}``=]{450}1r
\put(310,-50){\fbox{$ (k,K)$}}

 \putmorphism(1170,-250)(1,0)[\phantom{A\ot B}`(B'', A') ` (k', A')]{650}1a
\putmorphism(450,-250)(0,-1)[\phantom{Y_2}``=]{450}1r
\putmorphism(1750,-250)(0,-1)[\phantom{Y_2}``=]{450}1r

 \putmorphism(480,-700)(1,0)[(B, A') `(B'', A'') `(k'k, A')]{1320}1a
\put(920,-470){\fbox{$ (-,A')_{k'k}$}}
\efig}
=\quad
\scalebox{0.78}{
\bfig
 \putmorphism(-150,500)(1,0)[(B, A) `(B', A)`(k,A)]{580}1a
 \putmorphism(450,500)(1,0)[\phantom{(B, A)} `(B'', A) `(k',A)]{660}1a
\putmorphism(-180,500)(0,-1)[\phantom{Y_2}``=]{450}1r
\put(240,270){\fbox{$ (-,A)_{k'k}$}}

 \putmorphism(-150,50)(1,0)[(B,A)` `(k'k,A)]{1080}1a
 \putmorphism(1080,50)(1,0)[(B'',A)`(B'', A') ` (B'',K)]{680}1a

\putmorphism(1050,500)(0,-1)[\phantom{Y_2}``=]{450}1r

 \putmorphism(-150,-400)(1,0)[(B, A)`(B, A') `(B,K)]{640}1a
 \putmorphism(530,-400)(1,0)[\phantom{Y_2X}`(B'', A') `(k'k,A')]{1220}1a
\put(570,-200){\fbox{$ (k'k,K)$}}

\putmorphism(-180,50)(0,-1)[\phantom{Y_2}``=]{450}1l
\putmorphism(1700,50)(0,-1)[\phantom{Y_3}``=]{450}1r
\efig}
$$ 
where $(-,A)_{k'k}$ is the 2-cell from the laxity of $(-,A)$ 

\noindent $\bullet$ \quad  \axiom{($k,K'K$)}  
$$
\scalebox{0.78}{
\bfig
 \putmorphism(450,500)(1,0)[(B', A) `(B', A') `(B', K)]{680}1a
 \putmorphism(1140,500)(1,0)[\phantom{A\ot B}`(B', A'') ` (B', K')]{680}1a

 \putmorphism(-150,50)(1,0)[(B, A) `(B', A)`(k,A)]{600}1a
 \putmorphism(450,50)(1,0)[\phantom{A\ot B}`(B', A'') `(B', K'K)]{1350}1a

\putmorphism(450,500)(0,-1)[\phantom{Y_2}``=]{450}1r
\putmorphism(1750,500)(0,-1)[\phantom{Y_2}``=]{450}1r
\put(880,290){\fbox{$ (B',-)_{K'K}$}}

 \putmorphism(-150,-400)(1,0)[(B, A)`(B, A'') `(B, K'K)]{980}1a
 \putmorphism(780,-400)(1,0)[\phantom{A'\ot B'}`(B', A'') `(k,A'')]{980}1a

\putmorphism(-180,50)(0,-1)[\phantom{Y_2}``=]{450}1l
\putmorphism(1750,50)(0,-1)[\phantom{Y_3}``=]{450}1r
\put(560,-200){\fbox{$ (k,K'K)$}}

\efig}
=
\scalebox{0.78}{
\bfig
 \putmorphism(-150,450)(1,0)[(B,A)`(B',A)`(k,A)]{600}1a
 \putmorphism(450,450)(1,0)[\phantom{A\ot B}`(B', A') `(B',K)]{680}1a

 \putmorphism(-150,0)(1,0)[(B,A)`(B,A')`(B,K)]{600}1a
 \putmorphism(450,0)(1,0)[\phantom{A\ot B}`(B', A') `(k,A')]{680}1a
 \putmorphism(1120,0)(1,0)[\phantom{A'\ot B'}`(B', A'') `(B', K')]{660}1a

\putmorphism(-180,450)(0,-1)[\phantom{Y_2}``=]{450}1r
\putmorphism(1100,450)(0,-1)[\phantom{Y_2}``=]{450}1r
\put(350,210){\fbox{$(k,K)$}}
\put(1000,-250){\fbox{$(k,K')$}}

 \putmorphism(450,-450)(1,0)[\phantom{A''\ot B'}` (B, A'') `(B, K')]{680}1a
 \putmorphism(1100,-450)(1,0)[\phantom{A''\ot B'}`(B', A'') ` (k, A'')]{660}1a

\putmorphism(450,0)(0,-1)[\phantom{Y_2}``=]{450}1l
\putmorphism(1750,0)(0,-1)[\phantom{Y_2}``=]{450}1r
 \putmorphism(-150,-450)(1,0)[(B,A)`(B,A')`(B,K)]{600}1a
\putmorphism(-180,-450)(0,-1)[\phantom{Y_2}``=]{450}1r
\putmorphism(1100,-450)(0,-1)[\phantom{Y_2}``=]{450}1r
 \putmorphism(-150,-900)(1,0)[(B,A)`(B,A'')`(B, K'K)]{1280}1a
\put(260,-670){\fbox{$(B,-)_{K'K}$}}
\efig}
$$
where $(B,-)_{K'K}$ is the 2-cell from the laxity of $(B,-)$ 

\noindent $\bullet$ \quad   \axiom{($u, K'K$)}  
$$
\scalebox{0.86}{
\bfig
\putmorphism(-150,500)(1,0)[(B,A)`(B,A')`(B,K)]{600}1a
 \putmorphism(470,500)(1,0)[\phantom{F(A)}`(B,A'') `(B, K')]{600}1a
 \putmorphism(-150,50)(1,0)[(B,A)`(B,A'')`(B,K'K)]{1220}1a

\putmorphism(-180,500)(0,-1)[\phantom{Y_2}``=]{450}1r
\putmorphism(1080,500)(0,-1)[\phantom{Y_2}``=]{450}1r
\put(250,290){\fbox{$(B,-)_{K'K}$}}

\putmorphism(-150,-400)(1,0)[(\tilde B,A)`(\tilde B,A'') `(\tilde B,K'K)]{1200}1a

\putmorphism(-180,50)(0,-1)[\phantom{Y_2}``(u,A)]{450}1l
\putmorphism(1080,50)(0,-1)[\phantom{Y_3}``(u,A'')]{450}1r
\put(250,-160){\fbox{$(u,K'K)$}} 
\efig}
=
\scalebox{0.86}{
\bfig
\putmorphism(-150,500)(1,0)[(B,A)`(B,A')`(B,K)]{600}1a
 \putmorphism(470,500)(1,0)[\phantom{F(A)}`(B,A'') `(B, K')]{600}1a

 \putmorphism(-150,50)(1,0)[(\tilde B, A)`(\tilde B,A')`(\tilde B,K)]{600}1a
 \putmorphism(470,50)(1,0)[\phantom{F(A)}`(B'',\tilde A) `(\tilde B,K')]{620}1a

\putmorphism(-180,500)(0,-1)[\phantom{Y_2}``(u,A)]{450}1l
\putmorphism(450,500)(0,-1)[\phantom{Y_2}``]{450}1r
\putmorphism(300,500)(0,-1)[\phantom{Y_2}``(u,A')]{450}0r
\putmorphism(1080,500)(0,-1)[\phantom{Y_2}``(u,A'')]{450}1r
\put(-40,280){\fbox{$(u,K)$}}
\put(620,280){\fbox{$(u,K')$}}

\putmorphism(-150,-400)(1,0)[(\tilde B, A)`(\tilde B,A'') `(\tilde B,K'K)]{1200}1a

\putmorphism(-180,50)(0,-1)[\phantom{Y_2}``=]{450}1l
\putmorphism(1080,50)(0,-1)[\phantom{Y_3}``=]{450}1r
\put(260,-160){\fbox{$(\tilde B,-)_{K'K}$}}

\efig}
$$

\noindent $\bullet$ \quad  \axiom{($k'k, U$)} 
$$
\scalebox{0.86}{
\bfig
\putmorphism(-150,500)(1,0)[(B,A)`(B',A)`(k,A)]{600}1a
 \putmorphism(450,500)(1,0)[\phantom{F(A)}`(B'',A) `(k',A)]{620}1a

 \putmorphism(-150,50)(1,0)[(B,\tilde A)`(B',\tilde A)`(k,\tilde A)]{600}1a
 \putmorphism(470,50)(1,0)[\phantom{F(A)}`(B'',\tilde A) `(k',\tilde A)]{620}1a

\putmorphism(-180,500)(0,-1)[\phantom{Y_2}``(B,U)]{450}1l
\putmorphism(450,500)(0,-1)[\phantom{Y_2}``]{450}1r
\putmorphism(240,500)(0,-1)[\phantom{Y_2}``(B',U)]{450}0r
\putmorphism(1080,500)(0,-1)[\phantom{Y_2}``(B'',U)]{450}1r
\put(-40,280){\fbox{$(k,U)$}}
\put(620,280){\fbox{$(k',U)$}}

\putmorphism(-150,-400)(1,0)[(B,\tilde A)`(B'',\tilde A) `(k'k,\tilde A)]{1200}1a

\putmorphism(-180,50)(0,-1)[\phantom{Y_2}``=]{450}1l
\putmorphism(1080,50)(0,-1)[\phantom{Y_3}``=]{450}1r
\put(260,-160){\fbox{$(-,\tilde A)_{k'k}$}}
\efig}
=
\scalebox{0.86}{
\bfig
\putmorphism(-150,500)(1,0)[(B,A)`(B',A)`(k,A)]{600}1a
 \putmorphism(450,500)(1,0)[\phantom{F(A)}`(B'',A) `(k',A)]{620}1a
 \putmorphism(-150,50)(1,0)[(B,A)`(B'',A)`(k'k, A)]{1220}1a

\putmorphism(-180,500)(0,-1)[\phantom{Y_2}``=]{450}1r
\putmorphism(1080,500)(0,-1)[\phantom{Y_2}``=]{450}1r
\put(260,290){\fbox{$(-,A)_{k'k}$}}

\putmorphism(-150,-400)(1,0)[(B,\tilde A)`(B'',\tilde A) `(k'k,\tilde A)]{1200}1a

\putmorphism(-180,50)(0,-1)[\phantom{Y_2}``(B,U)]{450}1l
\putmorphism(1080,50)(0,-1)[\phantom{Y_3}``(B'',U)]{450}1r
\put(300,-180){\fbox{$(k'k,U)$}} 
\efig}
$$

\noindent $\bullet$ \quad   \axiom{($\frac{u}{u'},K$)} \qquad $(\frac{u}{u'}, K)=\frac{(u,K)}{(u', K)}$  \qquad\text{and}\qquad 
$\bullet$ \quad  \axiom{($k,\frac{U}{U'}$)}  \qquad $(k,\frac{U}{U'})=\frac{(k,U)}{(k,U')}$ 

\noindent $\bullet$ \quad  \axiom{($u,\frac{U}{U'}$)}
$$(u,\frac{U}{U'})=
\scalebox{0.86}{
\bfig
 \putmorphism(-150,500)(1,0)[(B,A)`(B,A) `=]{600}1a
\putmorphism(-180,500)(0,-1)[\phantom{Y_2}`(B, \tilde A) `(B,U)]{450}1l
\put(0,50){\fbox{$(u,U)$}}
\putmorphism(-150,-400)(1,0)[(\tilde B, \tilde A)`(\tilde B, \tilde A) `=]{640}1a
\putmorphism(-180,50)(0,-1)[\phantom{Y_2}``(u,\tilde A)]{450}1l
\putmorphism(450,50)(0,-1)[\phantom{Y_2}``(\tilde B, U)]{450}1r
\putmorphism(450,500)(0,-1)[\phantom{Y_2}`(\tilde B, A) `(u,A)]{450}1r
\putmorphism(-820,50)(1,0)[(B, \tilde A)``=]{520}1a
\putmorphism(-820,50)(0,-1)[\phantom{(B, \tilde A')}``(B,U')]{450}1l
\putmorphism(-820,-400)(0,-1)[(B, \tilde A')`(\tilde B, \tilde A')`(u,\tilde A')]{450}1l
\putmorphism(-820,-850)(1,0)[\phantom{(B, \tilde A)}``=]{520}1a
\putmorphism(-150,-400)(0,-1)[(\tilde B, \tilde A)`(\tilde B, \tilde A') `(\tilde B, U')]{450}1r
\put(-650,-630){\fbox{$(u,U')$}}
\efig}
$$

\noindent $\bullet$ \quad  \axiom{($\frac{u}{u'},U$)} 
$$(\frac{u}{u'},U)=
\scalebox{0.86}{
\bfig
 \putmorphism(-150,500)(1,0)[(B,A)`(B,A) `=]{600}1a
\putmorphism(-180,500)(0,-1)[\phantom{Y_2}`(B, \tilde A) `(B,U)]{450}1l
\put(0,50){\fbox{$(u,U)$}}
\putmorphism(-150,-400)(1,0)[(\tilde B, \tilde A)` `=]{500}1a
\putmorphism(-180,50)(0,-1)[\phantom{Y_2}``(u,\tilde A)]{450}1l
\putmorphism(450,50)(0,-1)[\phantom{Y_2}`(\tilde B, \tilde A)`(\tilde B, U)]{450}1r
\putmorphism(450,500)(0,-1)[\phantom{Y_2}`(\tilde B, A) `(u,A)]{450}1r
\putmorphism(450,50)(1,0)[\phantom{(B, \tilde A)}`(\tilde B, A)`=]{620}1a
\putmorphism(1070,50)(0,-1)[\phantom{(B, \tilde A')}``(u',A)]{450}1r
\putmorphism(1070,-400)(0,-1)[(\tilde B', A)`(\tilde B', \tilde A)`(\tilde B', U)]{450}1r
\putmorphism(450,-850)(1,0)[\phantom{(B, \tilde A)}``=]{520}1a
\putmorphism(450,-400)(0,-1)[\phantom{(B, \tilde A)}`(\tilde B', \tilde A) `(U',\tilde A)]{450}1l
\put(600,-630){\fbox{$(u',U)$}}
\efig}
$$

\noindent $\bullet$ \quad  \axiom{$(k,K)$-l-nat}  
$$
\scalebox{0.86}{
\bfig
 \putmorphism(-150,500)(1,0)[(B,A)`(B', A)`(k, A)]{600}1a
 \putmorphism(450,500)(1,0)[\phantom{A\ot B}`(B', A') `(B', K)]{680}1a
 \putmorphism(-150,50)(1,0)[(B,A)`(B,A')`(B,K)]{600}1a
 \putmorphism(450,50)(1,0)[\phantom{A\ot B}`(B', A') `(k, A')]{680}1a

\putmorphism(-180,500)(0,-1)[\phantom{Y_2}``=]{450}1r
\putmorphism(1100,500)(0,-1)[\phantom{Y_2}``=]{450}1r
\put(350,260){\fbox{$(k,K)$}}

\putmorphism(-150,-400)(1,0)[(\tilde B,A)`(\tilde B,A') `(\tilde B, K)]{640}1a
 \putmorphism(450,-400)(1,0)[\phantom{A'\ot B'}` (\tilde B',A') `(l,A')]{680}1a

\putmorphism(-180,50)(0,-1)[\phantom{Y_2}``(u,A)]{450}1l
\putmorphism(450,50)(0,-1)[\phantom{Y_2}``]{450}1r
\putmorphism(300,50)(0,-1)[\phantom{Y_2}``(u,A')]{450}0r
\putmorphism(1100,50)(0,-1)[\phantom{Y_2}``]{450}1r
\putmorphism(1080,50)(0,-1)[\phantom{Y_2}``(v,A')]{450}0r
\put(-20,-180){\fbox{$(u,K)$}}
\put(660,-180){\fbox{$(\omega,A')$}}

\efig}
\quad=\quad
\scalebox{0.86}{
\bfig
 \putmorphism(-150,500)(1,0)[(B,A)`(B', A)`(k, A)]{600}1a
 \putmorphism(450,500)(1,0)[\phantom{A\ot B}`(B', A') `(B', K)]{680}1a

 \putmorphism(-150,50)(1,0)[(\tilde B,A)`(\tilde B',A)`(l, A)]{600}1a
 \putmorphism(450,50)(1,0)[\phantom{A\ot B}`(\tilde B',A') `(\tilde B',K)]{680}1a
\putmorphism(-180,500)(0,-1)[\phantom{Y_2}``]{450}1l
\putmorphism(-160,500)(0,-1)[\phantom{Y_2}``(u,A)]{450}0l
\putmorphism(450,500)(0,-1)[\phantom{Y_2}``]{450}1l
\putmorphism(610,500)(0,-1)[\phantom{Y_2}``(v,A)]{450}0l 
\putmorphism(1120,500)(0,-1)[\phantom{Y_3}``(v,A')]{450}1r
\put(-40,270){\fbox{$(\omega,A)$}} 
\put(650,270){\fbox{$(v,K)$}}
\putmorphism(-150,-400)(1,0)[(\tilde B,A)`(\tilde B,A') `(\tilde B, K)]{640}1a
 \putmorphism(450,-400)(1,0)[\phantom{A'\ot B'}` (\tilde B',A') `(l,A')]{680}1a

\putmorphism(-180,50)(0,-1)[\phantom{Y_2}``=]{450}1l
\putmorphism(1120,50)(0,-1)[\phantom{Y_3}``=]{450}1r
\put(300,-200){\fbox{$(l,K)$}}

\efig}
$$

\noindent $\bullet$ \quad  \axiom{$(k,K)$-r-nat}  
$$
\scalebox{0.86}{
\bfig
 \putmorphism(-150,500)(1,0)[(B,A)`(B', A)`(k, A)]{600}1a
 \putmorphism(450,500)(1,0)[\phantom{A\ot B}`(B', A') `(B', K)]{680}1a
 \putmorphism(-150,50)(1,0)[(B,A)`(B,A')`(B,K)]{600}1a
 \putmorphism(450,50)(1,0)[\phantom{A\ot B}`(B', A') `(k, A')]{680}1a

\putmorphism(-180,500)(0,-1)[\phantom{Y_2}``=]{450}1r
\putmorphism(1100,500)(0,-1)[\phantom{Y_2}``=]{450}1r
\put(350,260){\fbox{$(k,K)$}}

\putmorphism(-180,50)(0,-1)[\phantom{Y_2}``(B,U)]{450}1l
\putmorphism(450,50)(0,-1)[\phantom{Y_2}``]{450}1r
\putmorphism(300,50)(0,-1)[\phantom{Y_2}``(B,V)]{450}0r
\putmorphism(1100,50)(0,-1)[\phantom{Y_2}``]{450}1r
\putmorphism(1080,50)(0,-1)[\phantom{Y_2}``(B',V)]{450}0r
\put(-20,-180){\fbox{$(B,\zeta)$}}
\put(660,-180){\fbox{$(k,V)$}}

\putmorphism(-150,-400)(1,0)[(B,\tilde A)`(B,\tilde A') `(B,L)]{640}1a
 \putmorphism(450,-400)(1,0)[\phantom{A'\ot B'}` (B',\tilde A') `(k,\tilde A')]{680}1a
\efig}
\quad=\quad
\scalebox{0.86}{
\bfig
 \putmorphism(-150,500)(1,0)[(B,A)`(B', A)`(k, A)]{600}1a
 \putmorphism(450,500)(1,0)[\phantom{A\ot B}`(B', A') `(B', K)]{680}1a

 \putmorphism(-150,50)(1,0)[(B,\tilde A)`(B',\tilde A)`(k,\tilde A)]{600}1a
 \putmorphism(450,50)(1,0)[\phantom{A\ot B}`(B',\tilde A') `(B',L)]{680}1a

\putmorphism(-180,500)(0,-1)[\phantom{Y_2}``]{450}1l
\putmorphism(-160,500)(0,-1)[\phantom{Y_2}``(B,U)]{450}0l
\putmorphism(450,500)(0,-1)[\phantom{Y_2}``]{450}1l
\putmorphism(610,500)(0,-1)[\phantom{Y_2}``(B',U)]{450}0l 
\putmorphism(1120,500)(0,-1)[\phantom{Y_3}``(B',V)]{450}1r
\put(-40,270){\fbox{$(k,U)$}} 
\put(650,270){\fbox{$(B', \zeta)$}}

\putmorphism(-180,50)(0,-1)[\phantom{Y_2}``=]{450}1l
\putmorphism(1120,50)(0,-1)[\phantom{Y_3}``=]{450}1r
\put(300,-200){\fbox{$(k,L)$}}

\putmorphism(-150,-400)(1,0)[(B,\tilde A)`(B,\tilde A') `(B,L)]{640}1a
 \putmorphism(450,-400)(1,0)[\phantom{A'\ot B'}` (B',\tilde A') `(k,\tilde A')]{680}1a
\efig}
$$

\noindent $\bullet$ \quad  \axiom{$(u,U)$-l-nat} 
$$
\scalebox{0.86}{
\bfig
 \putmorphism(-150,500)(1,0)[(B,A)`(B,A) `=]{600}1a
 \putmorphism(450,500)(1,0)[(B,A)` `(k,A)]{450}1a
\putmorphism(-180,500)(0,-1)[\phantom{Y_2}`(B, \tilde A) `(B,U)]{450}1l
\put(-20,50){\fbox{$(u,U)$}}
\putmorphism(-170,-400)(1,0)[(\tilde B, \tilde A)` `=]{480}1a
\putmorphism(-180,50)(0,-1)[\phantom{Y_2}``(u,\tilde A)]{450}1l
\putmorphism(450,50)(0,-1)[\phantom{Y_2}`(\tilde B, \tilde A)`(\tilde B, U)]{450}1l
\putmorphism(450,500)(0,-1)[\phantom{Y_2}`(\tilde B, A) `(u,A)]{450}1l
\put(600,260){\fbox{$(\omega,A)$}}
\putmorphism(430,50)(1,0)[\phantom{(B, \tilde A)}``(l, A)]{500}1a
\putmorphism(1070,50)(0,-1)[\phantom{(B, A')}`(\tilde B', \tilde A)`]{450}1r
\putmorphism(1050,50)(0,-1)[``(\tilde B',U)]{450}0r
\putmorphism(1070,500)(0,-1)[(B', A)`(\tilde B', A)`]{450}1r
\putmorphism(1050,500)(0,-1)[``(v,A)]{450}0r
\putmorphism(450,-400)(1,0)[\phantom{(B, \tilde A)}``(l, \tilde A)]{500}1a
\put(600,-170){\fbox{$(l,U)$}}
\efig}
\quad=\quad
\scalebox{0.86}{
\bfig
 \putmorphism(-150,500)(1,0)[(B,A)`(B',A) `(k,A)]{600}1a
 \putmorphism(450,500)(1,0)[\phantom{(B,A)}` `=]{450}1a
\putmorphism(-180,500)(0,-1)[\phantom{Y_2}`(B, \tilde A) `]{450}1l
\putmorphism(-160,500)(0,-1)[` `(B,U)]{450}0l
\put(620,50){\fbox{$(v,U)$}}
\putmorphism(-170,-400)(1,0)[(\tilde B, \tilde A)` `(l, \tilde A)]{470}1a
\putmorphism(-180,50)(0,-1)[\phantom{Y_2}``]{450}1l
\putmorphism(-160,50)(0,-1)[\phantom{Y_2}``(u,\tilde A)]{450}0l
\putmorphism(450,50)(0,-1)[\phantom{Y_2}`(\tilde B', \tilde A)`(v,\tilde A)]{450}1r
\putmorphism(450,500)(0,-1)[\phantom{Y_2}`(B', \tilde A) `(B',U)]{450}1r
\put(0,260){\fbox{$(k,U)$}}
\putmorphism(-190,50)(1,0)[\phantom{(B, \tilde A)}``(k, \tilde A)]{500}1a
\putmorphism(1070,50)(0,-1)[\phantom{(B, A')}`(\tilde B', \tilde A)`(\tilde B',U)]{450}1r
\putmorphism(1070,500)(0,-1)[(B', A)``(v,A)]{450}1r
\putmorphism(1100,500)(0,-1)[`(\tilde B', A)`]{450}0r
\putmorphism(450,-400)(1,0)[\phantom{(B, \tilde A)}``=]{500}1b
\put(0,-170){\fbox{$(\omega, \tilde A)$}}
\efig}
$$

\noindent $\bullet$ \quad \axiom{$(u,U)$-r-nat} 
$$
\scalebox{0.86}{
\bfig
 \putmorphism(-150,500)(1,0)[(B,A)`(B,A) `=]{600}1a
 \putmorphism(550,500)(1,0)[` `(B,K)]{380}1a
\putmorphism(-180,500)(0,-1)[\phantom{Y_2}`(B, \tilde A) `(B,U)]{450}1l
\put(-20,50){\fbox{$(u,U)$}}
\putmorphism(-170,-400)(1,0)[(\tilde B, \tilde A)` `=]{480}1a
\putmorphism(-180,50)(0,-1)[\phantom{Y_2}``(u,\tilde A)]{450}1l
\putmorphism(450,50)(0,-1)[\phantom{Y_2}`(\tilde B, \tilde A)`(\tilde B, U)]{450}1l
\putmorphism(450,500)(0,-1)[\phantom{Y_2}`(\tilde B, A) `(u,A)]{450}1l
\put(620,280){\fbox{$(u,K)$}}
\putmorphism(430,50)(1,0)[\phantom{(B, \tilde A)}``(\tilde B,K)]{500}1a
\putmorphism(1070,50)(0,-1)[\phantom{(B, A')}`(\tilde B, \tilde A')`]{450}1r
\putmorphism(1090,50)(0,-1)[\phantom{(B, A')}``(\tilde B,V)]{450}0r
\putmorphism(1070,500)(0,-1)[(B, A')`(\tilde B, A')`]{450}1r
\putmorphism(1090,500)(0,-1)[``(u,A')]{450}0r
\putmorphism(450,-400)(1,0)[\phantom{(B, \tilde A)}``(\tilde B, L)]{500}1a
\put(620,-170){\fbox{$ (\tilde{B},\zeta)$ } } 
\efig}
\quad=\quad
\scalebox{0.86}{
\bfig
 \putmorphism(-150,500)(1,0)[(B,A)`(B,A') `(B,K)]{600}1a
 \putmorphism(450,500)(1,0)[\phantom{(B,A)}` `=]{450}1a
\putmorphism(-180,500)(0,-1)[\phantom{Y_2}`(B, \tilde A) `]{450}1l
\putmorphism(-160,500)(0,-1)[\phantom{Y_2}` `(B,U)]{450}0l
\put(620,50){\fbox{$(u,V)$}}
\putmorphism(-170,-400)(1,0)[(\tilde B, \tilde A)` `(\tilde B, L)]{500}1a
\putmorphism(-180,50)(0,-1)[\phantom{Y_2}``]{450}1l
\putmorphism(-160,50)(0,-1)[\phantom{Y_2}``(u,\tilde A)]{450}0l
\putmorphism(450,50)(0,-1)[\phantom{Y_2}`(\tilde B, \tilde A')`(u,\tilde A')]{450}1r
\putmorphism(450,500)(0,-1)[\phantom{Y_2}`(B, \tilde A') `(B,V)]{450}1r
\put(0,260){\fbox{$(B,\zeta)$}}
\putmorphism(-190,50)(1,0)[\phantom{(B, \tilde A)}``(B,L)]{500}1a
\putmorphism(1070,50)(0,-1)[\phantom{(B, A')}`(\tilde B, \tilde A')`(\tilde B,V)]{450}1r
\putmorphism(1070,500)(0,-1)[(B, A')``(u,A')]{450}1r
\putmorphism(1110,500)(0,-1)[`(\tilde B, A')`]{450}0r
\putmorphism(450,-400)(1,0)[\phantom{(B, \tilde A)}``=]{500}1b
\put(0,-170){\fbox{$(u,L)$}}
\efig}
$$
for any 2-cells 
\begin{equation} \eqlabel{omega-zeta}
\scalebox{0.86}{
\bfig
\putmorphism(-150,175)(1,0)[B` B'`k]{450}1a
\putmorphism(-150,-175)(1,0)[\tilde B`\tilde B' `l]{440}1b
\putmorphism(-170,175)(0,-1)[\phantom{Y_2}``u]{350}1l
\putmorphism(280,175)(0,-1)[\phantom{Y_2}``v]{350}1r
\put(0,-15){\fbox{$\omega$}}
\efig}
\quad\text{and}\quad
\scalebox{0.86}{
\bfig
\putmorphism(-150,175)(1,0)[A` A'`K]{450}1a
\putmorphism(-150,-175)(1,0)[\tilde A`\tilde A' `L]{440}1b
\putmorphism(-170,175)(0,-1)[\phantom{Y_2}``U]{350}1l
\putmorphism(280,175)(0,-1)[\phantom{Y_2}``V]{350}1r
\put(0,-15){\fbox{$\zeta$}}
\efig}
\end{equation}
in $\Bb$, respectively $\Aa$. 
\end{prop}

\bigskip

{\bf  Appendix B.2}

\vspace{1,3cm}

{\bf Table 4 from the Appendix of \cite{Fem1}: }

\begin{table}[H]
\begin{center}
\begin{tabular}{ c c c } 
New axiom & \hspace{0,3cm} Origin from $\F\colon\Aa\to\llbracket\Bb,\Cc\rrbracket$ & \hspace{0,3cm}  additional interpretation\\ [0.5ex]
\hline
2-cell $(k,K)$ & part 3 of $(-,K)$ being a h.o.t. 
\\ [1ex]   
2-cell $(u,K)$ & part 2 of $(-,K)$ being a h.o.t.  \\ [1ex]
2-cell $(k,U)$ & part 2 of $(-,U)$ being a v.l.t.  \\ [1ex]
2-cell $(u,U)$ & part 3 of $(-,U)$ being a v.l.t.  \\ [1ex]
\hline
\axiomref{($1_B,K$)}  & \axiomref{h.o.t.-2} of $(-,K)$ \\ [1ex]    
\axiomref{($k,1_A$)}  & \axiomref{m.ho-vl.-1} of unitor \\ [1ex] 
  & $\F_A\colon \Id_{(-,A)}\Rrightarrow (-,1_A)$ & \axiomref{h.o.t.-2} of $(k,-)$ \\ [1ex]   
\axiomref{($1^B,K$)}  & \axiomref{h.o.t.-4} of $(-,K)$ \\ [1ex]    
\axiomref{($u,1_A$)}  & \axiomref{m.ho-vl.-2} of unitor \\ [1ex] 
  & $\F_A\colon \Id_{(-,A)}\Rrightarrow (-,1_A)$ & \axiomref{v.l.t.\x 2} of $(u,-)$ \\ [1ex]   
\axiomref{($1_B,U$)}  & \axiomref{v.l.t.\x 2} of $(-,U)$ \\ [1ex]   
\axiomref{($k,1^A$)}  & \axiomref{lx.f.v2} of $\F$ (is an equality of v.l.t.) \\ [1ex] 
  & evaluated at $k$ & \axiomref{h.o.t.-4} of $(k,-)$ \\ [1ex]    
\axiomref{($1^B,U$)}  & \axiomref{v.l.t.\x 4} of $(-,U)$ \\ [1ex]    
\axiomref{($u,1^A$)}  & \axiomref{lx.f.v2} of $\F$ (is an equality of v.l.t.) \\ [1ex] 
  & evaluated at $u$ & \axiomref{v.l.t.\x 2} of $(u,-)$ \\ [1ex]  
\hline
\axiomref{($k'k,K$)}  & \axiomref{h.o.t.-1} of $(-,K)$ \\ [1ex]    
\axiomref{($k,K'K$)}  & \axiomref{m.ho-vl.-1} of compositor \\ [1ex] 
  & $\F_{LK}\colon (-,L)(-,K)\Rrightarrow (-,LK)$ & \axiomref{h.o.t.-1} of $(k,-)$ \\ [1ex]   
\axiomref{($\frac{u}{u'},K$)}  & \axiomref{h.o.t.-3} of $(-,K)$ \\ [1ex]    
\axiomref{($u,K'K$)}  & \axiomref{m.ho-vl.-2} of compositor \\ [1ex] 
  & $\F_{LK}\colon (-,L)(-,K)\Rrightarrow (-,LK)$ & \axiomref{v.l.t.\x 1} of $(u,-)$ \\ [1ex]   
\axiomref{($k'k,U$)}  & \axiomref{v.l.t.\x 1} of $(-,U)$ \\ [1ex]   
\axiomref{($k,\frac{U}{U'}$)}  & \axiomref{lx.f.v1} of $\F$ (is an equality of v.l.t.) \\ [1ex] 
  & evaluated at $k$ {h.o.t.-3} of $(k,-)$ \\ [1ex]  
\axiomref{($u,\frac{U}{U'}$)}  & \axiomref{lx.f.v1} of $\F$ (is an equality of v.l.t.) \\ [1ex] 
  & evaluated at $u$ & \axiomref{v.l.t.\x 3} of $(u,-)$ \\ [1ex]
\axiomref{($\frac{u}{u'},U$)}  & \axiomref{v.l.t.\x 3} of \axiomref{$-,U$} \\ [1ex]   
\hline
\axiomref{$(k,K)$-l-nat}  & \axiomref{h.o.t.-5} of \axiomref{$-,K$} & \axiomref{m.hl-vo.-1} of \axiomref{$\omega,-$} \\ [1ex]   
\axiomref{$(k,K)$-r-nat}  & \axiomref{m.ho-vl.-1} of \axiomref{$-,\zeta$} & \axiomref{h.o.t.-5} of $(k,-)$ \\ [1ex]   
\axiomref{$(u,U)$-l-nat}  & \axiomref{v.l.t.\x 5} of \axiomref{$-,U$} & \axiomref{m.hl-vo.-2} of \axiomref{$\omega,-$} \\ [1ex]   
\axiomref{$(u,U)$-r-nat}  & \axiomref{m.ho-vl.-2} of \axiomref{$-,\zeta$} & \axiomref{v.l.t.\x 5} of $(u,-)$ \\ [1ex]   
\end{tabular}
\caption{Generation of a lax double quasi-functor $\Aa\times\Bb\protect\to\Cc$}
\label{table:1}
\end{center}
\end{table}

{\bf Appendix B.3}

\bigskip

The following two definitions are from \cite{Fem}. 

\begin{defn} \delabel{lax v tr cubical}
A vertical lax transformation $\theta_0\colon (-,-)_1\Rightarrow (-,-)_2$ between lax double quasi-functors $(-,-)_1,(-,-)_2\colon 
\Aa\times\Bb\to\Cc$ is given by: for each $A\in\Aa$ a vertical lax transformation $\theta_0^A\colon (-,A)_1\Rightarrow(-,A)_2$ and 
for each $B\in\Bb$ a vertical lax transformation $\theta_0^B\colon (B,-)_1\Rightarrow(B,-)_2$, both of lax double functors, such that 
$(\theta_0^A)_B=(\theta_0^B)_A$ and such that 

\medskip
\noindent \axiom{$VLT^q_1$} \vspace{-0,8cm}

$$
\scalebox{0.8}{
\bfig
 \putmorphism(-170,500)(1,0)[(B,A)_1`(B,A)_1 `=]{500}1a
\putmorphism(-180,500)(0,-1)[\phantom{Y_2}`(B, A)_2 `(\theta_0^A)_B]{450}1l
\put(-50,250){\fbox{$(\theta_0^A)^u$}}
\putmorphism(-150,-400)(1,0)[(\tilde B, A)_2`(\tilde B, A)_2 `=]{500}1a
\putmorphism(-180,50)(0,-1)[\phantom{Y_2}``(u,A)_2]{450}1l
\putmorphism(300,50)(0,-1)[\phantom{Y_2}``(\theta_0^A)_{\tilde B}]{450}1l
\putmorphism(300,500)(0,-1)[\phantom{Y_2}`(\tilde B, A)_1 `(u,A)_1]{450}1r
\putmorphism(-720,50)(1,0)[(B, A)_2``=]{400}1a
\putmorphism(-720,50)(0,-1)[\phantom{(B, \tilde A')}``(B,U)_2]{450}1l
\putmorphism(-720,-400)(0,-1)[(B, \tilde A)_2`(\tilde B, \tilde A)_2`(u,\tilde A)_2]{450}1l
\putmorphism(-720,-850)(1,0)[\phantom{(B, \tilde A)}``=]{370}1a 
\putmorphism(-180,-400)(0,-1)[`(\tilde B, \tilde A)_2 `(\tilde B, U)_2]{450}1r
\put(-610,-630){\fbox{$(u,U)_2$}}

\putmorphism(430,50)(1,0)[`(\tilde B, A)_1`=]{450}1a
\putmorphism(870,50)(0,-1)[\phantom{Y_2}`(\tilde B, \tilde A)_1`(\tilde B, U)_1]{450}1r
\putmorphism(300,-400)(0,-1)[`(\tilde B, \tilde A)_2 `(\tilde B,U)_2]{450}1r
\putmorphism(870,-400)(0,-1)[`(\tilde B, \tilde A)_2 `(\theta_0^{\tilde B})_{\tilde A}]{450}1r
\putmorphism(450,-850)(1,0)[` `=]{260}1a 
\put(430,-200){\fbox{$(\theta_0^{\tilde B})^U$}}
\efig}
=
\scalebox{0.8}{
\bfig
 \putmorphism(0,500)(1,0)[(B,A)_1`(B,A)_1 `=]{500}1a
\putmorphism(0,500)(0,-1)[\phantom{Y_2}`(B, A)_2 `(\theta^A_0)_B]{450}1l
\put(70,-200){\fbox{$(\theta^B_0)^U$}}
\putmorphism(-40,-400)(1,0)[(B, \tilde A)_2` `=]{330}1a %
\putmorphism(0,50)(0,-1)[\phantom{Y_2}``(B,U)_2]{450}1l
\putmorphism(450,50)(0,-1)[\phantom{Y_2}`(B, \tilde A)_2`(\theta^B_0)_{\tilde A}]{450}1r
\putmorphism(450,500)(0,-1)[\phantom{Y_2}`(B, \tilde A)_1 `(B,U)_1]{450}1l
\putmorphism(450,50)(1,0)[\phantom{(B, \tilde A)}`(B, \tilde A)_1`=]{620}1a %
\putmorphism(970,50)(0,-1)[\phantom{(B, \tilde A')}``(u,\tilde A)_1]{450}1r
\putmorphism(970,-400)(0,-1)[(\tilde B, \tilde A)_1`(\tilde B, \tilde A)_2`(\theta^{\tilde A}_0)_{\tilde B}]{450}1r
\putmorphism(450,-850)(1,0)[\phantom{(B, \tilde A)}``=]{350}1a 
\putmorphism(450,-400)(0,-1)[\phantom{(B, \tilde A)}`(\tilde B, \tilde A) `(u,\tilde A)_2]{450}1l
\put(600,-630){\fbox{$(\theta^{\tilde A}_0)^u$}}

 \putmorphism(970,500)(1,0)[(B,A)_1`(B,A)_1 `=]{500}1a
\putmorphism(970,500)(0,-1)[\phantom{Y_2}` `(B,U)_1]{450}1l
\putmorphism(1500,500)(0,-1)[\phantom{Y_2}`(\tilde B, A)_1 `(u,A)_1]{450}1r
\putmorphism(1500,50)(0,-1)[\phantom{Y_2}`(\tilde B, \tilde A)_1 `(\tilde B,U)_1]{450}1r
\putmorphism(1110,-400)(1,0)[` `=]{260}1a
\put(1050,200){\fbox{$(u,U)_1$}}
\efig}
$$
for every 1v-cells $U\colon A\to \tilde A$ and $u\colon B\to\tilde B$; 



\noindent \axiom{$VLT^q_2$} \vspace{-0,3cm}
$$\scalebox{0.86}{
\bfig
 \putmorphism(-150,500)(1,0)[(B,A)_1`(B,A)_1 `=]{600}1a
 \putmorphism(450,500)(1,0)[\phantom{(B,A)_1}` `(B,K)_1]{450}1a
\putmorphism(-200,500)(0,-1)[\phantom{Y_2}`(B, A)_2 `(\theta_0^A)_B]{450}1l
\put(-20,50){\fbox{$(\theta^A_0)^u$}}
\putmorphism(-170,-400)(1,0)[(\tilde B, A)_2` `=]{500}1a
\putmorphism(-200,50)(0,-1)[\phantom{Y_2}``(u,A)_2]{450}1l
\putmorphism(450,50)(0,-1)[\phantom{Y_2}`(\tilde B, A)_2 `(\theta_0^A)_{\tilde B}]{450}1l
\putmorphism(450,500)(0,-1)[\phantom{Y_2}`(\tilde B, A)_1 `(u,A)_1]{450}1l
\put(600,260){\fbox{$(u,K)_1$}}
\putmorphism(450,50)(1,0)[\phantom{(B, \tilde A)}``(\tilde B,K)_1]{500}1a
\putmorphism(1070,50)(0,-1)[\phantom{(B, A')}`(\tilde B, A')_2`(\theta^{\tilde B}_0)_{A'}]{450}1r
\putmorphism(1070,500)(0,-1)[(B, A')_1`(\tilde B, A')_1`(u,A')_1]{450}1r
\putmorphism(450,-400)(1,0)[\phantom{(B, \tilde A)}``(\tilde B,K)_2]{500}1a
\put(600,-170){\fbox{$(\theta_0^{\tilde B})_K$}}
\efig}
\quad=\quad
\scalebox{0.86}{
\bfig
 \putmorphism(-150,500)(1,0)[(B,A)_1`(B,A')_1 `(B,K)_1]{600}1a
 \putmorphism(450,500)(1,0)[\phantom{(B,A)}` `=]{540}1a
\putmorphism(-180,500)(0,-1)[\phantom{Y_2}`(B, A)_2 `(\theta^A_0)_B]{450}1l
\put(0,280){\fbox{$(\theta_0^B)_K$}}
\putmorphism(-180,-400)(1,0)[(\tilde B,A)_2` `(\tilde B,K)_2]{500}1a
\putmorphism(-180,50)(0,-1)[\phantom{Y_2}``(u,A)_2]{450}1l
\putmorphism(450,50)(0,-1)[\phantom{Y_2}`(\tilde B, A')_2`(u, A')_2]{450}1r
\putmorphism(450,500)(0,-1)[\phantom{Y_2}`(B, A')_2 `(\theta^{A'}_0)_B]{450}1r
\put(620,50){\fbox{$(\theta_0^{A'})^u$}}
\putmorphism(-180,50)(1,0)[\phantom{(B, \tilde A)}``(B,K)_2]{500}1a
\putmorphism(1130,50)(0,-1)[\phantom{(B, A')}`(\tilde B, A')_2`(\theta_0^{A'})_{\tilde B}]{450}1r
\putmorphism(1130,500)(0,-1)[(B, A')_1`(\tilde B, A')_1`(u,A')_1]{450}1r
\putmorphism(450,-400)(1,0)[\phantom{(B, \tilde A)}``=]{520}1b
\put(0,-170){\fbox{$(u,K)_2$}}
\efig}
$$
for every 1h-cell $K\colon A\to A'$ and 1v-cell $u\colon B\to \tilde B$, 

\medskip
\noindent \axiom{$VLT^q_3$} \vspace{-0,6cm}

$$\scalebox{0.86}{
\bfig
 \putmorphism(-150,500)(1,0)[(B,A)_1`(B,A)_1 `=]{600}1a
 \putmorphism(450,500)(1,0)[\phantom{(B',A)_1}` `(k,A)_1]{450}1a
\putmorphism(-200,500)(0,-1)[\phantom{Y_2}`(B, A)_2 `(\theta_0^B)_A]{450}1l
\put(-20,50){\fbox{$(\theta^B_0)^U$}}
\putmorphism(-170,-400)(1,0)[(B, \tilde A)_2` `=]{500}1a
\putmorphism(-200,50)(0,-1)[\phantom{Y_2}``(B,U)_2]{450}1l
\putmorphism(450,50)(0,-1)[\phantom{Y_2}`(B, \tilde A)_2 `(\theta_0^B)_{\tilde A}]{450}1l 
\putmorphism(450,500)(0,-1)[\phantom{Y_2}`(B, \tilde A)_1 `(B,U)_1]{450}1l
\put(600,260){\fbox{$(k,U)_1$}}
\putmorphism(450,50)(1,0)[\phantom{(B, \tilde A)}``(k,\tilde A)_1]{500}1a %
\putmorphism(1070,50)(0,-1)[\phantom{(B, A')}`(B', \tilde A)_2`(\theta^{\tilde A}_0)_{B'}]{450}1r
\putmorphism(1070,500)(0,-1)[(B, A')_1`(B', \tilde A)_1`(B',U)_1]{450}1r
\putmorphism(450,-400)(1,0)[\phantom{(B, \tilde A)}``(k,\tilde A)_2]{500}1a
\put(600,-170){\fbox{$(\theta_0^{\tilde A})_k$}}
\efig}
\quad=\quad
\scalebox{0.86}{
\bfig
 \putmorphism(-150,500)(1,0)[(B,A)_1`(B',A)_1 `(k,A)_1]{600}1a
 \putmorphism(450,500)(1,0)[\phantom{(B,A)}` `=]{540}1a
\putmorphism(-180,500)(0,-1)[\phantom{Y_2}`(B, A)_2 `(\theta^A_0)_B]{450}1l
\put(0,280){\fbox{$(\theta_0^A)_k$}}
\putmorphism(-180,-400)(1,0)[(\tilde B,A)_2` `(k,\tilde A)_2]{500}1a
\putmorphism(-180,50)(0,-1)[\phantom{Y_2}``(B,U)_2]{450}1l
\putmorphism(450,50)(0,-1)[\phantom{Y_2}`(B', \tilde A)_2`(B', U)_2]{450}1r
\putmorphism(450,500)(0,-1)[\phantom{Y_2}`(B', A)_2 `(\theta^{B'}_0)_A]{450}1r
\put(640,50){\fbox{$(\theta_0^{B'})^U$}}
\putmorphism(-180,50)(1,0)[\phantom{(B, \tilde A)}``(k,A)_2]{500}1a
\putmorphism(1150,50)(0,-1)[\phantom{(B, A')}`(B',\tilde A)_2`(\theta_0^{\tilde A})_{B'}]{450}1r
\putmorphism(1150,500)(0,-1)[(B', A)_1`(B',\tilde A)_1`(B', U)_1]{450}1r
\putmorphism(450,-400)(1,0)[\phantom{(B, \tilde A)}``=]{520}1b
\put(0,-170){\fbox{$(k,U)_2$}}
\efig}
$$
for every 1v-cell $U\colon A\to \tilde A$ and 1h-cell $k\colon B\to B'$, and 


\medskip
\noindent \axiom{$VLT^q_4$} \vspace{-0,3cm}
$$\scalebox{0.86}{
\bfig
\putmorphism(-150,500)(1,0)[(B,A)_1`(B',A)_1`(k,A)_1]{600}1a
 \putmorphism(480,500)(1,0)[\phantom{F(A)}`(B',A')_1 `(B',K)_1]{640}1a
 \putmorphism(-150,50)(1,0)[(B,A)_1`(B,A')_1`(B,K)_1]{600}1a
 \putmorphism(470,50)(1,0)[\phantom{F(A)}`(B',A')_1 `(k,A')_1]{660}1a

\putmorphism(-180,500)(0,-1)[\phantom{Y_2}``=]{450}1r
\putmorphism(1100,500)(0,-1)[\phantom{Y_2}``=]{450}1r
\put(320,280){\fbox{$(k,K)_1$}}

\putmorphism(-150,-400)(1,0)[(B,A)_2`(B,A')_2 `(B,K)_2]{640}1a
 \putmorphism(490,-400)(1,0)[\phantom{F(A')}` (B',A')_2 `(k,A')_2]{640}1a

\putmorphism(-180,50)(0,-1)[\phantom{Y_2}``(\theta_0^A)_B]{450}1l %
\putmorphism(450,50)(0,-1)[\phantom{Y_2}``]{450}1l
\putmorphism(610,50)(0,-1)[\phantom{Y_2}``(\theta_0^B)_{A'}]{450}0l 
\putmorphism(1120,50)(0,-1)[\phantom{Y_3}``(\theta_0^{A'})_{B'}]{450}1r
\put(-40,-180){\fbox{$(\theta^B_0)_K$}} 
\put(620,-180){\fbox{$(\theta_0^{A'})_k$}}
\efig}
\quad
=
\quad
\scalebox{0.86}{
\bfig
\putmorphism(-150,500)(1,0)[(B,A)_1`(B',A)_1`(k,A)_1]{600}1a
 \putmorphism(480,500)(1,0)[\phantom{F(A)}`(B',A')_1 `(B',K)_1]{640}1a

 \putmorphism(-150,50)(1,0)[(B,A)_2`(B',A)_2`(k,A)_2]{600}1a
 \putmorphism(450,50)(1,0)[\phantom{F(A)}`(B',K)_2 `(B', K)_2]{640}1a

\putmorphism(-180,500)(0,-1)[\phantom{Y_2}``(\theta_0^A)_B]{450}1l
\putmorphism(450,500)(0,-1)[\phantom{Y_2}``]{450}1r
\putmorphism(300,500)(0,-1)[\phantom{Y_2}``(\theta_0^{B'})_A]{450}0r
\putmorphism(1100,500)(0,-1)[\phantom{Y_2}``(\theta_0^{B'})_{A'}]{450}1r
\put(-40,280){\fbox{$(\theta_0^{A})_k$}}
\put(700,280){\fbox{$(\theta_0^{B'})_K$}}

\putmorphism(-150,-400)(1,0)[(B,A)_2`(B,A')_2 `(B,K)_2]{640}1a
 \putmorphism(490,-400)(1,0)[\phantom{F(A')}` (B',A')_2 `(k,A')_2]{640}1a

\putmorphism(-180,50)(0,-1)[\phantom{Y_2}``=]{450}1l
\putmorphism(1120,50)(0,-1)[\phantom{Y_3}``=]{450}1r
\put(320,-200){\fbox{$(k,K)_2$}}
\efig}
$$
for every 1h-cells $K\colon A\to A'$ and $k\colon B\to B'$. 
\end{defn}

\medskip

\begin{defn} \delabel{modif btw horiz}
Let horizontal oplax transformations $\theta, \theta'$ and vertical lax transformations $\theta_0, \theta'_0$ acting between lax double quasi-functors 
$H_1, H_2, H_3, H_4\colon \Aa\times\Bb\to\Cc$ be given as in the left diagram below. Denote by $(-,A)_i\colon\Bb\to\Cc, (B,-)_i\colon\Aa\to\Cc, i=1,2,3,4$ the pairs of 
lax double functors corresponding to $H_1, H_2, H_3, H_4$, respectively. A modification $\Theta$ (on the left below) is given by a pair of modifications 
$\tau^A, \tau^B$ acting between transformations among lax double functors:
\begin{equation} \eqlabel{q-modif}
\scalebox{0.86}{
\bfig
\putmorphism(-150,50)(1,0)[H_1` H_2` \theta]{400}1a
\putmorphism(-150,-270)(1,0)[H_3 ` H_4 ` \theta' ]{400}1b
\putmorphism(-170,50)(0,-1)[\phantom{Y_2}``\theta_0]{320}1l
\putmorphism(250,50)(0,-1)[\phantom{Y_2}``\theta_0']{320}1r
\put(-30,-140){\fbox{$\tau$}}
\efig}
\qquad\qquad
\scalebox{0.86}{
\bfig
\putmorphism(-180,50)(1,0)[(-,A)_1` (-,A)_2`\theta^A]{550}1a
\putmorphism(-180,-270)(1,0)[(-,A)_3`(-,A)_4 `\theta'^A]{550}1b
\putmorphism(-170,50)(0,-1)[\phantom{Y_2}``\theta_0^A]{320}1l
\putmorphism(350,50)(0,-1)[\phantom{Y_2}``\theta_0'^A]{320}1r
\put(0,-140){\fbox{$\tau^A$}}
\efig}
\qquad\qquad
\scalebox{0.86}{
\bfig
\putmorphism(-180,50)(1,0)[(B,-)_1` (B,-)_2`\theta^B]{550}1a
\putmorphism(-180,-270)(1,0)[(B,-)_3`(B,-)_4 `\theta'^B]{550}1b
\putmorphism(-170,50)(0,-1)[\phantom{Y_2}``\theta_0^B]{320}1l
\putmorphism(350,50)(0,-1)[\phantom{Y_2}``\theta_0'^B]{320}1r
\put(0,-140){\fbox{$\tau^B$}}
\efig}
\end{equation}
such that $\tau^A_B=\tau^B_A$ for every $A\in\Aa, B\in\Bb$. 
\end{defn}

\bigskip

{\bf\large{Appendix C}}

\bigskip

{\bf Appendix C.1}

\bigskip

{\bf Transformations, modifications and perturbations of Gray-categories}

\medskip

A lax transformation $\alpha:F\to G$  for each 1-cell $f:A\to A'$ has a 2-cell $\alpha_f:[\alpha(A)\vert G(f)]\Rightarrow
[F(f)\vert\alpha(A')]$, and for composable 1-cells $A\stackrel{f}\to A'\stackrel{f'}\to A''$ and each 2-cell $a:f\rightarrow g$ it has 
invertible 3-cells $\alpha^2_{f',f}: \threefrac{\alpha_f}{\alpha_{f'}}{G^2_{f',f}}\Rrightarrow\frac{F^2_{f',f}}{\alpha_{f'f}}$ (which we call {\em cocycles}) and 
$\alpha_a:\frac{\alpha_f}{F(a)}\Rrightarrow \frac{G(a)}{\alpha_g}$ obeying the axioms \\ 
$$(i)  \quad \alpha_{1_a}=\Id_{\alpha(A)}, \qquad 
(ii) \quad \frac{\alpha_a}{G(\xi)}=\frac{F(\xi)}{\alpha_b}, \qquad 
(iii) \quad \frac{\alpha_a}{\alpha_{a'}}=\alpha_\frac{a}{a'}, \qquad 
(iv) \quad \alpha_{id_f}=\Id_{\alpha_f}$$

$(v) \quad \frac{\alpha^2_{f',f}}{\alpha^2_{f'',f'f}}=\frac{\alpha^2_{f'',f'}}{\alpha^2_{f''f',f}}, \qquad
(vi) \quad \alpha^2_{1_B,f}=\Id_{\alpha_f}=\alpha^2_{f,1_A}$

$(vii) \quad (L) \quad \frac{\alpha_\psi}{\alpha^2_{g',f}}=\frac{\alpha^2_{f',f}}{\alpha^2_{[id_f\vert\psi]}}, \qquad
(vii) \quad (R) \quad \frac{\alpha_\phi}{\alpha^2_{f',g}}=\frac{\alpha^2_{f',f}}{\alpha^2_{[\phi\vert id_{f'}]}}$ 
\smallskip

\noindent for every 3-cell $\xi:a\Rrightarrow b$, vertically composable 2-cells $a,a'$, composable 1-cells $f,f',f''$, and 2-cells $\phi:f\Rightarrow g, \psi:f'\Rightarrow g'$. 

\bigskip

A lax modification $\Gamma: \alpha\Rightarrow\beta:F\to G$ among lax transformations for each object $A$ has a 2-cell $\Gamma(A):\alpha(A)\Rightarrow\beta(A)$, and at every 1-cell $f:A\to B$ it has a 3-cell $\Gamma_f:\frac{\Gamma(A)}{\beta_f}\Rrightarrow\frac{\alpha_f}{\Gamma(B)}$ satisfying $\Gamma_{id_A}=\Id_{\Gamma(A)}$ and for composable 1-cells $A\stackrel{f}\to A'\stackrel{f'}\to A''$ and each 2-cell $a:f\rightarrow g$ the axioms 
$$\threefrac{\Gamma_f}{\Gamma_{f'}}{\alpha^2_{f'f}}=\frac{\beta^2_{f'f}}{\Gamma_{f'f}}\quad\text{and}\quad 
\frac{\Gamma_f}{\alpha_a}=\frac{\beta_a}{\Gamma_g}.$$

\bigskip

A lax perturbation $\Pi:\Gamma\Rrightarrow\Lambda: \alpha\Rightarrow\beta:F\to G$ among lax modifications at every object $A$ has a 3-cell $\Pi_A$ obeying the axiom 
$$\frac{\Pi_A}{\Lambda_f}=\frac{\Gamma_f}{\Pi_{A'}}$$ 
for every 1-cell $f:A\to A'$. 

\bigskip
\bigskip

An oplax transformation $\alpha:F\to G$  for each 1-cell $f:A\to A'$ has a 2-cell $\alpha_f:[F(f)\vert\alpha(A')]\Rightarrow
[\alpha(A)\vert G(f)]$, and for composable 1-cells $A\stackrel{f}\to A'\stackrel{f'}\to A''$ and each 2-cell $a:f\rightarrow g$ it has 
invertible 3-cells $\alpha^2_{f',f}: \threefrac{\alpha_{f'}}{\alpha_f}{G^2_{f',f}}\Rrightarrow\frac{F^2_{f',f}}{\alpha_{f'f}}$ (which we call {\em cocycles}) and 
$\alpha_a:\frac{F(a)}{\alpha_g}\Rrightarrow \frac{\alpha_f}{G(a)}$ obeying analogous axioms as a lax transformation above.

\bigskip

An oplax modification $\Gamma: \alpha\Rightarrow\beta:F\to G$ among oplax transformations for each object $A$ has a 2-cell 
$\Gamma(A):\alpha(A)\Rightarrow\beta(A)$, and at every 1-cell $f:A\to B$ it has a 3-cell $\Gamma_f:\frac{\alpha_f}{\Gamma(A)}\Rrightarrow\frac{\Gamma(B)}{\beta_f}$ satisfying $\Gamma_{id_A}=\Id_{\Gamma(A)}$ and for composable 1-cells $A\stackrel{f}\to A'\stackrel{f'}\to A''$ and each 2-cell $a:f\rightarrow g$ the axioms 
$$\threefrac{\Gamma_f}{\Gamma_{f'}}{\beta^2_{f'f}}=\frac{\alpha^2_{f'f}}{\Gamma_{f'f}}\quad\text{and}\quad 
\frac{\alpha_a}{\Gamma_f}=\frac{\Gamma_g}{\beta_a}.$$

\bigskip

An oplax perturbation $\Pi:\Gamma\Rrightarrow\Lambda: \alpha\Rightarrow\beta:F\to G$ among oplax modifications at every object $A$ has a 3-cell $\Pi_A$ obeying the axiom 
\begin{equation}\eqlabel{op-pert}
\frac{\Gamma_f}{\Pi_{A'}}=\frac{\Pi_A}{\Lambda_f}
\end{equation}
for every 1-cell $f:A\to A'$. (Note that in our short fraction notation the axioms of a lax and an oplax perturbation look the same, though their domain and codomain 2-cells differ, so the full forms of the axioms are indeed different.)

\bigskip

{\bf Appendix C.2}

\medskip


The first six axioms for lax incubators (see \ssref{incub}) in our notation using fractions for transversal composition of 3-cells: 
$$\axiom{$A_1, A_2:B:C$}, \,\ \axiom{$A:B_1, B_2:C$}, \,\ \axiom{$A:B:C_1,C_2$}$$ 
$$\axiom{$\alpha:B:C$}, \,\ \axiom{$A:\beta: C$}, \,\ \axiom{$A:B:\gamma$}$$ 

\axiom{$f'f,g,h$} \qquad $\fourfrac{(-,g,C)^2_{f',f}}{(-,B',h)^2_{f',f}}{(f,g,h)}{(f',g,h)}=
\threefrac{(f'f,g,h)}{(-,B,h)^2_{f',f}}{(-,g,C')^2_{f',f}}$ \qquad\qquad 
\axiom{$a,g,h$} \qquad $\threefrac{(-,g,C)_a}{(-,B,h)_a}{(f,g,h)}=\threefrac{(f',g,h)}{(-,B,h)_a}{(-,g,C')_a}$

\axiom{$f, g'g,h$} \qquad $\threefrac{(f,-,C)^2_{g',g}}{(f,g'g,h)}{(A',-,h)^2_{g',g}}=
\fourfrac{(A,-,h)^2_{g',g}}{(f,g',h)}{(f,g,h)}{(f,-,C')^2_{g',g}}$ \qquad\qquad 
\axiom{$f,b,h$} \qquad $\threefrac{(f,-,C)_b}{(f,g',h)}{(A',-,h)_b}=\threefrac{(A,-,h)_b}{(f,g,h)}{(f,-,C')_b}$

\axiom{$f,g, h'h$}  \qquad $\fourfrac{(f,g,h)}{(f,g,h')}{(f,B,-)^2_{h',h}}{(A',g,-)^2_{h',h}}=
\threefrac{(A,g,-)^2_{h',h}}{(f,B',-)^2_{h',h}}{(f,g,h'h)}$ \qquad\qquad 
\axiom{$f,g, c$} \qquad $\threefrac{(f,g,h)}{(f,B,-)_c}{(A',g,-)_c}=\threefrac{(A,g,-)_c}{(f,B',-)_c}{(f,g,h')}$.

\end{document}